\documentclass[12pt,fullpage]{article}

\usepackage{amsfonts}
\usepackage{xr}
\usepackage{graphicx}
\usepackage{amsmath}
\usepackage{epsfig}

\def\R{\mathbb{R}}
\def\N{\mathbb{N}}

\def\epsilon{\varepsilon}
\def\hat{\widehat}
\def\tilde{\widetilde}
\def\div{\mbox{div}}
\def\Id{\ {\rm{Id}}}
\def\trait (#1) (#2) (#3){\vrule width #1pt height #2pt depth #3pt}
\def\fin{\hfill\trait (0.1) (5) (0) \trait (5) (0.1) (0) \kern-5pt 
\trait (5) (5) (-4.9) \trait (0.1) (5) (0)}

\newcommand{\SE}{\setcounter{equation}{0} \section}
\newcommand{\be}{\begin{equation}}
\newcommand{\ee}{\end{equation}}
\newcommand{\baa}{\begin{array}}
\newcommand{\eaa}{\end{array}}
\newcommand{\ba}{\begin{eqnarray}}
\newcommand{\ea}{\end{eqnarray}}

\newtheorem{theo}{\bf Theorem}[section]
\newtheorem{lem}[theo]{\bf Lemma}
\newtheorem{pro}[theo]{\bf Proposition}
\newtheorem{cor}[theo]{\bf Corollary}

\newtheorem{rem}[theo]{\bf Remark}

\begin{document}
\date{}
\title{\bf{Rearrangement inequalities and applications to isoperimetric problems for eigenvalues}}
\author{Fran\c cois Hamel$^{\hbox{\small{ a}}}$, Nikolai 
Nadirashvili$^{\hbox{\small{ b}}}$ and Emmanuel Russ$^{\hbox{\small{ 
a}}}$\\
\\
\footnotesize{$^{\hbox{a }}$Universit\'e Aix-Marseille III, LATP, 
Facult\'e des Sciences et Techniques, Case cour A}\\
\footnotesize{Avenue Escadrille Normandie-Niemen, F-13397 Marseille 
Cedex 20, France}\\
\footnotesize{$^{\hbox{b }}$ CNRS, LATP, CMI, 39 rue F. Joliot-Curie, 
F-13453 Marseille Cedex 13, France}}
\maketitle

\begin{abstract}
Let $\Omega$ be a bounded $C^{2}$ domain in $\R^n$, where $n$ is any positive integer, and let $\Omega^{\ast}$ be the Euclidean ball centered at $0$ and having the same Lebesgue measure as $\Omega$. Consider the operator $L=-\div(A\nabla)+v\cdot \nabla +V$ on $\Omega$ with Dirichlet boundary condition, where the symmetric matrix field $A$ is in $W^{1,\infty}(\Omega)$, the vector field $v$ is in $L^{\infty}(\Omega,\R^n)$ and $V$ is a continuous function in $\overline{\Omega}$. We prove that minimizing the principal eigenvalue of $L$ when the Lebesgue measure of $\Omega$ is fixed and when $A$, $v$ and $V$ vary under some constraints is the same as minimizing the principal eigenvalue of some operators $L^*$ in the ball $\Omega^*$ with smooth and radially symmetric coefficients. The constraints which are satisfied by the original coefficients in $\Omega$ and the new ones in $\Omega^*$ are expressed in terms of some distribution functions or some integral, pointwise or geometric quantities. Some strict comparisons are also established when $\Omega$ is not a ball. To these purposes, we associate to the principal eigenfunction $\varphi$ of $L$ a new symmetric rearrangement defined on $\Omega^*$, which is different from the classical Schwarz symmetrization, and which preserves the integral of $\div(A\nabla\varphi)$ on suitable equi-measurable sets. A substantial part of the paper is devoted to the proofs of pointwise and integral inequalities of independent interest which are satisfied by this rearrangement. The comparisons for the eigenvalues hold for general operators of the type $L$ and they are new even for symmetric operators. Furthermore they generalize, in particular, and provide an alternative proof of the well-known Rayleigh-Faber-Krahn isoperimetric inequality about the principal eigenvalue of the Laplacian under Dirichlet boundary condition on a domain with fixed Lebesgue measure.
\end{abstract}

\tableofcontents


\SE{Introduction}\label{intro}

Throughout all the paper, we fix an integer $n\geq 1$ and denote by $\alpha_n=\pi^{n/2}/\Gamma(n/2+1)$ the Lebesgue measure of the Euclidean unit ball in $\R^n$. By ``domain'', we mean a non-empty open connected subset of $\R^n$, and we denote by ${\mathcal C}$ the set of all bounded domains of $\R^n$ which are of class $C^{2}$. Throughout all the paper, unless otherwise specified, $\Omega$ will always be in the class $\mathcal{C}$. For any measurable subset $A\subset\R^n$, $\left\vert A\right\vert$ stands for the Lebesgue measure of $A$. If $\Omega\in {\mathcal C}$, $\Omega^{\ast}$ will denote the 
Euclidean ball centered at $0$ such that
$$\left\vert\Omega^{\ast}\right\vert=\left\vert \Omega\right\vert.$$
Define also $C(\overline{\Omega})$ (resp. $C(\overline{\Omega},\R^n)$) the space of real-valued (resp. $\R^n$-valued) continuous functions on $\overline{\Omega}$. For all $x\in \R^n\setminus\left\{0\right\}$, set 
\begin{equation} \label{radial}
e_r(x)=\frac{x}{\left\vert x\right\vert},
\end{equation}
where $\left\vert x\right\vert$ denotes the Euclidean norm of $x$. Finally, if $\Omega\in {\mathcal C}$, if $v:\Omega\rightarrow \R^n$ is measurable and if $1\leq p\leq +\infty$, we say that $v\in L^{p}(\Omega,\R^n)$ if $\left\vert v\right\vert\in L^{p}(\Omega)$, and write (somewhat abusively) $\left\Vert v\right\Vert_{p}$ or $\|v\|_{L^p(\Omega,\R^n)}$ instead of $\left\Vert\left\vert v\right\vert \right\Vert_{L^p(\Omega,\R)}$. 

Various rearrangement techniques for functions defined on $\Omega$ were considered in the literature. The most famous one is the Schwarz symmetrization. Let us briefly recall what the idea of this symmetrization is. For any function $u\in L^1(\Omega)$, denote by $\mu_u$ the 
distribution function of $u$, given by
\[
\mu_u(t)=\left\vert \left\{x\in \Omega;\ u(x)>t\right\}\right\vert
\]
for all $t\in\R$. Note that $\mu$ is right-continuous, non-increasing and $\mu_u(t)\to 0$ (resp. $\mu_u(t)\to|\Omega|$) as $t\to+\infty$ (resp. $t\to-\infty$). For all  $x\in\Omega^{\ast}\backslash\{0\}$, define
\[
u^{\ast}(x)=\sup\left\{t\in\R;\ \mu_u(t)\geq\alpha_n\left\vert x\right\vert^n\right\}.
\]
The function $u^{\ast}$ is clearly radially symmetric, non-increasing in the variable $|x|$ and it satisfies
$$\left|\left\{x\in\Omega,\ u(x)>\zeta\right\}\right|=\left|\left\{x\in\Omega^{\ast},\ u^{\ast}(x)>\zeta\right\}\right|$$
for all $\zeta\in\R$. An essential property of the Schwarz symmetrization is the following one: if $u\in H^1_0(\Omega)$, then $|u|^{\ast}\in H^1_0(\Omega^{\ast})$ and (see \cite{polyaszego})
\begin{equation} \label{Schwarzproperty}
\left\Vert |u|^{\ast}\right\Vert_{L^2(\Omega^{\ast})}=\left\Vert u\right\Vert_{L^2(\Omega)}\mbox{ and } \left\Vert \nabla |u|^{\ast}\right\Vert_{L^2(\Omega^{\ast})}\leq \left\Vert \nabla u\right\Vert_{L^2(\Omega)}.
\end{equation}
One of the main applications of this rearrangement technique is the resolution of optimization problems for the eigenvalues of some second-order elliptic operators on $\Omega$. Let us briefly recall some of these problems. If $\lambda_1(\Omega)$ denotes the first eigenvalue of the Laplace operator in $\Omega$ with Dirichlet boundary condition, it is well-known that $\lambda_1(\Omega)\geq \lambda_1(\Omega^{\ast})$ and that the inequality is strict unless $\Omega$ is a ball (remember that $\Omega$ is always assumed to be in the class ${\mathcal C}$). Since $\lambda_1(\Omega^{\ast})$ can be explicitly computed, this result provides the classical Rayleigh-Faber-Krahn inequality, which states that
\begin{equation} \label{RFK}
\lambda_1(\Omega)\geq\lambda_1(\Omega^*)=\left(\frac 1{\left\vert 
\Omega\right\vert}\right)^{2/n}\alpha_n^{2/n}(j_{n/2-1,1})^2,
\end{equation}
where $j_{m,1}$ the first positive zero of the Bessel function 
$J_m$. Moreover, equality in 
(\ref{RFK}) is attained if and only if $\Omega$ is a ball. This result was first 
conjectured by Rayleigh for $n=2$ (\cite{Rayleigh} pp. 339-340), and proved independently by Faber (\cite{Faber}) and Krahn (\cite{Krahndim2}) for $n=2$, and by Krahn for all $n$ in \cite{Krahndimn} (see \cite{Krahndimnenglish} for the English translation). The proof of the inequality $\lambda_1(\Omega)\geq \lambda_1(\Omega^{\ast})$ is an immediate consequence of the following variational formula for $\lambda_1(\Omega)$:
\begin{equation} \label{variational}
\lambda_1(\Omega)=\min_{v\in H^1_0(\Omega)\setminus\left\{0\right\}} \frac{\displaystyle{\int_{\Omega} 
\left\vert \nabla v(x)\right\vert^2dx}}{\displaystyle{\int_{\Omega} \left\vert 
v(x)\right\vert^2dx}},
\end{equation}
and of the properties (\ref{Schwarzproperty}) of the Schwarz symmetrization.

Lots of optimization results involving other eigenvalues of the Laplacian (or more general elliptic symmetric operators of the form $-\hbox{div}(A\nabla)$) on $\Omega$ under Dirichlet boundary condition have also been established. For instance, the minimum of $\lambda_2(\Omega)$ (the second eigenvalue of the Laplace operator in $\Omega$ under Dirichlet boundary condition) among bounded open sets of $\R^n$ with given Lebesgue measure is achieved by the union of two identical balls (this result is attributed to Szeg\"o, see \cite{polya}). Very few things seem to be known about optimization problems for the other eigenvalues, see \cite{bucurhenrot,henrot,polya,polya2,wk}. Various optimization results are also known for functions of the eigenvalues. For instance, it is proved in \cite{ashben} that $\lambda_2(\Omega)/\lambda_1(\Omega)\leq\lambda_2(\Omega^{\ast})/\lambda_1(\Omega^{\ast})$, and the equality is attained if and only if $\Omega$ is a ball. The same result was also extended in \cite{ashben} to elliptic operators in divergence form with definite weight. We also refer to \cite{ab2,ab3,ahs,barnes,co,karaa,karaa2,ly,m,o,ppw,polya2} for further bounds or other optimization results for some eigenvalues or some functions of the eigenvalues in fixed or varying domains of $\R^n$ (or of manifolds).

Other boundary conditions may also be considered. For instance, if $\mu_2(\Omega)$ is the first non-trivial eigenvalue of $-\Delta$ under the Neumann boundary condition, then $\mu_2(\Omega)\leq\mu_2(\Omega^{\ast})$ and the equality is attained if, and only if, $\Omega$ is a ball (see \cite{sz} in dimension $n=2$, and \cite{wein} in any dimension). Bounds or optimization results for other eigenvalues of the Laplacian under Neumann boundary condition (\cite{polya2,sz,wein}, see also \cite{bandle1} for inhomogeneous problems), for Robin boundary condition (\cite{bossel}) or for the Stekloff eigenvalue problem (\cite{brock}) have also been established. We also mention another Rayleigh conjecture for the lowest eigenvalue of the clamped plate. If $\Omega\subset\R^2$, denote by $\Lambda_1(\Omega)$ the lowest eigenvalue of the operator $\Delta^2$, so that $\Delta^2u_1=\Lambda_1(\Omega)u_1$ in $\Omega$ with $u_1=\nu\cdot\nabla u_1=0$ on $\partial\Omega$, where $u_1$ denotes the principal eigenfunction and $\nu$ denotes the outward unit normal on $\partial\Omega$. The second author proved in \cite{nad} that $\Lambda_1(\Omega)\geq \Lambda_1(\Omega^{\ast})$ and that equality holds if and only if $\Omega$ is a ball, that is a disk in dimension $n=2$. The analogous result was also established in $\R^3$ in \cite{ashbenclamped}, while the problem is still open in higher dimensions. Much more complete surveys of all these topics can be found in \cite{bandle2,henrot,henrot2}.

It is important to observe that the variational formula (\ref{variational}) relies heavily on the fact that $-\Delta$ is symmetric on $L^2(\Omega)$. More generally, all the optimization problems considered hitherto concern symmetric operators, and their resolution relies on a ``Rayleigh'' quotient (that is, a variational formula similar to (\ref{variational})) and the Schwarz symmetrization. Before going further, let us recall that other rearrangement techniques than the Schwarz symmetrization can be found in the literature. For instance, even if this kind of problem is quite different from the ones we are interested in for the present paper, the Steiner symmetrization is the key tool to show that, among all triangles with fixed area, the principal eigenvalue of the Laplacian with Dirichlet boundary condition is minimal for the equilateral triangle (see \cite{polyaszego}). Steiner symmetrization is indeed relevant to take into account the polygonal geometry of the domain.

A natural question then arises: can inequalities on eigenvalues of non-symmetric operators be obtained~? In view of what we have just explained, such problems require different rearrangement techniques.

Actually, even for symmetric operators, some optimization problems cannot be solved by means of the Schwarz symmetrization, and other rearrangements have to be used. For instance, consider an operator $L=-\div(A\nabla)$ on a domain $\Omega$ under Dirichlet boundary condition. Assume that $A(x)\geq \Lambda(x)\hbox{Id}$ on $\Omega$ in the sense of quadratic forms (see below for precise definitions; $\hbox{Id}$ denotes the $n\times n$ identity matrix) for some positive function $\Lambda$, and that the $L^1$ norm of $\Lambda^{-1}$ is given. Then, what can be said about the infimum of the principal eigenvalue of $L$ under this constraint ? In particular, is this infimum greater than the corresponding one on $\Omega^{\ast}$, which is a natural conjecture in view of all the previous results ? Solving such a problem, which is one of our results in the present paper, does not seem to be possible by means of a variational formula for $\lambda_1$ (although the operator $L$ is symmetric in $L^2(\Omega)$) and the Schwarz or Steiner symmetrizations. 

More general constraints (given distribution functions; integral, pointwise or geometric constraints) on the coefficients $A$, $v$ and $V$ of non-symmetric operators $L$ of the type $L=-\div(A\nabla)+v\cdot\nabla+V$ under Dirichlet boundary condition will also be investigated. In general, the operator $L$ is non-symmetric, and there is no simple variational formulation of its first eigenvalue such as (\ref{variational}) --min-max formulations of the pointwise type (see \cite{bnv}) or of the integral type (see \cite{holland}) certainly hold, but they do not help in our context. 

The purpose of the present paper is twofold. First, we present a new rearrangement technique and we show some properties of the rearranged function. The inequalities we obtain between the function in $\Omega$ and its symmetrization in $\Omega^*$ are of independent interest. Then, we show how this technique can be used to cope with new comparisons between the principal eigenvalues of general non-symmetric elliptic operators of the type $-\hbox{div}(A\nabla)+v\cdot\nabla+V$ in $\Omega$ and of some symmetrized operators in $\Omega^*$. Actually, the comparisons we establish are new even when the operators are symmetric or one-dimensional.


\SE{Main results}\label{results}

Let us now give precise statements. We are interested in operators of the form
$$L=-\div(A\nabla)+v\cdot\nabla+V$$
in $\Omega\in{\mathcal{C}}$ under Dirichlet boundary condition.

Throughout the paper, we denote by ${\mathcal S}_n(\R)$ the set of $n\times n$ symmetric matrices with real entries. We always assume that $A:\Omega\rightarrow {\mathcal S}_n(\R)$ is in $W^{1,\infty}(\Omega)$. This assumption will be denoted by $A=(a_{i,j})_{1\le i,j\le n}\in W^{1,\infty}(\Omega,{\mathcal S}_n(\R))$: all the components $a_{i,j}$ are in $W^{1,\infty}(\Omega)$ and they can therefore be assumed to be continuous in $\overline{\Omega}$ up to a modification on a zero-measure set. We set
$$\|A\|_{W^{1,\infty}(\Omega,{\mathcal{S}}_n(\R))}=\max_{1\le i,j\le n}\|a_{i,j}\|_{W^{1,\infty}(\Omega)},$$
where
$$\|a_{i,j}\|_{W^{1,\infty}(\Omega)}=\|a_{i,j}\|_{L^{\infty}(\Omega)}+\sum_{1\le k\le n}\left\Vert\frac{\partial a_{i,j}}{\partial x_k}\right\Vert_{L^{\infty}(\Omega)}.$$
We always assume that $A$ is uniformly elliptic on $\overline{\Omega}$, which means that there exists $\delta>0$ such that, for all $x\in\overline{\Omega}$ and for all $\xi\in \R^n$,
$$A(x)\xi\cdot \xi\geq \delta \left\vert \xi\right\vert^2.$$
For $B=(b_{i,j})_{1\le i,j\le n}\in{\mathcal{S}}_n(\R)$, $\xi=(\xi_1,\ldots,\xi_n)\in\R^n$ and $\xi'=(\xi'_1,\ldots,\xi'_n)\in\R^n$, we denote $B\xi\cdot\xi'=\sum_{1\le i,j\le n}b_{i,j}\xi_j\xi'_i$. Actually, in some statements we compare the matrix field $A$ with a matrix field of the type $x\mapsto\Lambda(x)\hbox{Id}$. We call
$$L^{\infty}_+(\Omega)=\{\Lambda\in L^{\infty}(\Omega),\ \mathop{\hbox{ess inf}}_{\Omega}\Lambda>0\},$$
and, for $A\in W^{1,\infty}(\Omega,{\mathcal S}_n(\R))$ and $\Lambda\in L^{\infty}_+(\Omega)$, we say that $A\ge\Lambda\ \hbox{Id}$ almost everywhere (a.e.) in $\Omega$ if, for almost every $x\in\Omega$,
$$\forall\ \xi\in\R^n,\ A(x)\xi\cdot\xi\ge\Lambda(x)|\xi|^2.$$
For instance, if, for each $x\in\overline{\Omega}$, $\Lambda[A](x)$ denotes the smallest eigenvalue of the matrix $A(x)$, then $\Lambda[A]\in L^{\infty}_+(\Omega)$ and there holds $A(x)\ge\Lambda[A](x)\hbox{Id}$ (this inequality is actually satisfied for all $x\in\overline{\Omega}$).
 
We also always assume that the vector field $v$ is in $L^{\infty}(\Omega,\R^n)$ 
and that the potential $V$ is in $L^{\infty}(\Omega)$. In some statements, $V$ will be in the space $C(\overline{\Omega})$ of continuous functions on $\overline{\Omega}$.

Denote by $\lambda_1(\Omega,A,v,V)$ the principal eigenvalue of $L=-\div(A\nabla)+v\cdot\nabla+V$ with Dirichlet boundary condition on $\Omega$, and $\varphi_{\Omega,A,v,V}$ the corresponding (unique) nonnegative eigenfunction with $L^{\infty}$-norm equal to $1$. Recall that the following properties hold for $\varphi_{\Omega,A,v,V}$ (see \cite{bnv}):
\begin{equation} \label{eq}
\left\{
\begin{array}{l}
-\div\left(A\nabla\varphi_{\Omega,A,v,V}\right)+v\cdot\nabla\varphi_{\Omega,A,v,V}+V\varphi_{\Omega,A,v,V}=\lambda_1(\Omega,A,v,V)\varphi_{\Omega,A,v,V}\mbox{ 
in }\Omega,\\
\\
\varphi_{\Omega,A,v,V}>0\mbox{ in }\Omega,\ \varphi_{\Omega,A,v,V}=0\mbox{ on }\partial\Omega,\ 
\left\Vert \varphi_{\Omega,A,v,V}\right\Vert_{L^{\infty}(\Omega)}=1.
\end{array}
\right.
\end{equation}
and $\varphi_{\Omega,A,v,V}\in W^{2,p}(\Omega)$ for all $1\leq p<+\infty$ by standard elliptic estimates, whence $\varphi_{\Omega,A,v,V}\in C^{1,\alpha}(\overline{\Omega})$ for all $0\leq\alpha<1$. Recall also that $\lambda_1(\Omega,A,v,V)>0$ if and only if the operator $L$ satisfies the maximum principle in $\Omega$, and that the inequality
$$\lambda_1(\Omega,A,v,V)>\displaystyle{\mathop{\hbox{ess inf}}_{\Omega}}\ V$$
always holds (see \cite{bnv} for details and further results).

We are interested in optimization problems for $\lambda_1(\Omega,A,v,V)$ when $\Omega,A,v$ and $V$ vary and satisfy some constraints. Our goal is to compare $\lambda_1(\Omega,A,v,V)$ with the principal eigenvalue $\lambda_1(\Omega^*,A^*,v^*,V^*)$ for some fields $A^*$, $v^*$ and $V^*$ which are defined in the ball $\Omega^*$ and satisfy the same constraints as $A$, $v$ and $V$. The constraints may be of different types: integral type, $L^{\infty}$ type, given distribution function of $V^-$, or bounds on the determinant of $A$ and on another symmetric function of the eigenvalues of $A$. Throughout the paper, we denote
$$s^-=\max(-s,0)\hbox{ and }s^+=\max(s,0)\hbox{ for all }s\in\R.$$


\subsection{Constraints on the distribution function of $V^-$ and on some integrals involving $\Lambda$ and $v$}

We fix here the $L^1$ norms of $\Lambda^{-1}$ and $|v|^2\Lambda^{-1}$, some $L^{\infty}$ bounds on $\Lambda$ and $v$, as well as the distribution function of the negative part of $V$, under the condition that $\lambda_1(\Omega,A,v,V)\ge 0$. Then we can associate some fields $A^*$, $v^*$ and $V^*$ satisfying the same constraints in $\Omega^*$, and for which the principal eigenvalue is not too much larger that $\lambda_1(\Omega,A,v,V)$, with the extra property that $A^*$, $|v^*|$ and $V^*$ are smooth and radially symmetric.

\begin{theo} \label{th1}
Let $\Omega\in {\mathcal C}$, $A\in W^{1,\infty}(\Omega,{\mathcal S}_n(\R))$, $\Lambda\in L^{\infty}_+(\Omega)$, $v\in L^{\infty}(\Omega,\R^n)$ and $V\in C(\overline{\Omega})$. Assume that $A\geq\Lambda\Id$ a.e. in $\Omega$, and that $\lambda_1(\Omega,A,v,V)\ge 0$. Then, for all $\varepsilon>0$, there exist three radially symmetric $C^{\infty}(\overline{\Omega^*})$ fields $\Lambda^*>0$, $\omega^*\ge 0$ and $\overline{V}^*\le 0$ such that, for $v^*=\omega^*e_r$ in $\overline{\Omega^*}\backslash\{0\}$,
\be\label{bounds1}\left\{\baa{l}
\displaystyle{\mathop{\rm{ess}\ \rm{inf}}_{\Omega}}\ \Lambda \le \displaystyle{\mathop{\min}_{\overline{\Omega^*}}}\ \Lambda^*\le \displaystyle{\mathop{\max}_{\overline{\Omega^*}}}\ \Lambda^*\le \displaystyle{\mathop{\rm{ess}\ \rm{sup}}_{\Omega}}\ \Lambda,\ \|(\Lambda^*)^{-1}\|_{L^1(\Omega^*)}=\|\Lambda^{-1}\|_{L^1(\Omega)},\\
\|v^*\|_{L^{\infty}(\Omega^*,\R^n)}\le\|v\|_{L^{\infty}(\Omega,\R^n)},\ \|\ |v^*|^2(\Lambda^*)^{-1}\|_{L^1(\Omega^*)}=\|\ |v|^2\Lambda^{-1}\|_{L^1(\Omega)},\\
\mu_{|\overline{V}^*|}=\mu_{(\overline{V}^*)^-}\le\mu_{V^-},\eaa\right.
\ee
and
\be\label{lambdacomparaison}
\lambda_1(\Omega^*,\Lambda^*{\rm{Id}},v^*,\overline{V}^*)\le\lambda_1(\Omega,A,v,V)+\epsilon.
\ee
There also exists a nonpositive radially symmetric $L^{\infty}(\Omega^*)$ field $V^*$ such that $\mu_{V^*}=\mu_{-V^-}$, $V^*\le\overline{V}^*\le 0$ in $\Omega^*$ and $\lambda_1(\Omega^*,\Lambda^*{\rm{Id}},v^*,V^*)\le\lambda_1(\Omega^*,\Lambda^*{\rm{Id}},v^*,\overline{V}^*)\le\lambda_1(\Omega,A,v,V)+\epsilon$.\par
If one further assumes that $\Lambda$ is equal to a constant $\gamma>0$ in $\Omega$, then there exist two radially symmetric bounded functions $\omega^*_0\ge 0$ and $V^*_0\le 0$ in $\Omega^*$ such that, for $v^*_0=\omega^*_0e_r$,
\be\label{bounds1bis}\left\{\baa{l}
\|v^*_0\|_{L^{\infty}(\Omega^*,\R^n)}\le\|v\|_{L^{\infty}(\Omega,\R^n)},\ \|v^*_0\|_{L^2(\Omega^*,\R^n)}\le\|v\|_{L^2(\Omega,\R^n)},\\
-\displaystyle{\mathop{\max}_{\overline{\Omega}}}\ V^-\le V^*_0\le 0\hbox{ a.e. in }\Omega^*,\ \|V^*_0\|_{L^p(\Omega^*)}\le\|V^-\|_{L^p(\Omega)}\hbox{ for all }1\le p\le+\infty,\eaa\right.
\ee
and
\be\label{ineqv0}
\lambda_1(\Omega^*,\gamma{\rm{Id}},v^*_0,V^*_0)\le
\lambda_1(\Omega,A,v,V).
\ee
\end{theo}

Remember (see \cite{bnv}) that the inequality $\lambda_1(\Omega,A,v,V)\ge\lambda_1(\Omega,A,v,-V^-)$ always holds. This is the reason why, in order to decrease $\lambda_1(\Omega,A,v,V)$, the rearranged potentials had better be nonpositive in $\Omega^*$, and only the negative part of $V$ plays a role. Notice that the quantities such as the integral of $\Lambda^{-1}$, which are preserved here after symmetrization, also appear in other contexts, like in homogenization of elliptic or parabolic equations.

In the case when $\Lambda$ is a constant, then the number $\epsilon$ can be dropped in (\ref{lambdacomparaison}). The price to pay is that the new fields in $\Omega^*$ may not be smooth anymore and the distribution function of the new potential $V^*_0$ in $\Omega^*$ is no longer equal to that of $-V^-$.

However, in the general case, neither $\Lambda$ is constant in $\Omega$ nor $\Lambda^*$ is constant in $\Omega^*$ (see Remark~\ref{remLambdanoncst} for details). For instance, as already underlined, an admissible $\Lambda$ is the continuous positive function $\Lambda[A]$, which is not constant in general. Actually, even in the case of operators $L$ which are written in a self-adjoint form (that is, with $v=0$), the comparison result stated in Theorem~\ref{th1} is new.

An optimization result follows immediately from Theorem~\ref{th1}. To state it, we need a few notations. Given
\be\label{constraints}
m>0,\ \overline{M}_{\Lambda}\ge\underline{m}_{\Lambda}>0,\ \alpha\in\left[\frac{m}{\overline{M}_{\Lambda}},\frac{m}{\underline{m}_{\Lambda}}\right],\ \overline{M}_v\ge 0,\ \tau\in\left[0,\alpha\overline{M}_v^2\right],\ \overline{M}_V\ge 0
\ee
and
$$\baa{rcl}
\mu\ \in\ {\mathcal{F}}_{0,\overline{M}_V}(m) & := & \{\rho:\R\to[0,m],\ \rho\hbox{ is right-continuous, non-increasing},\\
& & \quad\rho=m\hbox{ on }(-\infty,0),\ \rho=0\hbox{ on }[\overline{M}_V,+\infty)\},\eaa$$
we set, for all open set $\Omega\in{\mathcal{C}}$ such that $|\Omega|=m$,
$$\baa{rcl}
\displaystyle{\mathcal G}_{\overline{M}_{\Lambda},\underline{m}_{\Lambda},\alpha,\overline{M}_v,\tau,\overline{M}_V,\mu}(\Omega) & = & \displaystyle\left\{(A,v,V)\in W^{1,\infty}(\Omega,{\mathcal S}_n(\R))\times L^{\infty}(\Omega,\R^n)\times C(\overline{\Omega});\right.\\
& & \quad\exists\ \Lambda\in L^{\infty}_+(\Omega),\ A\geq \Lambda \Id\mbox{ a.e. in }\Omega,\\ & & \quad\displaystyle\underline{m}_{\Lambda}\le\displaystyle{\mathop{\hbox{ess inf}}_{\Omega}}\ \Lambda\le\displaystyle{\mathop{\hbox{ess sup}}_{\Omega}}\ \Lambda\le\overline{M}_{\Lambda},\ \left\Vert \Lambda^{-1}\right\Vert_{L^1(\Omega)}=\alpha,\\
& & \quad\left.\|v\|_{L^{\infty}(\Omega,\R^n)}\le\overline{M}_v,\ \left\Vert |v|^2\Lambda^{-1}\right\Vert_{L^1(\Omega)}=\tau\hbox{ and }\mu_{V^-}\le\mu\right\}\eaa$$
and
\be\label{defminlambda}
\underline{\lambda}_{\overline{M}_{\Lambda},\underline{m}_{\Lambda},\alpha,\overline{M}_v,\tau,\overline{M}_V,\mu}(\Omega)=\inf_{(A,v,V)\in{\mathcal{G}}_{\overline{M}_{\Lambda},\underline{m}_{\Lambda},\alpha,\overline{M}_v,\tau,\overline{M}_V,\mu}(\Omega)}\lambda_1(\Omega,A,v,V).
\ee
Notice that, given $\mu\in{\mathcal F}_{0,\overline{M}_V}(m)$ and $\Omega\in{\mathcal{C}}$ such that $|\Omega|=m$, there exists $V\in C(\overline{\Omega})$ such that $\mu_{V^-}\le\mu$ (for instance, $V=0$ is admissible; furthermore, there is $V\in L^{\infty}(\Omega)$ such that $\mu_{V^-}=\mu$, see Appendix~\ref{distribution}), and, necessarily, $V\ge-\overline{M}_V$ in $\overline{\Omega}$. It is immediate to see that, under the conditions (\ref{constraints}), each set ${\mathcal{G}}_{\overline{M}_{\Lambda},\underline{m}_{\Lambda},\alpha,\overline{M}_v,\tau,\overline{M}_V,\mu}(\Omega)$ is not empty.

\begin{cor}\label{cor1} Let $m$, $\overline{M}_{\Lambda}$, $\underline{m}_{\Lambda}$, $\alpha$, $\overline{M}_v$, $\tau$, $\overline{M}_V$ be as in $(\ref{constraints})$,  $\mu\in{\mathcal{F}}_{0,\overline{M}_V}(m)$ and $\Omega^*$ be the Euclidean ball centered at the origin such that $|\Omega^*|=m$. If $\underline{\lambda}_{\overline{M}_{\Lambda},\underline{m}_{\Lambda},\alpha,\overline{M}_v,\tau,\overline{M}_V,\mu}(\Omega)\ge 0$ for all $\Omega\in{\mathcal{C}}$ such that $|\Omega|=m$, then
\be\label{minlambda}
\min_{\Omega\in{\mathcal{C}},\ |\Omega|=m}\underline{\lambda}_{\overline{M}_{\Lambda},\underline{m}_{\Lambda},\alpha,\overline{M}_v,\tau,\overline{M}_V,\mu}(\Omega)=\underline{\lambda}_{\overline{M}_{\Lambda},\underline{m}_{\Lambda},\alpha,\overline{M}_v,\tau,\overline{M}_V,\mu}(\Omega^*).
\ee
Furthermore, in the definition of $\underline{\lambda}_{\overline{M}_{\Lambda},\underline{m}_{\Lambda},\alpha,\overline{M}_v,\tau,\overline{M}_V,\mu}(\Omega^*)$ in (\ref{defminlambda}), the data $A$, $v$ and $V$ can be assumed to be such that $A=\Lambda\ {\rm{Id}}$, $v=\omega e_r=|v|e_r$ and $V\le 0$ in $\Omega^*$, where $\Lambda$, $\omega$ and $V$ are $C^{\infty}(\overline{\Omega^*})$ and radially symmetric.
\end{cor}

Let us now discuss about the non-negativity condition $\lambda_1(\Omega,A,v,V)\ge 0$ in Theorem~\ref{th1}, as well as that of Corollary~\ref{cor1}. We recall (see \cite{bnv}) that
$$\lambda_1(\Omega,A,v,V)>\min_{\overline{\Omega}}V.$$
Therefore, the condition $\lambda_1(\Omega,A,v,V)\ge 0$ is satisfied in particular if $V\ge 0$ in $\overline{\Omega}$, and the condition $\underline{\lambda}_{\overline{M}_{\Lambda},\underline{m}_{\Lambda},\alpha,\overline{M}_v,\tau,\overline{M}_V,\mu}(\Omega)\ge 0$ in Corollary~\ref{cor1} is satisfied if $\overline{M}_V=0$. Another more complex condition which also involves $A$ and $v$ can be derived. To do so, assume $A\ge\Lambda\Id$ a.e. in $\Omega$ with $m_{\Lambda}:=\hbox{ess inf}_{\Omega}\ \Lambda>0$, and call $M_v=\|v\|_{\infty}$ and $m_V=\min_{\overline{\Omega}}V$. Multiply by $\varphi=\varphi_{\Omega,A,v,V}$ the equation (\ref{eq}) and integrate by parts over $\Omega$. It follows that, for all $\beta\in(0,1]$,
$$\baa{rcl}
\lambda_1(\Omega,A,v,V)\displaystyle{\int_{\Omega}}\varphi^2 & \ge & \displaystyle{\int_{\Omega}}\Lambda|\nabla\varphi|^2-\displaystyle{\int_{\Omega}}|v|\ |\nabla\varphi|\ \varphi+m_V\displaystyle{\int_{\Omega}}\varphi^2\\
& \ge & (1-\beta)\displaystyle{\int_{\Omega}}\Lambda|\nabla\varphi|^2+m_V\displaystyle{\int_{\Omega}}\varphi^2-\displaystyle{\frac{1}{4\beta}}\displaystyle{\int_{\Omega}}|v|^2\Lambda^{-1}\varphi^2\\
& \ge & [(1-\beta)m_{\Lambda}\lambda_1(\Omega)+m_V-(4\beta m_{\Lambda})^{-1}M_v^2]\displaystyle{\int_{\Omega}}\varphi^2,\eaa$$
where $\lambda_1(\Omega)=\lambda_1(\Omega,\hbox{Id},0,0)=\min_{\phi\in H^1_0(\Omega),\ \|\phi\|_2=1}\int_{\Omega}|\nabla\phi|^2$. If $M_v>0$ and $m_{\Lambda}\sqrt{\lambda_1(\Omega)}\ge M_v$, then the value $\beta=M_v/(2m_{\Lambda}\sqrt{\lambda_1(\Omega)})\in(0,1]$ gives the best inequality, that is $\lambda_1(\Omega,A,v,V)\ge m_V+\sqrt{\lambda_1(\Omega)}(m_{\Lambda}\sqrt{\lambda_1(\Omega)}-M_v)$. The same inequality also holds from the previous calculations if $M_v=0$. Therefore, the following inequality always holds:
$$\lambda_1(\Omega,A,v,V)\ge m_V+\sqrt{\lambda_1(\Omega)}\times\max(0,m_{\Lambda}\sqrt{\lambda_1(\Omega)}-M_v).$$
As a consequence, under the notations of Corollary~\ref{cor1}, it follows from (\ref{RFK}) that
$$\underline{\lambda}_{\overline{M}_{\Lambda},\underline{m}_{\Lambda},\alpha,\overline{M}_v,\tau,\overline{M}_V,\mu}(\Omega)\ge-\overline{M}_V+m^{-1/n}\alpha_n^{1/n}j_{n/2-1,1}\times\max(0,\underline{m}_{\Lambda}m^{-1/n}\alpha_n^{1/n}j_{n/2-1,1}-\overline{M}_v)$$
for all $\Omega\in{\mathcal{C}}$ such that $|\Omega|=m$. The conclusion of Corollary~\ref{cor1} is then true if the right-hand side of the above inequality is nonnegative. In particular, for given $n\in\N\backslash\{0\}$, $\underline{m}_{\Lambda}>0$, $\overline{M}_v\ge 0$ and $\overline{M}_V\ge 0$, this holds if $m>0$ is small enough.

To complete this section, we now give a more precise version of Theorem~\ref{th1} when $\Omega$ is not a ball.

\begin{theo}\label{th1bis} 
Under the notation of Theorem~\ref{th1}, assume that $\Omega\in{\mathcal{C}}$ is not a ball and let $\overline{M}_A>0$, $\underline{m}_{\Lambda}>0$, $\overline{M}_v\ge 0$ and $\overline{M}_V\ge 0$ be such that
\be\label{bornes}
\|A\|_{W^{1,\infty}(\Omega,{\mathcal{S}}_n(\R))}\le\overline{M}_A,\ \mathop{{\rm{ess}}\ {\rm{inf}}}_{\Omega}\Lambda\ge\underline{m}_{\Lambda},\ \|v\|_{L^{\infty}(\Omega,\R^n)}\le\overline{M}_v\hbox{ and }\|V\|_{L^{\infty}(\Omega,\R)}\le\overline{M}_V.
\ee
Then there exists a positive constant $\theta=\theta(\Omega,n,\overline{M}_A,\underline{m}_{\Lambda},\overline{M}_v,\overline{M}_V)>0$ depending only on $\Omega$, $n$, $\overline{M}_A$, $\underline{m}_{\Lambda}$, $\overline{M}_v$ and $\overline{M}_V$, such that if $\lambda_1(\Omega,A,v,V)>0$, then there exist three radially symmetric $C^{\infty}(\overline{\Omega^*})$ fields $\Lambda^*>0$, $\omega^*\ge 0$,  $\overline{V}^*\le 0$ and a nonpositive radially symmetric $L^{\infty}(\Omega^*)$ field $V^*$, which satisfy (\ref{bounds1}), $\mu_{V^*}=\mu_{-V^-}$, $V^*\le\overline{V}^*\le 0$ and are such that
$$\lambda_1(\Omega^*,\Lambda^*{\rm{Id}},v^*,V^*)\le\lambda_1(\Omega^*,\Lambda^*{\rm{Id}},v^*,\overline{V}^*)\le\frac{\lambda_1(\Omega,A,v,V)}{1+\theta},$$
where $v^*=\omega^*e_r$ in $\overline{\Omega^*}\backslash\{0\}$.
\end{theo}

Notice that the assumption $A\ge\Lambda\Id$ a.e. in $\Omega$ and the bounds (\ref{bornes}) imply necessarily that $\overline{M}_A\ge\underline{m}_{\Lambda}$.


\subsection{Constraints on the determinant and another symmetric function of the eigenvalues of $A$}

For our second type of comparison result, we keep the same constraints on $v$ and $V$ as in Theorem \ref{th1} but we modify the one on $A$: we now prescribe some conditions on the determinant and another symmetric function of the eigenvalues of $A$. {\bf{We assume in this subsection that $n\ge 2$}}. If $A\in{\mathcal{S}}_n(\R)$, if $p\in\{1,\ldots,n-1\}$ and if $\lambda_1[A]\le\cdots\le\lambda_n[A]$ denote the eigenvalues of $A$, then we call
$$\sigma_p(A)=\sum_{I\subset\{1,\ldots,n\},\ \hbox{card}(I)=p}\left(\prod_{i\in I}\lambda_i[A]\right).$$
Throughout the paper, the notation $\hbox{card}(I)$ means the cardinal of a finite set $I$. If $A$ is nonnegative, it follows from the arithmetico-geometrical inequality that  $\hbox{C}^p_n\times(\hbox{det}(A))^{p/n}\le\sigma_p(A)$, where $\hbox{C}^p_n$ is the binomial coefficient $\hbox{C}^p_n=n!/(p!\times(n-p)!)$. 

Our third result is as follows:

\begin{theo} \label{th2}
Assume $n\ge 2$. Let $\Omega\in {\mathcal C}$, $A\in W^{1,\infty}(\Omega,{\mathcal S}_n(\R))$, $v\in L^{\infty}(\Omega,\R^n)$, $V\in C(\overline{\Omega})$ and let $p\in\{1,\ldots,n-1\}$, $\omega>0$ and $\sigma>0$ be given. Assume that $A\ge\gamma{\rm{Id}}$ in $\overline{\Omega}$ for some constant $\gamma>0$, that
\be\label{detsigma}
{\rm{det}}(A(x))\ge\omega,\ \sigma_p(A(x))\le\sigma\hbox{ for all }x\in\overline{\Omega},
\ee
and that $\lambda_1(\Omega,A,v,V)\ge 0$. Then, there are two positive numbers $0<a_1\le a_2$ which only depend on $n$, $p$, $\omega$ and $\sigma$, such that, for all $\varepsilon>0$, there exist a matrix field $A^*\in C^{\infty}(\overline{\Omega^*}\backslash\{0\},{\mathcal{S}}_n(\R))$, two radially symmetric $C^{\infty}(\overline{\Omega^*})$ fields $\omega^*\ge 0$ and $\overline{V}^*\le 0$, and a nonpo\-si\-tive radially symmetric $L^{\infty}(\Omega^*)$ field $V^*$, such that, for $v^*=\omega^*e_r$ in $\overline{\Omega^*}\backslash\{0\}$,
\be\label{bounds2}\left\{\baa{l}
A\ge a_1{\rm{Id}}\hbox{ in }\Omega,\ A^*\ge a_1{\rm{Id}}\hbox{ in }\Omega^*,\\
{\rm{det}}(A^*(x))=\omega,\ \sigma_p(A^*(x))=\sigma\hbox{ for all }x\in\overline{\Omega^*}\backslash\{0\},\\
\|v^*\|_{L^{\infty}(\Omega^*,\R^n)}\le\|v\|_{L^{\infty}(\Omega,\R^n)},\ \|v^*\|_{L^2(\Omega^*,\R^n)}=\|v\|_{L^2(\Omega,\R^n)},\\
\mu_{|\overline{V}^*|}\le\mu_{V^-},\ \mu_{V^*}=\mu_{-V^-},\ V^*\le\overline{V}^*\le 0\hbox{ in }\Omega^*\eaa\right.
\ee
and
$$\lambda_1(\Omega^*,A^*,v^*,V^*)\le\lambda_1(\Omega^*,A^*,v^*,\overline{V}^*)\le\lambda_1(\Omega,A,v,V)+\epsilon.$$
Furthermore, the matrix field $A^*$ is defined, for all $x\in\overline{\Omega^*}\backslash\{0\}$, by:
$$A^*(x)x\cdot x=a_1|x|^2\hbox{ and }A^*(x)y\cdot y=a_2|y|^2\hbox{ for all }y\perp x.$$
Lastly, there exist two radially symmetric bounded functions $\omega^*_0\ge 0$ and $V^*_0\le 0$ in $\Omega^*$ satisfying (\ref{bounds1bis}) and $\lambda_1(\Omega^*,A^*,v^*_0,V^*_0)\le\lambda_1(\Omega,A,v,V)$, where $v^*_0=\omega^*_0e_r$ in $\Omega^*$.
\end{theo}

\begin{rem}\label{remdet}{\rm 
Notice that the assumptions of Theorem~\ref{th2} imply necessarily that $\hbox{C}^p_n\omega^{p/n}\le\sigma$. Actually, the matrix field $A^*$ cannot be extended by continuity at $0$, unless $a_1=a_2$, namely $\hbox{C}^p_n\omega^{p/n}=\sigma$. As a consequence, $A^*$ is not in $W^{1,\infty}(\Omega^*,{\mathcal{S}}_n(\R))$ if $\hbox{C}^p_n\omega^{p/n}\neq\sigma$, but we can still define $\lambda_1(\Omega^*,A^*,v^*,V^*)$. Indeed, for $\widetilde{A}^*=a_1\hbox{Id}$ in $\overline{\Omega^*}$, the principal eigenfunction $\overline{\varphi}^*$ (resp. $\varphi^*$) of the operator $-\hbox{div}(\widetilde{A}^*\nabla)+v^*\cdot\nabla+\overline{V}^*$ (resp. $-\hbox{div}(\widetilde{A}^*\nabla)+v^*\cdot\nabla+V^*$) is radially symmetric and belongs to all $W^{2,p}(\Omega^*)$ spaces for all $1\le p<+\infty$. Hence,
$$A^*\nabla\overline{\varphi}^*=\widetilde{A}^*\nabla\overline{\varphi}^*=a_1\nabla\overline{\varphi}^*$$
(resp. $A^*\nabla\varphi^*=\widetilde{A}^*\nabla\varphi^*=a_1\nabla\varphi^*$). With a slight abuse of notation, we say that $\overline{\varphi}^*$ (resp. $\varphi^*$) is the principal eigenfunction of $-\hbox{div}(A^*\nabla)+v^*\cdot\nabla+\overline{V}^*$ (resp. $-\hbox{div}(A^*\nabla)+v^*\cdot\nabla+V^*$) and we call
$$\lambda_1(\Omega^*,A^*,v^*,\overline{V}^*)=\lambda_1(\Omega^*,\widetilde{A}^*,v^*,\overline{V}^*)$$
(resp. $\lambda_1(\Omega^*,A^*,v^*,V^*)=\lambda_1(\Omega^*,\widetilde{A}^*,v^*,V^*)$).}
\end{rem}

An interpretation of the conditions (\ref{detsigma}) is that they provide some bounds for the local deformations induced by the matrices $A(x)$, uniformly with respect to $x\in\overline{\Omega}$. Notice that these constraints are saturated for the matrix field $A^*$ in the ball $\Omega^*$.

As for Theorem~\ref{th1}, an optimization result follows immediately from Theorem~\ref{th2}.

\begin{cor} Assume $n\ge 2$. Given $m>0$, $p\in\{1,\ldots,n-1\}$, $\omega>0$, $\sigma\ge{\rm{C}}_n^p\omega^{p/n}$, $\overline{M}_v\ge 0$, $\tau\in\left[0,\sqrt{m}\times\overline{M}_v\right]$, $\overline{M}_V\ge 0$ and $\mu\in{\mathcal{F}}_{0,\overline{M}_V}(m)$, we set, for all $\Omega\in{\mathcal{C}}$ such that $|\Omega|=m$,
$$\baa{rcl}
\displaystyle{\mathcal G}'_{p,\omega,\sigma,\overline{M}_v,\tau,\overline{M}_V,\mu}(\Omega) & = & \displaystyle\left\{(A,v,V)\in W^{1,\infty}(\Omega,{\mathcal S}_n(\R))\times L^{\infty}(\Omega,\R^n)\times C(\overline{\Omega});\right.\\
& & \quad\exists\ \gamma>0,\ A(x)\ge\gamma{\rm{Id}}\hbox{ for all }x\in\overline{\Omega},\\
& & \quad{\rm{det}}(A(x))\ge\omega,\ \sigma_p(A(x))\le\sigma\hbox{ for all }x\in\overline{\Omega},\\
& & \quad\left.\|v\|_{L^{\infty}(\Omega,\R^n)}\le\overline{M}_v,\ \left\Vert v\right\Vert_{L^2(\Omega,\R^n)}=\tau\hbox{ and }\mu_{V^-}\le\mu\right\}\eaa$$
and
$$\underline{\lambda}'_{p,\omega,\sigma,\overline{M}_v,\tau,\overline{M}_V,\mu}(\Omega)=\inf_{(A,v,V)\in{\mathcal{G}}'_{p,\omega,\sigma,\overline{M}_v,\tau,\overline{M}_V,\mu}(\Omega)}\lambda_1(\Omega,A,v,V).$$
If $\underline{\lambda}'_{p,\omega,\sigma,\overline{M}_v,\tau,\overline{M}_V,\mu}(\Omega)\ge 0$ for all $\Omega\in{\mathcal{C}}$ such that $|\Omega|=m$, then
$$\inf_{\Omega\in{\mathcal{C}},\ |\Omega|=m}\underline{\lambda}'_{p,\omega,\sigma,\overline{M}_v,\tau,\overline{M}_V,\mu}(\Omega)=\inf_{(v^*,V^*)\in{\mathcal{G}}^*_{\overline{M}_v,\tau,\overline{M}_V,\mu}}\lambda_1(\Omega^*,A^*,v^*,V^*),$$
where $\Omega^*$ is the ball centered at the origin such that $|\Omega^*|=m$, $A^*$  is given as in Theorem~\ref{th2} and
$$\baa{rcl}
{\mathcal{G}}^*_{\overline{M}_v,\tau,\overline{M}_V,\mu} & = & \left\{(v^*,V^*)\in L^{\infty}(\Omega^*,\R^n)\times C(\overline{\Omega}),\ v^*=|v^*|e_r,\ V^*\le 0,\right.\\
& & \quad |v^*|\hbox{ and }V^*\hbox{ are radially symmetric and }C^{\infty}(\overline{\Omega^*}),\\
& & \quad\left.\|v^*\|_{L^{\infty}(\Omega,\R^n)}\le\overline{M}_v,\ \left\Vert v^*\right\Vert_{L^2(\Omega,\R^n)}=\tau\hbox{ and }\mu_{(V^*)^-}\le\mu\right\}.\eaa$$
\end{cor}

Notice also that a sufficient condition for $\underline{\lambda}'_{p,\omega,\sigma,\overline{M}_v,\tau,\overline{M}_V,\mu}(\Omega)$ to be nonnegative for all $\Omega\in{\mathcal{C}}$ such that $|\Omega|=m$ is:
$$-\overline{M}_V+m^{-1/n}\alpha_n^{1/n}j_{n/2-1,1}\times\max(0,a_1m^{-1/n}\alpha_n^{1/n}j_{n/2-1,1}-\overline{M}_v)\ge 0,$$
where $a_1>0$ is the same as in Theorem~\ref{th2} and only depends on $n$, $p$, $\omega$ and $\sigma$ (see Lemma~\ref{algebra} for its definition). When $n$, $p$, $\omega$, $\sigma$, $\overline{M}_v$ and $\overline{M}_V$ are given, the above inequality is satisfied in particular if $m>0$ is small enough.

When $\Omega\in{\mathcal{C}}$ is not a ball, we can make Theorem~\ref{th2} more precise: under the same notations as in Theorem~\ref{th2}, if $\overline{M}_A>0$, $\overline{M}_v\ge 0$ and $\overline{M}_V\ge 0$ are such that $\|A\|_{W^{1,\infty}(\Omega,{\mathcal{S}}_n(\R))}\le\overline{M}_A$, $\|v\|_{L^{\infty}(\Omega,\R^n)}\le\overline{M}_v$ and $\|V\|_{L^{\infty}(\Omega,\R)}\le\overline{M}_V$, then there exists a positive constant
$$\theta'=\theta'(\Omega,n,p,\omega,\sigma,\overline{M}_A,\overline{M}_v,\overline{M}_V)>0$$
depending only on $\Omega$, $n$, $p$, $\omega$, $\sigma$, $\overline{M}_A$, $\overline{M}_v$ and $\overline{M}_V$, such that if $\lambda_1(\Omega,A,v,V)>0$, then there exist a matrix field $A^*\in C^{\infty}(\overline{\Omega^*}\backslash\{0\},{\mathcal{S}}_n(\R))$ (the same as in Theorem~\ref{th2}), two radially symmetric $C^{\infty}(\overline{\Omega^*})$ fields $\omega^*\ge 0$,  $\overline{V}^*\le 0$ and a nonpositive radially symmetric $L^{\infty}(\Omega^*)$ field $V^*$, which satisfy (\ref{bounds2}), $\mu_{V^*}=\mu_{-V^-}$, $V^*\le\overline{V}^*\le 0$ and are such that
$$\lambda_1(\Omega^*,A^*,v^*,V^*)\le\lambda_1(\Omega^*,A^*,v^*,\overline{V}^*)\le\frac{\lambda_1(\Omega,A,v,V)}{1+\theta'},$$
where $v^*=\omega^*e_r$ in $\overline{\Omega^*}\backslash\{0\}$. It is immediate to see that this fact is a consequence of Theorems~\ref{th1bis} and~\ref{th2} (notice in particular that the eigenvalues of $A(x)$ are between two positive constants which only depend on $n$, $p$, $\omega$ and $\sigma$).


\subsection{Faber-Krahn inequalities for non-symmetric operators}

An immediate corollary of Theorem~\ref{th1} is an optimization result, slightly different from Corollary \ref{cor1}, where the constraint over the potential $V$ is stated in terms of $L^p$ norms. Namely, given $m>0$, $\overline{M}_{\Lambda}\geq \underline{m}_{\Lambda}>0$, $\alpha\in\left[\frac{m}{\overline{M}_{\Lambda}},\frac{m}{\underline{m}_{\Lambda}}\right]$, $\overline{M}_v\geq 0$, $\tau\in\left[0,\alpha\overline{M}_v^2\right]$, $\tau_V\geq 0$, $1\le p\le+\infty$ and $\Omega\in {\mathcal C}$ such that $|\Omega|=m$, set
$$\baa{rcl}
\displaystyle{\mathcal H}_{\overline{M}_{\Lambda},\underline{m}_{\Lambda},\alpha,\overline{M}_v,\tau,\tau_V,p}(\Omega) & = & \displaystyle\left\{(A,v,V)\in W^{1,\infty}(\Omega,{\mathcal S}_n(\R))\times L^{\infty}(\Omega,\R^n)\times C(\overline{\Omega});\right.\\
& & \quad\exists\ \Lambda\in L^{\infty}_+(\Omega),\ A\geq \Lambda \Id\mbox{ a.e. in }\Omega,\\ & & \quad\displaystyle\underline{m}_{\Lambda}\le\displaystyle{\mathop{\hbox{ess inf}}_{\Omega}}\ \Lambda\le\displaystyle{\mathop{\hbox{ess sup}}_{\Omega}}\ \Lambda\le\overline{M}_{\Lambda},\ \left\Vert \Lambda^{-1}\right\Vert_{L^1(\Omega)}=\alpha,\\
& & \quad\left.\|v\|_{L^{\infty}(\Omega,\R^n)}\le\overline{M}_v,\ \left\Vert |v|^2\Lambda^{-1}\right\Vert_{L^1(\Omega)}=\tau\hbox{ and }\left\Vert V^{-}\right\Vert_{L^p(\Omega)}\leq \tau_V\right\}\eaa$$
and
$$\underline{\underline{\lambda}}_{\overline{M}_{\Lambda},\underline{m}_{\Lambda},\alpha,\overline{M}_v,\tau,\tau_V,p}(\Omega)=\inf_{(A,v,V)\in{\mathcal{H}_{\overline{M}_{\Lambda},\underline{m}_{\Lambda},\alpha,\overline{M}_v,\tau,\tau_V,p}(\Omega)}}\lambda_1(\Omega,A,v,V).$$
Since, in Theorem \ref{th1}, the $L^p$ norm of $\overline{V}^{\ast}$ is smaller than the one of $V^{-}$ (because the distribution functions of their absolute values are ordered this way), it follows from Theorem \ref{th1} that
\[
\min_{\Omega\in{\mathcal{C}},\ |\Omega|=m}\underline{\underline{\lambda}}_{\overline{M}_{\Lambda},\underline{m}_{\Lambda},\alpha,\overline{M}_v,\tau,\tau_V,p}(\Omega)=\underline{\underline{\lambda}}_{\overline{M}_{\Lambda},\underline{m}_{\Lambda},\alpha,\overline{M}_v,\tau,\tau_V,p}(\Omega^{\ast}),
\]
assuming that $\underline{\underline{\lambda}}_{\overline{M}_{\Lambda},\underline{m}_{\Lambda},\alpha,\overline{M}_v,\tau,\tau_V,p}(\Omega)\geq 0$ for all $\Omega\in {\mathcal C}$ such that $\left\vert \Omega\right\vert=m$. In other words, the infimum of $\lambda_1(\Omega,A,v,V)$ over all the previous constraints when $\Omega$ varies but still satisfies $\left\vert \Omega\right\vert=m$ is the same as the infimum in the ball $\Omega^{\ast}$. Observe that we do not know in general if this infimum is actually a minimum. However, specializing to the case of $L^{\infty}$ constraints for $v$ and $V$, we can solve a slightly different optimization problem and establish, as an application of Theorems~\ref{th2} and~\ref{fixedlinftyball} (see Section~\ref{sec5} below), a generalization of the classical Rayleigh-Faber-Krahn inequality for the principal eigenvalue of the Laplace operator.

\begin{theo} \label{faberkrahn}
Let $\Omega\in {\mathcal C}$, $\overline{M}_A>0$, $\underline{m}_{\Lambda}>0$, $\tau_1\ge 0$ and $\tau_2\ge 0$ be given. Assume that $\Omega$ is not a ball. Consider $A\in W^{1,\infty}(\Omega,{\mathcal S}_n(\R))$, $\Lambda\in L^{\infty}_+(\Omega)$, $v\in L^{\infty}(\Omega,\R^n)$ and $V\in L^{\infty}(\Omega)$ satisfying
$$\left\{\baa{l}
A\geq \Lambda\Id\mbox{ a.e. in }\Omega,\ \|A\|_{W^{1,\infty}(\Omega,{\mathcal{S}}_n(\R))}\le\overline{M}_A,\ \displaystyle{\mathop{{\rm{ess}}\ {\rm{inf}}}_{\Omega}}\ \Lambda\ge\underline{m}_{\Lambda},\\
\|v\|_{L^{\infty}(\Omega,\R^n)}\le\tau_1\hbox{ and }\|V\|_{L^{\infty}(\Omega)}\le\tau_2.\eaa\right.$$
Then there exists a positive constant $\eta=\eta(\Omega,n,\overline{M}_A,\underline{m}_{\Lambda},\tau_1)>0$ depending only on $\Omega$, $n$, $\overline{M}_A$, $\underline{m}_{\Lambda}$ and $\tau_1$, and there exists a radially symmetric $C^{\infty}(\overline{\Omega^*})$ field $\Lambda^*>0$ such that
\begin{equation} \label{lambdabounds}
\displaystyle{\mathop{\rm{ess}\ \rm{inf}}_{\Omega}}\ \Lambda \le \displaystyle{\mathop{\min}_{\overline{\Omega^*}}}\ \Lambda^*\le \displaystyle{\mathop{\max}_{\overline{\Omega^*}}}\ \Lambda^*\le \displaystyle{\mathop{\rm{ess}\ \rm{sup}}_{\Omega}}\ \Lambda,\ \|(\Lambda^*)^{-1}\|_{L^1(\Omega^*)}=\|\Lambda^{-1}\|_{L^1(\Omega)},
\end{equation}
and
\be\label{ineqfklambda}
\lambda_1(\Omega^{\ast},\Lambda^{\ast}{\rm{Id}},\tau_1e_r,-\tau_2)\le\lambda_1(\Omega,A,v,V)-\eta.
\ee
\end{theo}

Notice that, as in Theorem~\ref{th1bis}, the assumptions of Theorem~\ref{faberkrahn} imply necessarily that $\overline{M}_A\ge\underline{m}_{\Lambda}$. Notice also that, in Theorem~\ref{faberkrahn}, contrary to our other results, we do not assume that $\lambda_1(\Omega,A,v,V)\geq 0$. In the previous results, we imposed a constraint on the distribution function of the negative part of the potential and we needed the nonnegativity of $\lambda_1(\Omega,A,v,V)$. Here, we first write
$$\lambda_1(\Omega,A,v,V)\ge\lambda_1(\Omega,A,v,-\tau_2)=-\tau_2+\lambda_1(\Omega,A,v,0)$$
and we apply Theorem~\ref{th1bis} to $\lambda_1(\Omega,A,v,0)$, which is positive. We complete the proof with further results which are established in Section~\ref{sec5}.

Observe also that, in the inequality (\ref{ineqfklambda}), the constraints $\tau_1$ and $\tau_2$ on the $L^{\infty}$ norms of the drift and the potential are saturated in the ball $\Omega^*$.

Actually, in Theorem~\ref{faberkrahn}, if we replace the assumption $\|V\|_{L^{\infty}(\Omega)}\le\tau_2$ by $\hbox{ess inf}_{\Omega}\ V\ge\tau_3$ (where $\tau_3\in\R$), then inequality (\ref{ineqfklambda}) is changed into
$$\lambda_1(\Omega^{\ast},\Lambda^{\ast}{\rm{Id}},\tau_1e_r,\tau_3)\le\lambda_1(\Omega,A,v,V)-\eta.$$
Since $\lambda_1(\Omega^{\ast},\Lambda^{\ast}{\rm{Id}},\tau_1e_r,\tau)=\lambda_1(\Omega^{\ast},\Lambda^{\ast}{\rm{Id}},\tau_1e_r,0)+\tau$ for all $\tau\in\R$, the previous inequality is better than (\ref{ineqfklambda}). In the following corollary, we choose to compare directly $V$ with $\hbox{ess inf}_{\Omega}\ V$.

\begin{cor}\label{corFK}
Let $\Omega\in {\mathcal C}$, $A\in W^{1,\infty}(\Omega,{\mathcal S}_n(\R))$, $v\in L^{\infty}(\Omega,\R^n)$ and $V\in L^{\infty}(\Omega)$. Call $\Lambda[A](x)$ the smallest eigenvalue of the matrix $A(x)$ at each point $x\in\overline{\Omega}$ and assume that $\gamma_A=\min_{\overline{\Omega}}\Lambda[A]>0$. Then
\be\label{fk}
\lambda_1(\Omega,A,v,V)\ge F_n(|\Omega|,\min_{\overline{\Omega}}\Lambda[A],\|v\|_{L^{\infty}(\Omega,\R^n)},\mathop{{\rm{ess}}\ {\rm{inf}}}_{\Omega}\ V),
\ee
where $F_n:(0,+\infty)\times(0,+\infty)\times[0,+\infty)\times\R\to\R$ is defined by
$$F_n(m,\gamma,\alpha,\beta)=\lambda_1(B_{(m/\alpha_n)^{1/n}}^n,\gamma{\rm{Id}},\alpha\ e_r,\beta)$$
for all $(m,\gamma,\alpha,\beta)\in(0,+\infty)\times(0,+\infty)\times[0,+\infty)\times\R$, and $B_{(m/\alpha_n)^{1/n}}^n$ denotes the Euclidean ball of $\R^n$ with center $0$ and radius $(m/\alpha_n)^{1/n}$. Furthermore, the inequality (\ref{fk}) is strict if $\Omega$ is not a ball.
\end{cor}

In Corollary~\ref{corFK}, formula (\ref{fk}) reduces to (\ref{RFK}) when $A=\hbox{Id}$ and $v=0$, $V=0$. Theo\-rem~\ref{faberkrahn} can then be viewed as a natural extension of the first Rayleigh conjecture to more general elliptic operators with potential, drift and general diffusion. We refer to Remark~\ref{interpretation} for further comments on these results.


\subsection{Some comparisons with results in the literature}

If in Theorem~\ref{th1}, the function $\Lambda$ is identically equal to a constant $\gamma>0$ in $\Omega$, and if $V\ge 0$, then inequality (\ref{ineqv0}) could also be derived implicitely from Theorem~1 by Talenti \cite{talenti}. In \cite{talenti}, Talenti's argument relies on the Schwarz symmetrization and one of the key inequalities which is used in \cite{talenti} is
$$\int_{\Omega}-\hbox{div}(A\nabla\varphi)\times\varphi=\int_{\Omega}A\nabla\varphi\cdot\nabla\varphi\ge\gamma\int_{\Omega}|\nabla\varphi|^2.$$
This kind of inequality cannot be used directly for our purpose since it does not take into account the fact that $A\ge\Lambda\ \hbox{Id}$ a.e. in $\Omega$, where the function $\Lambda$ may not be constant. The proofs of the present paper use a completely different rearrangement technique which has its own interest, and which allows us to take into account any non-constant function $\Lambda\in L^{\infty}_+(\Omega)$. Actually, paper \cite{talenti} was not concerned with eigenvalue problems, but with various comparison results for solutions of elliptic problems (see also \cite{at1, at2, atlm, tv}). Even in the case when $\Lambda$ is constant and $V\ge 0$, proving the inequality (\ref{ineqv0}) between the principal eigenvalues of the initial and rearranged operators by means of Talenti's results requires several extra arguments, some of them using results contained in Section~\ref{sec5} of the present paper. We also refer to Section~\ref{sec62} for additional comments in the case when $\Lambda$ is constant.

But, once again, besides the own interest and the novelty of the tools we use in the present paper, one of the main features in Theorem~\ref{th1} (and in Theorems~\ref{th1bis} and~\ref{faberkrahn}) is that the ellipticity function $\Lambda$ and its symmetrization $\Lambda^*$ are not constant in general (see Remark~\ref{remLambdanoncst}). Optimizing with non-constant coefficients in the second-order terms creates additional and substantial difficulties. In particular, the conclusion of Theorem~\ref{th1} does not follow from previous works, even implicitely and even if the lower-order terms are zero. More generally speaking, all the comparison results of the present paper are new even when $v=0$, namely when the operator $L$ is symmetric. Moreover, all the results are new also when the operators are one-dimensional (except Theorem~\ref{th2} the statement of which does make sense only when $n\ge 2$).

The improved version of Theorem \ref{th1} when $\Omega$ is not a ball, namely Theorem \ref{th1bis}, is also new and does not follow from earlier results.\par

As far as Theorem \ref{th2} is concerned, optimization problems for eigenvalues when the constraint on $A$ is expressed in terms of the determinant 
and the trace, or more general symmetric functions of the eigenvalues of $A$, have not been considered hitherto. 

Let us now focus on Theorem \ref{faberkrahn} and Corollary~\ref{corFK}. In a previous work (\cite{hnr,hnrpreprint}), we proved a somewhat more complete version of this Faber-Krahn inequality in the case of the Laplace operator with a drift term. Namely, let $\Omega$ be a $C^{2,\alpha}$ non empty bounded domain of $\R^n$ for some $0<\alpha<1$. For any vector field $v\in L^{\infty}(\Omega,\R^n)$, denote by
\be\label{lambda1Omegav}
\lambda_1(\Omega,v)=\lambda_1(\Omega,\hbox{Id},v,0)
\ee
the principal eigenvalue of $-\Delta+v\cdot \nabla$ in $\Omega$ under Dirichlet boundary condition. Then, the following Faber-Krahn type inequality holds:

\begin{theo} \label{FKlaplace}{\rm{\cite{hnr,hnrpreprint}}}
Let $\Omega$ be a $C^{2,\alpha}$ non-empty bounded connected open subset of $\R^n$ for some $0<\alpha<1$, let $\tau\geq 0$ and $v\in L^{\infty}(\Omega,\R^n)$ be such that $\left\Vert v\right\Vert_{L^{\infty}(\Omega,\R^n)}\leq \tau$. Then
\begin{equation} \label{ineqlaplace}
\lambda_1(\Omega,v)\geq \lambda_1(\Omega^{\ast},\tau e_r),
\end{equation}
and the equality holds if and only if, up to translation, $\Omega=\Omega^{\ast}$ and $v=\tau e_r$.
\end{theo}

\begin{rem}{\rm Here we quote exactly the statement of \cite{hnr,hnrpreprint}, but actually it is enough to assume that $\Omega$ is of class $C^2$.}
\end{rem}

Notice that we can recover Theorem \ref{FKlaplace} from the results of the present paper. Indeed, when $\Omega$ is not ball, the strict inequality in (\ref{ineqlaplace}) follows at once from Theorem \ref{faberkrahn}, and when $\Omega$ is a ball (say, with center $0$) and $v\neq \tau e_r$, this strict inequality will follow from Theorem~\ref{fixedlinftyball} (see Section~\ref{sec5} below). Strictly speaking, the inequality (\ref{ineqlaplace}) could also be derived from Theorem~2 in \cite{talenti} (see also \cite{at1,at2}) and from extra arguments similar to the ones used in Section~\ref{sec61}. But the case of equality is new, while Theorem~\ref{faberkrahn} is entirely new. Indeed, an important feature in Theorem~\ref{faberkrahn} is the fact that the diffusion $A$ is assumed to be bounded from below by $\Lambda\Id$ where $\Lambda$ is a possibly non-constant function, and that $\lambda_1(\Omega,A,v,V)$ is compared with $\lambda_1(\Omega^{\ast},\Lambda^{\ast}\hbox{Id},\|v\|_{\infty}e_r,-\|V\|_{\infty})$, where $\Lambda^{\ast}$ is also possibly nonconstant (in other words, the operator $\div(\Lambda^{\ast}\nabla)$ is not necessarily equal to a constant times the Laplace operator). Furthermore, another novelty in Theorem~\ref{faberkrahn} is that, when $\Omega$ is not a ball, the difference $\lambda_1(\Omega,A,v,V)-\lambda_1(\Omega^{\ast},\Lambda^{\ast}\hbox{Id},\|v\|_{\infty}e_r,-\|V\|_{\infty})$ is estimated from below by a positive quantity depending only on $\Omega$, $n$ and on some structural constants of the operator. All these observations imply that Theorem~\ref{faberkrahn} is definitely more general than Theorem~\ref{FKlaplace} and is not implicit in \cite{talenti}, or even in more recent works in the same spirit (like \cite{atlm}, for instance).

When the vector field $v$ is divergence free (in the sense of distributions), then $\lambda_1(\Omega,v)\ge\lambda_1(\Omega)$ (multiply $-\Delta\varphi_{\Omega,{\small{\hbox{Id}}},v,0}+v\cdot\nabla\varphi_{\Omega,{\small{\hbox{Id}}},v,0}=\lambda_1(\Omega,v)\varphi_{\Omega,{\small{\hbox{Id}}},v,0}$ by $\varphi_{\Omega,{\small{\hbox{Id}}},v,0}$ and integrate by parts over $\Omega$).\footnote{We refer to \cite{bhn} for a detailed analysis of the behavior of $\lambda_1(\Omega,A,Bv,V)$ when $B\to+\infty$ and $v$ is a fixed divergence free vector field in $L^{\infty}(\Omega)$.} Thus, minimizing $\lambda_1(\Omega,v)$ when $|\Omega|=m$ and $v$ is divergence free and satisfies $\|v\|_{L^{\infty}(\Omega,\R^n)}\le\tau$ (with given $m>0$ and $\tau\ge 0$), is the same as minimizing $\lambda_1(\Omega)$ in the Rayleigh conjecture. We also refer to \cite{hnr} and \cite{hnrpreprint} for further optimization results for $\lambda_1(\Omega,v)$ with $L^{\infty}$ constraints on the drifts.

\begin{rem}\label{nonsmooth}{\rm For non-empty connected and possibly unbounded open sets $\Omega$ with finite measure, the principal eigenvalue $\lambda_1(\Omega,A,v,V)$ of the operator $L=-\hbox{div}(A\nabla)+v\cdot\nabla+V$ can be defined as
$$\lambda_1(\Omega,A,v,V)=\sup\ \{\lambda\in\R,\ \exists\ \phi\in C^2(\Omega),\ \phi>0\hbox{ in }\Omega,\ (-L+\lambda)\phi\le 0\hbox{ in }\Omega\}.$$
When $\Omega$ is bounded, this definition is taken from \cite{bnv} (see also \cite{a,np}), and it coincides with the characterization (\ref{eq}) when $\Omega\in{\mathcal{C}}$. It follows from the arguments of Chapter~2 of \cite{bnv} that
\be\label{lambda1inf1}
\lambda_1(\Omega,A,v,V)=\displaystyle{\mathop{\inf}_{\Omega'\subset\subset\Omega,\ \Omega'\in{\mathcal{C}}}} \lambda_1(\Omega',A|_{\Omega'},v|_{\Omega'},V|_{\Omega'}),
\ee
where $A|_{\Omega'}$, $v|_{\Omega'}$, $V|_{\Omega'}$ denote the restrictions of the fields $A$, $v$ and $V$ to $\Omega'$. When $\Omega$ is a general non-empty open set with finite measure, we then define
\be\label{lambda1inf2}
\lambda_1(\Omega,A,v,V)=\inf_{j\in J}\lambda_1(\Omega_j,A|_{\Omega_j},v|_{\Omega_j},V|_{\Omega_j}),
\ee
where the $\Omega_j$'s are the connected components of $\Omega$.\par
Some of the comparison results which are stated in the previous subsections can then be extended to the class of general open sets $\Omega$ with finite measure (see Remarks~\ref{nonsmooth1},~\ref{nonsmooth2} and~\ref{nonsmooth3}).}
\end{rem}


\subsection{Main tools: a new type of symmetrization}

As already underlined, the proofs of Theorems \ref{th1}, \ref{th1bis}, \ref{th2} and~\ref{faberkrahn} do not use the usual Schwarz symmetrization. The key tool in the proofs is a new (up to our knowledge) rearrangement technique for some functions on $\Omega$, which can take into account non-constant ellipticity functions $\Lambda$. Roughly speaking, given $\Omega$, $A$, $v$ and $V$ such that $A\geq \Lambda \Id$, if $\varphi=\varphi_{\Omega,A,v,V}$ denotes the principal eigenfunction of the operator $-\hbox{div}(A\nabla)+v\cdot\nabla+V$ in $\Omega$ under Dirichlet boundary condition (that is, $\varphi$ solves (\ref{eq})), we associate to $\varphi$, $\Lambda$, $v$ and $V$ some rearranged functions or vector fields, which are called $\widetilde{\varphi}$, $\widehat{\Lambda}$, $\widehat{v}$ and $\widehat{V}$. They are defined on $\Omega^{\ast}$ and are built so that some quantities are preserved. The precise definitions will be given in Section \ref{rearrangement}, but let us quickly explain how the function $\widetilde{\varphi}$ is defined. Denote by $R$ the radius of 
$\Omega^{\ast}$. For all $0\leq a<1$, define
\[
\Omega_a=\left\{x\in \Omega,\ a<\varphi(x)\le1\right\}
\]
and define $\rho(a)\in \left(0,R\right]$ such that $\left\vert 
\Omega_a\right\vert=\left\vert B_{\rho(a)}\right\vert$, where $B_s$ denotes the open Euclidean ball of radius $s>0$ and centre $0$. Define also $\rho(1)=0$. The function $\rho:\left[0,1\right]\rightarrow \left[0,R\right]$ is decreasing, continuous, one-to-one and onto. Then, the rearrangement of $\varphi$ is the radially symmetric decreasing function $\widetilde{\varphi}:\overline{\Omega^{\ast}}\rightarrow \R$ vanishing on $\partial\Omega^{\ast}$ such that, for all $0\leq a<1$,
\[
\int_{\Omega_a} \div(A\nabla\varphi)(x)dx=\int_{B_{\rho(a)}} 
\div(\widehat{\Lambda}\nabla\widetilde{\varphi})(x)dx
\]
(we do not wish to give the explicit expression of the function $\widehat{\Lambda}$ right now). The fundamental inequality satisfied by $\widetilde{\varphi}$ is the fact that, for all $x\in \overline{\Omega^{\ast}}$, 
\begin{equation} \label{fund}
\widetilde{\varphi}(x)\geq \rho^{-1}(\left\vert x\right\vert)
\end{equation}
(see Corollary \ref{cor36} below, and Lemma~\ref{lem311} for strict inequalities when $\Omega$ is not a ball).\par
This symmetrization is definitely different from the Schwarz symmetrization since the distribution functions of $\varphi$ and $\widetilde{\varphi}$ are not the same in general. Moreover, the $L^1$ norm of the gradient of $\widetilde{\varphi}$ on $\Omega^*$ is larger than or equal to that of $\varphi$ on $\Omega$, and, when $A=\gamma\hbox{Id}$ (for a positive constant $\gamma$), the $L^2$ norm of the gradient of $\widetilde{\varphi}$ on $\Omega^*$ is larger than or equal to that of $\varphi$ on $\Omega$ (see Remark~\ref{remgradients} below).\par
Actually, the function $\varphi$ is not regular enough for this construction to be correct, and we have to deal with suitable approximations of $\varphi$. We refer to Section~\ref{rearrangement} and the following ones for exact and complete statements and proofs. Let us just mention that the proof of (\ref{fund}) relies, apart from the definition of $\widetilde{\varphi}$, on the usual isoperimetric inequality on $\R^n$.\par
Notice that the tools which are developed in this paper not only give new comparison results for symmetric and non-symmetric second-order operators with non-constant coefficients, but they also provide an alternative proof of the Rayleigh-Faber-Krahn isoperimetric inequa\-li\-ty (\ref{RFK}) for the Dirichlet Laplacian.\par
Finally, the new rearrangement we introduce in this paper is likely to be used in other problems involving elliptic partial differential equations.\hfill\break

\noindent{\bf Let us give a few open problems related to our results.} In all our results, several minimization problems for the principal eigenvalue of a second-order elliptic operator in a domain $\Omega$ under some constraints have been reduced to the same problems on the ball $\Omega^*$ centered at $0$ with the same Lebesgue measure and for operators with radially symmetric coefficients. However, even in the case of the ball and for operators with radially symmetric coefficients, some of these optimization problems remain open. For instance, in Corollary~\ref{cor1}, is it possible to compute explicitly the right-hand side of (\ref{minlambda}) ? An analogous question may be asked for the other theorems, corresponding to different constraints (even for Theorem \ref{faberkrahn}). 

When we combine Theorems~\ref{th1} and~\ref{th1bis}, it follows that the inequality (\ref{ineqv0}) is strict when $\Omega$ is not a ball and $\lambda_1(\Omega,A,v,V)>0$. But in Theorem \ref{th1}, when $\Omega$ is a ball, for which $A$, $v$ and $V$ does the case of equality occur in (\ref{ineqv0})~? Does this require that the initial data should be all radially symmetric~? The same question can be asked in Theorem~\ref{th2} as well. An answer to these questions would provide a complete analogue of Theorem~\ref{FKlaplace} for general second-order elliptic operators in divergence form. Furthermore, in Theorem \ref{th1}, in the general case when $\Lambda$ is not constant and even if $\Omega$ is a ball, can one state a result without $\varepsilon$ but with still keeping the constraints (\ref{bounds1})~?

When $\Omega=\Omega^{\ast}$, $\Lambda^{\ast}$ is fixed and $v$ and $V$ vary with some constraints on their $L^{\infty}$ norms, we prove in Section \ref{sec5} that there exist a unique $v$ and a unique $V$ minimizing $\lambda_1(\Omega^{\ast},\Lambda^{\ast}\hbox{Id},v,V)$. In particular, if $\Lambda^{\ast}$ is radially symmetric, then we show that $v$ and $V$ are given by inequality (\ref{ineqfklambda}) of Theorem~\ref{faberkrahn}. Many other optimization results in the ball can be asked if some of the fields $\Lambda^{\ast},v^{\ast}$ and $V^{\ast}$ are fixed while the others vary under some constraints. We intend to come back to all these issues in a forthcoming paper.

Here are some other open problems. In Theorem \ref{th2}, can one replace the determinant of $A$ by more general functions of the eigenvalues of $A$, namely $\sigma_q(A)$ with $p<q\leq n-1$ ? \par

It would also be very interesting to obtain results similar to ours for general second-order elliptic operators of the form
\[
-\sum_{i,j} a_{i,j}\partial_{i,j} + \sum_i b_i\partial_i +c,
\]
where the $a_{i,j}$'s are continuous in $\overline{\Omega}$ (but do not necessarily belong to $W^{1,\infty}(\Omega)$), and the $b_{i}$'s and $c$ are bounded in $\Omega$ (recall that such operators still have a real principal eigenvalue, see \cite{bnv}), and to consider other boundary conditions (Neumann, Robin, Stekloff problems...)\hfill\break

\noindent{\bf Outline of the paper.} The paper is organized as follows. Section \ref{rearrangement} is devoted to the precise definitions of the rearranged function and the proof of the inequalities satisfied by this rearrangement, whereas improved inequalities are obtained in Section \ref{equality} when $\Omega$ is not a ball. The proofs of Theorems \ref{th1}, \ref{th1bis} and \ref{th2} are given in Section \ref{appl}, while the Faber-Krahn inequalities (Theorem \ref{faberkrahn} and Corollary~\ref{corFK}) are established in Section \ref{sec5}. Some optimization results in a fixed domain, which are interesting in their own right and are also required for the proof of Theorem~\ref{faberkrahn}, are also proved in Section~\ref{sec5}. Finally, the appendix contains the proof of a technical approximation result (which is used in the proofs of Section \ref{appl}), a short remark about distribution functions and some useful asymptotics of $\lambda_1(\Omega^*,\tau e_r)=\lambda_1(\Omega^*,\hbox{Id},\tau e_r,0)$ when $\tau\to+\infty$.\hfill\break

\noindent{\bf Acknowledgements. }The authors thank C. Bandle for pointing out to us reference \cite{talenti}, and L.~Roques for valuable discussions.


\SE{Inequalities for the rearranged functions} \label{rearrangement}

In this section, we present a new spherical rearrangement of functions and we prove some pointwise and integral inequalities for the rearranged data. The results are of independent interest and this is the reason why we present them in a separate section.


\subsection{General framework, definitions of the rearrangements and basic properties}\label{sec21}

In this subsection, we give some assumptions which will remain valid throughout all Section~\ref{rearrangement}. Fix $\Omega\in{\mathcal{C}}$, 
$A_{\Omega}\in C^1(\overline{\Omega},{\mathcal S}_n(\R))$, $\Lambda_{\Omega}\in C^1(\overline{\Omega})$, $\omega\in C(\overline{\Omega})$ and $V\in C(\overline{\Omega})$. Assume that
\be\label{ALambdaOmega}
A_{\Omega}(x)\geq\Lambda_{\Omega}(x)\hbox{Id}\hbox{ for all }x\in\overline{\Omega},
\ee
and that there exists $\gamma>0$ such that
$$\Lambda_{\Omega}(x)\geq \gamma\hbox{ for all }x\in \overline{\Omega}.$$\par
Let $\psi$ be a $C^1(\overline{\Omega})$ function, analytic and positive in $\Omega$, such that $\psi=0$ on $\partial\Omega$ and
$$\nabla\psi(x)\neq 0\hbox{ for all }x\in\partial\Omega,$$
so that $\nu\cdot\nabla\psi<0$ on 
$\partial\Omega$, where $\nu$ denotes the outward unit normal to 
$\partial\Omega$. We always assume throughout this section that
$$f:=-\div (A_{\Omega}\nabla\psi)\hbox{ in }\Omega$$
is a non-zero polynomial, so that $\psi\in W^{2,p}(\Omega)$ for all $1\leq p<+\infty$ and $\psi\in C^{1,\alpha}(\overline{\Omega})$ for all $0\le\alpha<1$.\par
Set
$$M=\max_{x\in \overline{\Omega}} \psi(x).$$
For all $a\in \left[0,M\right)$, define
\[
\Omega_a=\left\{x\in \Omega,\ \psi(x)>a\right\}
\]
and, for all $a\in \left[0,M\right]$,
\[
\Sigma_a=\left\{x\in \overline{\Omega},\ \psi(x)=a\right\}.
\]
The set $\{x\in\overline{\Omega},\ \nabla\psi(x)=0\}$ is included in some compact set $K\subset\Omega$, which implies that the set
$$Z=\{a\in[0,M],\ \exists\ x\in\Sigma_{a},\ \nabla\psi(x)=0\}$$
of the critical values of $\psi$ is finite (\cite{soucek}) and can then be written as
$$Z=\{a_{1},\cdots,a_{m}\}$$
for some $m\in\N^*=\N\backslash\{0\}$. Observe also that $M\in Z$ and that $0\not\in Z$. One can then assume without loss of generality that
$$0<a_{1}<\cdots<a_{m}=M.$$
The set $Y=[0,M]\backslash Z$ of the non critical values of $\psi$ is open relatively to $[0,M]$ and can be written as
$$Y\ =\ [0,M]\backslash Z\ =\ [0,a_{1})\cup(a_{1},a_{2})\cup\cdots\cup(a_{m-1},M).$$
For all $a\in Y$, the hypersurface $\Sigma_{a}$ is of class $C^2$ (notice also that $\Sigma_0=\partial\Omega$ is of class $C^2$ by assumption) and $|\nabla\psi|$ does not vanish on 
$\Sigma_{a}$. Therefore, the functions defined on $Y$ by
\begin{equation}\label{ghi}\left\{\baa{l}
g\ :\ Y\ni a\mapsto\displaystyle{\int_{\Sigma_{a}}}|\nabla\psi(y)|^{-1}d\sigma_{a}(y)\\
h\ :\ Y\ni a\mapsto\displaystyle{\int_{\Sigma_{a}}}f(y)|\nabla\psi(y)|^{-1}d\sigma_{a}(y)\\
i\ :\ Y\ni a\mapsto\displaystyle{\int_{\Sigma_{a}}}d\sigma_{a}(y)\eaa\right.\end{equation}
are (at least) continuous in $Y$ and $C^1$ in $Y\setminus\left\{0\right\}$, where $d\sigma_{a}$ denotes the surface measure on $\Sigma_{a}$ for $a\in Y$.\par
Denote by $R$ the radius of $\Omega^{\ast}$ (the open Euclidean ball centered at the origin and such that $|\Omega^*|=|\Omega|$, that is $\Omega^*=B_R$). For all $a\in[0,M)$, let
$\rho(a)\in(0,R]$ be defined so that
$$|\Omega_{a}|=|B_{\rho(a)}|=\alpha_n\rho(a)^n.$$
Recall that $\alpha_n$ is the volume of the unit ball $B_1$. The function $\rho$ is extended at $M$ by
$$\rho(M)=0.$$

\begin{lem}\label{lem31} The function $\rho$ is a continuous decreasing map from $[0,M]$ onto $[0,R]$.
\end{lem}

\noindent{\bf{Proof.}} The function $\rho\ :\ [0,M]\to[0,R]$ is clearly decreasing since
$$\left|\{x\in\Omega,\ a<\psi(x)\le b\}\right|>0$$
for all $0\le a<b\le M$. Fix now any $a\in(0,M]$. Since $\psi\in W^{2,p}(\Omega)$ (actually, for all $1\le p<+\infty$), one has
$$\frac{\partial^2\psi}{\partial x_i\partial x_j}\times{\bf{1}}_{\{\psi=a\}}=\frac{\partial\psi}{\partial x_i}\times{\bf{1}}_{\{\psi=a\}}=0\hbox{ almost everywhere in }\Omega$$
for all $1\le i,j\le n$, where ${\bf{1}}_E$ denotes the characteristic function of a set $E$. Therefore, $f\times{\bf{1}}_{\{\psi=a\}}=0$ almost everywhere in $\Omega$. Since $f$ is a nonzero polynomial, one gets that
$$|\Sigma_{a}|=0\hbox{ for all }a\in(0,M].$$
Notice that $|\Sigma_{0}|=|\partial\Omega|=0$ as well. Lastly, $\rho(0)=R$ and $\rho(M)=0$. As a conclusion, the function $\rho$ is continuous on $[0,M]$ and is a one-to-one and onto map from $[0,M]$ to $[0,R]$.\hfill\fin

\begin{lem}\label{lem32} The function $\rho$ is of class $C^1$ in $Y$ and
$$\forall\ a\in Y,\quad\rho'(a)=-(n\alpha_n\rho(a)^{n-1})^{-1}g(a)=-(n\alpha_n\rho(a)^{n-1})^{-1}\int_{\Sigma_a}|\nabla\psi(y)|^{-1}d\sigma_a(y).$$
\end{lem}

\noindent{\bf{Proof.}} Fix $a\in Y$. Let $\eta>0$ be such that $[a,a+\eta]\subset Y$. For $t\in(0,\eta)$,
$$\baa{rcl}
\alpha_n[\rho(a+t)^n-\rho(a)^n]=|\Omega_{a+t}|-|\Omega_{a}| & = & -\displaystyle{\int_{\{a<\psi(x)\le a+t\}}}dx\\
& = & -\displaystyle{\int_a^{a+t}}\left(\displaystyle{\int_{\Sigma_{b}}}|\nabla\psi(y)|^{-1}d\sigma_{b}(y)\right)db\eaa$$
from the co-area formula. Hence,
$$\frac{\alpha_n[\rho(a+t)^n-\rho(a)^n]}{t}\to -g(a)\ \hbox{ as }t\to 0^+$$
for all $a\in Y$, due to the continuity of $g$ on $Y$. Similarly, one has that
$$\frac{\alpha_n[\rho(a+t)^n-\rho(a)^n]}{t}\to -g(a)\ \hbox{ as }t\to 0^-$$
for all $a\in Y\backslash\{0\}$. The conclusion of the lemma follows since $Y\subset[0,M)$, whence $\rho(a)\neq 0$ for all $a\in Y$.\hfill\fin\break

We now define the function $\widetilde{\psi}$ in $\overline{\Omega^*}$, which is a spherical rearrangement of $\psi$ by means of a new type of symmetrization. The definition of $\widetilde{\psi}$ involves the rearrangement of the datum~$\Lambda_{\Omega}$.\par
First, call
$$E=\{x\in\overline{\Omega^*},\ |x|\in\rho(Y)\}.$$
The set $E$ is a finite union of spherical shells and, from Lemma \ref{lem31}, it is open relatively to $\overline{\Omega^*}$ and can be written as
$$E=\{x\in\R^n,\ |x|\in(0,\rho(a_{m-1}))\cup\cdots\cup(\rho(a_{2}),\rho(a_{1}))\cup(\rho(a_{1}),R]\}.$$
with
$$0=\rho(a_{m})=\rho(M)<\rho(a_{m-1})<\cdots<\rho(a_{1})<R.$$
Notice that $0\not\in E$.\par
Next, for all $r\in \rho(Y)$, set
\begin{equation} \label{G}
G(r)=\frac{\displaystyle{\int_{\Sigma_{\rho^{-1}(r)}} \left\vert \nabla\psi(y)\right\vert^{-1}d\sigma_{\rho^{-1}(r)}}}{\displaystyle{\int_{\Sigma_{\rho^{-1}(r)}} \Lambda_{\Omega}(y)^{-1}\left\vert \nabla\psi(y)\right\vert^{-1}d\sigma_{\rho^{-1}(r)}}}>0,
\end{equation}
where $\rho^{-1}\ :\ [0,R]\to[0,M]$ denotes the reciprocal of the function $\rho$. For all $x\in E$, define
\be\label{defhatLambda}
\widehat{\Lambda}(x)=G(|x|).
\ee
The function $\widehat{\Lambda}$ is then defined almost everywhere in $\overline{\Omega^*}$. By the obervations above and since $\Lambda_{\Omega}$ is positive and $C^1(\overline{\Omega})$, the function $\widehat{\Lambda}$ is continuous on $E$ and $C^1$ on $E\cap\Omega^*$. Furthermore, $\widehat{\Lambda}\in L^{\infty}(\Omega^*)$ and
\be\label{infiniLambda}
0<\min_{\overline{\Omega}}\Lambda_{\Omega}\le\mathop{\hbox{ ess inf}}_{\Omega^*}\widehat{\Lambda}\le\mathop{\hbox{ ess sup}}_{\Omega^*}\widehat{\Lambda}\le\max_{\overline{\Omega}}\Lambda_{\Omega}.
\ee
For any two real numbers $a<b$ such that $[a,b]\subset Y$, the co-area formula gives
$$\baa{rcl}
\displaystyle{\int_{\Omega_a\backslash\Omega_b}}\Lambda_{\Omega}(y)^{-1}dy & = & \displaystyle{\int_a^b}\left(\int_{\Sigma_s}\Lambda_{\Omega}(y)^{-1}|\nabla\psi(y)|^{-1}d\sigma_s(y)\right)ds\\
& = & \displaystyle{\int_{\rho(b)}^{\rho(a)}}\left(\frac{\displaystyle{\int_{\Sigma_{\rho^{-1}(t)}}}\Lambda_{\Omega}(y)^{-1}|\nabla\psi(y)|^{-1}d\sigma_{\rho^{-1}(t)}(y)}{\displaystyle{\int_{\Sigma_{\rho^{-1}(t)}}}|\nabla\psi(y)|^{-1}d\sigma_{\rho^{-1}(t)}(y)}\right)n\alpha_nt^{n-1}dt.\eaa$$
The last equality is obtained from Lemma~\ref{lem32} after the change of variables $s=\rho^{-1}(t)$. Since $\widehat{\Lambda}$ is radially symmetric, it follows by (\ref{G}-\ref{defhatLambda}) that
$$\int_{\Omega_a\backslash\Omega_b}\Lambda_{\Omega}(y)^{-1}dy=\int_{S_{\rho(b),\rho(a)}}\widehat{\Lambda}(x)^{-1}dx,$$
where, for any $0\le s<s'$, $S_{s,s'}$ denotes
$$S_{s,s'}=\{x\in\R^n,\ s<|x|<s'\}.$$
Lebesgue's dominated convergence theorem then implies that
\be\label{intLambda}
\int_{\Omega}\Lambda_{\Omega}(y)^{-1}dy=\int_{\Omega^*}\widehat{\Lambda}(x)^{-1}dx.
\ee\par
Lastly, set $F(0)=0$ and, for all $r\in \rho(Y)$, set
\be\label{defF}
F(r)=\frac{1}{n\alpha_nr^{n-1}G(r)}\int_{\Omega_{\rho^{-1}(r)}}\div(A_{\Omega}\nabla\psi)(x)dx.
\ee
The function $F$ is then defined almost everywhere in $[0,R]$.

\begin{lem}\label{lem33} The function $F$ belongs to $L^{\infty}(\left[0,R\right])$ and is continuous on $\rho(Y)\cup \left\{0\right\}$. Moreover, $F<0$ on $\rho(Y)$.
\end{lem}

\noindent{\bf{Proof.}} The continuity of $F$ on $\rho(Y)$ is a consequence 
of Lemma \ref{lem31}, of the continuity of $\widehat{\Lambda}$ on $E$ and of the fact that $\div(A_{\Omega}\nabla\psi)=f$ in $\Omega$, with $f$ continuous and thus bounded in $\overline{\Omega}$.\par
Observe that, since $\Lambda_{\Omega}(x)\geq \gamma>0$ for all $x\in \overline{\Omega}$, one has $\widehat{\Lambda}(x)\geq \gamma$ for all $x\in E$. For $0<r\le R$ with $r\in \rho(Y)$ ($\supset(0,\rho(a_{m-1}))$), one has
$$|F(r)|\le (n\alpha_nr^{n-1}\gamma)^{-1}\|f\|_{L^{\infty}(\Omega)}\ \alpha_nr^n=(n\gamma)^{-1}\|f\|_{L^{\infty}(\Omega)}\ r,$$
thus $F$ is continuous at $0$ as well and belongs to $L^{\infty}\left(\left[0,R\right]\right)$. Finally, for all $r\in \rho(Y)$, since $\psi(y)=\rho^{-1}(r)$ for all $y\in \Sigma_{\rho^{-1}(r)}$ and since $\psi>\rho^{-1}(r)$ in $\Omega_{\rho^{-1}(r)}$ and $\left\vert \nabla\psi(y)\right\vert\neq 0$ for all $y\in \Sigma_{\rho^{-1}(r)}$, one has
$$\nu_{\rho^{-1}(r)}\cdot\nabla\psi<0\hbox{ on }\Sigma_{\rho^{-1}(r)},$$
where, for any $a\in Y$, $\nu_a$ denotes the outward unit normal on $\partial\Omega_a$. Therefore
$$\nabla\psi(y)=-\left\vert \nabla\psi(y)\right\vert \nu_{\rho^{-1}(r)}(y)\hbox{ for all }r\in\rho(Y)\hbox{ and }y\in\Sigma_{\rho^{-1}(r)}.$$
As a consequence, the Green-Riemann formula yields that, for all $r\in \rho(Y)$,
$$
\begin{array}{lll}
\displaystyle \int_{\Omega_{\rho^{-1}(r)}} \div(A_{\Omega}\nabla\psi)(y)dy&=&\displaystyle \int_{\Sigma_{\rho^{-1}(r)}} A_{\Omega}(y)\nabla\psi(y)\cdot \nu_{\rho^{-1}(r)}(y)d\sigma_{\rho^{-1}(r)}(y)\\
\\&=&\displaystyle -\int_{\Sigma_{\rho^{-1}(r)}}  A_{\Omega}(y)\nu_{\rho^{-1}(r)}(y)\cdot \nu_{\rho^{-1}(r)}(y) \left\vert \nabla\psi(y)\right\vert d\sigma_{\rho^{-1}(r)}(y)<0,
\end{array}$$
which ends the proof.  \hfill\fin\break

For all $x\in\overline{\Omega^*}$, set
\be\label{defpsi}
\widetilde{\psi}(x)=-\int_{|x|}^RF(r)dr.
\ee
The function $\widetilde{\psi}$ is then radially symmetric and it vanishes on $\partial\Omega^*=\partial B_R$. From Lemma~\ref{lem33},
$$\widetilde{\psi}>0\hbox{ in }\Omega^*,$$
$\widetilde{\psi}$ is continuous on $\overline{\Omega^*}$, decreasing with respect to $\left\vert x\right\vert$ in $\overline{\Omega^*}$, and $C^1$ on $E\cup\left\{0\right\}$ (remember that $F(0)=0$). Note that
$$\widetilde{\psi}\in H^1_0(\Omega^{\ast})\cap W^{1,\infty}(\Omega^*).$$
Moreover, the following statement holds true:

\begin{lem}\label{lem34} The function $\widetilde{\psi}$ is of class $C^2$ in $E\cap\Omega^*$.
\end{lem}

\noindent{\bf{Proof.}} By definition of $\widetilde{\psi}$ and since $\widehat{\Lambda}$ is $C^1$ in $E\cap\Omega^*$, it is enough to prove that the function
$$z\ :\ r\mapsto\int_{\Omega_{\rho^{-1}(r)}}\div (A_{\Omega}\nabla\psi)(x)dx=-\int_{\Omega_{\rho^{-1}(r)}}f(x)dx$$
is of class $C^1$ on $\rho(Y)$. It would actually be enough to prove that $z$ is $C^1$ on $\rho(Y)\backslash\{R\}$.\par
Let $r$ be fixed in $\rho(Y)=(0,\rho(a_{m-1}))\cup\cdots\cup(\rho(a_{2}),\rho(a_{1}))\cup(\rho(a_{1}),R]$ and let $\eta>0$ be such that $[r-\eta,r]\subset\rho(Y)$. For $t\in(0,\eta)$, one has
$$\baa{rcl}
z(r-t)-z(r) & = & \displaystyle{\int_{\{\rho^{-1}(r)<\psi(x)\le\rho^{-1}(r-t)\}}}f(x)dx\\
& = & \displaystyle{\int_{\rho^{-1}(r)}^{\rho^{-1}(r-t)}}\left(\displaystyle{\int_{\Sigma_{a}}}f(y)|\nabla\psi(y)|^{-1}d\sigma_{a}(y)\right)da=\displaystyle{\int_{\rho^{-1}(r)}^{\rho^{-1}(r-t)}}h(a)da,\eaa$$
where $h$ is defined in (\ref{ghi}). Since $\rho^{-1}$ is of class $C^1$ on $\rho(Y)$ from Lemma \ref{lem32} and since $h$ is continuous on $Y$, it follows that
$$\frac{z(r-t)-z(r)}{-t}\to h(\rho^{-1}(r))(\rho^{-1})'(r)=-\frac{n\alpha_nr^{n-1}h(\rho^{-1}(r))}{g(\rho^{-1}(r))}\hbox{ as }t\to 0^+.$$
The same limit holds as $t\to 0^-$ for all $r\in\rho(Y)\backslash\{R\}$. Therefore, the function $z$ is differentiable on $\rho(Y)$ and
$$z'(r)=-\frac{n\alpha_nr^{n-1}h(\rho^{-1}(r))}{g(\rho^{-1}(r))}\ \hbox{ for all }r\in\rho(Y).$$
Since $\rho^{-1}$ is continuous on $[0,R]$, and $g$ and $h$ are continuous on $Y$, the function $z$ is of class $C^1$ on $\rho(Y)$. That completes the proof of Lemma \ref{lem34}.\hfill\fin\break\par

We now define a rearranged drift $\widehat{v}$ and a rearranged potential $\widehat{V}$. For all $x\in E$, define
\be\label{defhatv}
\widehat{v}(x)=\left(\frac{\displaystyle{\int_{\Sigma_{\rho^{-1}(|x|)}} \omega(y)^{2}\Lambda_{\Omega}(y)^{-1}\vert\nabla\psi(y)\vert^{-1}d\sigma_{\rho^{-1}(|x|)}(y)}}{\displaystyle{\int_{\Sigma_{\rho^{-1}(|x|)}} \Lambda_{\Omega}(y)^{-1}\left\vert \nabla\psi(y)\right\vert^{-1}d\sigma_{\rho^{-1}(|x|)}(y)}}\right)^{1/2}e_r(x),
\ee
(remember that $e_r$ is defined by (\ref{radial})). The vector field $\widehat{v}$ is then defined almost everywhere in $\overline{\Omega^*}$. Notice also that $|\widehat{v}|$ is radially symmetric, that $\widehat{v}(x)$ points in the direction $e_r(x)$ at each point $x\in E$, that $\widehat{v}$ belongs to $L^{\infty}(\Omega^*)$ and that
\be\label{infiniv}
\mathop{\hbox{ess inf}}_{\Omega}|\omega|\le\mathop{\hbox{ess inf}}_{\Omega^*}|\hat{v}|\le\mathop{\hbox{ess sup}}_{\Omega^*}|\hat{v}|\le\mathop{\hbox{ess sup}}_{\Omega}|\omega|=\|\omega\|_{L^{\infty}(\Omega)}.
\ee
Furthermore, since $\Lambda_{\Omega}^{-1}$ and $\omega$ are continuous in $\overline{\Omega}$, the vector field $\hat{v}$ is continuous in $E$, and, as it was done for (\ref{intLambda}), it is easy to check that
\be\label{intv2}
\int_{\Omega}\omega(y)^2\Lambda_{\Omega}(y)^{-1}dy=\int_{\Omega^*}|\widehat{v}(x)|^2\widehat{\Lambda}(x)^{-1}dx.
\ee\par
Lastly, for all $x\in E$, define
\be\label{defhatV}
\widehat{V}(x)=\frac{-\displaystyle{\int_{\Sigma_{\rho^{-1}(|x|)}} V^-(y) \left\vert \nabla\psi(y)\right\vert^{-1}d\sigma_{\rho^{-1}(|x|)}(y)}}{\displaystyle{\int_{\Sigma_{\rho^{-1}(|x|)}} \left\vert \nabla\psi(y)\right\vert^{-1}d\sigma_{\rho^{-1}(|x|)}(y)}},
\ee
where $V^-(y)$ denotes the negative part of $V(y)$, that is $V^-(y)=\max(0,-V(y))$. The function $\hat{V}$ is then defined almost everywhere in $\overline{\Omega^*}$. Observe that $\widehat{V}$ is radially symmetric, nonpositive, belongs to $L^{\infty}(\Omega^*)$, is continuous in $E$, and that
\be\label{infiniV}
-\|V\|_{L^{\infty}(\Omega)}\le\min_{\overline{\Omega}}(-V^-)\le\mathop{\hbox{ess inf}}_{\Omega^*}\widehat{V}\le\mathop{\hbox{ess sup}}_{\Omega^*}\widehat{V}\le 0.
\ee


\subsection{Pointwise comparison between $\psi$ and $\widetilde{\psi}$}

The first interest of the spherical rearrangement which was defined in the previous subsection is that the functions $\psi$ and $\tilde{\psi}$ can be compared on the sets $\Sigma_a$ and $\partial B_{\rho(a)}$. Namely, the function $\widetilde{\psi}$ satisfies the following key inequalities, which are summarized in Proposition~\ref{pro35} and Corollary~\ref{cor36}:

\begin{pro}\label{pro35} For any unit vector $e$ of $\R^n$, the function
$$\baa{rcl}
\widetilde{\Psi}\ :\ [0,M] & \to & \R_+\\
a & \mapsto & \widetilde{\psi}(\rho(a)e)\eaa$$
is continuous on $[0,M]$, differentiable on $Y$, and
\begin{equation}\label{eq2pro35}
\forall a\in Y,\quad \widetilde{\Psi}'(a)\ge 1.
\end{equation}
\end{pro}

Before giving the proof of Proposition \ref{pro35}, let us first establish the following important corollary.

\begin{cor}\label{cor36} For all $x\in\overline{\Omega^*}$,
$$\widetilde{\psi}(x)\ge \rho^{-1}(|x|).$$
\end{cor}

\noindent{\bf{Proof.}} Since $\widetilde{\Psi}$ is continuous on $\left[0,M\right]$ and differentiable on $\left[0,M\right]$ except on a finite set of points and since $\widetilde{\Psi}(0)=0$, the mean-value theorem and (\ref{eq2pro35}) show that $\widetilde{\Psi}(a)\geq a$ for all $a\in \left[0,M\right]$, which means that $\widetilde{\psi}(\rho(a)e)\geq a$ for all $a\in \left[0,M\right]$ and all unit vector $e$. Since $\widetilde{\psi}$ is radially symmetric, Corollary~\ref{cor36} follows from Lemma~\ref{lem31}.\hfill\fin\break

\noindent{\bf{Proof of Proposition \ref{pro35}.}} Let us first observe that the function $\widetilde{\Psi}$ is differentiable on $Y$, from Lemma \ref{lem32} and the fact that $\widetilde{\psi}$ is $C^1$ in $E$ (and even in $E\cup\{0\}$). Furthermore, since $\widetilde{\psi}$ is radially symmetric, 
and decreasing with respect to the variable $|x|$ and since $\rho$ is itself decreasing, 
it is enough to prove that
\be\label{rho'psi}
\forall\ x\in E,\quad|\rho'(\rho^{-1}(|x|))|\times|\nabla \widetilde{\psi}(x)|\ge 1.
\ee\par
We will make use of the following inequality:
\begin{equation} \label{ineq}
\forall\ x\in E,\ \frac{\displaystyle{\int_{\Sigma_{\rho^{-1}(|x|)}}A_{\Omega}(y)\nu_{\rho^{-1}(|x|)}(y)\cdot\nu_{\rho^{-1}(|x|)}(y) \left\vert \nabla\psi(y)\right\vert d\sigma_{\rho^{-1}(|x|)}(y)}}{\displaystyle{\int_{\Sigma_{\rho^{-1}(|x|)}} \left\vert \nabla\psi(y)\right\vert^{-1} d\sigma_{\rho^{-1}(|x|)}(y)}} \leq \widehat{\Lambda}(x)\ \vert \nabla \widetilde{\psi}(x)\vert^2,
\end{equation}
where one recalls that $\nu_{\rho^{-1}(|x|)}$ denotes the outward unit normal on $\partial\Omega_{\rho^{-1}(|x|)}$. We postpone the proof of (\ref{ineq}) to the end of this subsection and go on in the proof of Proposition \ref{pro35}.\par
Fix $x\in E$ and set $r=\left\vert x\right\vert$. Since $\rho^{-1}(r)\in Y$, there exists $\eta>0$ such that $\rho^{-1}(r-t)\in Y$ for all $t\in[0,\eta]$. For $t\in(0,\eta]$, the Cauchy-Schwarz inequality gives
\begin{equation} \label{cauchy}
\begin{array}{lll}
\left(\frac{\displaystyle{\int_{\Omega_{\rho^{-1}(r)}\backslash\Omega_{\rho^{-1}(r-t)}} \left\vert \nabla\psi(y)\right\vert dy}}{\displaystyle{\left\vert \Omega_{\rho^{-1}(r)}\backslash\Omega_{\rho^{-1}(r-t)} \right\vert}} \right)^2  &\leq & \frac{\displaystyle{\int_{\Omega_{\rho^{-1}(r)}\backslash\Omega_{\rho^{-1}(r-t)}} \Lambda_{\Omega}(y)^{-1}dy}}{\displaystyle{\left\vert \Omega_{\rho^{-1}(r)}\backslash\Omega_{\rho^{-1}(r-t)} \right\vert}} \\
& & \times\frac{\displaystyle{\int_{\Omega_{\rho^{-1}(r)}\backslash\Omega_{\rho^{-1}(r-t)}} \Lambda_{\Omega}(y)\left\vert \nabla\psi(y)\right\vert^2dy}}{\displaystyle{\left\vert \Omega_{\rho^{-1}(r)}\backslash\Omega_{\rho^{-1}(r-t)} \right\vert}}.
\end{array}
\end{equation}
The left-hand side of (\ref{cauchy}) is equal to
\[
\left(\frac{\displaystyle{\int_{\Omega_{\rho^{-1}(r)}\backslash\Omega_{\rho^{-1}(r-t)}} \left\vert \nabla\psi(y)\right\vert dy}}{\displaystyle{\left\vert \Omega_{\rho^{-1}(r)}\backslash\Omega_{\rho^{-1}(r-t)} \right\vert}} \right)^2= 
\left(\frac{\displaystyle{\int_{\Omega_{\rho^{-1}(r)}\backslash\Omega_{\rho^{-1}(r-t)}} \left\vert \nabla\psi(y)\right\vert dy}}{\displaystyle{\rho^{-1}(r-t)-\rho^{-1}(r)}}\displaystyle \times\frac{\rho^{-1}(r-t)-\rho^{-1}(r)}{\left\vert \Omega_{\rho^{-1}(r)}\backslash\Omega_{\rho^{-1}(r-t)} \right\vert}\right)^2.
\]
By the co-area formula, 
\[
\lim_{t\rightarrow 0^{+}} \frac{\displaystyle{\int_{\Omega_{\rho^{-1}(r)}\backslash\Omega_{\rho^{-1}(r-t)}} \left\vert \nabla\psi(y)\right\vert dy}}{\rho^{-1}(r-t)-\rho^{-1}(r)}=\int_{\Sigma_{\rho^{-1}(r)}} d\sigma_{\rho^{-1}(r)}(y)=i(\rho^{-1}(r)),
\]
and
\[
\lim_{t\rightarrow 0^{+}} \frac{\rho^{-1}(r-t)-\rho^{-1}(r)}{\left\vert \Omega_{\rho^{-1}(r)}\backslash\Omega_{\rho^{-1}(r-t)} \right\vert}=\frac{1}{\displaystyle{\int_{\Sigma_{\rho^{-1}(r)}}}|\nabla\psi(y)|^{-1}d\sigma_{\rho^{-1}(r)}(y)}=\frac 1{n\alpha_nr^{n-1}\left\vert \rho^{\prime}(\rho^{-1}(r))\right\vert}
\]
from Lemma~\ref{lem32}. By the isoperimetric inequality applied to $\Sigma_{\rho^{-1}(r)}=\partial\Omega_{\rho^{-1}(r)}$ and $\partial B_r$, one has
\begin{equation} \label{isop}
0<n\alpha_nr^{n-1}\le i(\rho^{-1}(r))=\int_{\Sigma_{\rho^{-1}(r)}} d\sigma_{\rho^{-1}(r)}(y),
\end{equation}
Therefore, one obtains
\begin{equation} \label{cauchy1}\baa{rcl}
\displaystyle{\mathop{\lim}_{t\rightarrow 0^{+}}} \left(\frac{\displaystyle{\int_{\Omega_{\rho^{-1}(r)}\backslash\Omega_{\rho^{-1}(r-t)}} \left\vert \nabla\psi(y)\right\vert dy}}{\left\vert \Omega_{\rho^{-1}(r)}\backslash\Omega_{\rho^{-1}(r-t)} \right\vert} \right)^2 & \ge & \left(\displaystyle{\frac{i(\rho^{-1}(r))}{n\alpha_nr^{n-1}}}\right)^2\times\displaystyle\frac 1{\left\vert \rho^{\prime}(\rho^{-1}(r))\right\vert^2}\\
& \geq & \displaystyle{\frac 1{\left\vert \rho^{\prime}(\rho^{-1}(r))\right\vert^2}}.\eaa
\end{equation}
The first factor of the right-hand side of (\ref{cauchy}) is equal to
\[
\frac{\displaystyle{\int_{\Omega_{\rho^{-1}(r)}\backslash\Omega_{\rho^{-1}(r-t)}} \Lambda_{\Omega}(y)^{-1}dy}}{\displaystyle{\left\vert \Omega_{\rho^{-1}(r)}\backslash\Omega_{\rho^{-1}(r-t)} \right\vert}}= \frac{\displaystyle{\int_{\Omega_{\rho^{-1}(r)}\backslash\Omega_{\rho^{-1}(r-t)}} \Lambda_{\Omega}(y)^{-1}dy}}{\displaystyle{\rho^{-1}(r-t)-\rho^{-1}(r)}} \displaystyle \times \frac{\rho^{-1}(r-t)-\rho^{-1}(r)}{\left\vert \Omega_{\rho^{-1}(r)}\backslash\Omega_{\rho^{-1}(r-t)} \right\vert},
\]
and the co-area formula therefore shows that
\begin{equation}\label{cauchy2}
\lim_{t\rightarrow 0^{+}} \frac{\displaystyle{\int_{\Omega_{\rho^{-1}(r)}\backslash\Omega_{\rho^{-1}(r-t)}} \Lambda_{\Omega}(y)^{-1}dy}}{\left\vert \Omega_{\rho^{-1}(r)}\backslash\Omega_{\rho^{-1}(r-t)} \right\vert} =\widehat{\Lambda}(x)^{-1}
\end{equation}
from (\ref{G}) and (\ref{defhatLambda}). Finally, the coarea formula again implies that
\begin{equation} \label{cauchy3}\baa{l}
\displaystyle{\mathop{\lim}_{t\rightarrow 0^{+}}} \frac{\displaystyle{\int_{\Omega_{\rho^{-1}(r)}\backslash\Omega_{\rho^{-1}(r-t)}} \Lambda_{\Omega}(y)\left\vert \nabla\psi(y)\right\vert^2dy}}{\left\vert \Omega_{\rho^{-1}(r)}\backslash\Omega_{\rho^{-1}(r-t)} \right\vert}=\displaystyle{\frac{\displaystyle{\int_{\Sigma_{\rho^{-1}(r)}}} \Lambda_{\Omega}(y)|\nabla\psi(y)|d\sigma_{\rho^{-1}(r)}(y)}
{\displaystyle{\int_{\Sigma_{\rho^{-1}(r)}}}|\nabla\psi(y)|^{-1}d\sigma_{\rho^{-1}(r)}(y)}}\\
\qquad\qquad\le\displaystyle{\frac{\displaystyle{\int_{\Sigma_{\rho^{-1}(r)}}} A_{\Omega}(y)\nu_{\rho^{-1}(r)}(y)\cdot\nu_{\rho^{-1}(r)}(y)|\nabla\psi(y)|d\sigma_{\rho^{-1}(r)}(y)}
{\displaystyle{\int_{\Sigma_{\rho^{-1}(r)}}}|\nabla\psi(y)|^{-1}d\sigma_{\rho^{-1}(r)}(y)}}\\
\qquad\qquad\le\hat{\Lambda}(x)|\nabla\tilde{\psi}(x)|^2\eaa
\ee
by (\ref{ineq}).\par
Finally, (\ref{cauchy}), (\ref{cauchy1}), (\ref{cauchy2}) and (\ref{cauchy3}) imply that
\be\label{rho'psibis}
\frac{1}{|\rho'(\rho^{-1}(r))|^2}\le\left(\frac{i(\rho^{-1}(r))}{n\alpha_nr^{n-1}}\right)^2\times\frac{1}{|\rho'(\rho^{-1}(r))|^2}\le|\nabla\tilde{\psi}(x)|^2.
\ee
Therefore, inequality (\ref{rho'psi}) holds and so does inequality (\ref{eq2pro35}). \hfill\fin

\begin{rem} Observe that (\ref{cauchy}), (\ref{cauchy2}) and (\ref{cauchy3}) together with the co-area formula yield
\be\label{gradient}
\displaystyle{\frac{\displaystyle{\int_{\Sigma_{\rho^{-1}(|x|)}}}d\sigma_{\rho^{-1}(|x|)}(y)}
{\displaystyle{\int_{\Sigma_{\rho^{-1}(|x|)}}}|\nabla\psi(y)|^{-1}d\sigma_{\rho^{-1}(|x|)}(y)}}=\mathop{\lim}_{t\rightarrow 0^{+}} \frac{\displaystyle{\int_{\Omega_{\rho^{-1}(|x|)}\backslash\Omega_{\rho^{-1}(|x|-t)}}\left\vert \nabla\psi(y)\right\vert dy}}{\left\vert \Omega_{\rho^{-1}(|x|)}\backslash\Omega_{\rho^{-1}(|x|-t)} \right\vert} \leq|\nabla\widetilde{\psi}(x)|
\ee
for all $x\in E$.
\end{rem}

We now give the\hfill\break
\noindent{\bf{Proof of (\ref{ineq}).}} Fix $x\in E$ and call $r=|x|$. Notice first that, as was already observed, for all $y\in  \partial\Omega_{\rho^{-1}(r)}$, 
\[
\nabla\psi(y)=-\left\vert \nabla\psi(y)\right\vert \nu_{\rho^{-1}(r)}(y).
\]
The Green-Riemann formula and the choice of $\widetilde{\psi}$ therefore yield
\begin{equation}\label{eqsym1}\baa{rcl}
\displaystyle{\int_{\Sigma_{\rho^{-1}(r)}}}A_{\Omega}(y)\nu_{\rho^{-1}(r)}(y)
\cdot \nu_{\rho^{-1}(r)}(y)|\nabla\psi(y)|d\sigma_{\rho^{-1}(r)}(y) & = & -\displaystyle{\int_{\Omega_{\rho^{-1}(r)}}}\div(A_{\Omega}\nabla\psi)(y)dy\\
& = & -n\alpha_nr^{n-1}\widehat{\Lambda}(x)F(r)\\
& = & n\alpha_nr^{n-1}\widehat{\Lambda}(x) |\nabla \widetilde{\psi}(x)|.\eaa
\end{equation}
By Cauchy-Schwarz,
\[
\begin{array}{rcl}
\displaystyle i(\rho^{-1}(r))^2 & = & \displaystyle \left(\int_{\Sigma_{\rho^{-1}(r)}} d\sigma_{\rho^{-1}(r)}(y)\right)^2\\
& \leq & \displaystyle\int_{\Sigma_{\rho^{-1}(r)}}A_{\Omega}(y)\nu_{\rho^{-1}(r)}(y)\cdot\nu_{\rho^{-1}(r)}(y) \left\vert \nabla\psi(y)\right\vert d\sigma_{\rho^{-1}(r)}(y)\\
& & \displaystyle \times\int_{\Sigma_{\rho^{-1}(r)}}(A_{\Omega}(y)\nu_{\rho^{-1}(r)}(y)\cdot\nu_{\rho^{-1}(r)}(y))^{-1} \left\vert \nabla\psi(y)\right\vert^{-1} d\sigma_{\rho^{-1}(r)}(y)\\
& \leq & \displaystyle\int_{\Sigma_{\rho^{-1}(r)}}A_{\Omega}(y)\nu_{\rho^{-1}(r)}(y)\cdot\nu_{\rho^{-1}(r)}(y) \left\vert \nabla\psi(y)\right\vert d\sigma_{\rho^{-1}(r)}(y)\\
& & \displaystyle \times\int_{\Sigma_{\rho^{-1}(r)}} \Lambda_{\Omega}(y)^{-1}\left\vert \nabla\psi(y)\right\vert^{-1} d\sigma_{\rho^{-1}(r)}(y)\\
& = & \displaystyle \widehat{\Lambda}(x)^{-1}\times\int_{\Sigma_{\rho^{-1}(r)}}A_{\Omega}(y)\nu_{\rho^{-1}(r)}(y)\cdot\nu_{\rho^{-1}(r)}(y) \left\vert \nabla\psi(y)\right\vert d\sigma_{\rho^{-1}(r)}(y)\\
& & \displaystyle \times \int_{\Sigma_{\rho^{-1}(r)}} \left\vert \nabla\psi(y)\right\vert^{-1}d\sigma_{\rho^{-1}(r)}(y).
\end{array}
\]
In other words,
\[
\begin{array}{l}
\displaystyle\frac{\displaystyle{\int_{\Sigma_{\rho^{-1}(r)}}A_{\Omega}(y)\nu_{\rho^{-1}(r)}(y)\cdot\nu_{\rho^{-1}(r)}(y) \left\vert \nabla\psi(y)\right\vert d\sigma_{\rho^{-1}(r)}(y)}}{\displaystyle{\int_{\Sigma_{\rho^{-1}(r)}} \left\vert \nabla\psi(y)\right\vert^{-1} d\sigma_{\rho^{-1}(r)}(y)}}\\
\qquad\qquad\le\ \hat{\Lambda}(x)^{-1}\times\displaystyle  \left(\frac{\displaystyle{\int_{\Sigma_{\rho^{-1}(r)}}A_{\Omega}(y)\nu_{\rho^{-1}(r)}(y)\cdot\nu_{\rho^{-1}(r)}(y) \left\vert \nabla\psi(y)\right\vert d\sigma_{\rho^{-1}(r)}(y)}}{i(\rho^{-1}(r))}\right)^2\\
\qquad\qquad=\ \left(\displaystyle{\frac{n\alpha_nr^{n-1}}{i(\rho^{-1}(r))}}\right)^2 \widehat{\Lambda}(x)\vert \nabla \widetilde{\psi}(x)\vert^2\eaa
\]
by (\ref{eqsym1}). The isoperimetric inequality (\ref{isop}) ends the proof of (\ref{ineq}).\hfill\fin


\subsection{A pointwise differential inequality for the rearranged data}

In the previous subsection, we could compare the values of $\psi$ and of its symmetrized function $\tilde{\psi}$. Here, we prove a partial differential inequality involving $\psi$ and $\tilde{\psi}$, as well as the rearranged data $\hat{\Lambda}$, $\hat{v}$ and $\hat{V}$.

\begin{pro} \label{comparison}
Let $\omega_0\in\R_+$ and $x\in E\cap\Omega^*$. Then, there exists $y\in \Omega$ such that $\psi(y)=\rho^{-1}(\left\vert x\right\vert)$, that is $y\in\Sigma_{\rho^{-1}(|x|)}$, and
\[
\begin{array}{ll}
& -\displaystyle{\rm{div}}(\widehat{\Lambda}\nabla\widetilde{\psi})(x) +\widehat{v}(x)\cdot\nabla\widetilde{\psi}(x)-\omega_0\vert\nabla\widetilde{\psi}(x)\vert+\widehat{V}(x)\widetilde{\psi}(x)\\
\leq & -\displaystyle{\rm{div}}(A_{\Omega}\nabla\psi)(y)-|\omega(y)|\times|\nabla\psi(y)|-\omega_0|\nabla\psi(y)|+V(y)\psi(y).
\end{array}
\]
Notice that $\widehat{v}(x)\cdot\nabla\widetilde{\psi}(x)=-|\widehat{v}(x)|\times|\nabla\widetilde{\psi}(x)|$.
\end{pro}

\noindent{\bf Proof. }Let $x\in E\cap\Omega^*$, $r=\left\vert x\right\vert$ and $\eta>0$ such that $\overline{S_{r-\eta,r}}\subset E\cap\Omega^*$. As done in the proof of Proposition~\ref{pro35}, the co-area formula and Cauchy-Schwarz inequality yield
\be\label{limbis}\begin{array}{l}
\left(\displaystyle{\mathop{\lim}_{t\to 0^+}} \frac{\displaystyle{\int_{\Omega_{\rho^{-1}(r)}\backslash\Omega_{\rho^{-1}(r-t)}} |\omega(y)|\times|\nabla\psi(y)|\ dy}}{\displaystyle{\left\vert \Omega_{\rho^{-1}(r)}\backslash\Omega_{\rho^{-1}(r-t)}\right\vert}} \right)^2\\
=\ \left(\frac{\displaystyle{\int_{\Sigma_{\rho^{-1}(r)}}\vert \omega(y)|d\sigma_{\rho^{-1}(r)}(y)}}{\displaystyle{\int_{\Sigma_{\rho^{-1}(r)}}|\nabla\psi(y)|^{-1}d\sigma_{\rho^{-1}(r)}(y)}} \right)^2\\
\le\ \frac{\displaystyle{\int_{\Sigma_{\rho^{-1}(r)}}\Lambda_{\Omega}(y)^{-1}\omega(y)^2|\nabla\psi(y)|^{-1}d\sigma_{\rho^{-1}(r)}(y)}}{\displaystyle{\int_{\Sigma_{\rho^{-1}(r)}}|\nabla\psi(y)|^{-1}d\sigma_{\rho^{-1}(r)}(y)}}\times\frac{\displaystyle{\int_{\Sigma_{\rho^{-1}(r)}}\Lambda_{\Omega}(y)|\nabla\psi(y)|d\sigma_{\rho^{-1}(r)}(y)}}{\displaystyle{\int_{\Sigma_{\rho^{-1}(r)}}|\nabla\psi(y)|^{-1}d\sigma_{\rho^{-1}(r)}(y)}}.\eaa
\ee
Using (\ref{ineq}), one obtains
\be\label{ineqbis}\baa{l}
\frac{\displaystyle{\int_{\Sigma_{\rho^{-1}(r)}}\Lambda_{\Omega}(y)|\nabla\psi(y)|d\sigma_{\rho^{-1}(r)}(y)}}{\displaystyle{\int_{\Sigma_{\rho^{-1}(r)}}|\nabla\psi(y)|^{-1}d\sigma_{\rho^{-1}(r)}(y)}}\\
\qquad\qquad\le\ \frac{\displaystyle{\int_{\Sigma_{\rho^{-1}(r)}}A_{\Omega}(y)\nu_{\rho^{-1}(r)}(y)\cdot\nu_{\rho^{-1}(r)}(y)|\nabla\psi(y)|d\sigma_{\rho^{-1}(r)}(y)}}{\displaystyle{\int_{\Sigma_{\rho^{-1}(r)}}|\nabla\psi(y)|^{-1}d\sigma_{\rho^{-1}(r)}(y)}}\\
\qquad\qquad\le\ \hat{\Lambda}(x)|\nabla\tilde{\psi}(x)|^2.\eaa
\ee
Finally, (\ref{limbis}) and (\ref{ineqbis}), together with definitions (\ref{G}-\ref{defhatLambda}) and (\ref{defhatv}), give that
\be\label{lim1}
\displaystyle{\mathop{\lim}_{t\to 0^+}} \frac{\displaystyle{\int_{\Omega_{\rho^{-1}(r)}\backslash\Omega_{\rho^{-1}(r-t)}} \vert\omega(y)|\times|\nabla\psi(y)\vert\ dy}}{\displaystyle{\left\vert \Omega_{\rho^{-1}(r)}\backslash\Omega_{\rho^{-1}(r-t)}\right\vert}}\le|\hat{v}(x)|\times|\nabla\tilde{\psi}(x)|=-\hat{v}(x)\cdot\nabla\tilde{\psi}(x).
\ee
The last equality follows also from (\ref{defpsi}) and Lemma~\ref{lem33}.\par
Remember also from (\ref{gradient}) that
\be\label{lim2}
\mathop{\lim}_{t\rightarrow 0^+} \frac{\displaystyle{\int_{\Omega_{\rho^{-1}(r)}\backslash\Omega_{\rho^{-1}(r-t)}} \left\vert \nabla\psi(y)\right\vert dy}}{\displaystyle{\left\vert \Omega_{\rho^{-1}(r)}\backslash\Omega_{\rho^{-1}(r-t)} \right\vert}} \leq \vert \nabla \widetilde{\psi}(x)\vert=-e_r(x)\cdot \nabla \widetilde{\psi}(x).
\ee\par
As far as $V$ is concerned, for any fixed unit vector $e$ in $\R^n$ and for any $t\in(0,\eta)$, it follows from (\ref{defhatV}) and Lemma~\ref{lem32} that
\[
\begin{array}{l}
\displaystyle \int_{\Omega_{\rho^{-1}(r)}\backslash\Omega_{\rho^{-1}(r-t)}} V(y)\psi(y)dy\\
\qquad\qquad\ge\ \displaystyle \int_{\rho^{-1}(r)}^{\rho^{-1}(r-t)} \left(\int_{\Sigma_a}(-V^-(y))\psi(y)\left\vert \nabla\psi(y)\right\vert^{-1}d\sigma_a(y)\right)da\\
\qquad\qquad=\ -\displaystyle \int_{\rho^{-1}(r)}^{\rho^{-1}(r-t)}a\left(\int_{\Sigma_a}V^-(y)\left\vert\nabla\psi(y)\right\vert^{-1}d\sigma_a(y)\right)da\\
\qquad\qquad=\ -\displaystyle \int_{r-t}^{r} \left(\int_{\Sigma_{\rho^{-1}(s)}} 
V^-(y)\left\vert \nabla\psi(y)\right\vert^{-1}d\sigma_{\rho^{-1}(s)}(y)\right)
\frac{\rho^{-1}(s)ds}{\left\vert \rho^{\prime}(\rho^{-1}(s))\right\vert}\\
\qquad\qquad=\ -\displaystyle n\alpha_n \int_{r-t}^{r} s^{n-1} \rho^{-1}(s) 
\left(\frac{\displaystyle{\int_{\Sigma_{\rho^{-1}(s)}}}V^-(y)\left\vert 
\nabla\psi(y)\right\vert^{-1}d\sigma_{\rho^{-1}(s)}(y)}{\displaystyle{\int_{\Sigma_{\rho^{-1}(s)}}} \left\vert \nabla\psi(y)\right\vert^{-1}d\sigma_{\rho^{-1}(s)}(y)}\right)ds\\
\qquad\qquad=\ \displaystyle n\alpha_n \int_{r-t}^{r} s^{n-1} \rho^{-1}(s) \widehat{V}(se)ds.
\end{array}
\]
Moreover, the radial symmetry of $\widehat{V}$ and $\widetilde{\psi}$ yields
\[
\displaystyle \int_{S_{r-t,r}} \widehat{V}(y)\widetilde{\psi}(y)dy  =  \displaystyle n\alpha_n\int_{r-t}^r s^{n-1}\widehat{V}(se) \widetilde{\psi}(se)ds.
\]
Corollary \ref{cor36} and the facts that $\left\vert \Omega_{\rho^{-1}(r)}\backslash\Omega_{\rho^{-1}(r-t)} \right\vert=\left\vert S_{r-t,r}\right\vert$ and that $\widehat{V}\le 0$ therefore show that
\[
\frac{\displaystyle{\int_{\Omega_{\rho^{-1}(r)}\backslash\Omega_{\rho^{-1}(r-t)}} V(y)\psi(y)dy}}{\displaystyle{\left\vert \Omega_{\rho^{-1}(r)}\backslash\Omega_{\rho^{-1}(r-t)} \right\vert}} \geq \frac{\displaystyle{\int_{S_{r-t,r}} \widehat{V}(y)\widetilde{\psi}(y)dy}}{\displaystyle{\left\vert S_{r-t,r}\right\vert}}.  
\]
Since $\widehat{V}$ and $\widetilde{\psi}$ are continuous in $E$ and radially symmetric, one therefore obtains, together with the co-area formula,
\be\label{lim3}\baa{rcl}
\frac{\displaystyle{\int_{\Sigma_{\rho^{-1}(r)}}V(y)\psi(y)|\nabla\psi(y)|^{-1}d\sigma_{\rho^{-1}(r)}(y)}}{\displaystyle{\int_{\Sigma_{\rho^{-1}(r)}}|\nabla\psi(y)|^{-1}d\sigma_{\rho^{-1}(r)}(y)}} & = & \displaystyle{\mathop{\lim}_{t\rightarrow 0^{+}}} \frac{\displaystyle{\int_{\Omega_{\rho^{-1}(r)}\backslash\Omega_{\rho^{-1}(r-t)}} V(y)\psi(y)dy}}{\displaystyle{\left\vert \Omega_{\rho^{-1}(r)}\backslash\Omega_{\rho^{-1}(r-t)}\right\vert}}\\
& \geq & \widehat{V}(x)\widetilde{\psi}(x).\eaa
\ee\par
Let now $t$ be any real number in $(0,\eta)$. Since $\tilde{\psi}$ (resp. $\hat{\Lambda}$) is radially symmetric, and $C^2$ (resp. $C^1$) on $\overline{S_{r-t,r}} \subset E\cap\Omega^*$, the Green Riemann formula gives
\be\label{green}\begin{array}{lll}
\displaystyle \int_{S_{r-t,r}}\div(\widehat{\Lambda}\nabla\widetilde{\psi})(y)dy & = & \displaystyle\int_{\partial S_{r-t,r}} \widehat{\Lambda}(y)\nabla \widetilde{\psi}(y)\cdot \nu(y) d\sigma(y)\\
& =& \displaystyle n\alpha_n[r^{n-1}G(r)F(r)-(r-t)^{n-1} G(r-t)F(r-t)],\end{array}
\ee
where $d\sigma$ and $\nu$ here denote the superficial measure on $\partial S_{r-t,r}$ and the outward unit normal to $S_{r-t,r}$, and $G$ and $F$ were defined in (\ref{G}) and (\ref{defF}). By definition of $F$, one gets that
$$\int_{S_{r-t,r}}\div(\widehat{\Lambda}\nabla\widetilde{\psi})(y)dy=\int_{\Omega_{\rho^{-1}(r)}\backslash\Omega_{\rho^{-1}(r-t)}}\div(A_{\Omega}\nabla\psi)(y)dy,$$
whence
\be\label{lim4}\baa{l}
\frac{\displaystyle{\int_{\Sigma_{\rho^{-1}(r)}}\div(A_{\Omega}\nabla\psi)(y)|\nabla\psi(y)|^{-1}d\sigma_{\rho^{-1}(r)}(y)}}{\displaystyle{\int_{\Sigma_{\rho^{-1}(r)}}|\nabla\psi(y)|^{-1}d\sigma_{\rho^{-1}(r)}(y)}}\\
\qquad\qquad=\ \displaystyle{\mathop{\lim}_{t\to 0^+}}\displaystyle\frac{\displaystyle{\int_{\Omega_{\rho^{-1}(r)}\backslash\Omega_{\rho^{-1}(r-t)}}}\div(A_{\Omega}\nabla\psi)(y)dy}{\left\vert \Omega_{\rho^{-1}(r)}\backslash\Omega_{\rho^{-1}(r-t)}\right\vert}\ =\ \div(\hat{\Lambda}\nabla\tilde{\psi})(x)\eaa
\ee
since $\left\vert S_{r-t,r}\right\vert=\left\vert \Omega_{\rho^{-1}(r)}\backslash\Omega_{\rho^{-1}(r-t)}\right\vert$.\par
It follows from (\ref{lim1}), (\ref{lim2}), (\ref{lim3}) and (\ref{lim4}) that
\[
\begin{array}{l}
\displaystyle\mathop{\lim}_{t\rightarrow 0^+} \frac{\displaystyle{\int_{\Omega_{\rho^{-1}(r)}\backslash\Omega_{\rho^{-1}(r-t)}} \left[\div(A_{\Omega}\nabla\psi)(y)+|\omega(y)|\times|\nabla\psi(y)|+\omega_0 |\nabla\psi(y)|-V(y)\psi(y)\right]dy}}{\displaystyle{\left\vert  \Omega_{\rho^{-1}(r)}\backslash\Omega_{\rho^{-1}(r-t)} \right\vert}}\\
\\
\quad\le\displaystyle \div(\widehat{\Lambda}\nabla\widetilde{\psi})(x)- \widehat{v}(x)\cdot \nabla \widetilde{\psi}(x)+\omega_0|\nabla\widetilde{\psi}(x)|-\widehat{V}(x)\widetilde{\psi}(x).
\end{array}
\]
To finish the proof, pick any sequence of positive numbers $(\epsilon_l)_{l\in\N}$ such that $\epsilon_l\rightarrow 0$ as $l\to+\infty$. Since $\psi$ is $C^2$ in $\Omega$, since $A_{\Omega}$ is $C^1$ in $\Omega$ (and even in $\overline{\Omega}$) and $\omega$ and $V$ are continuous in $\Omega$ (and even in $\overline{\Omega}$), the previous inequality provides the existence of a sequence of positive numbers $(t_l)_{l\in\N}\in(0,\eta)$ such that $t_l\rightarrow 0$ as $l\to+\infty$, and a sequence of points
$$y_l\in\overline{\Omega_{\rho^{-1}(r)}\backslash\Omega_{\rho^{-1}(r-t_l)}}\subset\overline{\Omega_{\rho^{-1}(r)}}\subset\Omega$$
such that
$$\baa{rl}
& \div(A_{\Omega}\nabla\psi)(y_l)+|\omega(y_l)|\times|\nabla\psi(y_l)|+\omega_0 
|\nabla\psi(y_l)|-V(y_l)\psi(y_l)\\
\leq &
\div(\widehat{\Lambda} \nabla\widetilde{\psi})(x)-\widehat{v}(x)\cdot \nabla \widetilde{\psi}(x)+\omega_0|\nabla\widetilde{\psi}(x)|-\widehat{V}(x)\widetilde{\psi}(x)+\epsilon_l.\eaa$$
Since $\rho^{-1}(r)\le\psi(y_l)\le\rho^{-1}(r-t_l)$ and $\rho^{-1}$ is continuous, the points $y_l$ converge, up to the extraction of some subsequence, to a point $y\in\Sigma_{\rho^{-1}(r)}$ such that
$$\baa{rl}
& \div(A_{\Omega}\nabla\psi)(y)+|\omega(y)|\times|\nabla\psi(y)|+\omega_0 |\nabla\psi(y)|-V(y)\psi(y)\\
\leq & 
\div(\widehat{\Lambda}\nabla\widetilde{\psi})(x)- \widehat{v}(x)\cdot \nabla \widetilde{\psi}(x)+\omega_0|\nabla\widetilde{\psi}(x)|-\widehat{V}(x)\widetilde{\psi}(x),\eaa$$
which is the conclusion of Proposition~\ref{comparison}. \hfill\fin

\begin{cor}\label{corpointwise}
If there are $\omega_0\ge 0$ and $\mu\ge 0$ such that
$$-{\rm{div}}(A_{\Omega}\nabla\psi)(y)-|\omega(y)|\times|\nabla\psi(y)|-\omega_0|\nabla\psi(y)|+V(y)\psi(y)\le\mu\psi(y)\ \hbox{ for all }y\in\Omega,$$
then
$$-\displaystyle{\rm{div}}(\widehat{\Lambda}\nabla\widetilde{\psi})(x) +\widehat{v}(x)\cdot\nabla \widetilde{\psi}(x)-\omega_0|\nabla\widetilde{\psi}(x)|+\widehat{V}(x)\widetilde{\psi}(x)\le\mu\widetilde{\psi}(x)\ \hbox{ for all }x\in E\cap\Omega^*.$$
\end{cor}

\noindent{\bf{Proof.}} It follows immediately from Corollary~\ref{cor36} and Proposition~\ref{comparison}.\hfill\fin


\subsection{An integral inequality for the rearranged data}

A consequence of the pointwise comparisons which were established in the previous subsections is the following integral comparison result:

\begin{pro} \label{integralcompar}
With the previous notations, assume that, for some $(\omega_0,\mu)\in\R^2$,
\begin{equation} \label{ptwiseineq}
-{\rm{div}}(\widehat{\Lambda}\nabla\widetilde{\psi})(x)+\widehat{v}(x)\cdot\nabla\widetilde{\psi}(x)-\omega_0\vert \nabla\widetilde{\psi}(x)\vert+\widehat{V}(x)\widetilde{\psi}(x)\leq\mu\widetilde{\psi}(x)\ \hbox{ for all }x\in E\cap\Omega^*.
\end{equation}
Fix a unit vector $e\in \R^n$. For all $r\in \left[0,R\right]$, define
\be\label{H}
H(r)=\int_0^r \left\vert \widehat{v}(re)\right\vert \widehat{\Lambda}(re)^{-1}dr
\ee
and, for all $x\in\overline{\Omega^*}$, let
\be\label{U}
U(x)=H(\left\vert x\right\vert).
\ee
Then, the following integral inequality is valid:
\be\label{ineqintegral}
\int_{\Omega^{\ast}} \left[\widehat{\Lambda}(x)\vert \nabla\widetilde{\psi}(x)\vert^2-\omega_0\vert \nabla\widetilde{\psi}(x)\vert \widetilde{\psi}(x)+\widehat{V}(x)\widetilde{\psi}(x)^2\right]e^{-U(x)}dx\leq \mu\int_{\Omega^{\ast}} \widetilde{\psi}(x)^2e^{-U(x)}dx.
\ee
\end{pro}

\noindent{\bf{Proof.}} Note first that, since $|\widehat{v}|$ and $\widehat{\Lambda}$ are radially symmetric and since $|\widehat{v}|\in L^{\infty}(\Omega^{\ast})$ and $\widehat{\Lambda}$ satisfies (\ref{infiniLambda}), the function $H$ is well-defined and continuous in $\left[0,R\right]$. Furthermore, its definition is independent from the choice of $e$. The radially symmetric function $U$ is then continuous in $\overline{\Omega^*}$ and, since the radially symmetric functions $\hat{v}=|\hat{v}|e_r$ and $1/\hat{\Lambda}$ are (at least) continuous in $E$, the function $U$ is of class $C^1$ in $E$ and
\be\label{gradientU}
\nabla U(x)=\hat{\Lambda}(x)^{-1}\hat{v}(x)\hbox{ for all }x\in E.
\ee
Observe also that the integrals in (\ref{ineqintegral}) are all well-defined since $\tilde{\psi}\in H^1_0(\Omega^*)$ and $\hat{\Lambda}$, $\hat{V}$, $U\in L^{\infty}(\Omega^*)$ (even, $U\in C(\overline{\Omega^*})$).\par
Now, recall that the set of critical values of $\psi$ is $Z=\left\{a_1,\ldots,a_m\right\}$ with
$$0<a_1<\cdots<a_m=M$$
and remember that the function $\rho$ defined in Subsection~\ref{sec21} is continuous and decreasing from $[0,M]$ onto $[0,R]$, from Lemma~\ref{lem31}. Fix $j\in \left\{1,\ldots,m-1\right\}$ and $r,r^{\prime}$ such that
$$0\le\rho(a_{j+1})<r<r^{\prime}<\rho(a_j)<R.$$
Multiplying (\ref{ptwiseineq}) by the nonnegative function $\widetilde{\psi}e^{-U}$ and integrating over $S_{r,r^{\prime}}$ yields
\be\label{int1}
\begin{array}{l}
\displaystyle \int_{S_{r,r'}} \left[-\div(\widehat{\Lambda}\nabla\widetilde{\psi})(x)+\widehat{v}(x)\cdot\nabla\widetilde{\psi}(x)-\omega_0\vert \nabla\widetilde{\psi}(x)\vert+\widehat{V}(x)\widetilde{\psi}(x)\right]\widetilde{\psi}(x)e^{-U(x)}dx\\
\qquad \leq\ 
\displaystyle \mu\int_{S_{r,r'}} \widetilde{\psi}(x)^2e^{-U(x)}dx.
\end{array}
\ee
Notice that all integrals above are well-defined since $\tilde{\psi}$ in $C^2$ in $E\cap\Omega^*$, $\hat{\Lambda}$ is $C^1$ in $E\cap\Omega^*$, $\hat{v}$, $\hat{V}$ are continuous in $E$, $U$ is continuous in $\overline{\Omega^*}$ and $\overline{S_{r,r'}}\subset E\cap\Omega^*$. Furthermore, as in (\ref{green}), the Green-Riemann formula yields
$$\begin{array}{l}
\displaystyle \int_{S_{r,r'}} -\div(\widehat{\Lambda}\nabla\widetilde{\psi})(x)\ \widetilde{\psi}(x)\ e^{-U(x)} dx\\
\qquad\qquad =\ \displaystyle \int_{S_{r,r'}} \widehat{\Lambda}(x)\vert \nabla\widetilde{\psi}(x)\vert^2 e^{-U(x)}dx-  \int_{S_{r,r'}} \widehat{\Lambda}(x) \widetilde{\psi}(x) \nabla\widetilde{\psi}(x)\cdot \nabla U(x) e^{-U(x)}dx \\
\qquad\qquad\ \ \ \ \displaystyle - n\alpha_n (r^{\prime})^{n-1}G(r^{\prime})F(r^{\prime}) \widetilde{\psi}(r^{\prime}e)e^{-H(r^{\prime})}+ n\alpha_n r^{n-1}G(r)F(r) \widetilde{\psi}(re)e^{-H(r)}.
\end{array}$$
By (\ref{gradientU}), it follows then that
\be\label{int2}\begin{array}{l}
\displaystyle \int_{S_{r,r'}} \left[-\div(\widehat{\Lambda}\nabla\widetilde{\psi})(x)+\widehat{v}(x)\cdot\nabla\widetilde{\psi}(x)-\omega_0\vert \nabla\widetilde{\psi}(x)\vert+\widehat{V}(x)\widetilde{\psi}(x)\right]\widetilde{\psi}(x)e^{-U(x)}dx  \\
\qquad\qquad\displaystyle =\int_{S_{r,r'}}\left[\widehat{\Lambda}(x)\vert \nabla\widetilde{\psi}(x)\vert^2-\omega_0\vert \nabla\widetilde{\psi}(x)\vert\tilde{\psi}(x)+\widehat{V}(x)\widetilde{\psi}(x)^2\right]e^{-U(x)}dx\\
\qquad\qquad\ \ \ \ \displaystyle -n\alpha_n (r^{\prime})^{n-1}G(r^{\prime})F(r^{\prime}) \widetilde{\psi}(r^{\prime}e)e^{-H(r^{\prime})} + n\alpha_n r^{n-1}G(r)F(r) \widetilde{\psi}(re)e^{-H(r)}.\end{array}
\ee\par
On the other hand, for all $s\in \rho(Y)$,
\[
n\alpha_ns^{n-1}F(s)G(s)=\int_{\Omega_{\rho^{-1}(s)}} \div(A_{\Omega}\nabla\psi)(x)dx,
\]
by (\ref{defF}). The function
$$s\mapsto I(s)=n\alpha_ns^{n-1}F(s)G(s),$$
which was a priori defined only in $\rho(Y)$, can then be extended continuously in $[0,R]$ from the results in Lemma~\ref{lem31} and since $\div(A_{\Omega}\nabla\psi)=-f$ is bounded in $\Omega$. The continuous extension of $I$ in $[0,R]$ is still called $I$. Passing to the limit as $r\to\rho(a_{j+1})^+$ and $r'\to\rho(a_j)^-$ in (\ref{int1}) and (\ref{int2}) yields, for each $j\in\{1,\ldots,m-1\}$,
\begin{equation} \label{sum1}\begin{array}{ll}
& \displaystyle\int_{S_{\rho(a_{j+1}),\rho(a_j)}}\left[\widehat{\Lambda}(x)\vert \nabla\widetilde{\psi}(x)\vert^2-\omega_0\vert \nabla\widetilde{\psi}(x)\vert \widetilde{\psi}(x)+\widehat{V}(x)\widetilde{\psi}(x)^2\right]e^{-U(x)}dx\\
& -I(\rho(a_j))\ \tilde{\psi}(\rho(a_j)e)\ e^{-H(\rho(a_j))}+I(\rho(a_{j+1}))\ \tilde{\psi}(\rho(a_{j+1})e)\ e^{-H(\rho(a_{j+1}))}\\
\le & \mu\displaystyle\int_{S_{\rho(a_{j+1}),\rho(a_j)}}\tilde{\psi}(x)^2e^{-U(x)}dx.\eaa
\ee
Once again, all integrals above are well-defined. Arguing similarly in the spherical shell $S_{\rho(a_1),R}$ and since $\psi(Re)=0$, one obtains
\begin{equation} \label{sum2}\begin{array}{ll}
& \displaystyle\int_{S_{\rho(a_1),R}}\left[\widehat{\Lambda}(x)\vert \nabla\widetilde{\psi}(x)\vert^2-\omega_0\vert \nabla\widetilde{\psi}(x)\vert \widetilde{\psi}(x)+\widehat{V}(x)\widetilde{\psi}(x)^2\right]e^{-U(x)}dx\\
& +I(\rho(a_1))\ \tilde{\psi}(\rho(a_1)e)\ e^{-H(\rho(a_1))}\\
\le & \mu\displaystyle\int_{S_{\rho(a_1),R}}\tilde{\psi}(x)^2e^{-U(x)}dx.\eaa
\ee
Summing up (\ref{sum1}) for all $1\leq j\leq m-1$ and (\ref{sum2}) and using the fact that $I(\rho(a_m))=I(0)=0$ yield (\ref{ineqintegral}).\hfill\fin

\begin{cor}\label{corintegral}
If there are $\omega_0\ge 0$ and $\mu\ge 0$ such that
$$-{\rm{div}}(A_{\Omega}\nabla\psi)(y)-|\omega(y)|\times|\nabla\psi(y)|-\omega_0|\nabla\psi(y)|+V(y)\psi(y)\le\mu\psi(y)\ \hbox{ for all }y\in\Omega,$$
then, under the notation of Proposition~\ref{integralcompar},
$$\int_{\Omega^{\ast}} \left[\widehat{\Lambda}(x)\vert \nabla\widetilde{\psi}(x)\vert^2-\omega_0\vert \nabla\widetilde{\psi}(x)\vert\widetilde{\psi}(x)+\widehat{V}(x)\widetilde{\psi}(x)^2\right]e^{-U(x)}dx\leq \mu\int_{\Omega^{\ast}} \widetilde{\psi}(x)^2e^{-U(x)}dx.$$
\end{cor}

\noindent{\bf{Proof.}} It follows immediately from Corollary~\ref{corpointwise} and Proposition~\ref{integralcompar}.\hfill\fin\break

We complete this section by two remarks which proceed from the previous results and provide comparisons between some norms of the function $\psi$ and its symmetrization $\tilde{\psi}$.

\begin{rem}{\rm The calculations of the previous subsections (see in particular the proof of Proposition~\ref{comparison}) and Corollary~\ref{cor36} imply that, for any nondecreasing function $\Theta:[0,+\infty)\to[0,+\infty)$,
$$\int_{\Omega_{\rho^{-1}(s)}\backslash\Omega_{\rho^{-1}(r)}}\Theta(\psi(y))dy=\int_{S_{r,s}}\Theta(\rho^{-1}(|x|))dx\le\int_{S_{r,s}}\Theta(\tilde{\psi}(x))dx$$
for all $0<r<s\le R$ such that $[r,s]\subset\rho(Y)$, and then for all $0\le r<s\le R$ from Lebesgue's dominated convergence theorem. In particular,
$$\int_{\Omega}(\psi(y))^pdy\le\int_{\Omega^*}(\tilde{\psi}(x))^pdx$$
for all $0\le p<+\infty$. Remember also (as an immediate consequence of Corollary~\ref{cor36}) that
$$\max_{\overline{\Omega^*}}\tilde{\psi}=\tilde{\psi}(0)\ge\rho^{-1}(0)=\max_{\overline{\Omega}}\psi.$$}
\end{rem}

\begin{rem}\label{remgradients}{\rm For any $0<r<s\le R$ such that $[r,s]\subset\rho(Y)$, it follows from the co-area formula and a change of variables that
$$\baa{l}
\displaystyle{\int_{\Omega_{\rho^{-1}(s)}\backslash\Omega_{\rho^{-1}(r)}}}A_{\Omega}(y)\nabla\psi(y)\cdot\nabla\psi(y)dy\\
\qquad=n\alpha_n\displaystyle{\int_r^s}t^{n-1}\times\left(\displaystyle{\frac{\displaystyle{\int_{\Sigma_{\rho^{-1}(t)}}}A_{\Omega}(y)\nu_{\rho^{-1}(t)}(y)\cdot\nu_{\rho^{-1}(t)}(y)|\nabla\psi(y)|d\sigma_{\rho^{-1}(t)}(y)}{\displaystyle{\int_{\Sigma_{\rho^{-1}(t)}}}|\nabla\psi(y)|^{-1}d\sigma_{\rho^{-1}(t)}(y)}}\right)dt\\
\qquad\le\displaystyle{\int_{S_{r,s}}}\hat{\Lambda}(x)|\nabla\tilde{\psi}(x)|^2dx,\eaa$$
where the last inequality is due to (\ref{ineq}). As usual, Lebesgue's dominated convergence theorem then yields
$$\int_{\Omega_{\rho^{-1}(s)}\backslash\Omega_{\rho^{-1}(r)}}A_{\Omega}(y)\nabla\psi(y)\cdot\nabla\psi(y)dy\le\int_{S_{r,s}}\hat{\Lambda}(x)|\nabla\tilde{\psi}(x)|^2dx$$
for all $0\le r<s\le R$, whence
$$\int_{\Omega}A_{\Omega}(y)\nabla\psi(y)\cdot\nabla\psi(y)dy\le\int_{\Omega^*}\hat{\Lambda}(x)|\nabla\tilde{\psi}(x)|^2dx.$$
As a consequence,
$$\|\nabla\psi\|_{L^2(\Omega,\R^n)}\le\sqrt{\frac{M_{\Lambda}}{m_{\Lambda}}}\times\|\nabla\tilde{\psi}\|_{L^2(\Omega^*,\R^n)}$$
from (\ref{ALambdaOmega}) and (\ref{infiniLambda}), where $M_{\Lambda}=\max_{\overline{\Omega}}\Lambda_{\Omega}$ and $m_{\Lambda}=\min_{\overline{\Omega}}\Lambda_{\Omega}$. In particular, $\|\nabla\psi\|_{L^2(\Omega,\R^n)}\le\|\nabla\tilde{\psi}\|_{L^2(\Omega^*,\R^n)}$
if $\Lambda$ is constant in $\Omega$.\par
With the same notations as above, it follows from (\ref{gradient}) that
$$\baa{rcl}
\displaystyle{\int_{\Omega_{\rho^{-1}(s)}\backslash\Omega_{\rho^{-1}(r)}}}|\nabla\psi(y)|dy & = & n\alpha_n\displaystyle{\int_r^s}t^{n-1}\times\left(\displaystyle{\frac{\displaystyle{\int_{\Sigma_{\rho^{-1}(t)}}}d\sigma_{\rho^{-1}(t)}(y)}{\displaystyle{\int_{\Sigma_{\rho^{-1}(t)}}}|\nabla\psi(y)|^{-1}d\sigma_{\rho^{-1}(t)}(y)}}\right)dt\\
& \le & \displaystyle{\int_{S_{r,s}}}|\nabla\tilde{\psi}(x)|dx.\eaa$$
The inequality then holds for all $0\le r<s\le R$, whence
$$\|\nabla\psi\|_{L^1(\Omega,\R^n)}\le\|\nabla\tilde{\psi}\|_{L^1(\Omega^*,\R^n)}.$$}
\end{rem}


\SE{Improved inequalities when $\Omega$ is not a ball} \label{equality}

Throughout this section, we assume that $\Omega\in {\mathcal C}$ is not a ball. Fix real numbers $\alpha\in(0,1)$, $N>0$ and $\delta>0$. Denote by $E_{\alpha,N,\delta}(\Omega)$ the set of all functions $\psi\in C^{1,\alpha}(\overline{\Omega})$, positive in $\Omega$, vanishing on $\partial\Omega$, such that
\begin{equation} \label{emNalpha}
\left\Vert \psi\right\Vert_{C^{1,\alpha}(\overline{\Omega})}\leq N\mbox{ and }\psi(x)\geq\delta\times d(x,\partial\Omega)\hbox{ for all }x\in\overline{\Omega},
\end{equation}
where $d(x,\partial\Omega)=\min_{y\in\partial\Omega}|x-y|$ denotes the distance between $x$ and $\partial\Omega$, and we set
$$\left\Vert \psi\right\Vert_{C^{1,\alpha}(\overline{\Omega})}=\|\psi\|_{L^{\infty}(\Omega)}+\|\nabla\psi\|_{L^{\infty}(\Omega,\R^n)}+\sup_{z\neq z'\in\Omega}\frac{|\nabla\psi(z)-\nabla\psi(z')|}{|z-z'|^{\alpha}}.$$
Notice that, for each $\psi\in E_{\alpha,N,\delta}(\Omega)$, one has
$$\nabla\psi(y)\cdot\nu(y)=-|\nabla\psi(y)|\le-\delta\hbox{ for all }y\in\partial\Omega,$$
where $\nu$ denotes the outward unit normal on $\partial\Omega$.\par
Our goal here is to prove stronger versions of Corollary~\ref{cor36} and Corollary~\ref{corintegral}, using the fact that $\Omega$ is not a ball. In the sequel, unless explicitly mentioned, all the constants only depend on some of the data $\Omega$, $n$, $\alpha$, $N$ and $\delta$.\par
Denote again by $R$ the radius of $\Omega^{\ast}$, so that $\Omega^*=B_R$. First, the isoperimetric inequality yields the existence of $\beta=\beta(\Omega,n)>0$ such that
\be\label{areabeta}
\hbox{area}(\partial\Omega)=\int_{\partial\Omega}d\sigma_{\partial\Omega}(y)\ge(1+\beta)n\alpha_nR^{n-1},
\ee
where the left-hand side is the $(n-1)$-dimensional measure of $\partial\Omega$.\par
For all $\gamma>0$, define
$$U_{\gamma}=\{x\in\overline{\Omega},\ d(x,\partial\Omega)\le\gamma\}.$$
Since $\partial\Omega$ is of class $C^2$, there exists $\gamma_1=\gamma_1(\Omega)>0$ such that $\overline{\Omega\backslash U_{\gamma_1}}\neq\emptyset$ and the segments $[y,y-\gamma_1\nu(y)]$ are included in
$\overline{\Omega}$ and pairwise disjoint when $y$ describes $\partial\Omega$. Thus, for all $\gamma\in(0,\gamma_1]$, the segments $[y,y-\gamma\nu(y)]$ (resp. $(y,y-\gamma\nu(y)]$) describe the set $U_{\gamma}$ (resp. $\{x\in\Omega$, $d(x,\partial\Omega)\le\gamma\}$) as $y$ describes $\partial\Omega$.

\begin{lem}\label{lem39} There exists a constant
$$\gamma_2=\gamma_2(\Omega,\alpha,N,\delta)\in(0,\gamma_1]$$ such that, for all $\psi\in E_{\alpha,N,\delta}(\Omega)$, one has:
$$|\nabla\psi|\ge\frac{\delta}{2}\hbox{ in }U_{\gamma_2},$$
$$\nabla\psi(y-r\nu(y))\cdot\nu(y)\le-\frac{\delta}{2}\hbox{ for all }y\in\partial\Omega\hbox{ and }r\in[0,\gamma_2],$$
and
$$\psi\ge\frac{\gamma\delta}{2}\hbox{ in }\overline{\Omega\backslash U_{\gamma}}\hbox{ for all }\gamma\in[0,\gamma_2].$$
\end{lem}

\noindent{\bf{Proof.}} Assume that the conclusion of the lemma does not hold. Then there exists a sequence of positive numbers $(\gamma^l)_{l\in\N}\to
0$ and a sequence of functions $(\psi^l)_{l\in \N}\in E_{\alpha,N,\delta}(\Omega)$ such that one of the three following cases occur:\par
1) either for each $l\in\N$, there is a point $x^l\in U_{\gamma^l}$ such
that $|\nabla\psi^l(x^l)|<\delta/2$;\par
2) or for each $l\in\N$, there are a point $y^l\in\partial\Omega$ and a number $r^l\in[0,\gamma^l]$ such that $\nabla\psi^l(y^l-r^l\nu(y^l))\cdot\nu(y^l)>-\delta/2$;\par
3) or for each $l\in\N$, there is a point
$x^l\in\overline{\Omega\backslash U_{\gamma^l}}$ such that $\psi^l(x^l)<\gamma^l\delta/2$.\par
Observe first that, by Ascoli theorem, up to the extraction of a subsequence, 
there exists a function $\psi\in E_{\alpha,N,\delta}(\Omega)$ such that
$$\psi^l\rightarrow \psi\hbox{ in }C^{1}(\overline{\Omega})\hbox{ as }l\to+\infty.$$
If case~1) occurs, then, up to some subsequence, one can assume without loss of generality that $x^l\rightarrow x\in \partial\Omega$ as $l\to+\infty$, and one obtains $|\nabla\psi(x)|\le\delta/2$, which is impossible since $\psi\in E_{\alpha,N,\delta}(\Omega)$ and $\delta>0$. Similarly, if case~2) occurs, $y^l\rightarrow y\in \partial\Omega$ as $l\to+\infty$ up to a subsequence, and one has $-|\nabla\psi(y)|=\nabla\psi(y)\cdot \nu(y)\geq-\delta/2$, which is also impossible.\par
Therefore, only case~3) can occur. For each $l\in\N$, let $y^l\in\partial\Omega$ be such that
$$d^l:=|x^l-y^l|=d(x^l,\partial\Omega)\ge\gamma^l>0.$$
Up to extraction of some subsequence, one has
$$x^l\to x\in\overline{\Omega}\hbox{ as }l\to+\infty\hbox{ and }\psi(x)\le 0$$
by passing to the limit as $l\to+\infty$ in the inequality $\psi^l(x^l)<\gamma^l\delta/2$. Since
$$\psi(x)\ge\delta\ d(x,\partial\Omega)>0\hbox{ for all }x\in\Omega,$$
it follows that $x\in\partial\Omega$, whence $|x^l-y^l|\to 0$ and $y^l\to x$ as $l\to+\infty$. On the one hand, the mean value theorem implies that
$$\frac{\psi^l(x^l)-\psi^l(y^l)}{|x^l-y^l|}=\nabla\psi^l(z^l)\cdot\frac{x^l-y^l}{|x^l-y^l|}\to-\nabla\psi(x)\cdot\nu(x)=|\nabla\psi(x)|\
\hbox{ as }l\to+\infty,$$
where $z^l$ is a point lying on the segment between $x^l$ and $y^l$ (whence, $z^l\to x$ as $l\to+\infty$). On the other hand,
since $\psi^l=0$ on $\partial\Omega$,
$$\frac{\psi^l(x^l)-\psi^l(y^l)}{|x^l-y^l|}=\frac{\psi^l(x^l)}{d^l}<\frac{\gamma^l\delta}{2\gamma^l}=\frac{\delta}{2}.$$
Hence, $|\nabla\psi(x)|\le\delta/2$ at the limit as $l\to+\infty$, which contradicts the positivity of $\delta$ and the fact that $\psi\in E_{\alpha,N,\delta}(\Omega)$.\par
Therefore, case~3) is ruled out too and the proof of Lemma~\ref{lem39} is complete.\hfill\fin\break

In the sequel, we assume that $\psi\in E_{\alpha,N,\delta}(\Omega)$ and is analytic in $\Omega$. The data $A_{\Omega}$, $\Lambda_{\Omega}$, $\omega$ and $V$ are as in Section~\ref{rearrangement}. We assume that
$$\div(A_{\Omega}\nabla\psi)=-f\hbox{ in }\Omega,$$
where $f$ is a non zero polynomial, and we use the same sets $Z$, $Y$, $E$, $\Omega_{a}$, $\Sigma_{a}$ and the same functions $\rho$, $\widetilde{\psi}$, $\hat{\Lambda}$, $\hat{v}$, $\hat{V}$ and $U$ as in Section~\ref{rearrangement}.

\begin{lem}\label{lem310} Assume that $\Omega\in{\mathcal{C}}$ is not a ball. Then there exists a constant $a_0=a_0(\Omega,n,\alpha,N,\delta)>0$ only depending on $\Omega$, $n$, $\alpha$, $N$ and $\delta$ such that $[0,a_0]\subset Y$ and 
$$i(a)=\int_{\Sigma_{a}}d\sigma_{a}(y)={\rm{area}}(\Sigma_{a})\ge\left(1+\frac{\beta}{2}\right)n\alpha_nR^{n-1}$$
for all $a\in[0,a_0]$, where $\beta=\beta(\Omega,n)>0$ was given in (\ref{areabeta}).
\end{lem}

\noindent{\bf{Proof.}} Let $\gamma_2=\gamma_2(\Omega,\alpha,N,\delta)>0$ be as in Lemma~\ref{lem39}. Since $0<\gamma_2\le\gamma_1$, it follows that $\overline{\Omega\backslash U_{\gamma_2}}\neq\emptyset$, and
$$\psi\ge\frac{\gamma_2\delta}{2}>0\hbox{ in }\overline{\Omega\backslash U_{\gamma_2}}$$
from Lemma~\ref{lem39}. Therefore,
$$M=\max_{\overline{\Omega}}\psi\ge\frac{\gamma_2\delta}{2}$$
and
\begin{equation}\label{sigmaa}
\Sigma_{a}=\{x\in\overline{\Omega},\ \psi(x)=a\}\subset U_{\gamma_2}\hbox{ for all }a\in\left[0,\frac{\gamma_2\delta}{4}\right].
\end{equation}
Call
$$a'_0=a'_0(\Omega,\alpha,N,\delta)=\frac{\gamma_2\delta}{4}>0.$$
From Lemma~\ref{lem39} and the assumptions made on $\psi$ and $\partial\Omega$, one has that $|\nabla\psi|\neq 0$ everywhere on the $C^2$ hypersurface $\Sigma_{a}$ for all $a\in[0,a'_0]$. Thus,
\be\label{a'0}
[0,a'_0]\subset Y.
\ee\par
On the other hand, for all $y\in\partial\Omega$, the segment $[y,y-\gamma_2\nu(y)]$ is included in $\overline{\Omega}$
and there exists $\theta\in[0,1]$ such that
$$\psi(y-\gamma_2\nu(y))=\underbrace{\psi(y)}_{=0}-\gamma_2\nu(y)\cdot\nabla\psi(y-\theta\gamma_2\nu(y))\ge \frac{\gamma_2\delta}{2},$$
again from Lemma~\ref{lem39}. Actually, more precisely, for each $y\in\partial\Omega$, the function
$$\kappa\ :\ [0,\gamma_2]\to\R,\ s\mapsto\psi(y-s\nu(y))$$
is differentiable and $\kappa'(s)\ge\delta/2$ for all $s\in[0,\gamma_2]$. It follows that, for all $a\in[0,a'_0]$ and $y\in\partial\Omega$, there exists a unique
point
$$\phi_a(y)\in[y,y-\gamma_2\nu(y)]\cap\Sigma_a.$$
Moreover, for such a choice of $a$, the map $\phi_a$ is one-to-one since the segments $[y,y-\gamma_2\nu(y)]$ are pairwise
disjoint (and describe $U_{\gamma_2}$) when $y$ describes $\partial\Omega$ (because $0<\gamma_2\le\gamma_1$). Lastly,
\be\label{onto}
\Sigma_a=\{\phi_a(y),\ y\in\partial\Omega\}
\ee
from (\ref{sigmaa}).\par
Let us now prove that the area of $\Sigma_a$ is close to that of $\partial\Omega$ for $a\ge 0$ small enough. To do so, call
$${\mathcal{B}}=\{x'=(x_1,\ldots,x_{n-1}),\ |x'|<1\}$$ 
and represent $\partial\Omega$ by a finite number of $C^2$ maps $y^1,\ldots,y^p$ (for some $p=p(\Omega)\in\N^*$) defined in $\overline{\mathcal{B}}$, depending only on $\Omega$, and for which
$$\partial_1y^j(x')\times\cdots\times\partial_{n-1}y^j(x')\neq 0\ \hbox{ for all }1\le j\le p\hbox{ and }x'\in\overline{\mathcal{B}}.$$
Here,
$$\partial_iy^j(x')=(\partial_{x_i}y^j_1(x'),\ldots,\partial_{x_i}y^j_n(x'))\in\R^n$$
for each $1\le i\le n-1$, $1\le j\le p$ and $x'\in\overline{\mathcal{B}}$, where $y^j(x')=(y^j_1(x'),\ldots,y^j_n(x'))$. The maps $y^j$ are chosen so that
$$\partial\Omega=\{y^j(x'),\ 1\le j\le p,\ x'\in\overline{\mathcal{B}}\}.$$\par
For each $a\in[0,a'_0]$ and for each $1\le j\le p$, there exists then a map $t^j_{a}:\overline{\mathcal{B}}\to[0,\gamma_2]$ such that
\be\label{tja}
\psi(y^j(x')-t^j_{a}(x')\nu(y^j(x')))=a\hbox{ for all }x'\in\overline{\mathcal{B}},
\ee
and
$$\Sigma_a=\{y^j(x')-t^j_{a}(x')\nu(y^j(x')),\ 1\le j\le p,\ x'\in\overline{\mathcal{B}}\}.$$
Namely,
$$y^j(x')-t^j_a(x')\nu(y^j(x'))=\phi_a(y^j(x'))\hbox{ for all }1\le j\le p\hbox{ and }x'\in\overline{\mathcal{B}}.$$
From the arguments above, each real number $t^j_{a}(x')$ is then uniquely determined, and $t^j_{0}(x')=0$. Since the functions $\psi$, $y^j$ and $\nu\circ y^j$ (for all $1\le j\le p$) are (at least) of class $C^1$ (respectively in $\overline{\Omega}$, $\overline{\mathcal{B}}$ and $\overline{\mathcal{B}}$), it follows from the implicit function theorem and Lemma~\ref{lem39} that the functions $t_{a}^j$ (for all $a\in[0,a'_0]$ and $1\le j\le p$) are of class $C^1(\overline{\mathcal{B}})$ and that the functions
$$h^j_{x'}:[0,a_1]\to[0,\gamma_2],\ a\mapsto t_{a}^j(x')$$
(for all $1\le j\le p$ and $x'\in\overline{\mathcal{B}}$) are of class $C^1([0,a'_0])$. From the chain rule applied to (\ref{tja}), it is straightforward to check that, for all $a\in[0,a'_0]$, $1\le j\le p$ and $x'\in\overline{\mathcal{B}}$,
$$(h^j_{x'})'(a)=\frac{-1}{\nu(y^j(x'))\ \cdot\ \nabla\psi(y^j(x')-t^j_{a}(x')\nu(y^j(x')))}\in(0,2\delta^{-1}]\ \hbox{ (from Lemma~\ref{lem39})},$$
whence
\be\label{tja2}
0\le t^j_{a}(x')=h^j_{x'}(a)\le 2\delta^{-1}a
\ee
because $h^j_{x'}(0)=t^j_{0}(x')=0$. Similarly,
\be\label{chain}
\partial_{x_i}t^j_{a}(x')=\frac{[\partial_iy^j(x')-t^j_{a}(x')\partial_i(\nu\circ y^j)(x')]\ \cdot\ \nabla\psi(y^j(x')-t^j_{a}(x')\nu(y^j(x')))}{\nu(y^j(x'))\ \cdot\ \nabla\psi(y^j(x')-t^j_{a}(x')\nu(y^j(x')))}
\ee
for all $a\in[0,a'_0]$, $1\le i\le n-1$, $1\le j\le p$ and $x'\in\overline{\mathcal{B}}$. For all $1\le j\le p$ and $x'\in\overline{\mathcal{B}}$, one has $\psi(y^j(x'))=0$, whence $\partial_iy^j(x')\cdot\nabla\psi(y^j(x'))=0$ (for all $1\le i\le n-1$). As a consequence, 
$$|\partial_iy^j(x')\cdot\nabla\psi(y^j(x')-t^j_{a}(x')\nu(y^j(x')))|\le C_1\times(t^j_{a}(x'))^{\alpha}$$
for all $a\in[0,a'_0]$, $1\le i\le n-1$, $1\le j\le p$, $x'\in\overline{\mathcal{B}}$, and for some constant $C_1$ defined by
$$C_1=\max_{1\le i'\le n-1,\ 1\le j'\le p,\ \xi\in\overline{\mathcal{B}}}|\partial_{i'}y^{j'}(\xi)|\ \times\ \sup_{z\neq z'\in\Omega}\frac{|\nabla\psi(z)-\nabla\psi(z')|}{|z-z'|^{\alpha}}\ <\ +\infty.$$
Observe that $C_1=C_1(\Omega,N)$ only depends on $\Omega$ and $N$ (remember that $\|\psi\|_{C^{1,\alpha}(\overline{\Omega})}\le N$). Call now
$$C_2=C_2(\Omega,N)=\max_{1\le i'\le n-1,\ 1\le j'\le p,\ \xi\in\overline{\mathcal{B}}}|\partial_{i'}(\nu\circ y^{j'})(\xi))|\ \times\ \sup_{z\in\overline{\Omega}}|\nabla\psi(z)|\ <\ +\infty,$$
which also depends on $\Omega$ and $N$ only. Together with (\ref{chain}) and Lemma~\ref{lem39}, the above arguments imply that
$$\baa{rcl}
|\partial_{x_i}t^j_{a}(x')| & \le & 2\delta^{-1}[C_1\times(t^j_{a}(x'))^{\alpha}+C_2\times t^j_{a}(x')]\\
& \le & 2\delta^{-1}[C_1\times(2\delta^{-1}a)^{\alpha}+2C_2\times\delta^{-1}a]\ \hbox{ from (\ref{tja2})}\eaa$$
for all $a\in[0,a'_0]$, $1\le i\le n-1$, $1\le j\le p$ and $x'\in\overline{\mathcal{B}}$.\par
It follows that, for all $\eta>0$, there exists $a''_0=a''_0(\Omega,\alpha,N,\delta,\eta)\in(0,a'_0]$ such that, for all $a\in[0,a''_0]$,
\be\label{suptja}
\sup_{1\le i\le n-1,\ 1\le j\le p,\ x'\in\overline{\mathcal{B}}}|\partial_{x_i}(t^j_{a}\ \nu\circ y^j)(x')|\ \leq \eta.
\ee\par
Finally, there are some open sets $U^1,\ldots,U^p\subset{\mathcal{B}}$ such that
$$\hbox{area}(\partial\Omega)=\sum_{j=1}^p\int_{U^j}|\partial_1y^j(x')\times\cdots\times\partial_{n-1}y^j(x')|\ dx',$$
where the sets $\{y^j(x'),\ x'\in U^j\}$ for $j=1,\ldots,p$ are pairwise disjoint and, for any $\epsilon>0$, there are some measurable sets $V^1,\ldots,V^p\subset\overline{\mathcal{B}}$ such that
$$\partial\Omega=\{y^j(x'),\ 1\le j\le p,\ x'\in V^j\},$$
and $V^j\supset U^j$, $\displaystyle{\int_{\mathcal{B}}}{\bf{1}}_{V^j\backslash U^j}(x')dx'\le\epsilon$ for all $1\le j\le p$. Since all functions $y^j$ and $t^j_{a}\ \nu\circ y^j$ (for all $a\in[0,a'_0]$ and $1\le j\le p$) are of class $C^1(\overline{\mathcal{B}})$, since each function $\phi_a$ is one-to-one and since (\ref{onto}) holds, it follows that
$$\hbox{area}(\Sigma_a)=\sum_{j=1}^p\int_{U^j}|\partial_1(y^j-t^j_{a}\nu\circ y^j)(x')\times\cdots\times\partial_{n-1}(y^j-t^j_{a}\nu\circ y^j)(x')|\ dx'$$
for all $a\in[0,a'_0]$.\par
One concludes from (\ref{suptja}) that, since $\beta=\beta(\Omega,n)$ in (\ref{areabeta}) is positive, there exists a positive constant $a_0=a_0(\Omega,n,\alpha,N,\delta)\in(0,a'_0]$ which only depends on $\Omega$, $n$, $\alpha$, $N$ and $\delta$, and which is such that
$$|\hbox{area}(\Sigma_a)-\hbox{area}(\partial\Omega)|\leq \frac{\beta}2n\alpha_nR^{n-1}\hbox{ for all }a\in[0,a_0].$$
As a consequence, one has
$$i(a)=\hbox{area}(\Sigma_a)\ge\left(1+\frac{\beta}{2}\right)n\alpha_nR^{n-1}$$
for all $a\in[0,a_0]$. The area of $\Sigma_a$ is well-defined for all $a\in[0,a_0]$ since $[0,a_0]\subset[0,a'_0]\subset Y$ because of (\ref{a'0}).\par
That completes the proof of Lemma~\ref{lem310}.\hfill\fin

\begin{lem}\label{lem311} Assume that $\Omega\in{\mathcal{C}}$ is not a ball. Then, with the notations of Section~\ref{rearrangement} and Lemma~\ref{lem310}, one has
$$\widetilde{\psi}(x)\ge\left(1+\frac{\beta}{2}\right)\rho^{-1}(|x|)$$
for all $x\in\overline{\Omega^*}$ such that $\rho(a_0)\le|x|\le R$.
\end{lem}

\noindent{\bf{Proof.}} From Lemma \ref{lem310}, one knows that
$$\overline{S_{\rho(a_0),R}}\subset E.$$
Notice that $0<\rho(a_0)<R$. Fix any
$x\in\overline{\Omega^*}$ such that $r=|x|\in[\rho(a_0),R]$ (that is $\rho^{-1}(r)\in[0,a_0]\subset Y$). Formula (\ref{rho'psibis}) of Section~\ref{rearrangement} implies that
$$\frac{i(\rho^{-1}(r))}{n\alpha_nr^{n-1}}\le|\rho'(\rho^{-1}(r))|\times|\nabla \widetilde{\psi}(x)|.$$
But 
$$\baa{rcl}
i(\rho^{-1}(r))=\displaystyle{\int_{\Sigma_{\rho^{-1}(r)}}}d\sigma_{\rho^{-1}(r)}(y)=\hbox{area}(\Sigma_{\rho^{-1}(r)})
& \ge & \left(1+\displaystyle{\frac{\beta}{2}}\right)n\alpha_nR^{n-1}\\
& \ge & \left(1+\displaystyle{\frac{\beta}{2}}\right)n\alpha_nr^{n-1}\eaa$$
from Lemma \ref{lem310}. Thus,
$$1+\frac{\beta}{2}\le|\rho'(\rho^{-1}(r))|\times|\nabla 
\widetilde{\psi}(x)|.$$
The conclusion of Lemma~\ref{lem311} follows from the above inequality, 
as in the proof of Corollary~\ref{cor36}.\hfill\fin\break

The improved version of Corollary~\ref{cor36} is the following:

\begin{cor}\label{corball1} Assume that $\Omega\in{\mathcal{C}}$ is not a ball. Then there exists a positive constant $\eta=\eta(\Omega,n,\alpha,N,\delta)>0$ depending only on $\Omega$, $n$, $\alpha$, $N$ and $\delta$, such that
$$\widetilde{\psi}(x)\ge(1+\eta)\ \rho^{-1}(|x|)$$
for all $x\in\overline{\Omega^*}$.
\end{cor}

\noindent{\bf{Proof.}} Let $e$ be any unit vector in $\R^n$. Let $\widetilde{\Phi}$ be the function defined in
$[0,R]$ by
$$\widetilde{\Phi}(r)=\widetilde{\psi}(re)\hbox{ for all }r\in[0,R].$$
This function is continuous in $[0,R]$, differentiable (except at a finite set of points) and decreasing in $[0,R]$. Furthermore,
Proposition \ref{pro35} and the fact that $\rho^{-1}$ is decreasing in $[0,R]$ imply that
$$-\widetilde{\Phi}'(r)\ge-\frac{d}{dr}(\rho^{-1}(r))\ge 0$$
for all $r\in\rho(Y)=(0,\rho(a_{m-1}))\cup\cdots\cup(\rho(a_{2}),\rho(a_{1}))\cup(\rho(a_{1}),R]$. As in the proof of Corollary~\ref{cor36}, the mean value theorem yields
$$\widetilde{\psi}(re)-\widetilde{\psi}(\rho(a_0)e)\ge\rho^{-1}(r)-a_0$$
for all $r\in[0,\rho(a_0)]$. For each such a $r$ in $[0,\rho(a_0)]$, one has $\rho^{-1}(r)\in[a_0,M]\subset(0,M]$, whence
$$\frac{\widetilde{\psi}(re)}{\rho^{-1}(r)}\ge 1+\frac{\widetilde{\psi}(\rho(a_0)e)-a_0}{\rho^{-1}(r)}\ge 1+\frac{\beta a_0}{2\rho^{-1}(r)}$$
from Lemma \ref{lem311}. Remind that $M$ denotes the 
maximum of $\psi$, so that $\rho^{-1}(r)\le M\leq N$. Hence, one obtains
$$\widetilde{\psi}(re)\ge\left(1+\frac{\beta a_0}{2N}\right)\rho^{-1}(r)\hbox{ for all }r\in[0,\rho(a_0)].$$\par
As in the proof of Corollary~\ref{cor36}, the conclusion of Corollary~\ref{corball1}
follows from the above inequality and from Lemma \ref{lem311}, with the choice
$$\eta=\eta(\Omega,n,\alpha,N,\delta)=\min\left(\frac{\beta}{2},\frac{\beta a_0}{2N}\right)=\frac{\beta a_0}{2N}>0$$
for instance (notice that $a_0\le M\le N$).\hfill\fin\break

Lastly, the following corollary is an improved version of Corollary~\ref{corintegral}.

\begin{cor}\label{corball2} Assume that $\Omega\in{\mathcal{C}}$ is not a ball. If there are $\omega_0\ge 0$ and $\mu\ge 0$ such that
$$-{\rm{div}}(A_{\Omega}\nabla\psi)(y)-|\omega(y)|\times|\nabla\psi(y)|-\omega_0|\nabla\psi(y)|+V(y)\psi(y)\le\mu\psi(y)\ \hbox{ for all }y\in\Omega,$$
then, under the notation of Section~\ref{rearrangement},
$$\int_{\Omega^{\ast}} \left[\widehat{\Lambda}(x)\vert \nabla\widetilde{\psi}(x)\vert^2-\omega_0\vert \nabla\widetilde{\psi}(x)\vert\widetilde{\psi}(x)+\widehat{V}(x)\widetilde{\psi}(x)^2\right]e^{-U(x)}dx\leq\frac{\mu}{1+\eta}\int_{\Omega^{\ast}} \widetilde{\psi}(x)^2e^{-U(x)}dx,$$
where the positive constant $\eta=\eta(\Omega,n,\alpha,N,\delta)>0$ is given in Corollary~\ref{corball1} and depends only on $\Omega$, $n$, $\alpha$, $N$ and $\delta$.
\end{cor}

\noindent{\bf{Proof.}} Under the assumptions of Corollary~\ref{corball2}, it follows from Proposition~\ref{comparison} and Corollary~\ref{corball1} that
$$-{\rm{div}}(\widehat{\Lambda}\nabla\widetilde{\psi})(x)+\widehat{v}(x)\cdot\nabla\widetilde{\psi}(x)-\omega_0\vert \nabla\widetilde{\psi}(x)\vert+\widehat{V}(x)\widetilde{\psi}(x)\leq\mu\rho^{-1}(|x|)\le\frac{\mu}{1+\eta}\tilde{\psi}(x)$$
for all $x\in E\cap\Omega^*$. The conclusion of Corollary~\ref{corball2} then follows from Proposition~\ref{integralcompar}.\hfill\fin


\SE{Application to eigenvalue problems} \label{appl}

The present section is devoted to the proofs of some of the main theorems which were stated in Section~\ref{results}. We apply the rearrangement inequalities of the previous two sections to get some comparison results for the principal eigenvalues of operators which are defined in $\Omega$ and in $\Omega^*$. Here, the data have given averages or given distribution functions, or satisfy other types of pointwise constraints.\par
We shall use a triple approximation process. First, we approximate the diffusion and the drift coefficients in $\Omega$ by smooth functions. Second, we approximate the principal eigenfunctions in $\Omega$ by analytic functions. Lastly, we approximate the symmetrized data in $\Omega^*$ by coefficients having the same distribution functions or satisfying the same constraints as the original data in $\Omega$. Section~\ref{sec4.1} is concerned with the latter approximation process. In Section~\ref{sec4.2}, we deal with the case of general non-symmetric operators for which the inverse of the lower bound $\Lambda$ of the diffusion matrix field $A$ has a given $L^1$ norm, the drift $v$ has a given $L^2$ norm with weight $\Lambda^{-1}$ and the negative part of the potential $V$ has a given distribution function. Lastly, in Section~\ref{sec4.3}, we consider diffusion matrix fields $A$ whose trace and determinant satisfy some pointwise constraints.


\subsection{Approximation of symmetrized fields by fields having given distribution functions}\label{sec4.1}

In this subsection, $\Omega$ denotes an open connected bounded non-empty $C^1$ subset of $\R^n$ and $\Omega^*$ is the open Euclidean ball which is centered at the origin and such that $|\Omega^*|=|\Omega|$. In this subsection we do not require $\Omega$ to be of class $C^2$. Call $R>0$ the radius of $\Omega^*$, that is $\Omega^*=B_R$. For $0\le s<s'\le R$, one recalls that $S_{s,s'}=\{z\in\R^n,\ s<|z|<s'\}$.\par
Let $\psi\ :\ \overline{\Omega}\to\R$ be a continuous function. Call
$$m=\min_{\overline{\Omega}}\psi\ \hbox{ and }\ M=\max_{\overline{\Omega}}\psi.$$
For all $a<b\in\R$, denote
$$\Omega_a=\{x\in\Omega,\ a<\psi(x)\},\ \Omega_{a,b}=\{x\in\Omega,\ a<\psi(x)\le b\}\hbox{ and }\Sigma_a=\{x\in\Omega,\ \psi(x)=a\}.$$
One assumes that
\be\label{hyppsi}
|\Sigma_a|=0\hbox{ for all }a\in\R.
\ee
It follows then that $m<M$. For each $a\in[m,M]$, set
$$\rho(a)=\left(\frac{|\Omega_a|}{\alpha_n}\right)^{1/n}.$$
The function $\rho\ :\ [m,M]\to\R$ is then continuous, decreasing and it ranges onto $[0,R]$.

\begin{lem}\label{lemapprox} Under assumption $(\ref{hyppsi})$, let $g$ be in $L^{\infty}(\Omega,\R)$ and $h$ in $L^{\infty}(\Omega^*,\R)$ and radially symmetric, and assume that
\be\label{gh}
\forall\ m\le a<b\le M, \quad\int_{\Omega_{a,b}}g=\int_{S_{\rho(b),\rho(a)}}h.
\ee
Then there exist a sequence of radially symmetric functions $(g_k)_{k\in\N}$ in $L^{\infty}(\Omega^*,\R)$ and two sequences of radially symmetric $C^{\infty}(\overline{\Omega^*},\R)$ functions $(\underline{g}_k)_{k\in\N}$ and $(\overline{g}_k)_{k\in\N}$ such that
$$g_k\rightharpoonup h,\ \ \underline{g}_k\rightharpoonup h\ \hbox{ and }\ \overline{g}_k\rightharpoonup h\ \hbox{ as }\ k\to+\infty$$
in $\sigma(L^p(\Omega^*),L^{p'}(\Omega^*))$ for all $1<p\le+\infty$ $($weak convergence, with $1/p+1/p'=1$ for $1<p<+\infty$, and weak-* convergence for $p=+\infty$ and $p'=1)$,\footnote{This convention is used throughout the paper.} and
$$\left\{\baa{ll}
|\{x\in\Omega,\ g(x)>t\}|=|\{x\in\Omega^*,\ g_k(x)>t\}| & \hbox{for all }t\in\R\hbox{ and }k\in\N,\\
\\
|\{x\in\Omega,\ g(x)\ge t\}|=|\{x\in\Omega^*,\ g_k(x)\ge t\}| & \hbox{for all }t\in\R\hbox{ and }k\in\N,\\
\\
\displaystyle{\mathop{\hbox{\rm{ess} \rm{inf}}}_{\Omega}}\ g\le\underline{g}_k\le g_k\le\overline{g}_k\le\displaystyle{\mathop{\hbox{\rm{ess} \rm{sup}}}_{\Omega}}\ g & \hbox{a.e. in }\Omega^*\hbox{ for all }k\in\N.\eaa\right.$$
\end{lem}

\begin{rem}{\rm The functions $g_k$ and $\underline{g}_k$ are actually constructed so that
$$g_k-\underline{g}_k\to 0\ \hbox{ and }\ g_k-\overline{g}_k\to 0\ \hbox{ as }\ k\to+\infty\hbox{ in }L^q(\Omega^*)$$
for all $q\in[1,+\infty)$ $($see the proof of Lemma~\ref{lemapprox}$)$.}
\end{rem}

\begin{rem}{\rm The functions $g_k$ then have the same distribution function as $g$ and the functions $\underline{g}_k$ $($resp. $\overline{g}_k)$ then have distribution functions which are less than $($resp. larger than$)$ or equal to that of~$g$. Moreover, if one further assumes that $g$ and $h$ are nonnegative almost everywhere in $\Omega$ and $\Omega^*$ respectively, then the functions $g_k$ $($resp. $\underline{g}_k$, $\overline{g}_k)$ are nonnegative almost everywhere in $\Omega^*$ $($resp. everywhere in $\overline{\Omega^*})$. The same property holds good with nonpositivity instead of nonnegativity.}
\end{rem}

In order not to lengthen the reading of the paper, the proof of this lemma is postponed in the Appendix. To finish this subsection, we just point out an immediate corollary of Lemma~\ref{lemapprox}.

\begin{cor}\label{corghat} In addition to $(\ref{hyppsi})$, assume that $\psi$ is in $C^1(\overline{\Omega})$, $\psi=m$ on $\partial\Omega$, $\psi>m$ in $\Omega$, and $\psi$ has a finite number $p$ of critical values $a_i$ with $m<a_1<\cdots<a_p=M$. Let $g\in L^{\infty}(\Omega,\R)$ and $\hat{g}\in L^{\infty}(\Omega^*,\R)$ be defined by
\be\label{defhatg}
\hat{g}(x)=\frac{\displaystyle{\int_{\Sigma_{\rho^{-1}(|x|)}}}g(y)|\nabla\psi(y)|^{-1}d\sigma_{\rho^{-1}(|x|)}(y)}{\displaystyle{\int_{\Sigma_{\rho^{-1}(|x|)}}}|\nabla\psi(y)|^{-1}d\sigma_{\rho^{-1}(|x|)}(y)}
\ee
for almost every $x\in\Omega^*$ such that $|x|\neq\rho(a_p),\ldots,\rho(a_1)$, where $d\sigma_a$ denotes the surface measure on $\Sigma_a$ for $a\in[m,M]$ which is not a critical value of $\psi$. Then the conclusion of Lemma~\ref{lemapprox} holds for $h=\hat{g}$.
\end{cor}

\noindent{\bf{Proof.}} It is enough to prove that the function $\hat{g}$ defined by (\ref{defhatg}) is indeed well-defined, bounded and radially symmetric and that property (\ref{gh}) is satisfied with $h=\hat{g}$. To do so, choose first any two real numbers $a<b$ in $[m,M]$ such that $a_i\not\in[a,b]$ for all $i=1,\ldots,p$. From the co-area formula and Fubini's theorem, one has
\be\label{coarea}
\int_{\Omega_{a,b}}|g(y)|dy=\int_a^b\left(\int_{\Sigma_s}|g(y)|\times|\nabla\psi(y)|^{-1}d\sigma_s(y)\right)ds<+\infty
\ee
and the quantity
$$\int_{\Sigma_s}|g(y)|\times|\nabla\psi(y)|^{-1}d\sigma_s(y)$$
is therefore finite for almost every $s\in[a,b]$. The quantity
$$\frac{\displaystyle{\int_{\Sigma_{\rho^{-1}(r)}}}g(y)|\nabla\psi(y)|^{-1} d\sigma_{\rho^{-1}(r)}(y)}{\displaystyle{\int_{\Sigma_{\rho^{-1}(r)}}}|\nabla\psi(y)|^{-1} d\sigma_{\rho^{-1}(r)}(y)}$$
is then finite for almost every $r$ in $[\rho^{-1}(b),\rho^{-1}(a)]$ and for all $m\le a<b\le M$ such that $\{a_1,\ldots,a_p\}\cap[a,b]=\emptyset$. In other words, $\hat{g}$ is well-defined for almost every $x\in\Omega^*$ such that $|x|\neq\rho(a_1),\ldots,\rho(a_p)$. Moreover, the function $\hat{g}$ is in $L^{\infty}(\Omega^*,\R)$, 
$$\mathop{\hbox{ess inf}}_{\Omega^*}\ g\le\mathop{\hbox{ess inf}}_{\Omega^*}\ \hat{g}\le\mathop{\hbox{ess sup}}_{\Omega^*}\ \hat{g}\le\mathop{\hbox{ess sup}}_{\Omega^*}\ g$$
and $\hat{g}$ is clearly radially symmetric.\par
On the other hand, the same calculations as the ones which were done in Lemma~\ref{lem32} in Section~\ref{rearrangement} imply that the function $\rho$ is actually differentiable at each value $a\not\in\{a_1,\ldots,a_p\}$, with
$$\rho'(a)=-(n\alpha_n\rho(a)^{n-1})^{-1}\int_{\Sigma_a}|\nabla\psi(y)|^{-1}d\sigma_a(y).$$
Coming back to (\ref{coarea}), the change of variables $s=\rho^{-1}(r)$ then yields
$$\int_{\Omega_{a,b}}g(y)dy=\int_{\rho(b)}^{\rho(a)}
\left(\int_{\Sigma_{\rho^{-1}(r)}}g(y)|\nabla\psi(y)|^{-1} d\sigma_{\rho^{-1}(r)}(y)\right)
\displaystyle{\frac{n\alpha_nr^{n-1}}{\displaystyle{\int_{\Sigma_{\rho^{-1}(r)}}} |\nabla\psi(y)|^{-1}d\sigma_{\rho^{-1}(r)}(y)}}dr,$$
that is
$$\int_{\Omega_{a,b}}g(y)dy=\displaystyle{\int_{S_{\rho(b),\rho(a)}}}\hat{g}(x)dx$$
for all $a<b$ in $[m,M]$ such that $\{a_1,\ldots,a_p\}\cap[a,b]=\emptyset$. Since $|\Sigma_c|=0$ for all $c\in[m,M]$ (and then $\rho$ is continuous), one gets from Lebesgue's dominated convergence theorem that
$$\int_{\Omega_{a,b}}g=\int_{S_{\rho(b),\rho(a)}}\hat{g}$$
for all $a<b$ in $[m,M]$ and then the conclusion of Corollary~\ref{corghat} 
follows from Lemma~\ref{lemapprox}.\hfill\fin


\subsection{Operators whose coefficients have given averages or given distribution functions}\label{sec4.2}

In this subsection, we consider operators for which $A\ge\Lambda\Id$, and $\Lambda$ and $v$ satisfy some integral constraints while the negative part of the potential $V$ will be fixed.\hfill\break

\noindent{\bf{Proof of Theorem~\ref{th1}.}} It will be divided into several steps. Throughout the proof, we fix
$$\Omega\in{\mathcal{C}},\ A\in W^{1,\infty}(\Omega,{\mathcal{S}}_n(\R)),\ \Lambda\in L^{\infty}_+(\Omega),\ v\in L^{\infty}(\Omega,\R^n)\hbox{ and }V\in C(\overline{\Omega}).$$
Call
\be\label{mMLambda}
0<m_{\Lambda}=\mathop{\hbox{ess inf}}_{\Omega}\Lambda\le\mathop{\hbox{ess sup}}_{\Omega}\Lambda=M_{\Lambda}<+\infty
\ee
and assume that
\be\label{ALambda}
A\ge\Lambda\Id\hbox{ a.e. in }\Omega.
\ee
and that
$$\lambda_1(\Omega,A,v,V)\ge 0.$$\hfill\break\par

{\underline{Step 1: Approximation of $\Lambda$ in $\Omega$}}. Write first
$$A(x)=(a_{i,j}(x))_{1\le i,j\le n}.$$
Each function $a_{i,j}$ is in $W^{1,\infty}(\Omega)$ and is therefore continuous in $\overline{\Omega}$ (up to a modification on a zero-measure set). For each $x\in\overline{\Omega}$, call $\Lambda[A](x)$ the lowest eigenvalue of $A(x)$, that is:
$$\forall\ x\in\overline{\Omega},\quad\Lambda[A](x)=\min_{\xi\in\R^n,\ |\xi|=1}A(x)\xi\cdot\xi.$$
The function $\Lambda[A]$ is then continuous in $\overline{\Omega}$ and
$$\Lambda[A](x)\ge\Lambda(x)\hbox{ a.e. in }\Omega$$
because of (\ref{ALambda}). In particular, $\Lambda[A](x)\ge m_{\Lambda}$ for all $x\in\overline{\Omega}$, where $m_{\Lambda}>0$ has been defined in (\ref{mMLambda}). There exists then a continuous function $\overline{\Lambda}$ in $\R^n$ such that
$$\overline{\Lambda}(x)=\Lambda[A](x)\hbox{ for all }x\in\overline{\Omega},\hbox{ and }m_{\Lambda}\le\overline{\Lambda}(x)\le\|\Lambda[A]\|_{L^{\infty}(\Omega)}\hbox{ for all }x\in\R^n.$$\par
{\bf{We first consider the case when $m_{\Lambda}<M_{\Lambda}$}}. Thus,
$$\int_{\Omega}\Lambda(y)^{-1}(y)dy>M_{\Lambda}^{-1}|\Omega|,$$
whence
$$\overline{\epsilon}:=M_{\Lambda}-\|\Lambda^{-1}\|_{L^1(\Omega)}^{-1}|\Omega|\in(0,M_{\Lambda}-m_{\Lambda}].$$\par
Pick any $\epsilon$ in $(0,\overline{\epsilon})$. Let $J_{\epsilon}$ be the function defined in $[0,M_{\Lambda}-m_{\Lambda}-\epsilon]$ by
$$\forall\ \tau\in[0,M_{\Lambda}-m_{\Lambda}-\epsilon],\ J_{\epsilon}(\tau)=\int_{\Omega}
\max\left(\min(\Lambda(y),M_{\Lambda}-\epsilon),m_{\Lambda}+\tau\right)^{-1}dy.$$
This means that the function $\Lambda$ is truncated between $m_{\Lambda}+\tau$ and $M_{\lambda}-\epsilon$. The function $J_{\epsilon}$ is Lipschitz-continuous and nonincreasing in $[0,M_{\Lambda}-m_{\Lambda}-\epsilon]$. Furthermore,
$$J_{\epsilon}(0)=\int_{\Omega}\min(\Lambda(y),M_{\Lambda}-\epsilon)^{-1}dy>\int_{\Omega}\Lambda(y)^{-1}dy$$
since $M_{\Lambda}-\epsilon<M_{\Lambda}=\hbox{ess sup}_{\Omega}\ \Lambda$, while
$$J_{\epsilon}(M_{\Lambda}-m_{\Lambda}-\epsilon)=(M_{\Lambda}-\epsilon)^{-1}|\Omega|<(M_{\Lambda}-\overline{\epsilon})^{-1}|\Omega|=\int_{\Omega}\Lambda(y)^{-1}dy$$
owing to the definition of $\overline{\epsilon}$. Therefore,
$$\tau(\epsilon):=\min\left\{\tau\in[0,M_{\Lambda}-m_{\Lambda}-\epsilon],\ J_{\epsilon}(\tau)=\int_{\Omega}\Lambda(y)^{-1}dy\right\}$$
is well-defined and $0<\tau(\epsilon)<M_{\Lambda}-m_{\lambda}-\epsilon$.\par
Moreover, we claim that
\be\label{tauepsilon}
\tau(\epsilon)\to 0^+\hbox{ as }\epsilon\to 0^+.
\ee
If not, there exist $\tau_{\infty}\in(0,M_{\Lambda}-m_{\Lambda}]$ and a sequence $(\epsilon_p)_{p\in\N}\in(0,\overline{\epsilon})$ such that $\epsilon_p\to 0$ and $\tau(\epsilon_p)\to\tau_{\infty}$ as $p\to+\infty$. Then,
$$\baa{rcl}
\displaystyle{\int_{\Omega}}\Lambda(y)^{-1}dy\ =\ J_{\epsilon_p}(\tau(\epsilon_p)) & = & \displaystyle{\int_{\Omega}}\max\left(\min(\Lambda(y),M_{\Lambda}-\epsilon_p),m_{\Lambda}+\tau(\epsilon_p)\right)^{-1}dy\\
& \displaystyle{\mathop{\to}_{p\to+\infty}} & \displaystyle{\int_{\Omega}}\max\left(\min(\Lambda(y),M_{\Lambda}),m_{\Lambda}+\tau_{\infty}\right)^{-1}dy,\eaa$$
whence
$$\int_{\Omega}\Lambda(y)^{-1}dy=\int_{\Omega}\max(\Lambda(y),m_{\Lambda}+\tau_{\infty})^{-1}dy$$
and $\Lambda\ge m_{\Lambda}+\tau_{\infty}>m_{\Lambda}$ a.e. in $\Omega$, which is impossible.\par
Choose a sequence $(\epsilon_k)_{k\in\N}$ of real numbers such that
$$0<\epsilon_k<\overline{\epsilon}\hbox{ for all }k\in\N,\hbox{ and }\epsilon_k\to 0\hbox{ as }k\to+\infty.$$
For each $k\in\N$, call
$$\tau_k=\tau(\epsilon_k)\in(0,M_{\Lambda}-m_{\Lambda}-\epsilon_k).$$
It then follows from (\ref{tauepsilon}) that
$$\tau_k\to 0\hbox{ as }k\to+\infty.$$
Then, for each $k\in\N$, denote
\be\label{lambdaomegak}
\Lambda_{\Omega,k}=\max\left(\min(\Lambda,M_{\Lambda}-\epsilon_k),m_{\Lambda}+\tau_k\right)\hbox{ a.e. in }\Omega
\ee
and define the function $\underline{\Lambda}_k$ in $\R^n$ (almost everywhere) by
$$\underline{\Lambda}_k(y)=\left\{\baa{ll}\Lambda_{\Omega,k}(y) & \hbox{if }y\in\Omega,\\ m_{\Lambda}+\tau_k & \hbox{if }y\in\R^n\backslash\Omega,\eaa\right.$$
and the continuous function $\overline{\Lambda}_k$ by
$$\overline{\Lambda}_k(y)=\max(\overline{\Lambda}(y),m_{\Lambda}+\tau_k)\hbox{ for all }y\in\R^n.$$
Notice that
\be\label{lambdak-1}
\int_{\Omega}\underline{\Lambda}_k(y)^{-1}dy=\int_{\Omega}\Lambda_{\Omega,k}(y)^{-1}dy=J_{\epsilon_k}(\tau_k)=\int_{\Omega}\Lambda(y)^{-1}dy
\ee
for each $k\in\N$. Observe also that
$$0<m_{\Lambda}<m_{\Lambda}+\tau_k\le\underline{\Lambda}_k\le M_{\Lambda}-\epsilon_k<M_{\Lambda},\quad\underline{\Lambda}_k\le\overline{\Lambda}_k\hbox{ a.e. in }\R^n,$$
and that
$$\overline{\Lambda}_k\le\max(\|\Lambda[A]\|_{L^{\infty}(\Omega)},m_{\Lambda}+\tau_k)=\|\Lambda[A]\|_{L^{\infty}(\Omega)}\hbox{ in }\R^n,$$
because $\|\Lambda[A]\|_{L^{\infty}(\Omega)}\ge\hbox{ess sup}_{\Omega}\ \Lambda=M_{\Lambda}\ge m_{\Lambda}+\tau_k$.\par
Let $(\rho_{k'})_{k'\in\N}$ be a sequence of mollifiers in $\R^n$. For each $(k,k')\in\N^2$, call
$$\overline{\Lambda}_{k,k'}=\left(\rho_{k'}*\overline{\Lambda}_k^{-1}\right)^{-1}\hbox{ and }\underline{\Lambda}_{k,k'}=\left(\rho_{k'}*\underline{\Lambda}_k^{-1}\right)^{-1}$$
that is
$$\overline{\Lambda}_{k,k'}(y)^{-1}=\int_{\R^n}\rho_{k'}(z)\overline{\Lambda}_k(y-z)^{-1}dz\hbox{ and }\underline{\Lambda}_{k,k'}(y)^{-1}=\int_{\R^n}\rho_{k'}(z)\underline{\Lambda}_k(y-y)^{-1}dz$$
for all $y\in\R^n$. The functions $\overline{\Lambda}_{k,k'}$ and $\underline{\Lambda}_{k,k'}$ are of class $C^{\infty}(\R^n)$ and they satisfy
$$0<m_{\Lambda}<m_{\Lambda}+\tau_k\le\underline{\Lambda}_{k,k'}(y)\le M_{\Lambda}-\epsilon_k<M_{\Lambda}$$
and
$$\underline{\Lambda}_{k,k'}(y)\le\overline{\Lambda}_{k,k'}(y)\le\|\Lambda[A]\|_{L^{\infty}(\Omega)}\hbox{ for all }y\in\R^n.$$
Furthermore, for each $k\in\N$,
\be\label{limitk'}
\underline{\Lambda}_{k,k'}^{-1}\to\underline{\Lambda}_k^{-1}\hbox{ in }L^p_{loc}(\R^n)\hbox{ for all }1\le p<+\infty\hbox{ as }k'\to+\infty
\ee
and
$$\overline{\Lambda}_{k,k'}\to\overline{\Lambda}_k=\max(\Lambda[A],m_{\Lambda}+\tau_k)\hbox{ uniformly in }\overline{\Omega}\hbox{ as }k'\to+\infty.$$
Actually, since $\overline{\Lambda}_k^{-1}\to\overline{\Lambda}^{-1}$ uniformly in $\R^n$ as $k\to+\infty$ (because $\tau_k\to 0$ as $k\to+\infty$ and $\overline{\Lambda}\ge m_{\Lambda}>0$ in $\R^n$), one even has that
\be\label{Lambdakk'}
\overline{\Lambda}_{k,k'}-\overline{\Lambda}_k\to 0\hbox{ uniformly in }\overline{\Omega}\hbox{ as }(k,k')\to(+\infty,+\infty).
\ee
Call, for each $(k,k')\in\N^2$,
$$\alpha_{k,k'}=\frac{\|\underline{\Lambda}_{k,k'}^{-1}\|_{L^1(\Omega)}}{\|\Lambda^{-1}\|_{L^1(\Omega)}}>0.$$
Define the function $\underline{\Lambda}$ almost everywhere in $\R^n$ by:
$$\underline{\Lambda}=\left\{\baa{ll}\Lambda & \hbox{in }\Omega,\\ m_{\Lambda} & \hbox{in }\R^n\backslash\Omega.\eaa\right.$$
Since $0<m_{\Lambda}\le\underline{\Lambda}\le M_{\Lambda}$ a.e. in $\R^n$ and $\epsilon_k,\ \tau_k\to 0$ as $k\to+\infty$, it follows that
$$\|\underline{\Lambda}_k-\underline{\Lambda}\|_{L^{\infty}(\R^n)}\to 0\hbox{ and }\|\underline{\Lambda}_k^{-1}-\underline{\Lambda}^{-1}\|_{L^{\infty}(\R^n)}\to 0\hbox{ as }k\to+\infty.$$
Thus,
$$\|\underline{\Lambda}_{k,k'}^{-1}-\Lambda^{-1}\|_{L^1(\Omega)}\le\|\rho_{k'}*(\underline{\Lambda}_k^{-1}-\underline{\Lambda}^{-1})\|_{L^1(\Omega)}+\|\rho_{k'}*\underline{\Lambda}^{-1}-\underline{\Lambda}^{-1}\|_{L^1(\Omega)}\to 0\hbox{ as }(k,k')\to(+\infty,+\infty),$$
whence
\be\label{alphakk'}
\alpha_{k,k'}-1=\frac{\|\underline{\Lambda}_{k,k'}^{-1}\|_{L^1(\Omega)}-\|\Lambda^{-1}\|_{L^1(\Omega)}}{\|\Lambda^{-1}\|_{L^1(\Omega)}}\to 0\hbox{ as }(k,k')\to(+\infty,+\infty).
\ee
Furthermore, because of (\ref{lambdak-1}) and (\ref{limitk'}), there holds
\be\label{alphakk'bis}
\alpha_{k,k'}-1=\frac{\|\underline{\Lambda}_{k,k'}^{-1}\|_{L^1(\Omega)}-\|\underline{\Lambda}_k^{-1}\|_{L^1(\Omega)}}{\|\Lambda^{-1}\|_{L^1(\Omega)}}\to 0\hbox{ as }k'\to+\infty,\hbox{ for each }k\in\N.
\ee\par
Define now
$$\Lambda_{k,k'}(y)=\alpha_{k,k'}\underline{\Lambda}_{k,k'}(y)\hbox{ for all }y\in\overline{\Omega}\hbox{ and }(k,k')\in\N^2.$$
The functions $\Lambda_{k,k'}$ are of class $C^{\infty}(\overline{\Omega})$ and they satisfy
\be\label{Lambdak1}
\int_{\Omega}\Lambda_{k,k'}(y)^{-1}dy=\int_{\Omega}\Lambda(y)^{-1}dy\hbox{ for all }(k,k')\in\N^2,
\ee
and
\be\label{Lambdak2}\left\{\baa{l}
0<\alpha_{k,k'}\times(m_{\Lambda}+\tau_k)\le\Lambda_{k,k'}\le\alpha_{k,k'}\times(M_{\Lambda}-\epsilon_k)\hbox{ in }\overline{\Omega}\hbox{ for all }(k,k')\in\N^2,\\
\alpha_{k,k'}\times(m_{\Lambda}+\tau_k)\displaystyle{\mathop{\to}_{k'\to+\infty}}m_{\Lambda}+\tau_k\hbox{ for all }k\in\N,\\
\alpha_{k,k'}\times(M_{\Lambda}-\epsilon_k)\displaystyle{\mathop{\to}_{k'\to+\infty}}M_{\Lambda}-\epsilon_k\hbox{ for all }k\in\N,\eaa\right.
\ee
together with
\be\label{Lambdak2bis}
\|\Lambda_{k,k'}^{-1}-\Lambda^{-1}\|_{L^1(\Omega)}\to0\hbox{ as }(k,k')\to(+\infty,+\infty).
\ee
Lastly,
\be\label{LambdakA}
\Lambda_{k,k'}=\alpha_{k,k'}\underline{\Lambda}_{k,k'}\le\alpha_{k,k'}\overline{\Lambda}_{k,k'}\hbox{ in }\overline{\Omega}\hbox{ for all }(k,k')\in\N^2
\ee
and
\be\label{LambdakAbis}
\alpha_{k,k'}\overline{\Lambda}_{k,k'}-\overline{\Lambda}_k\to 0\hbox{ uniformly in }\overline{\Omega}\hbox{ as }(k,k')\to(+\infty,+\infty)
\ee
from (\ref{Lambdakk'}) and (\ref{alphakk'}).\par
{\bf{In the case when $m_{\Lambda}=M_{\Lambda}$}}, namely when $\Lambda$ is equal to a constant (up to modification on a zero-measure set), then one sets $\Lambda_{k,k'}=\underline{\Lambda}_{k,k'}=\underline{\Lambda}_k=\Lambda_{\Omega,k}=\Lambda$, $\alpha_{k,k'}=1$, $\overline{\Lambda}_k=\overline{\Lambda}$, $\epsilon_k=\tau_k=0$ and properties (\ref{Lambdak1}), (\ref{Lambdak2}), (\ref{Lambdak2bis}), (\ref{LambdakA}) and (\ref{LambdakAbis}) hold immediately.\hfill\break\par

{\underline{Step~2: Approximation of $A$ in $\Omega$}}. Let us now approximate the $W^{1,\infty}(\Omega)$ matrix field $A=(a_{i,j})_{1\le i,j\le n}$. First, each function $a_{i,j}$ can be extended to a $W^{1,\infty}(\R^n)$ function $\overline{a}_{i,j}$ such that
$$\|\overline{a}_{i,j}\|_{L^{\infty}(\R^n)}=\|a_{i,j}\|_{L^{\infty}(\Omega)}\hbox{ and }\|\nabla \overline{a}_{i,j}\|_{L^{\infty}(\R^n)}=\|\nabla a_{i,j}\|_{L^{\infty}(\Omega)},$$
whence
$$\|\overline{a}_{i,j}\|_{W^{1,\infty}(\R^n)}=\|a_{i,j}\|_{W^{1,\infty}(\Omega)}\le\|A\|_{W^{1,\infty}(\Omega)},$$
where we recall that $\|A\|_{W^{1,\infty}(\Omega)}=\max_{1\le i,j\le n}\|a_{i,j}\|_{W^{1,\infty}(\Omega)}$. Since the matrix field $A=(a_{i,j})_{1\le i,j\le n}$ is symmetric, the matrix field $(\overline{a}_{i,j})_{1\le i,j\le n}$ can be assumed to be symmetric. For each $1\le i,j\le n$, the functions $\rho_{k'}*\overline{a}_{i,j}$ are of class $C^{\infty}(\R^n)$ and converge uniformly to $a_{i,j}$ in $\overline{\Omega}$ as $k'\to+\infty$. Furthermore,
$$\|\nabla(\rho_{k'}*\overline{a}_{i,j})\|_{L^{\infty}(\R^n)}\le\|\nabla\overline{a}_{i,j}\|_{L^{\infty}(\R^n)}=\|\nabla a_{i,j}\|_{L^{\infty}(\Omega)}\le\|a_{i,j}\|_{W^{1,\infty}(\Omega)}\le\|A\|_{W^{1,\infty}(\Omega)}$$
for all $k'\in\N$ and $1\le i,j\le n$. For each $k'\in\N$, the matrix field $(\rho_{k'}*\overline{a}_{i,j})_{1\le i,j\le n}$ can be approximated in $C^1(\overline{\Omega})$ norm by symmetric matrix fields with polynomial entries in $\overline{\Omega}$. Therefore, there exists a sequence of symmetric matrix fields $(A'_{k'})_{k'\in\N}=((a'_{k',i,j})_{1\le i,j\le n})_{k'\in\N}$ in $\overline{\Omega}$ with polynomial entries $a'_{k',i,j}$ such that, for all $1\le i,j\le n$,
\be\label{a'kij}
a'_{k',i,j}\to a_{i,j}\hbox{ uniformly in }\overline{\Omega}\hbox{ as }k'\to+\infty
\ee
and
$$\limsup_{k'\to+\infty}\|\nabla a'_{k',i,j}\|_{L^{\infty}(\Omega)}\le\|\nabla a_{i,j}\|_{L^{\infty}(\Omega)}\le\|A\|_{W^{1,\infty}(\Omega)}.$$
Call
$$\eta_{k,k'}=n^2\times\max_{1\le i,j\le n}\|a'_{k',i,j}-a_{i,j}\|_{L^{\infty}(\Omega)}\ +\ \|\overline{\Lambda}_k-\alpha_{k,k'}\overline{\Lambda}_{k,k'}\|_{L^{\infty}(\Omega)}+\tau_k.$$
Because of (\ref{LambdakAbis}), (\ref{a'kij}), and since $\tau_k\to 0$ as $k\to+\infty$, there holds
$$\eta_{k,k'}\to 0\hbox{ as }(k,k')\to(+\infty,+\infty).$$
The symmetric matrix fields
$$A_{k,k'}=A'_{k'}+\eta_{k,k'}\hbox{Id}=(a_{k,k',i,j})_{1\le i,j\le n}$$
with polynomial entries $a_{k,k',i,j}$ are such that
\be\label{Ak1}
a_{k,k',i,j}\to a_{i,j}\hbox{ uniformly in }\overline{\Omega}\hbox{ as }(k,k')\to(+\infty,+\infty)
\ee
and
\be\label{Ak2}
\sup_{k\in\N}\ \limsup_{k'\to+\infty}\ \|\nabla a_{k,k',i,j}\|_{L^{\infty}(\Omega)}\le\|\nabla a_{i,j}\|_{L^{\infty}(\Omega)}\le\|A\|_{W^{1,\infty}(\Omega)}.
\ee
Furthermore, for each $(k,k')\in\N^2$, $x\in\overline{\Omega}$ and $\xi=(\xi_1,\ldots,\xi_n)\in\R^n$ with $|\xi|=1$, there holds
$$\baa{l}
A_{k,k'}(x)\xi\cdot\xi-\Lambda_{k,k'}(x)\\
\qquad\qquad=\displaystyle{\mathop{\sum}_{1\le i,j\le n}}(a'_{k',i,j}(x)-a_{i,j}(x))\xi_i\xi_j\ +\ n^2\times\displaystyle{\mathop{\max}_{1\le i,j\le n}}\|a'_{k',i,j}-a_{i,j}\|_{L^{\infty}(\Omega)}\\
\qquad\qquad\quad+A(x)\xi\cdot\xi-\Lambda_{k,k'}(x)+\|\overline{\Lambda}_k-\alpha_{k,k'}\overline{\Lambda}_{k,k'}\|_{L^{\infty}(\Omega)}+\tau_k\\
\qquad\qquad\ge\Lambda[A](x)+\tau_k-\alpha_{k,k'}\overline{\Lambda}_{k,k'}(x)+\|\overline{\Lambda}_k-\alpha_{k,k'}\overline{\Lambda}_{k,k'}\|_{L^{\infty}(\Omega)}\hbox{ (because of (\ref{LambdakA}))}\\
\qquad\qquad\ge 0\eaa$$
because $\Lambda[A]+\tau_k\ge\max(\Lambda[A],m_{\Lambda}+\tau_k)=\overline{\Lambda}_k$ in $\overline{\Omega}$. In other words,
\be\label{Ak3}
A_{k,k'}(x)\ge \Lambda_{k,k'}(x)\hbox{Id}\hbox{ for all }(k,k')\in\N^2\hbox{ and }x\in\overline{\Omega},
\ee
in the sense of symmetric matrices.\hfill\break\par

{\underline{Step~3: Approximation of $v$ in $\Omega$}}. Call
$$0\le m_v=\mathop{\hbox{ess inf}}_{\Omega}|v|\le\mathop{\hbox{ess sup}}_{\Omega}|v|=\|v\|_{L^{\infty}(\Omega,\R^n)}=\|v\|_{\infty}=M_v<+\infty.$$\par
{\bf{We first consider the case when $m_v<M_v$ and $m_{\Lambda}<M_{\Lambda}$}}. In particular, it follows that $M_v=\|v\|_{\infty}>0$ and that there exists
\be\label{overepsilon'}
\overline{\epsilon}'\in(0,M_v-m_v)\hbox{ such that }(M_v-\overline{\epsilon}')^2\times\int_{\Omega}\Lambda(y)^{-1}dy>\int_{\Omega}|v(y)|^2\Lambda(y)^{-1}dy.
\ee
Call $K$ the function defined for all $\epsilon'\in[0,\overline{\epsilon}']$ by
$$K(\epsilon')=\int_{\Omega}\left[|v(y)|^2-\min\left(|v(y)|^2,(M_v-\epsilon')^2\right)\right]\times\Lambda(y)^{-1}dy.$$
The function $K$ is continuous and nondecreasing in $[0,\overline{\epsilon}']$, vanishes at $0$ and is positive in $(0,\overline{\epsilon}']$ due to the definition of $M_v$. Let the sequences of positive numbers $(\epsilon_k)_{k\in\N}$ and $(\tau_k)_{k\in\N}$ be as in Step~1. Since $\max(\epsilon_k,\tau_k)\to 0$ as $k\to+\infty$, one can assume without loss of generality that
$$\frac{|\Omega|M_v^2}{m_{\Lambda}^2}\times\max(\epsilon_k,\tau_k)<K(\overline{\epsilon}')\hbox{ for all }k\in\N.$$
For each $k\in\N$, call
$$\epsilon'_k=\min\left\{\epsilon'\in[0,\overline{\epsilon}'],\ \frac{|\Omega|M_v^2}{m_{\Lambda}^2}\times\max(\epsilon_k,\tau_k)=K(\epsilon')\right\}.$$
From the above remarks, $\epsilon'_k$ is well-defined and $0<\epsilon'_k<\overline{\epsilon}'$. Furthermore, $K(\epsilon'_k)=\max(\epsilon_k,\tau_k)\times|\Omega|M_v^2m_{\Lambda}^{-2}\to 0$ as $k\to+\infty$, whence
$$\epsilon'_k\to 0^+\hbox{ as }k\to+\infty.$$\par
Fix now any unit vector $e\in\R^n$. For each $k\in\N$, let $L_k$ be the function defined for all $\tau'\in[0,M_v-m_v-\epsilon'_k]$ by
\be\label{Lktau'}
L_k(\tau')=\int_{\Omega}|v_{\Omega,k,\tau'}(y)|^2\Lambda_{\Omega,k}(y)^{-1}dy,
\ee
where
\be\label{vomegaktau'}
v_{\Omega,k,\tau'}(y)=\left\{\baa{ll}
(M_v-\epsilon'_k)|v(y)|^{-1}v(y) & \hbox{if }|v(y)|>M_v-\epsilon'_k,\\
v(y) & \hbox{if }m_v+\tau'\le|v(y)|\le M_v-\epsilon'_k,\\
(m_v+\tau')|v(y)|^{-1}v(y) & \hbox{if }m_v<|v(y)|<m_v+\tau',\\
(m_v+\tau')|v(y)|^{-1}v(y) & \hbox{if }|v(y)|=m_v\hbox{ and }m_v>0,\\
(m_v+\tau')e & \hbox{if }|v(y)|=0\hbox{ and }m_v=0.\eaa\right.
\ee
Each function $L_k$ is Lipschitz-continuous in $[0,M_v-m_v-\epsilon'_k]$ and
$$L_k(0)=\int_{\Omega}\min\left(|v(y)|^2,(M_v-\epsilon'_k)^2\right)\Lambda_{\Omega,k}(y)^{-1}dy,$$
whence
$$\baa{rcl}
L_k(0)-\displaystyle{\int_{\Omega}}|v(y)|^2\Lambda(y)^{-1}dy & = & \displaystyle{\int_{\Omega}}\min\left(|v(y)|^2,(M_v-\epsilon'_k)^2\right)\cdot(\Lambda_{\Omega,k}(y)^{-1}-\Lambda(y)^{-1})dy-K(\epsilon'_k)\\
& \le & \displaystyle{\frac{|\Omega|M_v^2}{m_{\Lambda}^2}}\times\max(\epsilon_k,\tau_k)\ -\ K(\epsilon'_k)\\
& = & 0\eaa$$
due to the definitions of $\Lambda_{\Omega,k}$ (see (\ref{lambdaomegak})) and $m_{\Lambda}$, $M_{\Lambda}$, $M_v$ and $\epsilon'_k$. Furthermore,
$$L_k(M_v-m_v-\epsilon'_k)=(M_v-\epsilon'_k)^2\int_{\Omega}\Lambda_{\Omega,k}(y)^{-1}dy=(M_v-\epsilon'_k)^2\int_{\Omega}\Lambda(y)^{-1}dy>\int_{\Omega}|v(y)|^2\Lambda(y)^{-1}dy$$
from (\ref{lambdak-1}) and (\ref{overepsilon'}) (remember that $0<\epsilon'_k<\overline{\epsilon}'<M_v-m_v\le M_v$). Therefore, by continuity of $L_k$, the real number
\be\label{tau'k}
\tau'_k=\min\left\{\tau'\in[0,M_v-m_v-\epsilon'_k),\ L_k(\tau')=\int_{\Omega}|v(y)|^2\Lambda(y)^{-1}dy\right\}
\ee
is well-defined and $0\le\tau'_k<M_v-m_v-\epsilon'_k$. Moreover, $\tau'_k\to 0$ as $k\to+\infty$. Otherwise, up to extraction, there exists $\tau>0$ such that $\tau'_k\to\tau$ as $k\to+\infty$, and
$$L_k(\tau'_k)\mathop{\to}_{k\to+\infty}\int_{\Omega}\max\left(|v(y)|^2,(m_v+\tau)^2\right)\Lambda(y)^{-1}dy>\int_{\Omega}|v(y)|^2\Lambda(y)^{-1}dy$$
by definition of $m_v$. This is impossible by definition of $\tau'_k$.\par
Call now
$$v_{\Omega,k}=v_{\Omega,k,\tau'_k}.$$
Notice that, for each $k\in\N$, the vector field $v_{\Omega,k}$ is in $L^{\infty}(\Omega,\R^n)$ and it satisfies
$$m_v+\tau'_k\le\mathop{\hbox{ess inf}}_{\Omega}|v_{\Omega,k}|\le\mathop{\hbox{ess sup}}_{\Omega}|v_{\Omega,k}|\le M_v-\epsilon'_k$$
and
\be\label{vklambdak-1}
\int_{\Omega}|v_{\Omega,k}(y)|^2\Lambda_{\Omega,k}(y)^{-1}dy=L_k(\tau'_k)=\int_{\Omega}|v(y)|^2\Lambda(y)^{-1}dy.
\ee\par
Write now the vector fields $v$ and $v_{\Omega,k}$ as $v=(v_1,\ldots,v_n)$ and $v_{\Omega,k}=(v_{k,1},\ldots,v_{k,n})$, extend all functions $v_i$ and $v_{k,i}$ by $0$ in $\R^n\backslash\Omega$, call $\underline{v}_i$ and $\underline{v}_{k,i}$ these extensions. Set
$$\underline{v}=(\underline{v}_1,\ldots,\underline{v}_n)\hbox{ and }\underline{v}_k=(\underline{v}_{k,1},\ldots,\underline{v}_{k,n}).$$
One then has that
\be\label{underlinevk}
\|\underline{v}_k-\underline{v}\|_{L^{\infty}(\R^n,\R^n)}\to 0\hbox{ as }k\to+\infty
\ee
since $(\epsilon'_k,\tau'_k)\to(0,0)$ as $k\to+\infty$. For each $(k,k')\in\N^2$, denote
$$\underline{v}_{k,k'}=(\rho_{k'}*\underline{v}_{k,1},\ldots,\rho_{k'}*\underline{v}_{k,n}),$$
where $(\rho_{k'})_{k'\in\N}$ is the same sequence of mollifiers as in Step~1. The vector fields $\underline{v}_{k,k'}$ are then of class $C^{\infty}(\R^n,\R^n)$ and they satisfy
$$\|\underline{v}_{k,k'}\|_{L^{\infty}(\R^n,\R^n)}\le\|\underline{v}_k\|_{L^{\infty}(\R^n,\R^n)}\le M_v-\epsilon'_k<M_v=\|v\|_{\infty}\hbox{ for all }(k,k')\in\N^2.$$
For each fixed $k\in\N$, the fields $\underline{v}_{k,k'}$ converge to ${\underline{v}_k}_{|\Omega}=v_{\Omega,k}$ as $k'\to+\infty$ in all spaces $L^p(\Omega)$ for $1\le p<+\infty$. Actually, one also has that
$$\|\underline{v}_{k,k'}-v_{\Omega,k}\|_{L^p(\Omega,\R^n)}+\|\underline{v}_{k,k'}-v\|_{L^p(\Omega,\R^n)}\to 0\hbox{ as }(k,k')\to(+\infty,+\infty)$$
because of (\ref{underlinevk}). Together with (\ref{alphakk'}), (\ref{Lambdak2}) and (\ref{Lambdak2bis}), it follows that
$$\int_{\Omega}|\underline{v}_{k,k'}(y)|^2\Lambda_{k,k'}(y)^{-1}dy\to\int_{\Omega}|v(y)|^2\Lambda(y)^{-1}dy\hbox{ as }(k,k')\to(+\infty,+\infty).$$
Remember that the limit in the right-hand side is positive because $M_v>0$ here. Then the positive real numbers
$$\beta_{k,k'}=\left(\frac{\displaystyle{\int_{\Omega}}|v(y)|^2\Lambda(y)^{-1}dy}{\displaystyle{\int_{\Omega}}|\underline{v}_{k,k'}(y)|^2\Lambda_{k,k'}(y)^{-1}dy}\right)^{1/2}$$
are well-defined for $k$ and $k'$ large enough (one can then assume for all $(k,k')\in\N^2$ without loss of generality) and $\beta_{k,k'}\to 1$ as $(k,k')\to(+\infty,+\infty)$. Moreover, for each $k\in\N$,
$$\beta_{k,k'}^{-2}=\frac{\displaystyle{\int_{\Omega}}|\underline{v}_{k,k'}(y)|^2\Lambda_{k,k'}(y)^{-1}dy}{\displaystyle{\int_{\Omega}}|v(y)|^2\Lambda(y)^{-1}dy}\mathop{\to}_{k'\to+\infty}\frac{\displaystyle{\int_{\Omega}}|\underline{v}_k(y)|^2\underline{\Lambda}_k(y)^{-1}dy}{\displaystyle{\int_{\Omega}}|v(y)|^2\Lambda(y)^{-1}dy}=\frac{\displaystyle{\int_{\Omega}}|v_{\Omega,k}(y)|^2\Lambda_{\Omega,k}(y)^{-1}dy}{\displaystyle{\int_{\Omega}}|v(y)|^2\Lambda(y)^{-1}dy}=1$$
because of (\ref{limitk'}), (\ref{alphakk'bis}) and (\ref{vklambdak-1}).\par
Set
$$v_{k,k'}(y)=\beta_{k,k'}\underline{v}_{k,k'}(y)\hbox{ for all }y\in\overline{\Omega}\hbox{ and }(k,k')\in\N^2.$$
The vector fields $v_{k,k'}$ are in $C^{\infty}(\overline{\Omega},\R^n)$ and they satisfy
\be\label{vk1}
\int_{\Omega}|v_{k,k'}(y)|^2\Lambda_{k,k'}(y)^{-1}dy=\int_{\Omega}|v(y)|^2\Lambda(y)^{-1}dy\hbox{ for all }(k,k')\in\N^2
\ee
and
\be\label{vk2}\left\{\baa{l}
\|v_{k,k'}\|_{L^{\infty}(\Omega,\R^n)}\le\beta_{k,k'}\|\underline{v}_{k,k'}\|_{L^{\infty}(\R^n,\R^n)}\le\beta_{k,k'}\times(\|v\|_{\infty}-\epsilon'_k)\hbox{ for all }(k,k')\in\N^2,\\
\beta_{k,k'}\times(\|v\|_{\infty}-\epsilon'_k)\to\|v\|_{\infty}-\epsilon'_k\hbox{ as }k'\to+\infty\hbox{ for all }k\in\N,\\
\|v_{k,k'}-v\|_{L^p(\Omega,\R^n)}\to 0\hbox{ as }(k,k')\to(+\infty,+\infty)\hbox{ for all }1\le p<+\infty.\eaa\right.
\ee\par
{\bf{Consider now the case when $m_v<M_v$ and $m_{\Lambda}=M_{\Lambda}$.}} Namely, up to modification on a zero-measure set, $\Lambda$ is constant. Choose $\overline{\epsilon}'\in(0,M_v-m_v)$ such that (\ref{overepsilon'}) holds, namely
\be\label{overepsilon'bis}
(M_v-\overline{\epsilon}')^2|\Omega|>\int_{\Omega}|v(y)|^2dy.
\ee
Take any sequence $(\epsilon'_k)_{k\in\N}$ in $(0,\overline{\epsilon}')$ such that $\epsilon'_k\to 0$ as $k\to+\infty$. For each $k\in\N$ and $\tau'\in[0,M_v-m_v-\epsilon'_k]$, define $L_k(\tau')$ as in (\ref{Lktau'}) with $\Lambda_{\Omega,k}=\Lambda$. Each function $L_k$ is Lipschitz-continuous. Moreover,
$$L_k(0)=\Lambda^{-1}\int_{\Omega}\min\left(|v(y)|^2,(M_v-\epsilon'_k)^2\right)dy<\Lambda^{-1}\int_{\Omega}|v(y)|^2dy=\int_{\Omega}|v(y)|^2\Lambda(y)^{-1}dy$$
since $0\le M_v-\epsilon'_k<M_v$, and
$$\baa{rcl}
L_k(M_v-m_v-\epsilon'_k)=(M_v-\epsilon'_k)^2\Lambda^{-1}|\Omega| & > & (M_v-\overline{\epsilon}')^2\Lambda^{-1}|\Omega|\\
& > & \Lambda^{-1}\displaystyle{\int_{\Omega}}|v(y)|^2dy=\displaystyle{\int_{\Omega}}|v(y)|^2\Lambda(y)^{-1}dy\eaa$$
from (\ref{overepsilon'bis}). Therefore, the real numbers $\tau'_k$ given in (\ref{tau'k}) are well-defined and are such that $0<\tau'_k<M_v-m_v-\epsilon'_k$ for each $k\in\N$. We then keep the same definitions of $v_{\Omega,k}$, $\underline{v}$, $\underline{v}_k$, $\underline{v}_{k,k'}$, $\beta_{k,k'}$ and $v_{k,k'}$ as above and properties (\ref{vk1}) and (\ref{vk2}) hold.\par
{\bf{Lastly, in the case when $m_v=M_v=\|v\|_{\infty}$}}, namely when $|v|$ is equal to the constant $\|v\|_{\infty}$ almost everywhere, then $v$ is kept unchanged. We set $\underline{v}_{k,k'}=v_{k,k'}=v$, $\beta_{k,k'}=1$, $\epsilon'_k=0$, and properties (\ref{vk1}) and (\ref{vk2}) hold.\hfill\break\par

{\underline{Step 4: Approximation of the eigenvalue $\lambda_1(\Omega,A,v,V)$}}. {\bf{Consider first the case when $m_v<M_v=\|v\|_{\infty}$}}. For each $k\in\N$, it follows from (\ref{Lambdak1}), (\ref{Lambdak2}), (\ref{Ak3}), (\ref{vk1}) and (\ref{vk2}) that there exists an integer $k'(k)\ge k$ such that the $C^{\infty}$ fields
$$A_k=(a_{k,i,j})_{1\le i,j\le n}=A_{k,k'(k)},\ \Lambda_k=\Lambda_{k,k'(k)}\hbox{ and }v_k=v_{k,k'(k)}$$
satisfy, for all $k\in\N$,
\be\label{Akvk1}\left\{\baa{l}
\displaystyle{\int_{\Omega}}\Lambda_k(y)^{-1}dy=\displaystyle{\int_{\Omega}}\Lambda(y)^{-1}dy,\\
0<m_{\Lambda}\le m_{\Lambda}+\displaystyle{\frac{\tau_k}{2}}\le\Lambda_k\le M_{\Lambda}-\displaystyle{\frac{\epsilon_k}{2}}\le M_{\Lambda}\hbox{ in }\overline{\Omega},\\
A_k(y)\ge\Lambda_k(y)\hbox{Id}\hbox{ for all }y\in\overline{\Omega},\\
\displaystyle{\int_{\Omega}}|v_k(y)|^2\Lambda_k(y)^{-1}dy=\displaystyle{\int_{\Omega}}|v(y)|^2\Lambda(y)^{-1}dy,\\
\|v_k\|_{L^{\infty}(\Omega,\R^n)}\le\|v\|_{L^{\infty}(\Omega,\R^n)}-\displaystyle{\frac{\epsilon'_k}{2}}\le\|v\|_{L^{\infty}(\Omega,\R^n)}.\eaa\right.
\ee
Notice that this is possible in both cases $m_{\Lambda}<M_{\Lambda}$ (then, $\epsilon_k$ and $\tau_k$ are all positive) or $m_{\Lambda}=M_{\Lambda}$ (then, $\epsilon_k=\tau_k=0$ for all $k\in\N$). Notice also that $\epsilon'_k>0$ for all $k\in\N$ in this case when $m_v<M_v$. Furthermore, the matrix field $A_k$ is symmetric with polynomial entries $a_{k,i,j}$ in $\overline{\Omega}$, and, by (\ref{Ak1}), (\ref{Ak2}) and (\ref{vk2}),
\be\label{Akvk2}\left\{\baa{l}
a_{k,i,j}\to a_{i,j}\hbox{ uniformly in }\overline{\Omega}\hbox{ as }k\to+\infty,\\
\displaystyle{\mathop{\limsup}_{k\to+\infty}}\|\nabla a_{k,i,j}\|_{L^{\infty}(\Omega)}\le\|\nabla a_{i,j}\|_{L^{\infty}(\Omega)}\le\|A\|_{W^{1,\infty}(\Omega)}\hbox{ for all }1\le i,j\le n,\\
\|v_k-v\|_{L^p(\Omega,\R^n)}\to 0\hbox{ as }k\to+\infty\hbox{ for all }1\le p<+\infty.\eaa\right.
\ee
{\bf{In the case when $m_v=M_v=\|v\|_{\infty}$}}, then one sets $v_k=v$ for all $k\in\N$, and properties (\ref{Akvk1}) and (\ref{Akvk2}) still hold (with $\epsilon'_k=0$ in that case).\par
Let us now prove that
\be\label{claimlambdak}
\lambda_1(\Omega,A_k,v_k,V)\to\lambda_1(\Omega,A,v,V)\hbox{ as }k\to+\infty.
\ee
Notice first that each operator $-\div(A_k\nabla)+v_k\cdot\nabla+V$ is elliptic because of (\ref{Akvk1}). Fix an open non-empty ball $B$ such that $B\subset\Omega$. It follows from \cite{bnv} that
\be\label{OmegaB}
\min_{\overline{\Omega}}V<\lambda_1(\Omega,A_k,v_k,V)\le\lambda_1(B,A_k,v_k,V)\hbox{ for all }k\in\N.
\ee
Furthermore, properties (\ref{Akvk1}) and (\ref{Akvk2}) imply that the sequences of matrix fields $(A_k)_{k\in\N}$ and $(A_k^{-1})_{k\in\N}$ are bounded in $W^{1,\infty}(\Omega)$ (and then in $W^{1,\infty}(B)$), and that the sequence of vector fields $(v_k)_{k\in\N}$ is bounded in $L^{\infty}(\Omega)$, and then in $L^{\infty}(B)$. From Lemma~1.1 in \cite{bnv}, there exists then a constant $C$ independent from $k$ such that
$$\lambda_1(B,A_k,v_k,V)\le C\hbox{ for all }k\in\N.$$
Together with (\ref{OmegaB}), it resorts that the sequence $(\lambda_1(\Omega,A_k,v_k,V))_{k\in\N}$ is bounded. Thus, for a sequence of integers $n(k)\to+\infty$ as $k\to+\infty$, one has that
$$\lambda_1(\Omega,A_{n(k)},v_{n(k)},V)\to\lambda\in\R\hbox{ as }k\to+\infty.$$
For each $k\in\N$, call $\varphi_k$ the principal eigenfunction of the operator $-\hbox{div}(A_k\nabla)+v_k\cdot\nabla+V$ in $\Omega$ with Dirichlet boundary condition, such that $\max_{\overline{\Omega}}\varphi_k=1$. Namely, each function $\varphi_k$ satisfies
\be\label{varphik}\left\{\baa{rcll}
-\hbox{div}(A_k\nabla\varphi_k)+v_k\cdot\nabla\varphi_k+V\varphi_k & = & \lambda_1(\Omega,A_k,v_k,V)\varphi_k & \hbox{in }\Omega,\\
\varphi_k & > & 0 & \hbox{in }\Omega,\\
\|\varphi_k\|_{L^{\infty}(\Omega)} & = & 1, & \\
\varphi_k & = & 0 & \hbox{on }\partial\Omega.\eaa\right.
\ee
From standard elliptic estimates, each function $\varphi_k$ is in  $W^{2,p}(\Omega)$ for all $1\le p<+\infty$ and in $C^{1,\alpha}(\overline{\Omega})$ for all $0\le\alpha<1$. Furthermore, since the eigenvalues $\lambda_1(\Omega,A_k,v_k,V)$ are bounded and $\|\varphi_k\|_{L^{\infty}(\Omega)}= 1$, it follows from (\ref{Akvk1}) and (\ref{Akvk2}) that the sequence $(\varphi_k)_{k\in\N}$ is bounded in all $W^{2,p}(\Omega)$ and $C^{1,\alpha}(\overline{\Omega})$, for all $1\le p<+\infty$ and $0\le\alpha<1$. Up to extraction of another subsequence, one can assume that there exists $\varphi_{\infty}\in\displaystyle{\mathop{\bigcap}_{1\le p<+\infty}}W^{2,p}(\Omega)$ such that
$$\varphi_{n(k)}\to\varphi_{\infty}\hbox{ as }k\to+\infty,\hbox{ weakly in }W^{2,p}(\Omega)\hbox{ and strongly in }C^{1,\alpha}(\overline{\Omega})$$
for all $1\le p<+\infty$ and $0\le\alpha<1$. Notice that (\ref{Akvk2}) implies that $\partial_{x_{i'}}a_{k,i,j}\rightharpoonup\partial_{x_{i'}}a_{i,j}$ in $\sigma(L^{\infty}(\Omega),L^1(\Omega))$ as $k\to+\infty$, for all $1\le i',i,j\le n$. Multiplying (\ref{varphik}) for $n(k)$ by any test function in ${\mathcal{D}}(\Omega)$, integrating over $\Omega$ and passing to the limit as $k\to+\infty$ leads to, because of (\ref{Akvk1}) and (\ref{Akvk2}),
$$-\hbox{div}(A\nabla\varphi_{\infty})+v\cdot\nabla\varphi_{\infty}+V\varphi_{\infty}=\lambda\varphi_{\infty}\hbox{ in }\Omega,$$
together with
$$\varphi_{\infty}=0\hbox{ on }\partial\Omega,\quad\varphi_{\infty}\ge 0\hbox{ in }\overline{\Omega}\hbox{ and }\max_{\overline{\Omega}}\varphi_{\infty}=1.$$
The strong maximum principle and the characterization of the principal eigenvalue and eigenfunction thanks to Krein-Rutman theory imply that $\lambda=\lambda_1(\Omega,A,v,V)$ and $\varphi_{\infty}$ is the principal eigenfunction of the operator $-\hbox{div}(A\nabla)+v\cdot\nabla+V$ in $\Omega$ with Dirichlet boundary condition. The limiting function $\varphi_{\infty}$ is uniquely determined because of the normalization $\max_{\overline{\Omega}}\varphi_{\infty}=1$. Since the limits $\lambda=\lambda_1(\Omega,A,v,V)$ and $\varphi_{\infty}$ do not depend on any subsequence, one concludes that (\ref{claimlambdak}) holds and that the whole sequence $(\varphi_k)_{k\in\N}$ converges to $\varphi_{\infty}$ as $k\to+\infty$, weakly in $W^{2,p}(\Omega)$ and strongly in $C^{1,\alpha}(\overline{\Omega})$ for all $1\le p<+\infty$ and $0\le\alpha<1$.\hfill\break\par

{\underline{Step 5: Approximation of the principal eigenfunction of $\varphi_K$ for a large $K$}}. Choose any arbitrary $\epsilon>0$. Let $\epsilon'\in(0,1)$ be such that
\be\label{epsilon'}
\frac{\lambda_1(\Omega,A,v,V)+3\epsilon'+\epsilon'm_{\Lambda}^{-1}+\epsilon'\|V^-\|_{\infty}}{1-\epsilon'}\le\lambda_1(\Omega,A,v,V)+\epsilon.
\ee
Thanks to (\ref{claimlambdak}), there exists then $K\in\N$ such that
\be\label{lambdaK}
\lambda_1(\Omega,A_K,v_K,V)\le\lambda_1(\Omega,A,v,V)+\epsilon'.
\ee
Remember that (\ref{Akvk1}) holds with $k=K$. Let $\varphi_K$ be the (unique) solution (\ref{varphik}) with $k=K$. Notice that
$$\nu\cdot\nabla\varphi_K=-|\nabla\varphi_K|<0\hbox{ on }\partial\Omega$$
from Hopf lemma.\par
Call $\mathcal{F}$, $\underline{\mathcal{F}}$, $\overline{\mathcal{F}}$, $\underline{f}$ and $\overline{f}$ the functions defined in $\overline{\Omega}$ by
$$\left\{\baa{l}
{\mathcal{F}}=-V\varphi_K+\lambda_1(\Omega,A_K,v_K,V)\varphi_K,\\
\underline{\mathcal{F}}=-v_K\cdot\nabla\varphi_K,\ \overline{\mathcal{F}}=|v_K|\times|\nabla\varphi_K|,\\ \underline{f}=\underline{\mathcal{F}}+{\mathcal{F}}=-v_K\cdot\nabla\varphi_K-V\varphi_K+\lambda_1(\Omega,A_K,v_K,V)\varphi_K,\\
\overline{f}=\overline{\mathcal{F}}+{\mathcal{F}}=|v_K|\times|\nabla\varphi_K|-V\varphi_K+\lambda_1(\Omega,A_K,v_K,V)\varphi_K.\eaa\right.$$
The function $\mathcal{F}$ is continuous in $\overline{\Omega}$. There exists then a sequence $({\mathcal{F}}_l)_{l\in\N}$ of polynomials such that
$${\mathcal{F}}_l\to{\mathcal{F}}\hbox{ uniformly in }\overline{\Omega}\hbox{ as }l\to+\infty.$$
Observe also that the function $\underline{\mathcal{F}}$ is in $L^{\infty}(\Omega)$, and that the function $\overline{\mathcal{F}}$ is nonnegative and continuous in $\overline{\Omega}$: this is true if $m_v<M_v$ because $v_K$ is then actually of class $C^{\infty}(\overline{\Omega})$ and $\varphi_K\in C^1(\overline{\Omega})$; this is also true if $m_v=M_v$ because $v_K$ is then equal to $v$ and $|v_K|=M_v$ in $\overline{\Omega}$ up to a modification on a zero-measure set. Let $R_0>0$ be such that
$$\overline{\Omega}\subset B_{R_0},$$
where $B_{R_0}$ is the open Euclidean ball of radius $R_0$ and center $0$. Denote by $\overline{\mathcal{F}}$ (with a slight abuse of notation) a continuous extension of $\overline{\mathcal{F}}$ in $\R^n$ such that $\overline{\mathcal{F}}\ge 0$ in $\R^n$ and $\overline{\mathcal{F}}=0$ in $\R^n\backslash B_{R_0}$. Extend by $0$ the function $\underline{\mathcal{F}}$ in $\R^n\backslash\overline{\Omega}$ and still call $\underline{\mathcal{F}}$ this extension. Notice that
\be\label{mathcalF}
\underline{\mathcal{F}}\le\overline{\mathcal{F}}\hbox{ in }\R^n.
\ee
For each $l\in\N$, call $\zeta_l$ the function defined in $\overline{B_{2R_0}}$ by
$$\zeta_l(z)=\frac{\left[(2R_0)^2-|z|^2\right]^l}{\displaystyle{\int_{B_{2R_0}}}\left[(2R_0)^2-|z'|^2\right]^ldz'}.$$
Extend the functions $\zeta_l$ by $0$ outside $B_{2R_0}$ and still call $\zeta_l$ these extensions. In $\overline{\Omega}$, define the functions
$$\underline{\mathcal{F}}_l=\zeta_l*\underline{\mathcal{F}}\hbox{ and }\overline{\mathcal{F}}_l=\zeta_l*\overline{\mathcal{F}}\hbox{ for all }l\in\N.$$
Owing to the choices of $R_0$ and $\zeta_l$, the functions $\underline{\mathcal{F}}_l$ and $\overline{\mathcal{F}}_l$ are polynomials in $\overline{\Omega}$. Furthermore, since $\underline{\mathcal{F}}\in L^{\infty}(\R^n)$ and $\overline{\mathcal{F}}\in C(\R^n)$, there holds
$$\|\underline{\mathcal{F}}_l-\underline{\mathcal{F}}\|_{L^p(\Omega)}\mathop{\to}_{l\to+\infty}0\hbox{ for all }1\le p<+\infty\hbox{ and }\overline{\mathcal{F}}_l\mathop{\to}_{l\to+\infty}\overline{\mathcal{F}}\hbox{ uniformly in }\overline{\Omega}.$$
Furthermore,
$$\underline{\mathcal{F}}_l=\zeta_l*\underline{\mathcal{F}}\le\zeta_l*\overline{\mathcal{F}}=\overline{\mathcal{F}}_l\hbox{ in }\overline{\Omega}\hbox{ for all }l\in\N$$
because of (\ref{mathcalF}). It follows that
\be\label{fl1}
\underline{f}_l=\underline{\mathcal{F}}_l+{\mathcal{F}}_l\mathop{\to}_{l\to+\infty}\underline{f}\hbox{ in }L^p(\Omega)\hbox{ for all }1\le p<+\infty,\ \ \overline{f}_l=\overline{\mathcal{F}}_l+{\mathcal{F}}_l\mathop{\to}_{l\to+\infty}\overline{f}\hbox{ uniformly in }\overline{\Omega}
\ee
together with
\be\label{fl2}
\underline{f}_l\le\overline{f}_l\hbox{ in }\overline{\Omega}\hbox{ for all }l\in\N.
\ee
Remember that the function $\varphi_K$ satisfies
$$-\div(A_K\nabla\varphi_K)=\underline{f}\hbox{ in }\Omega,$$
with $\varphi_K=0$ on $\partial\Omega$. For each $l\in\N$, call $\psi_l$ the solution of
\be\label{psil}\left\{\baa{rcll}
-\hbox{div}(A_K\nabla\psi_l) & = & \underline{f}_l & \hbox{in }\Omega,\\
\psi_l & = & 0 & \hbox{on }\partial\Omega.\eaa\right.
\ee
Each function $\psi_l$ is then analytic in $\Omega$ (remember that $A_K$ is a field of symmetric positive definite matrices with polynomial entries, and that each $\underline{f}_l$ is a polynomial in $\overline{\Omega}$). From standard elliptic estimates, the functions $\psi_l$ converge to the function $\varphi_K$ as $l\to+\infty$ in $W^{2,p}(\Omega)$ and $C^{1,\alpha}(\overline{\Omega})$ for all $1\le p<+\infty$ and $0\le\alpha<1$. Since $\varphi_K>0$ in $\Omega$ and $|\nabla\varphi_K|>0=\varphi_K$ on $\partial\Omega$, one then has
\be\label{psil+}
\psi_l>0\hbox{ in }\Omega\hbox{ and }|\nabla\psi_l|>0\hbox{ on }\partial\Omega\hbox{ for }l\hbox{ large enough}.
\ee
Furthermore, from (\ref{fl2}), there holds
$$\baa{rl}
-\hbox{div}(A_K\nabla\psi_l)-|v_K|\times|\nabla\psi_l|-\epsilon'|\nabla\psi_l|+[V-\lambda_1(\Omega,A,v,V)-2\epsilon']\ \psi_l & \\
\qquad\qquad\le\overline{f}_l-|v_K|\times|\nabla\psi_l|-\epsilon'|\nabla\psi_l|+[V-\lambda_1(\Omega,A,v,V)-2\epsilon']\ \psi_l & \hbox{ in }\Omega.\eaa$$
From (\ref{fl1}) and the definition of $\overline{f}$, it follows that
$$\baa{l}
\overline{f}_l-|v_K|\times|\nabla\psi_l|-\epsilon'|\nabla\psi_l|+[V-\lambda_1(\Omega,A,v,V)-2\epsilon']\ \psi_l\\
\to[\lambda_1(\Omega,A_K,v_K,V)-\lambda_1(\Omega,A,v,V)-\epsilon']\ \varphi_K-\epsilon'(|\nabla\varphi_K|+\varphi_K)\eaa$$
as $l\to+\infty$, uniformly in $\overline{\Omega}$. Using (\ref{lambdaK}) and the properties of $\varphi_K$, it follows that the continuous function
$$[\lambda_1(\Omega,A_K,v_K,V)-\lambda_1(\Omega,A,v,V)-\epsilon']\ \varphi_K-\epsilon'(|\nabla\varphi_K|+\varphi_K)$$
is negative in $\overline{\Omega}$. Therefore, there is $L\in\N$ large enough so that (\ref{psil+}) holds with $l=L$ and
\be\label{psilbis}
-\hbox{div}(A_K\nabla\psi_L)-|v_K|\times|\nabla\psi_L|-\epsilon'|\nabla\psi_L|+[V-\lambda_1(\Omega,A,v,V)-2\epsilon']\ \psi_L<0\hbox{ in }\Omega.
\ee\par

{\underline{Step 6: An inequality for the rearranged fields in the ball $\Omega^*$}}. Apply then the results of Section~\ref{rearrangement} to the function
$$\psi=\psi_L$$
and to the data
$$(A_{\Omega},\Lambda_{\Omega},\omega,V)=(A_K,\Lambda_K,|v_K|,V).$$
From the previous steps, these fields satisfy all assumptions of Section~\ref{rearrangement}. Given $\psi$, one can then define the sets $Z$, $Y$, $E$ and the function $\rho$ as in Section~\ref{sec21}. Given $\psi$ and the data $A_{\Omega}=A_K$, $\Lambda_{\Omega}=\Lambda_K$, $\omega=|v_K|$ and $V$, one can also define the corresponding rearranged fields $\tilde{\psi}$, $\hat{\Lambda}$, $\hat{v}$, $\hat{V}$ and $U$ given by (\ref{G}), (\ref{defhatLambda}), (\ref{defF}), (\ref{defpsi}), (\ref{defhatv}), (\ref{defhatV}), (\ref{H}) and~(\ref{U}).\par
One recalls that $\lambda_1(\Omega,A,v,V)\ge 0$ by assumption. From (\ref{psilbis}) and Corollary~\ref{corintegral} applied with $\omega_0=\epsilon'\ge 0$ and
$$\mu=\lambda_1(\Omega,A,v,V)+2\epsilon'\ge 0,$$
it follows that
\be\label{tildepsiU}
\int_{\Omega^*}\left[\widehat{\Lambda}(x)|\nabla\widetilde{\psi}(x)|^2-\epsilon'|\nabla\widetilde{\psi}(x)|\widetilde{\psi}(x)+\widehat{V}(x)\widetilde{\psi}(x)^2\right]e^{-U(x)}dx\le\mu\int_{\Omega^*}\widetilde{\psi}(x)^2e^{-U(x)}dx.
\ee
Remember that the function $\widetilde{\psi}$ is radially symmetric, continuous and decreasing with respect to $|x|$ in $\overline{\Omega^*}$, and that $\widetilde{\psi}\in H^1_0(\Omega^*)$. The field $\widehat{\Lambda}$ is radially symmetric and belongs in $L^{\infty}(\Omega^*)$. Furthermore,
\be\label{Lambdak4}
\int_{\Omega^*}\widehat{\Lambda}(x)^{-1}dx=\int_{\Omega}\Lambda_K(y)^{-1}dy=\int_{\Omega}\Lambda(y)^{-1}dy
\ee
from (\ref{intLambda}) and (\ref{Akvk1}). On the other hand,
\be\label{Lambdak5}
0<m_{\Lambda}\le m_{\Lambda}+\frac{\tau_K}{2}\le\min_{\overline{\Omega}}\Lambda_K\le\mathop{\hbox{ess inf}}_{\Omega^*}\widehat{\Lambda}\le\mathop{\hbox{ess sup}}_{\Omega^*}\widehat{\Lambda}\le\max_{\overline{\Omega}}\Lambda_K\le M_{\Lambda}-\frac{\epsilon_K}{2}\le M_{\Lambda}
\ee
from (\ref{infiniLambda}) and (\ref{Akvk1}).\par
From (\ref{tildepsiU}) and (\ref{Lambdak5}), it follows that
$$(1-\epsilon')\displaystyle{\int_{\Omega^*}}\widehat{\Lambda}|\nabla\widetilde{\psi}|^2e^{-U}+\displaystyle{\int_{\Omega^*}}\widehat{V}\widetilde{\psi}^2e^{-U}\le \mu\displaystyle{\int_{\Omega^*}}\widetilde{\psi}^2e^{-U}+\displaystyle{\frac{\epsilon'}{4}}\displaystyle{\int_{\Omega^*}}\widehat{\Lambda}^{-1}\widetilde{\psi}^2e^{-U}\le(\mu+\epsilon' m_{\Lambda}^{-1})\displaystyle{\int_{\Omega^*}}\widetilde{\psi}^2e^{-U}.$$
On the other hand, the field $\widehat{V}$ is radially symmetric, it belongs to $L^{\infty}(\Omega^*)$, and
$$-\|V\|_{\infty}\le\min_{\overline{\Omega}}(-V^-)\le\mathop{\hbox{ess inf}}_{\Omega^*}\widehat{V}\le\mathop{\hbox{ess sup}}_{\Omega^*}\widehat{V}\le 0$$
(see (\ref{infiniV})). Therefore,
$$(1-\epsilon')\displaystyle{\int_{\Omega^*}}
\left(\widehat{\Lambda}|\nabla\widetilde{\psi}|^2+\widehat{V}\widetilde{\psi}^2\right)e^{-U}\le(\mu+\epsilon' m_{\Lambda}^{-1}+\epsilon'\|V^-\|_{\infty}) \displaystyle{\int_{\Omega^*}}\widetilde{\psi}^2e^{-U},$$
that is
\be\label{ineqpsi}
\displaystyle{\int_{\Omega^*}}
\left(\widehat{\Lambda}|\nabla\widetilde{\psi}|^2+\widehat{V}\widetilde{\psi}^2\right)e^{-U}\le\frac{\lambda_1(\Omega,A,v,V)+2\epsilon'+\epsilon' m_{\Lambda}^{-1}+\epsilon'\|V^-\|_{\infty}}{1-\epsilon'}\times\displaystyle{\int_{\Omega^*}}\widetilde{\psi}^2e^{-U}.
\ee\par

{\underline{Step 7: Approximation of the rearranged fields in $\Omega^*$}}. First, define the function $\underline{\hat{\Lambda}}$ almost everywhere in $\R^n$ by:
$$\underline{\hat{\Lambda}}(x)=\left\{\baa{ll}\hat{\Lambda}(x) & \hbox{if }x\in\Omega^*,\\ \displaystyle{\mathop{\min}_{\overline{\Omega}}}\ \Lambda_K & \hbox{if }x\in\R^n\backslash\Omega^*,\eaa\right.$$
and then, for each $m\in\N$,
$$\underline{\Lambda}^*_m=\left(\rho_m*\underline{\hat{\Lambda}}^{-1}\right)^{-1}\hbox{ in }\R^n,$$
where the mollifiers $\rho_m$ can be assumed to be radially symmetric for all $m\in\N$. Next, for every $m\in\N$, call
$$\gamma_m=\frac{\|(\underline{\Lambda}^*_m)^{-1}\|_{L^1(\Omega^*)}}{\|\Lambda^{-1}\|_{L^1(\Omega)}}>0$$
and
$$\Lambda^*_m(x)=\gamma_m\underline{\Lambda}^*_m(x)\hbox{ for all }x\in\overline{\Omega^*}.$$
As in Step~1, it follows from the above definitions and from (\ref{Lambdak4}) and (\ref{Lambdak5}) that each function $\Lambda^*_m$ is radially symmetric and of class $C^{\infty}(\overline{\Omega^*})$, that $\lim_{m\to+\infty}\gamma_m=1$, that
\be\label{Lambdam1}
\int_{\Omega^*}\Lambda^*_m(x)^{-1}dx=\int_{\Omega}\Lambda(y)^{-1}dy\hbox{ for all }m\in\N
\ee
and that
\be\label{Lambdam2}\left\{\baa{l}
0<\gamma_m\times\left(m_{\Lambda}+\displaystyle{\frac{\tau_K}{2}}\right)\le\displaystyle{\mathop{\min}_{\overline{\Omega^*}}}\ \Lambda^*_m\le\displaystyle{\mathop{\max}_{\overline{\Omega^*}}}\ \Lambda^*_m\le\gamma_m\times\left(M_{\Lambda}-\displaystyle{\frac{\epsilon_K}{2}}\right)\hbox{ for all }m\in\N,\\
\gamma_m\times\left(m_{\Lambda}+\displaystyle{\frac{\tau_K}{2}}\right)\displaystyle{\mathop{\to}_{m\to+\infty}}m_{\Lambda}+\displaystyle{\frac{\tau_K}{2}},\ \gamma_m\times\left(M_{\Lambda}-\displaystyle{\frac{\epsilon_K}{2}}\right)\displaystyle{\mathop{\to}_{m\to+\infty}}M_{\Lambda}-\displaystyle{\frac{\epsilon_K}{2}}.\eaa\right.
\ee
From reciprocal of Lebesgue's theorem, one can also assume without loss of generality (even if it means extracting a subsequence) that
$$\Lambda^*_m(x)\to\underline{\hat{\Lambda}}(x)=\hat{\Lambda}(x)\hbox{ as }m\to+\infty\hbox{ for a.e. }x\in\Omega^*.$$
Remember also that if $m_{\Lambda}<M_{\Lambda}$ then $\epsilon_K>0$, $\tau_K>0$ and notice that if $m_{\Lambda}=M_{\Lambda}$ then $\epsilon_K=\tau_K=0$, $\gamma_m=1$, $\Lambda^*_m=m_{\Lambda}$ in $\overline{\Omega^*}$ for all $m\in\N$ and properties (\ref{Lambdam1}) and (\ref{Lambdam2}) hold immediately.\par
Next, owing to the definition (\ref{defhatv}), the vector field $\hat{v}$ can be written as
$$\hat{v}(x)=|\hat{v}(x)|\ e_r(x)\hbox{ in }\Omega^*,$$
where $|\hat{v}|$ is radially symmetric. Furthermore, as in Step~3 and since $\omega=|v_K|$, it follows from (\ref{infiniv}), (\ref{intv2}) and (\ref{Akvk1}) that
\be\label{vk4}
\int_{\Omega^*}|\hat{v}(x)|^2\hat{\Lambda}(x)^{-1}dx=\int_{\Omega}|v_K(y)|^2\Lambda_K(y)^{-1}dy=\int_{\Omega}|v(y)|^2\Lambda(y)^{-1}dy
\ee
and
\be\label{vk5}
\|\hat{v}\|_{L^{\infty}(\Omega^*,\R^n)}\le\|\omega\|_{L^{\infty}(\Omega)}=\|v_K\|_{L^{\infty}(\Omega,\R^n)}\le\|v\|_{L^{\infty}(\Omega,\R^n)}-\frac{\epsilon'_K}{2}\le\|v\|_{L^{\infty}(\Omega,\R^n)}.
\ee
Call
$$\underline{\omega}(x)=\left\{\baa{ll}|\hat{v}(x)| & \hbox{if }x\in\Omega^*,\\ 0 & \hbox{if }x\in\R^n\backslash\Omega^*\eaa\right.$$
and, for each $m\in\N$, 
$$\underline{v}_m^*(x)=(\rho_m*\underline{\omega})(x)\ e_r(x)\hbox{ for all }x\in\R^n\backslash\{0\}.$$
One has $\|\underline{v}^*_m\|_{L^{\infty}(\R^n,\R^n)}\le\|\underline{\omega}\|_{L^{\infty}(\R^n,\R^n)}\le\|v\|_{L^{\infty}(\Omega,\R^n)}$ and one can assume without loss of generality that
$$\underline{v}_m^*(x)\to\hat{v}(x)\hbox{ as }m\to+\infty\hbox{ for a.e. }x\in\Omega^*.$$
Consider first the case when $\|v\|_{L^{\infty}(\Omega,\R^n)}>0$. Therefore, for $m$ large enough (one can then assume that this holds for all $m$ without loss of generality), the real numbers
\be\label{deltam}
\delta_m=\left(\frac{\displaystyle{\int_{\Omega}}|v(y)|^2\Lambda(y)^{-1}dy}{\displaystyle{\int_{\Omega^*}}|\underline{v}_m^*(x)|^2\Lambda_m^*(x)^{-1}dx}\right)^{1/2}=\left(\frac{\displaystyle{\int_{\Omega^*}}|\hat{v}(x)|^2\hat{\Lambda}(x)^{-1}dx}{\displaystyle{\int_{\Omega^*}}|\underline{v}_m^*(x)|^2\Lambda_m^*(x)^{-1}dx}\right)^{1/2}
\ee
are well-defined, positive, and they are such that $\delta_m\to 1$ as $m\to+\infty$. Therefore, the vector fields defined by
$$v^*_m(x)=\delta_m\underline{v}^*_m(x)\hbox{ for all }x\in\overline{\Omega^*}\backslash\{0\}$$
are of class $C^{\infty}(\overline{\Omega^*}\backslash\{0\})$ and converge to $\hat{v}(x)$ as $m\to+\infty$ for almost every $x\in\Omega^*$. These fields can be written as
$$v_m^*(x)=|v_m^*(x)|e_r(x)\hbox{ for all }x\in\overline{\Omega^*}\backslash\{0\}$$
and $|v^*_m|$ is radially symmetric, of class $C^{\infty}(\overline{\Omega^*}\backslash\{0\})$ and can actually be extended at $0$ to a $C^{\infty}(\overline{\Omega^*})$ function. Furthermore, it follows from (\ref{vk5}) and (\ref{deltam}) that
\be\label{vm1}
\int_{\Omega^*}|v^*_m(x)|^2\Lambda^*_m(x)^{-1}dx=\int_{\Omega}|v(y)|^2\Lambda(y)^{-1}dy
\ee
and
\be\label{vm2}
\|v^*_m\|_{L^{\infty}(\Omega^*,\R^n)}\le\delta_m\|\hat{v}\|_{L^{\infty}(\Omega^*,\R^n)}\le\delta_m\times\left(\|v\|_{L^{\infty}(\Omega,\R^n)}-\frac{\epsilon'_K}{2}\right)\displaystyle{\mathop{\to}_{m\to+\infty}} \|v\|_{L^{\infty}(\Omega,\R^n)}-\frac{\epsilon'_K}{2}.\ee
Lastly, remember that if $m_v<M_v$, then $\epsilon'_K>0$, and notice that if $m_v=M_v$ (this is the case if $\|v\|_{L^{\infty}(\Omega,\R^n)}=0$), then $\epsilon'_K=0$, $\omega=m_v$ in $\Omega$, $\hat{v}=v^*_m=m_ve_r$ in $\Omega^*$ and properties (\ref{vm1}) and (\ref{vm2}) hold immediately with $\delta_m=1$.\par
Fix now an arbitrary unit vector $e$ in $\R^n$ and define, for each $m\in\N$,
$$\forall\ x\in\overline{\Omega^*},\quad U^*_m(x)=\int_0^{|x|}|v^*_m(re)|\Lambda^*_m(re)^{-1}dr.$$
As in Proposition~\ref{integralcompar}, the definition of $U^*_m$ does not depend on the choice of $e$. Furthermore, each function $U^*_m$ is continuous in $\overline{\Omega^*}$,  radially symmetric, of class $C^{\infty}(\overline{\Omega^*}\backslash\{0\})$ and it satisfies
\be\label{nablaUm}
\nabla U^*_m(x)=\Lambda^*_m(x)^{-1}|v^*_m(x)|e_r(x)=\Lambda^*_m(x)^{-1}v^*_m(x)\hbox{ for all }x\in\overline{\Omega^*}\backslash\{0\}.
\ee
On the other hand, each function $U^*_m$ is nonnegative in $\overline{\Omega^*}$ and it follows from (\ref{Lambdam2}) and (\ref{vm2}) that
\be\label{Um1}
\|U^*_m\|_{L^{\infty}(\Omega^*)}\le\delta_mm_{\Lambda}^{-1}R\times\left(\|v\|_{L^{\infty}(\Omega,\R^n)}-\frac{\epsilon'_K}{2}\right)\mathop{\to}_{m\to+\infty}m_{\Lambda}^{-1}R\times\left(\|v\|_{L^{\infty}(\Omega,\R^n)}-\frac{\epsilon'_K}{2}\right).
\ee
Moreover, since all fields $|v^*_m|$, $|\hat{v}|$, $\Lambda^*_m$ and $\hat{\Lambda}$ are radially symmetric, it follows from the above estimates and Lebesgue's dominated convergence theorem that
$$U^*_m(x)\mathop{\to}_{m\to+\infty}U(x)\hbox{ for all }x\in\overline{\Omega^*},$$
where $U$ is given in (\ref{U}).\par
Lastly, Corollary~\ref{corghat} applied with $g=-V^-$ and $\hat{g}=\hat{V}$ provides the existence of two sequences of radially symmetric fields $(V^*_m)_{m\in\N}$ and $(\overline{V}^*_m)_{m\in\N}$ in $\Omega^*$ such that, for each $m\in\N$, $V^*_m\in L^{\infty}(\Omega^*)$, $\overline{V}^*_m\in C^{\infty}(\overline{\Omega^*})$,
\be\label{V*m}
-\|V\|_{L^{\infty}(\Omega)}\le\min_{\overline{\Omega}}(-V^-)\le V^*_m\le\overline{V}^*_m\le 0\hbox{ in }\Omega^*
\ee
and the distribution functions of $V^*_m$ and $\overline{V}^*_m$ satisfy
$$\mu_{V^*_m}=\mu_{-V^-}\hbox{ and }\mu_{|\overline{V}^*_m|}\le\mu_{V^-}\le\mu_{|V|}.$$
Furthermore, the fields $V^*_m$ and $\overline{V}^*_m$ are constructed so that
$$V^*_m,\ \overline{V}^*_m\mathop{\rightharpoonup}_{m\to+\infty}\hat{V}\hbox{ in }\sigma(L^p(\Omega^*),L^{p'}(\Omega^*))$$
for all $1<p\le+\infty$, with $1/p+1/(p')=1$.\hfill\break\par

{\underline{Step~8: An inequality for the eigenvalue $\lambda_1(\Omega^*,\Lambda^*_m{\rm{Id}},v^*_m,\overline{V}^*_m)$ for large $m$}}. Remember first that the function $\tilde{\psi}$ is continuous nonnegative in $\overline{\Omega^*}$ and that $\tilde{\psi}(0)\ge\max_{\overline{\Omega}}\psi_L>0$ because of  Corollary~\ref{cor36} and (\ref{psil+}). It also follows from Lebesgue's dominated convergence theorem and all estimates of Step~7 that
$$\left\{\baa{rcl}
\displaystyle{\int_{\Omega}}\left(\Lambda^*_m|\nabla\tilde{\psi}|^2+\overline{V}^*_m\tilde{\psi}^2\right)e^{-U^*_m} & \to & \displaystyle{\int_{\Omega}}\left(\hat{\Lambda}|\nabla\tilde{\psi}|^2+\hat{V}\tilde{\psi}^2\right)e^{-U},\\
\displaystyle{\int_{\Omega}}\tilde{\psi}^2e^{-U^*_m} & \to & \displaystyle{\int_{\Omega}}\tilde{\psi}^2e^{-U}>0,\eaa\right.$$
as $m\to+\infty$.
Therefore, from (\ref{ineqpsi}), there exists $M\in\N$ such that
\be\label{ineqpsi2}\baa{rcl}
\displaystyle{\int_{\Omega}}\left(\Lambda^*_M|\nabla\tilde{\psi}|^2+\overline{V}^*_M\tilde{\psi}^2\right)e^{-U^*_M} & \le & \displaystyle{\frac{\lambda_1(\Omega,A,v,V)+3\epsilon'+\epsilon' m_{\Lambda}^{-1}+\epsilon'\|V^-\|_{\infty}}{1-\epsilon'}}\\
& & \times\displaystyle{\int_{\Omega^*}}\widetilde{\psi}^2e^{-U^*_M}.\eaa
\ee
Remember that (\ref{Lambdam1}) and (\ref{vm1}) hold with $m=M$. Furthermore, because of (\ref{Lambdam2}) and (\ref{vm2}), one can choose $M$ large enough so that
$$m_{\Lambda}\le\Lambda^*_M\le M_{\Lambda}\hbox{ in }\overline{\Omega^*}\hbox{ and }\|v^*_M\|_{L^{\infty}(\Omega^*,\R^n)}\le\|v\|_{L^{\infty}(\Omega,\R^n)}.$$
Notice that this holds also when $m_{\Lambda}=M_{\Lambda}$ (in this case, $\Lambda^*_M=m_{\Lambda}$ in $\overline{\Omega^*}$) or when $m_v=M_v$ (in this case, $v^*_M=m_ve_r$ in $\overline{\Omega^*}\backslash\{0\}$).\par
Call now
$$I=\inf_{\phi\in H^1_0(\Omega^*)\backslash\{0\}}\frac{\displaystyle{\int_{\Omega^*}}\left(\Lambda^*_M|\nabla\phi|^2+\overline{V}^*_M\phi^2\right)e^{-U^*_M}}{\displaystyle{\int_{\Omega^*}}\phi^2e^{-U^*_M}}.$$
It follows from (\ref{ineqpsi2}) that
\be\label{ineqI}
I\le\frac{\lambda_1(\Omega,A,v,V)+3\epsilon'+\epsilon' m_{\Lambda}^{-1}+\epsilon'\|V^-\|_{\infty}}{1-\epsilon'}.
\ee
Furthermore, $I$ is clearly finite and $I\ge\min_{\overline{\Omega^*}}\overline{V}^*_M\ge-\|V^-\|_{\infty}$. It is classical to check that $I$ is actually a minimum, which is reached at a function $\varphi^*_M\in H^1_0(\Omega^*)\backslash\{0\}$ such that $\varphi^*_M\ge 0$ a.e. in $\Omega^*$ and
$$\int_{\Omega^*}\left(\Lambda^*_M\nabla\varphi^*_M\cdot\nabla\phi+\overline{V}^*_M\varphi^*_M\phi\right)e^{-U^*_M}=I\int_{\Omega^*}\varphi^*_M\phi e^{-U^*_M}$$
for all $\phi\in H^1_0(\Omega^*)$. Because of (\ref{nablaUm}) with $m=M$, the change of functions $\phi=\Phi e^{U^*_M}$ leads to
$$\int_{\Omega^*}\Lambda^*_M\nabla\varphi^*_M\cdot\nabla\Phi+v^*_M\cdot\nabla\varphi^*_M\ \Phi+\overline{V}^*_M\varphi^*_M\Phi=I\int_{\Omega^*}\varphi^*_M\Phi$$
for all $\Phi\in H^1_0(\Omega^*)$. From $H^2$ regularity and $W^{2,p}$ estimates, it then follows that $\varphi^*_M$ is actually in all $W^{2,p}(\Omega^*)$ and $C^{1,\alpha}(\overline{\Omega^*})$ for all $1\le p<+\infty$ and $0\le\alpha<1$ and that
$$\left\{\baa{rcll}
-\div(\Lambda^*_M\nabla\varphi^*_M)+v^*_M\cdot\nabla\varphi^*_M+\overline{V}^*_M\varphi^*_M & = & I\varphi^*_M & \hbox{in }\Omega^*,\\
\varphi^*_M & = & 0 & \hbox{on }\partial\Omega^*,\\
\varphi^*_M & \ge & 0 & \hbox{in }\Omega^*.\eaa\right.$$
Since $\varphi^*_M\neq 0$, one concludes from Krein-Rutman theory that $\varphi^*_M$ is --up to multiplication by a positive constant-- the principal eigenfunction of the operator $-\div(\Lambda^*_M\nabla)+v^*_M\cdot\nabla+\overline{V}^*_M$ with Dirichlet boundary condition, and that $I$ is the principal eigenvalue $I=\lambda_1(\Omega^*,\Lambda^*_M{\rm{Id}},v^*_M,\overline{V}^*_M)$. Together with (\ref{epsilon'}) and (\ref{ineqI}), it resorts that
$$\lambda_1(\Omega^*,\Lambda^*_M{\rm{Id}},v^*_M,\overline{V}^*_M)\le\frac{\lambda_1(\Omega,A,v,V)+3\epsilon'+\epsilon' m_{\Lambda}^{-1}+\epsilon'\|V^-\|_{\infty}}{1-\epsilon'}\le\lambda_1(\Omega,A,v,V)+\epsilon.$$\par
Lastly, since $V^*_M\le\overline{V}^*_M$ in $\Omega$ (see (\ref{V*m}) with $m=M$), if follows from \cite{bnv} that
$$\lambda_1(\Omega^*,\Lambda^*_M{\rm{Id}},v^*_M,V^*_M)\le\lambda_1(\Omega^*,\Lambda^*_M{\rm{Id}},v^*_M,\overline{V}^*_M)\le\lambda_1(\Omega,A,v,V)+\epsilon.$$\par

{\underline{Step~9: Conclusion}}. Since $\epsilon\in(0,1)$ was arbitrary, the proof of the first part of Theorem~\ref{th1} is complete with the choice
$$(\Lambda^*,\omega^*,v^*,V^*,\overline{V}^*)=(\Lambda^*_M,|v^*_M|,v^*_M,V^*_M,\overline{V}^*_M).$$\par

{\underline{Step~10: The case when $\Lambda$ is equal to a constant $\gamma>0$ in $\Omega$}}. It follows from the previous steps that $\Lambda^*$ is equal to the same constant $\gamma$ in $\Omega^*$. Furthermore, there exists a family $(v^*_{\epsilon},V^*_{\epsilon})_{\epsilon>0}$ of fields satisfying the same bounds (\ref{bounds1}) as $v^*$ and $V^*$, together with
$$\lambda_1(\Omega^*,\gamma\hbox{Id},v^*_{\epsilon},V^*_{\epsilon})\le\lambda_1(\Omega,A,v,V)+\epsilon$$
for all $\epsilon>0$. Furthermore, $v^*_{\epsilon}=|v^*_{\epsilon}|e_r$, and $|v^*_{\epsilon}|$ and $V^*_{\epsilon}$ are $C^{\infty}(\overline{\Omega^*})$ and radially symmetric.\par
Take any sequence $(\epsilon_k)_{k\in\N}$ of positive numbers such that $\epsilon_k\to 0$ as $k\to+\infty$. Up to extraction of a subsequence, there exist two radially symmetric functions $\omega^*_0\ge 0$ and $V^*_0\le 0$ such that $|v^*_{\epsilon_k}|\rightharpoonup\omega^*_0$ and $V^*_{\epsilon_k}\rightharpoonup V^*_0$ weakly in $L^p(\Omega^*)$ for all $1\le p<+\infty$ and weak-* in $L^{\infty}(\Omega^*)$. Furthermore, the fields $v^*_0=\omega^*_0e_r$ and $V^*_0$ satisfy the bounds (\ref{bounds1bis}). Since $-\max_{\overline{\Omega}}V^-<\lambda_1(\Omega^*,\gamma{\hbox{Id}},v^*_{\epsilon_k},V^*_{\epsilon_k})\le\lambda_1(\Omega,A,v,V)+\epsilon_k$ for all $k\in\N$, the sequence $(\lambda_1(\Omega^*,\gamma{\hbox{Id}},v^*_{\epsilon_k},V^*_{\epsilon_k}))_{k\in\N}$ converges, up to extraction of a subsequence, to a real number $\lambda^*_0\in[-\max_{\overline{\Omega}}V^-,\lambda_1(\Omega,A,v,V)]$. From standard elliptic estimates, the principal eigenfunctions $\varphi_k=\varphi_{\Omega^*,\gamma{\hbox{Id}},v^*_{\epsilon_k},V^*_{\epsilon_k}}$ are bounded independently of $k$ in all $W^{2,r}(\Omega^*)$ for $1\le r<\infty$. Up to extraction of a subsequence, they converge weakly in  $W^{2,r}(\Omega^*)$ for all $1\le r<\infty$ and strongly in $C^{1,\alpha}(\overline{\Omega^*})$ for all $0\le\alpha<1$ to a solution $\varphi^*_0$ of
$$-\gamma\Delta\varphi^*_0+v^*_0\cdot\nabla\varphi^*_0+V^*_0\varphi^*_0
=\lambda^*_0\varphi^*_0\mbox{ in }\Omega^*$$
such that $\varphi^*_0\ge 0$ in $\Omega^*$, $\varphi^*_0=0$ on $\partial\Omega^*$ and $\|\varphi^*_0\|_{L^{\infty}(\Omega^*)}=1$. By uniqueness, it resorts that $\lambda^*_0=\lambda_1(\Omega^*,\gamma\hbox{Id},v^*_0,V^*_0)$ and $\varphi^*_0=\varphi_{\Omega^*,\gamma\hbox{Id},v^*_0,V^*_0}$. Thus,
$$\lambda_1(\Omega^*,\gamma\hbox{Id},v^*_0,V^*_0)\le\lambda_1(\Omega,A,v,V).$$
Notice that the radially symmetric fields $|v^*_0|$ and $V^*_0$ satisfy the bounds (\ref{bounds1bis}), but may not be smooth anymore.\hfill\fin

\begin{rem}\label{remLambdanoncst}{\rm Assume now that $0<m_{\Lambda}<M_{\Lambda}$ and that $\Lambda=m_{\Lambda}$ in a neighbourhood of $\partial\Omega$, that is $\Lambda=m_{\Lambda}$ in the set $\{x\in\Omega,\ d(x,\partial\Omega)<\gamma\}$ for some $\gamma>0$. Then, under the notations of the previous proof, besides the aforementioned conditions of Steps~1 to 5, one can choose $K\in\N$ large enough so that $\|\Lambda^{-1}\|_{L^1(\Omega)}<
|\Omega|\times(m_{\Lambda}+4\tau_K)^{-1}$ and
$$\Lambda_K=\min_{\overline{\Omega}}\Lambda_K\in\left[m_{\Lambda}+\frac{\tau_K}{2},m_{\Lambda}+2\tau_K\right]$$
in a (smaller) neighbourhood of $\partial\Omega$. Then, owing to the definition of $\hat{\Lambda}$ in Step~6 (from Section~\ref{rearrangement}), it follows that $\hat{\Lambda}=\min_{\overline{\Omega}}\Lambda_K=\hbox{ess inf}_{\Omega^*}\hat{\Lambda}$ in a neighbourhood of $\partial\Omega^*$. Finally, besides the conditions of Steps~7 and~8, one can choose $M\in\N$ large enough so that
$$\Lambda^*_M=\min_{\overline{\Omega^*}}\Lambda^*_M\in\left[m_{\Lambda}+\frac{\tau_K}{4},m_{\Lambda}+4\tau_K\right]$$
in a neighbourhood of $\partial\Omega^*$. Since $\|\Lambda^{-1}\|_{L^1(\Omega)}<
|\Omega|\times(m_{\Lambda}+4\tau_K)^{-1}$ and $\|(\Lambda^*_M)^{-1}\|_{L^1(\Omega^*)}=\|\Lambda^{-1}\|_{L^1(\Omega)}$, one concludes that $\Lambda^*_M\ (=\Lambda^*)$ is not constant in $\overline{\Omega^*}$.\par
It then follows that, in Theorem~\ref{th1}, the functions $\Lambda$ and $\Lambda^*$ are not constant in general. Actually, for the same reason, the same observation is true for Theorems~\ref{th1bis} and~\ref{faberkrahn}.}
\end{rem}

\begin{rem}\label{nonsmooth1}{\rm Consider now the case when $\Omega$ is a general open subset of $\R^n$ with finite measure, and let $A\in W^{1,\infty}(\Omega,{\mathcal{S}}_n(\R))$, $v\in L^{\infty}(\Omega,\R^n)$, $V\in L^{\infty}(\Omega)\cap C(\overline{\Omega})$ be such that $A\ge\gamma\hbox{Id}$ in $\overline{\Omega}$ for some constant $\gamma>0$. Assume that $\lambda_1(\Omega,A,v,V)\ge 0$, where $\lambda_1(\Omega,A,v,V)$ has been characterized in Remark~\ref{nonsmooth}. We claim that we get a similar conclusion as in Theorem~\ref{th1} with $\Lambda=\gamma$. Indeed, for each $\epsilon>0$, it follows from (\ref{lambda1inf1}) and (\ref{lambda1inf2}) that there exists $\Omega'\in{\mathcal{C}}$ such that $\Omega'\subset\subset\Omega$ and
$$0\le\lambda_1(\Omega,A,v,V)\le\lambda_1(\Omega',A,v,V)\le\lambda_1(\Omega,A,v,V)+\epsilon.$$
To make notations simpler, we use the same symbols for the fields in $\Omega$ and their restrictions in $\Omega'$. Apply Theorem~\ref{th1} to $\Omega'\ (\subset\Omega)$ and to the fields $A$, $v$, $V$ with $\Lambda=\gamma$, and call $\Omega^*_{\epsilon}$ the ball centered at the origin with the same measure as $\Omega'$. There exist two radially symmetric bounded functions $\omega^*_{\epsilon}\ge 0$ and $V^*_{\epsilon}\le 0$ in $\Omega^*_{\epsilon}$ such that, for $v^*_{\epsilon}=\omega^*_{\epsilon}e_r$ in $\Omega^*_{\epsilon}$,
\be\label{omegaepsilon}
\lambda_1(\Omega^*_{\epsilon},\gamma\hbox{Id},v^*_{\epsilon},V^*_{\epsilon})\le\lambda_1(\Omega',A,v,V)\le\lambda_1(\Omega,A,v,V)+\epsilon
\ee
and
$$\left\{\baa{l}
\|v^*_{\epsilon}\|_{L^{\infty}(\Omega^*_{\epsilon},\R^n)}\le\|v\|_{L^{\infty}(\Omega',\R^n)}\le\|v\|_{L^{\infty}(\Omega,\R^n)},\ \|v^*_{\epsilon}\|_{L^2(\Omega^*_{\epsilon},\R^n)}\le\|v\|_{L^2(\Omega',\R^n)}\le\|v\|_{L^2(\Omega,\R^n)},\\
-\displaystyle{\mathop{\hbox{ess sup}}_{\Omega}}\ V^-\le-\displaystyle{\mathop{\max}_{\overline{\Omega'}}}\ V^-\le V^*_{\epsilon}\le 0\hbox{ a.e. in }\Omega^*_{\epsilon},\\
\|V^*_{\epsilon}\|_{L^p(\Omega^*_{\epsilon})}\le\|V^-\|_{L^p(\Omega')}\le\|V^-\|_{L^p(\Omega)}\hbox{ for all }1\le p\le+\infty.\eaa\right.$$
Extend $v^*_{\epsilon}$ by the vector $0$ in $\Omega^*\backslash\Omega^*_{\epsilon}$ and extend $V^*_{\epsilon}$ by $0$ in $\Omega^*\backslash\Omega^*_{\epsilon}$ too, and still call $v^*_{\epsilon}$ and $V^*_{\epsilon}$ these extensions, which are now defined in $\Omega^*\ (\supset\Omega^*_{\epsilon})$ and satisfy the same bounds as above in $\Omega^*$ (remember that $\Omega^*$ denotes the ball centered at the origin and having the same Lebesgue measure as $\Omega$). Furthermore,
$$\lambda_1(\Omega^*,\gamma\hbox{Id},v^*_{\epsilon},V^*_{\epsilon})\le
\lambda_1(\Omega^*_{\epsilon},\gamma\hbox{Id},v^*_{\epsilon},V^*_{\epsilon})
\le\lambda_1(\Omega,A,v,V)+\epsilon$$
because $\Omega^*_{\epsilon}\subset\Omega^*$ and because of (\ref{omegaepsilon}). One can then argue as in Step~10 of the proof of Theorem~\ref{th1} and one can then pass to the limit for a sequence $\epsilon_k\to 0^+$ as $k\to+\infty$. There exist then two radially symmetric bounded functions $\omega^*_0\ge 0$ and $V^*_0\le 0$ in $\Omega^*$ such that, for $v^*_0=\omega^*_0e_r$ in $\Omega^*$,
$$\lambda_1(\Omega^*,\gamma\hbox{Id},v^*_0,V^*_0)\le\lambda_1(\Omega,A,v,V)$$
and the bounds (\ref{bounds1bis}) are satisfied (with $-\hbox{ess sup}_{\Omega}V^-$ instead of $-\max_{\overline{\Omega}}V^-$ in (\ref{bounds1bis})).}
\end{rem}

Let us now turn to the proof of Theorem \ref{th1bis}. It shall use the results of Section~\ref{equality} and it follows the same scheme as the one of Theorem~\ref{th1}.\hfill\break

\noindent{\bf{Proof of Theorem \ref{th1bis}.}} Assume here that $\Omega\in {\mathcal C}$ is not a ball. Let $A\in W^{1,\infty}(\Omega,{\mathcal S}_n(\R))$, $\Lambda\in L^{\infty}_+(\Omega)$, $v\in L^{\infty}(\Omega,\R^n)$, $V\in C(\overline{\Omega})$, $\overline{M}_A>0$, $\underline{m}_{\Lambda}>0$, $\overline{M}_v\ge 0$ and $\overline{M}_V\ge 0$ be such that $A\geq \Lambda \Id$ a.e. in $\Omega$ and
$$\|A\|_{W^{1,\infty}(\Omega,{\mathcal{S}}_n(\R))}\le\overline{M}_A,\ \mathop{{\rm{ess}}\ {\rm{inf}}}_{\Omega}\Lambda\ge\underline{m}_{\Lambda},\ \|v\|_{L^{\infty}(\Omega,\R^n)}\le\overline{M}_v\hbox{ and }\|V\|_{L^{\infty}(\Omega)}\le\overline{M}_V.$$
Throughout the proof, the notation $\rho=\rho(\Omega,n,\overline{M}_A,\underline{m}_{\Lambda},\overline{M}_v,\overline{M}_V)$ denotes a constant $\rho$ which only depends on $\Omega$, $n$, $\overline{M}_A$, $\underline{m}_{\Lambda}$, $\overline{M}_v$ and $\overline{M}_V$.\par
Assume that $\lambda_1(\Omega,A,v,V)>0$ and call $\varphi$ the unique principal eigenfunction of the operator $-\div(A\nabla)+v\cdot\nabla+V$ in $\Omega$ with Dirichlet boundary condition, such that $\max_{\overline{\Omega}}\varphi=1$. Namely, the function $\varphi$ satisfies
\be\label{varphi}\left\{\baa{rcll}
-\hbox{div}(A\nabla\varphi)+v\cdot\nabla\varphi+V\varphi & = & \lambda_1(\Omega,A,v,V)\varphi & \hbox{in }\Omega,\\
\varphi & > & 0 & \hbox{in }\Omega,\\
\|\varphi\|_{L^{\infty}(\Omega)} & = & 1, & \\
\varphi & = & 0 & \hbox{on }\partial\Omega,\eaa\right.
\ee
and it is in $W^{2,p}(\Omega)$ and in $C^{1,\alpha}(\overline{\Omega})$ for all $1\le p<+\infty$ and $0\le\alpha<1$.\par
First, remember that
$$-\overline{M}_V\le\min_{\overline{\Omega}}V<\lambda_1(\Omega,A,v,V).$$
Then, let $B$ be an open ball included in $\Omega$. As observed in Step~4 of the proof of Theorem~\ref{th1}, it follows from Lemma~1.1 in \cite{bnv} that  there exists a constant $C>0$ only depending on $B$, $n$, $\overline{M}_A$, $\underline{m}_{\Lambda}$, $\overline{M}_v$ and $\overline{M}_V$, such that
$$\lambda_1(\Omega,A,v,V)\le\lambda_1(B,A,v,V)\le C.$$
From standard elliptic estimates, there exists then a constant
$$N'=N'(\Omega,n,\overline{M}_A,\underline{m}_{\Lambda},\overline{M}_v,
\overline{M}_V)>0$$
such that
$$\|\varphi\|_{C^{1,1/2}(\overline{\Omega})}\le N'.$$\par
We now claim the existence of a positive constant
$$\delta'=\delta'(\Omega,n,\overline{M}_A,\underline{m}_{\Lambda},\overline{M}_v,\overline{M}_V)>0$$
such that
$$\varphi(x)\ge\delta'\times d(x,\partial\Omega)\hbox{ for all }x\in\overline{\Omega}.$$
Assume not. Then there is a sequence of fields $(A_p,\Lambda_p,v_p,V_p)_{p\in\N}$ in $W^{1,\infty}(\Omega,{\mathcal S}_n(\R))\times L^{\infty}_+(\Omega)\times L^{\infty}(\Omega,\R^n)\times C(\overline{\Omega})$ such that $A_p\geq \Lambda_p \Id$ a.e. in $\Omega$,
$$\|A_p\|_{W^{1,\infty}(\Omega,{\mathcal{S}}_n(\R))}\le\overline{M}_A,\ \mathop{{\rm{ess}}\ {\rm{inf}}}_{\Omega}\Lambda_p\ge\underline{m}_{\Lambda},\ \|v_p\|_{L^{\infty}(\Omega,\R^n)}\le\overline{M}_v,\ \|V_p\|_{L^{\infty}(\Omega)}\le\overline{M}_V$$
for all $p\in\N$, and a sequence of points $(x_p)_{p\in\N}$ in $\overline{\Omega}$ such that
\be\label{xp}
0\le\varphi_p(x_p)<\frac{d(x_p,\partial\Omega)}{p+1}\hbox{ for all }p\in\N,
\ee
where one calls $(\lambda_1(\Omega,A_p,v_p,V_p),\varphi_p)$ the unique solution of
$$\left\{\baa{rcll}
-\hbox{div}(A_p\nabla\varphi_p)+v_p\cdot\nabla\varphi_p+V_p\varphi_p & = & \lambda_1(\Omega,A_p,v_p,V_p)\varphi_p & \hbox{in }\Omega,\\
\varphi_p & > & 0 & \hbox{in }\Omega,\\
\|\varphi_p\|_{L^{\infty}(\Omega)} & = & 1, & \\
\varphi_p & = & 0 & \hbox{on }\partial\Omega\eaa\right.$$
for each $p\in\N$. We have already noticed that the sequence $(\lambda_1(\Omega,A_p,v_p,V_p))_{p\in\N}$ is bounded. From standard elliptic estimates, the sequence $(\varphi_p)_{p\in\N}$ is also bounded in  $W^{2,q}(\Omega)$ and $C^{1,\alpha}(\overline{\Omega})$ for each $1\le q<+\infty$ and $0\le\alpha<1$. Up to a subsequence, one can assume without loss of generality that $A_p\rightarrow A_{\infty}$ (componentwise) uniformly in $\overline{\Omega}$, $\nabla A_p\rightharpoonup\nabla A_{\infty}$ (componentwise) in $\sigma(L^{\infty},L^1)$, $v_p\rightharpoonup v_{\infty}$ (componentwise) in $\sigma(L^{\infty},L^1)$, $V_p\rightharpoonup V_{\infty}$ in $\sigma(L^{\infty},L^1)$, $\varphi_p\to\varphi_{\infty}$ weakly in $W^{2,q}(\Omega)$ for all $1\le q<+\infty$ and strongly in $C^{1,\alpha}(\overline{\Omega})$ for all $0\leq \alpha<1$, $\lambda_1(\Omega,A_p,v_p,V_p)\rightarrow \lambda$ and $x_p\rightarrow x_{\infty}\in\overline{\Omega}$ as $p\to+\infty$. It follows that
$$\left\{\baa{rcll}
-\hbox{div}(A_{\infty}\nabla\varphi_{\infty})+v_{\infty}\cdot\nabla\varphi_{\infty}+V_{\infty}\varphi_{\infty} & = & \lambda\varphi_{\infty} & \hbox{in }\Omega,\\
\varphi_{\infty} & \ge & 0 & \hbox{in }\Omega,\\
\|\varphi_{\infty}\|_{L^{\infty}(\Omega)} & = & 1, & \\
\varphi_{\infty} & = & 0 & \hbox{on }\partial\Omega.\eaa\right.$$
and $\varphi_{\infty}(x_{\infty})=0$. Since $A_{\infty}\ge\underline{m}_{\Lambda}\hbox{Id}$ in $\Omega$, the strong maximum principle yields that $\varphi_{\infty}>0$ in $\Omega$, whence $x_{\infty}\in\partial\Omega$. On the other hand, (\ref{xp}) implies that $|\nabla\varphi_{\infty}(x_{\infty})|=0$, which is impossible from Hopf lemma. One has then reached a contradiction.\par
Therefore, coming back to the fields $(A,\Lambda,v,V)$ and to the function $\varphi$ solving (\ref{varphi}), we get the existence of $\delta'=\delta'(\Omega,n,\overline{M}_A,\underline{m}_{\Lambda},\overline{M}_v,
\overline{M}_V)>0$ such that
$$\varphi(x)\ge\delta'\times d(x,\partial\Omega)\hbox{ for all }x\in\overline{\Omega}.$$
In other words, the function $\varphi$ is in the set 
$E_{1/2,N',\delta'}(\Omega)$.\par
Call
$$N=3N'(\Omega,n,\overline{M}_A,\underline{m}_{\Lambda},
\overline{M}_v,\overline{M}_V)>0$$
and
$$\delta=\frac{\delta'(\Omega,n,\overline{M}_A,\underline{m}_{\Lambda},
\overline{M}_v,\overline{M}_V)}{3}>0.$$
Then, call $\eta>0$ the positive constant which is given in Corollary~\ref{corball2}, with the choice $\alpha=1/2$. It only depends on $\Omega$, $N$ and $\delta$, and therefore,
$$\eta=\eta(\Omega,n,\overline{M}_A,\underline{m}_{\Lambda},\overline{M}_v,
\overline{M}_V).$$
Call now
$$\theta=\theta(\Omega,n,\overline{M}_A,\underline{m}_{\Lambda},\overline{M}_v,
\overline{M}_V)=\frac{\eta(\Omega,n,\overline{M}_A,\underline{m}_{\Lambda},
\overline{M}_v,\overline{M}_V)}{2}>0.$$
Then, choose any $\epsilon' \in(0,1)$ such that
\be\label{epsilon'bis}
\frac{\displaystyle{\frac{\lambda_1(\Omega,A,v,V)+2\epsilon'}{1+\eta}}+\epsilon'+\epsilon'm_{\Lambda}^{-1}+\epsilon'\|V^-\|_{\infty}}{1-\epsilon'}\le\frac{\lambda_1(\Omega,A,v,V)}{1+\theta}.
\ee
It is indeed possible to choose such a $\epsilon'$ since $\lambda_1(\Omega,A,v,V)>0$ and $0<\theta<\eta$.\par
Under the notation of Step~4 of the proof of Theorem~\ref{th1}, there exists a sequence of $C^{\infty}$ fields $(A_k,\Lambda_k,v_k)_{k\in\N}$ satisfying (\ref{Akvk1}) and (\ref{Akvk2}), and such that the solutions $\varphi_k$ of (\ref{varphik}) converge to $\varphi$ as $k\to+\infty$ weakly in $W^{2,p}(\Omega)$ and strongly in $C^{1,\alpha}(\overline{\Omega})$ for all $1\le p<+\infty$ and $0\le\alpha<1$. Furthermore, it has been proved that $\lambda_1(\Omega,A_k,v_k,V)\to\lambda_1(\Omega,A,v,V)$ as $k\to+\infty$. Then there exists $K\in\N$ such that (\ref{lambdaK}) holds and
$$\varphi_K\in E_{1/2,2N',\delta'/2}(\Omega).$$
Under the same notation as in Step~5 of the proof of Theorem~\ref{th1}, the functions $\psi_l$ converge to $\varphi_K$ as $l\to+\infty$ in $W^{2,p}(\Omega)$ and $C^{1,\alpha}(\overline{\Omega})$ for all $1\le p<+\infty$ and $0\le\alpha<1$. Therefore, as in Step~5, there exists $L\in\N$ such that (\ref{psilbis}) holds and
$$\psi_L\in E_{1/2,3N',\delta'/3}(\Omega)=E_{1/2,N,\delta}(\Omega).$$
Therefore, all assumptions of Corollary~\ref{corball2} are satisfied with
$$(A_{\Omega},\Lambda_{\Omega},\omega,V,\psi,\omega_0,\mu)=(A_K,
\Lambda_K,|v_K|,V,\psi_L,\epsilon',\lambda_1(\Omega,A,v,V)+2\epsilon').$$
Notice especially that $\omega_0$ and $\mu$ are nonnegative. With the same notations as in Sections~\ref{rearrangement} and~\ref{equality}, it then follows from Corollary~\ref{corball2} that
$$\int_{\Omega^{\ast}} \left[\widehat{\Lambda}(x)\vert \nabla\widetilde{\psi}(x)\vert^2-\epsilon'\vert \nabla\widetilde{\psi}(x)\vert\widetilde{\psi}(x)+\widehat{V}(x)\widetilde{\psi}(x)^2\right]e^{-U(x)}dx\leq\frac{\mu}{1+\eta}\int_{\Omega^{\ast}} \widetilde{\psi}(x)^2e^{-U(x)}dx.$$
The same calculations as in Steps~6,~7 and~8 of the proof of Theorem~\ref{th1} can be carried out, where $\mu=\lambda_1(\Omega,A,v,V)+2\epsilon'$ is replaced by $(\lambda_1(\Omega,A,v,V)+2\epsilon')/(1+\eta)$ in (\ref{tildepsiU}). One then gets the existence of three radially symmetric $C^{\infty}(\overline{\Omega^*})$ fields $\Lambda^*>0$, $\omega^*\ge 0$,  $\overline{V}^*\le 0$ and a nonpositive radially symmetric $L^{\infty}(\Omega^*)$ field $V^*$, which satisfy (\ref{bounds1}) and are such that $\mu_{V^*}=\mu_{-V^-}$, $V^*\le\overline{V}^*$ and
$$\baa{rcl}
\lambda_1(\Omega^*,\Lambda^*{\rm{Id}},v^*,V^*) & \le & 
\lambda_1(\Omega^*,\Lambda^*{\rm{Id}},v^*,\overline{V}^*)\\
\\
& \le & \displaystyle{\frac{\displaystyle{\frac{\lambda_1(\Omega,A,v,V)+2\epsilon'}{1+\eta}}+\epsilon'+\epsilon'm_{\Lambda}^{-1}+\epsilon'\|V^-\|_{\infty}}{1-\epsilon'}},\eaa$$
where $v^*=\omega^*e_r$ in $\Omega^*$. As a consequence, there holds
$$\lambda_1(\Omega^*,\Lambda^*{\rm{Id}},v^*,V^*)\le\lambda_1(\Omega^*,\Lambda^*{\rm{Id}},v^*,\overline{V}^*)\le
\frac{\lambda_1(\Omega,A,v,V)}{1+\theta}$$
from the choice of $\epsilon'$ in (\ref{epsilon'bis}). The proof of Theorem~\ref{th1bis} is now complete.\hfill\fin


\subsection{Constraints on the eigenvalues of the matrix field $A$}\label{sec4.3}

We now give the proof of Theorem \ref{th2}. The following elementary lemma will be needed:

\begin{lem}\label{algebra}
Let $n\ge 2$, $p\in\{1,\ldots,n-1\}$, $\omega>0$, $\sigma>0$ and $A\in{\mathcal S}_n(\R)$ be positive definite such that ${\rm{det}}(A)\ge\omega$ and $\sigma_p(A)\le\sigma$. Then, there exist two positive numbers $a_1$, $a_2$ which only depend on $n$, $p$, $\omega$ and $\sigma$ such that ${\rm{det}}(D)=\omega$, $\sigma_p(D)=\sigma$, $A\ge a_1{\rm{Id}}$ and $D\ge a_1{\rm{Id}}$, where $D$ is the diagonal matrix $D={\rm{diag}}(a_1,a_2,\ldots,a_2)\in{\mathcal{S}}_n(\R)$.
\end{lem}

\noindent{\bf{Proof.}} Notice first that, as already underlined in the introduction,  the assumptions of the lemma imply that $\hbox{C}_n^p\omega^{p/n}\le\sigma$. Let $f(s)$ be defined for all $s>0$ as
$$f(s)=\omega\hbox{C}_{n-1}^{p-1}s^{p-n}+\hbox{C}_{n-1}^ps^p.$$
The function $s$ is continuous and strictly convex in $(0,+\infty)$. Furthermore, $f(0^+)=f(+\infty)=+\infty$ and elementary calculations give:
$$\min_{s\in(0,+\infty)}f(s)=f(\omega^{1/n})=\hbox{C}_n^p\omega^{p/n}\le\sigma.$$
Call
$$a_2=\max\left\{s>0,\ f(s)\le\sigma\right\},\ a_1=\frac{\omega}{(a_2)^{n-1}},\ D={\rm{diag}}(a_1,a_2,\ldots,a_2)\in{\mathcal{S}}_n(\R).$$
The real numbers $a_1$ and $a_2$ are well-defined and positive. They only depend on $n$, $p$, $\omega$ and $\sigma$. Furthermore, $a_1\le\omega^{1/n}\le a_2$, $D\in{\mathcal{S}}_n(\R)$ and
$$\hbox{det}(D)=\omega,\ \sigma_p(D)=\hbox{C}_{n-1}^{p-1}a_1(a_2)^{p-1}+\hbox{C}_{n-1}^p(a_2)^p=f(a_2)=\sigma.$$\par
Denote $0<\lambda_1\le\cdots\le\lambda_n$ the eigenvalues of $A$. Call
$$0<\widetilde{a}_2=\left(\prod_{2\leq i\leq n}\lambda_i\right)^{1/(n-1)}$$
and
\be\label{tildea1}
0<\widetilde{a}_1=\frac{\omega}{(\widetilde{a}_2)^{n-1}}\le\frac{\hbox{det}(A)}{\displaystyle{\mathop{\prod}_{2\leq i\leq n}}\lambda_i}=\lambda_1.
\ee
Since $\omega=\tilde{a}_1(\tilde{a}_2)^{n-1}\le\lambda_1(\tilde{a}_2)^{n-1}$, one has
$$f(\tilde{a}_2)=\omega\hbox{C}_{n-1}^{p-1}(\tilde{a}_2)^{p-n}+\hbox{C}_{n-1}^p(\tilde{a}_2)^p\le\hbox{C}_{n-1}^{p-1}\lambda_1\times\left(\prod_{2\leq i\leq n}\lambda_i\right)^{\frac{p-1}{n-1}}+\hbox{C}_{n-1}^p\times\left(\prod_{2\leq i\leq n}\lambda_i\right)^{\frac{p}{n-1}}.$$
For $q=p-1$ or $q=p$, call
$$\mathcal{J}_q=\{J\subset\{2,\ldots,n\},\ \hbox{card}(J)=q\}$$
and, for all $I\subset\{1,\ldots,n\}$,
$$\pi_I=\prod_{i\in I}\lambda_i.$$
Observe that $\hbox{card}(\mathcal{J}_q)=\hbox{C}_{n-1}^q$ and that
$$\prod_{2\le i\le n}\lambda_i=\left(\prod_{J\in\mathcal{J}_{p-1}}\pi_J\right)^{1/\hbox{C}_{n-2}^{p-2}}=\left(\prod_{J\in\mathcal{J}_p}\pi_J\right)^{1/\hbox{C}_{n-2}^{p-1}}.$$
Thus,
$$\baa{rcl}
f(\tilde{a}_2) & \le & \hbox{C}_{n-1}^{p-1}\lambda_1\times\left(\displaystyle{\mathop{\prod}_{J\in\mathcal{J}_{p-1}}}\pi_J\right)^{\frac{p-1}{(n-1)\hbox{C}_{n-2}^{p-2}}}+\hbox{C}_{n-1}^p\times\left(\displaystyle{\mathop{\prod}_{J\in\mathcal{J}_p}}\pi_J\right)^{\frac{p}{(n-1)\hbox{C}_{n-2}^{p-1}}}\\
& = & \hbox{C}_{n-1}^{p-1}\lambda_1\times\left(\displaystyle{\mathop{\prod}_{J\in\mathcal{J}_{p-1}}}\pi_J\right)^{1/\hbox{C}_{n-1}^{p-1}}+\hbox{C}_{n-1}^p\times\left(\displaystyle{\mathop{\prod}_{J\in\mathcal{J}_p}}\pi_J\right)^{1/\hbox{C}_{n-1}^p}\\
& \le & \lambda_1\times\displaystyle{\mathop{\sum}_{J\in\mathcal{J}_{p-1}}}\pi_J+\displaystyle{\mathop{\sum}_{J\in\mathcal{J}_p}}\pi_J\\
& = & \displaystyle{\mathop{\sum}_{I\subset\{1,\ldots,n\},\ \hbox{card}(I)=p}}\pi_I\\
& = & \sigma_p(A)\\
& \le & \sigma.\eaa$$
Hence, $\tilde{a}_2\le a_2$ and
$$\tilde{a}_1=\frac{\omega}{(\tilde{a}_2)^{n-1}}\ge\frac{\omega}{(a_2)^{n-1}}=a_1.$$
Together with (\ref{tildea1}), it follows that $\lambda_1\ge a_1$, whence $A\ge a_1\hbox{Id}$.\hfill\fin\break

\noindent{\bf{Proof of Theorem~\ref{th2}.}} Let $\Omega$, $A$, $v$, $V$, $p$, $\omega$, $\sigma$ be as in Theorem~\ref{th2}, and let $a_1>0$ and $a_2>0$ be given by Lemma~\ref{algebra}. Therefore,
$$A(x)\ge a_1\hbox{Id}\hbox{ for all }x\in\overline{\Omega}.$$
For all $x\in\overline{\Omega^*}\backslash\{0\}$, define now $A^*(x)$ as in Theorem~\ref{th2}. Thus, $A^*(x)\ge a_1\hbox{Id}$ and there is an invertible matrix $P(x)$ of size $n\times n$ such that
$$A^*(x)=P(x)DP(x)^{-1},$$
where $D=\hbox{diag}(a_1,a_2,\ldots,a_2)\in{\mathcal{S}}_n(\R)$. Hence,
$$\hbox{det}(A^*(x))=\hbox{det}(D)=\omega\hbox{ and }\sigma_p(A^*(x))=\sigma_p(D)=\sigma\hbox{ for all }x\in\overline{\Omega^*}\backslash\{0\}.$$\par
Let $\epsilon>0$ be given. From Theorem~\ref{th1} applied with $\Lambda$ equal to the positive constant $a_1$, it follows that there exist two radially symmetric $C^{\infty}(\overline{\Omega^*})$ fields $\omega^*\ge 0$ and $\overline{V}^*\le 0$, and a nonpositive radially symmetric $L^{\infty}(\Omega^*)$ field $V^*$, such that, for $v^*=\omega^*e_r$ in $\overline{\Omega^*}\backslash\{0\}$,
$$\left\{\baa{l}
\|v^*\|_{L^{\infty}(\Omega^*,\R^n)}\le\|v\|_{L^{\infty}(\Omega,\R^n)},\ \|v^*\|_{L^2(\Omega^*)}=\|v\|_{L^2(\Omega)},\\
\mu_{|\overline{V}^*|}\le\mu_{V^-},\ \mu_{V^*}=\mu_{-V^-},\ V^*\le\overline{V}^*\eaa\right.$$
and
$$\lambda_1(\Omega^*,a_1\hbox{Id},v^*,V^*)\le\lambda_1(\Omega^*,a_1\hbox{Id},v^*,\overline{V}^*)\le\lambda_1(\Omega,A,v,V)+\epsilon.$$
From Remark~\ref{remdet}, one concludes that
$$\lambda_1(\Omega^*,A^*,v^*,V^*)\le\lambda_1(\Omega^*,A^*,v^*,\overline{V}^*)\le\lambda_1(\Omega,A,v,V)+\epsilon.$$\par
Lastly, since $A\ge a_1\hbox{Id}$ in $\Omega$ and $A^*\nabla\phi=a_1\nabla\phi$ for any radially symmetric function $\phi$ in $\Omega^*$, the existence of two radially symmetric bounded functions $\omega^*_0\ge 0$ and $V^*_0\le 0$ in $\Omega^*$ satisfying (\ref{bounds1bis}) and
$$\lambda_1(\Omega^*,A^*,v^*_0,V^*_0)\le\lambda_1(\Omega,A,v,V),$$
where $v^*_0=\omega^*_0e_r$ in $\Omega^*$, can be done as in Step~10 of the proof of Theorem~\ref{th1}\par
This completes the proof of Theorem~\ref{th2}.\hfill\fin

\begin{rem}\label{nonsmooth2}{\rm Consider now the case when $\Omega$ is a general open subset of $\R^n$ with finite measure, and let $(A,v,V)$ be as in Theorem~\ref{th2}, with the extra assumption $V\in L^{\infty}(\Omega)$. Since $A\ge a_1\hbox{Id}$ in $\Omega$, it follows, as in Remark~\ref{nonsmooth1}, that there exist two radially symmetric bounded functions $\omega^*_0\ge 0$ and $V^*_0\le 0$ in $\Omega^*$ satisfying (\ref{bounds1bis}) (with $-\hbox{ess sup}_{\Omega}V^-$ instead of $-\max_{\overline{\Omega}}V^-$) and $\lambda_1(\Omega^*,A^*,v^*_0,V^*_0)\le\lambda_1(\Omega,A,v,V)$, where $v^*_0=\omega^*_0e_r$ in $\Omega^*$ and $A^*$ is the same as in Theorem~\ref{th2}.}
\end{rem}


\SE{The cases of $L^p$ constraints}\label{sec5}

In this section, we focus on optimization and comparison results for $\lambda_1(\Omega,A,v,V)$ when $\Omega$ has fixed Lebesgue measure, $A$ satisfies the same constraints as in Theorem~\ref{th1} and $v$ and $V$ satisfy $L^p$ constraints. We first give in Section~\ref{sec61} some optimization results when $\Omega$ is fixed. Then, relying on Theorem \ref{th1}, we derive in Section~\ref{sec62} some Faber-Krahn type inequalities.


\subsection{Optimization in fixed domains}\label{sec61}

In Theorem~\ref{th1}, we were able to compare the eigenvalue $\lambda_1(\Omega,A,v,V)$ to $\lambda_1(\Omega^*,\Lambda^*\hbox{Id},v^*,V^*)$ in $\Omega^*$, where the $L^1$ norms of $\Lambda^{-1}$ (with $A\ge\Lambda\ \hbox{Id}$) and $|v|^2\Lambda^{-1}$ and the distribution function of $V^-$ were the same as for the new fields in $\Omega^*$ (in particular, all $L^p$ norms of $V^-$ were preserved). In this section, we fix the domain $\Omega$ and we minimize and maximize $\lambda_1(\Omega,A,v,V)$ when $A\ge\Lambda\Id$ and $v$ and $V$ satisfy some $L^p\cap L^{\infty}$ constraints with some given weights. Furthermore, we prove the uniqueness of the optimal fields when $A$ is fixed and $v$ and $V$ satisfy given $L^{\infty}$ constraints.


\subsubsection{The case of $L^p$ constraints, $1<p\leq +\infty$}

Given $\Omega\in {\mathcal C}$, $M\ge 0$, $\Lambda\in L^{\infty}_+(\Omega)$, $1<p,q\le+\infty$, $w_{1,p}$, $w_{1,\infty}$, $w_{2,q}$, $w_{2,\infty}\in L^{\infty}_+(\Omega)$, $\tau_{1,p}$, $\tau_{1,\infty}$, $\tau_{2,q}$, $\tau_{2,\infty}\ge 0$ such that $M\ge{\hbox{ess sup}}_{\Omega}\ \Lambda$, define
\[
\begin{array}{l}
\displaystyle {\mathcal A}_{\Omega,M,\Lambda,p,q,w_{1,p},w_{1,\infty},w_{2,q},w_{2,\infty},\tau_{1,p},\tau_{1,\infty},\tau_{2,q},\tau_{2,\infty}}\\
\qquad\qquad=\displaystyle \left\{(A,v,V)\in W^{1,\infty}(\Omega,{\mathcal S}_n(\R))\times L^{\infty}(\Omega,\R^n) \times L^{\infty}(\Omega);\right.\\
\qquad\qquad\qquad\displaystyle \left. \left\Vert A\right\Vert_{W^{1,\infty}(\Omega,{\mathcal{S}}_n(\R))}\leq M,\ A\geq \Lambda\Id\mbox{ a.e. in }\Omega,\right.\\
\qquad\qquad\qquad\left\Vert w_{1,p}v\right\Vert_p\leq \tau_{1,p},\ \left\Vert w_{1,\infty}v\right\Vert_{\infty}\leq \tau_{1,\infty},\ \displaystyle \left. \left\Vert w_{2,q}V\right\Vert_q\leq \tau_{2,q},\ \left\Vert w_{2,\infty}V\right\Vert_{\infty}\leq \tau_{2,\infty}\right\},
\end{array}
\]
and
\be\label{optim}\left\{
\begin{array}{l}
\displaystyle \underline{\lambda}(\Omega,M,\Lambda,p,q,w_{1,p},w_{1,\infty},w_{2,q},w_{2,\infty},\tau_{1,p},\tau_{1,\infty},\tau_{2,q},\tau_{2,\infty})\\
\qquad\qquad\qquad\qquad=\displaystyle \inf_{(A,v,V)\in{\mathcal A}_{\Omega,M,\Lambda,p,q,w_{1,p},w_{1,\infty},w_{2,q},w_{2,\infty},\tau_{1,p},\tau_{1,\infty},\tau_{2,q},\tau_{2,\infty}}} \lambda_1(\Omega,A,v,V),\\
\displaystyle \overline{\lambda}(\Omega,M,\Lambda,p,q,w_{1,p},w_{1,\infty},w_{2,q},w_{2,\infty},\tau_{1,p},\tau_{1,\infty},\tau_{2,q},\tau_{2,\infty})\\
\qquad\qquad\qquad\qquad=\displaystyle \sup_{(A,v,V)\in {\mathcal A}_{\Omega,M,\Lambda,p,q,w_{1,p},w_{1,\infty},w_{2,q},w_{2,\infty},\tau_{1,p},\tau_{1,\infty},\tau_{2,q},\tau_{2,\infty}}} \lambda_1(\Omega,A,v,V).
\end{array}\right.
\ee
Observe that ${\mathcal A}_{\Omega,M,\Lambda,p,q,w_{1,p},w_{1,\infty},w_{2,q},w_{2,\infty},\tau_{1,p},\tau_{1,\infty},\tau_{2,q},\tau_{2,\infty}}\neq\emptyset$, that 
$$\underline{\lambda}(\Omega,M,\Lambda,p,q,w_{1,p},w_{1,\infty},w_{2,q},w_{2,\infty},\tau_{1,p},\tau_{1,\infty},\tau_{2,q},\tau_{2,\infty})\ge-\frac{\tau_{2,\infty}}{\hbox{ess inf}_{\Omega}w_{2,\infty}}$$
and that
$$\overline{\lambda}(\Omega,M,\Lambda,p,q,w_{1,p},w_{1,\infty},w_{2,q},w_{2,\infty},\tau_{1,p},\tau_{1,\infty},\tau_{2,q},\tau_{2,\infty})<+\infty$$
by \cite{bnv}.\par

Our first result deals with the optimization problem for (\ref{optim}):

\begin{theo} \label{fixedlp}
Let $\Omega\in {\mathcal C}$, $M\ge 0$, $\Lambda\in L^{\infty}_+(\Omega)$, $1<p,q\le+\infty$, $w_{1,p}$, $w_{1,\infty}$, $w_{2,q}$, $w_{2,\infty}\in L^{\infty}_+(\Omega)$ and $\tau_{1,p}$, $\tau_{1,\infty}$, $\tau_{2,q}$, $\tau_{2,\infty}\ge 0$ be given such that $M\ge{\rm{ess}}\ {\rm{sup}}_{\Omega}\ \Lambda$. Then, 
\begin{itemize}
\item[$(1)$]
there exists $(\underline{A},\underline{v},\underline{V})\in{\mathcal A}_{\Omega,M,\Lambda,p,q,w_{1,p},w_{1,\infty},w_{2,q},w_{2,\infty},\tau_{1,p},\tau_{1,\infty},\tau_{2,q},\tau_{2,\infty}}$ such that, if $\underline{\varphi}=\varphi_{\Omega,\underline{A},\underline{v},\underline{V}}$:
\begin{itemize}
\item[$(a)$]
$\underline{\lambda}(\Omega,M,\Lambda,p,q,w_{1,p},w_{1,\infty},w_{2,q},w_{2,\infty},\tau_{1,p},\tau_{1,\infty},\tau_{2,q},\tau_{2,\infty})=\lambda_1(\Omega,\underline{A},\underline{v},\underline{V})$,
\item[$(b)$]
$\underline{v}\cdot\nabla\underline{\varphi}=-\left\vert \underline{v}\right\vert\times\left\vert \nabla\underline{\varphi}\right\vert\mbox{a.e. in }\Omega$,
\item[$(c)$]
$\underline{V}(x)\leq 0\mbox{ a.e. in }\Omega$,
\item[$(d)$]
$\left\Vert w_{1,p}\underline{v}\right\Vert_p=\tau_{1,p}\mbox{ or }\left\Vert w_{1,\infty}\underline{v}\right\Vert_{\infty}=\tau_{1,\infty}$, $\left\Vert w_{2,q}\underline{V}\right\Vert_q=\tau_{2,q}\mbox{ or }\left\Vert w_{2,\infty}\underline{V}\right\Vert_{\infty}=\tau_{2,\infty}$.
\end{itemize}
Moreover, if $(A,v,V)\in {\mathcal A}_{\Omega,M,\Lambda,p,q,w_{1,p},w_{1,\infty},w_{2,q},w_{2,\infty},\tau_{1,p},\tau_{1,\infty},\tau_{2,q},\tau_{2,\infty}}$ is such that
$$\underline{\lambda}(\Omega,M,\Lambda,p,q,w_{1,p},w_{1,\infty},w_{2,q},w_{2,\infty},\tau_{1,p},\tau_{1,\infty},\tau_{2,q},\tau_{2,\infty})=\lambda_1(\Omega,A,v,V)$$
and if $\varphi=\varphi_{\Omega,A,v,V}$, then properties $(b),(c)$ and $(d)$ hold with $\varphi,v,V$ instead of $\underline{\varphi},\underline{v}$ and $\underline{V}$,
\item[$(2)$]
there exists $(\overline{A},\overline{v},\overline{V})\in {\mathcal A}_{\Omega,M,\Lambda,p,q,w_{1,p},w_{1,\infty},w_{2,q},w_{2,\infty},\tau_{1,p},\tau_{1,\infty},\tau_{2,q},\tau_{2,\infty}}$ such that, if $\overline{\varphi}=\varphi_{\Omega,\overline{A},\overline{v},\overline{V}}$:
\begin{itemize}
\item[$(a)$]
$\overline{\lambda}(\Omega,M,\Lambda,p,q,w_{1,p},w_{1,\infty},w_{2,q},w_{2,\infty},\tau_{1,p},\tau_{1,\infty},\tau_{2,q},\tau_{2,\infty})=\lambda_1(\Omega,\overline{A},\overline{v},\overline{V})$,
\item[$(b)$]
$\overline{v}\cdot\nabla\overline{\varphi}=\left\vert \overline{v}\right\vert\times\left\vert \nabla\overline{\varphi}\right\vert\mbox{a.e. in }\Omega$,
\item[$(c)$]
$\overline{V}(x)\geq 0\mbox{ a.e. in }\Omega$,
\item[$(d)$]
$\left\Vert w_{1,p}\overline{v}\right\Vert_p=\tau_{1,p}\mbox{ or }\left\Vert w_{1,\infty}\overline{v}\right\Vert_{\infty}=\tau_{1,\infty}$, $\left\Vert w_{2,q}\overline{V}\right\Vert_q=\tau_{2,q}\mbox{ or }\left\Vert w_{2,\infty}\overline{V}\right\Vert_{\infty}=\tau_{2,\infty}$.
\end{itemize}
Moreover, if $(A,v,V)\in {\mathcal A}_{\Omega,M,\Lambda,p,q,w_{1,p},w_{1,\infty},w_{2,q},w_{2,\infty},\tau_{1,p},\tau_{1,\infty},\tau_{2,q},\tau_{2,\infty}}$ is such that
$$\overline{\lambda}(\Omega,M,\Lambda,p,q,w_{1,p},w_{1,\infty},w_{2,q},w_{2,\infty},\tau_{1,p},\tau_{1,\infty},\tau_{2,q},\tau_{2,\infty})=\lambda_1(\Omega,A,v,V)$$
and if $\varphi=\varphi_{\Omega,A,v,V}$, then properties $(b),(c)$ and $(d)$ hold with $\varphi,v,V$ instead of $\overline{\varphi},\overline{v}$ and $\overline{V}$.
\end{itemize}
\end{theo}

We will use several times in the proof the following comparison result:

\begin{lem} \label{lemcomparison}
Let $\Omega\in{\mathcal{C}}$, $\mu\in\R$, $A\in W^{1,\infty}(\Omega,{\mathcal S}_n(\R))$ with $A\geq\gamma\Id$ a.e. in $\Omega$ for some $\gamma>0$, $v\in L^{\infty}(\Omega,\R^n)$ and $V\in L^{\infty}(\Omega)$. Assume that 
$\varphi$ and $\psi$ are functions in 
$W^{2,r}(\Omega)$ for all $1\leq r<+\infty$, satisfying $\left\Vert \varphi\right\Vert_{\infty}=
\left\Vert \psi\right\Vert_{\infty}$ and $\varphi=\psi=0$ on $\partial\Omega$. Assume also that $\varphi\geq 
0$ in $\Omega$, $\psi>0$ in $\Omega$ and
\[
\left\{
\begin{array}{l}
-{\rm{div}}(A\nabla\psi)+v.\nabla\psi+V\psi\geq \mu\psi\mbox{ a.e. in }\Omega,\\
\\
-{\rm{div}}(A\nabla\varphi)+v.\nabla \varphi+V\varphi\leq \mu\varphi\mbox{ a.e. in }\Omega.
\end{array}
\right.
\]
Then $\varphi=\psi$ in $\Omega$.
\end{lem}

\noindent{\bf Proof.} The proof uses a classical comparison method (this method was used for instance in \cite{bnv}). We give it here for the sake of completeness. Since $\psi>0$ in $\Omega$ and $\psi=0$ on $\partial\Omega$, the Hopf lemma yields $\nu(x)\cdot\nabla\psi(x)<0$ on $\partial\Omega$, where, for all $x\in \partial\Omega$, $\nu(x)$ denotes the outward normal unit vector at $x$. Since $\varphi\in C^{1,\beta}(\overline{\Omega})$ for all $0\leq \beta<1$, $\varphi\geq 0$ in $\Omega$ and $\varphi=0$ on $\partial\Omega$, it follows that there exists $\gamma>0$ such that $\gamma\psi>\varphi$ in $\Omega$. Define
\[
\gamma^{\ast}=\inf\left\{\gamma>0,\ \gamma\psi>\varphi\mbox{ in 
}\Omega\right\}.
\]
One clearly has $\gamma^{\ast}\psi\geq \varphi$ in $\Omega$, so that 
$\gamma^{\ast}>0$. Define $w=\gamma^{\ast}\psi-\varphi\ge 0$ and assume 
that $w>0$ everywhere in $\Omega$. Since 
\begin{equation} \label{diff}
-\div(A\nabla w)+v\cdot \nabla w+Vw-\mu w\geq 0\mbox{ a.e. in }\Omega
\end{equation}
and $w=0$ on $\partial\Omega$, the Hopf maximum principle implies that $\nu\cdot\nabla w<0$ on $\partial\Omega$. As above, this yields the existence of $\kappa>0$ such that $w>\kappa\varphi$ in $\Omega$, whence
$\gamma^{\ast}\psi/(1+\kappa)>\varphi$ in $\Omega$. This is a contradiction with the minimality of $\gamma^{\ast}$. Thus, there exists $x_0\in \Omega$ such that $w(x_0)=0$ (i.e. $\gamma^{\ast}\psi(x_0)=\varphi(x_0)$). Since $w\ge 0$ in $\Omega$, it follows from (\ref{diff}) and from the strong maximum principle, that $w=0$ in $\Omega$, which means that $\varphi$ and $\psi$ are proportional. Since they are non-negative in $\Omega$ and have the same $L^{\infty}$ norm in $\Omega$, one has $\varphi=\psi$, which ends the proof of Lemma \ref{lemcomparison}. \hfill\fin\break

For the proof of Theorem \ref{fixedlp}, we will treat the minimization problem only, the maximization problem being clearly analogous. It is plain to see that the result is a consequence of the two following lemmata:

\begin{lem} \label{step1}
Let $\Omega\in {\mathcal C}$, $M\ge 0$, $\Lambda\in L^{\infty}_+(\Omega)$, $1<p,q\le+\infty$, $w_{1,p}$, $w_{1,\infty}$, $w_{2,q}$, $w_{2,\infty}\in L^{\infty}_+(\Omega)$ and $\tau_{1,p}$, $\tau_{1,\infty}$, $\tau_{2,q}$, $\tau_{2,\infty}\ge 0$ be given such that $M\ge{\rm{ess}}\ {\rm{sup}}_{\Omega}\ \Lambda$. Assume that $(\underline{A},\underline{v},\underline{V})\in {\mathcal A}_{\Omega,M,\Lambda,p,q,w_{1,p},w_{1,\infty},w_{2,q},w_{2,\infty},\tau_{1,p},\tau_{1,\infty},\tau_{2,q},\tau_{2,\infty}}$ is such that
\be\label{lambdatau}
\underline{\lambda}(\Omega,M,\Lambda,p,q,w_{1,p},w_{1,\infty},w_{2,q},w_{2,\infty},\tau_{1,p},\tau_{1,\infty},\tau_{2,q},\tau_{2,\infty})=\lambda_1(\Omega,\underline{A},\underline{v},\underline{V}),
\ee
and let $\underline{\varphi}=\varphi_{\Omega,\underline{A},\underline{v},\underline{V}}$. Then, properties $(b),(c)$ and $(d)$ in Theorem \ref{fixedlp} hold.
\end{lem}

\begin{lem} \label{exist}
Let $\Omega\in {\mathcal C}$, $M\ge 0$, $\Lambda\in L^{\infty}_+(\Omega)$, $1<p,q\le+\infty$, $w_{1,p}$, $w_{1,\infty}$, $w_{2,q}$, $w_{2,\infty}\in L^{\infty}_+(\Omega)$ and $\tau_{1,p}$, $\tau_{1,\infty}$, $\tau_{2,q}$, $\tau_{2,\infty}\ge 0$ be given such that $M\ge{\rm{ess}}\ {\rm{sup}}_{\Omega}\ \Lambda$. Then, there exists $(\underline{A},\underline{v},\underline{V})\in {\mathcal A}_{\Omega,M,\Lambda,p,q,w_{1,p},w_{1,\infty},w_{2,q},w_{2,\infty},\tau_{1,p},\tau_{1,\infty},\tau_{2,q},\tau_{2,\infty}}$ such that (\ref{lambdatau}) holds.
\end{lem}

\noindent{\bf Proof of Lemma \ref{step1}.} Let $\underline{A},\underline{v},\underline{V}$ and $\underline{\varphi}$ as in Lemma \ref{step1}. Remember that $\underline{\varphi}\in C^1(\overline{\Omega})$. Define, for a.e. $x\in \Omega$,
\[
w(x)=\left\{
\begin{array}{ll}
\displaystyle -\left\vert \underline{v}(x)\right\vert \frac{\nabla\underline{\varphi}(x)}{\left\vert \nabla\underline{\varphi}(x)\right\vert} &\mbox{ if }\nabla\underline{\varphi}(x)\neq 0,\\
0 & \mbox{ if }\nabla\underline{\varphi}(x)=0
\end{array}
\right.
\]
so that $w\cdot \nabla\underline{\varphi}=-\left\vert \underline{v}\right\vert\times \left\vert \nabla\underline{\varphi}\right\vert\leq \underline{v}\cdot \nabla\underline{\varphi}$ a.e. in $\Omega$, and set
$$\mu=\lambda_1(\Omega,\underline{A},w,\underline{V})\hbox{ and }\psi=\varphi_{\Omega,\underline{A},w,\underline{V}}.$$
Notice that $|w|\le|v|$ a.e. in $\Omega$, whence $(\underline{A},w,\underline{V})\in {\mathcal A}_{\Omega,M,\Lambda,p,q,w_{1,p},w_{1,\infty},w_{2,q},w_{2,\infty},\tau_{1,p},\tau_{1,\infty},\tau_{2,q},\tau_{2,\infty}}$ and
$$\underline{\lambda}:=\underline{\lambda}(\Omega,M,\Lambda,p,q,w_{1,p},w_{1,\infty},w_{2,q},w_{2,\infty},\tau_{1,p},\tau_{1,\infty},\tau_{2,q},\tau_{2,\infty})\leq \mu.$$
Thus, one has
\[
\left\{\begin{array}{ll}
-\div(\underline{A}\nabla\underline{\varphi})+w\cdot\nabla\underline{\varphi}+\underline{V}\underline{\varphi}\leq -\div(\underline{A}\nabla\underline{\varphi})+\underline{v}\cdot\nabla\underline{\varphi}+\underline{V}\underline{\varphi}=\underline{\lambda}\underline{\varphi}\leq \mu\underline{\varphi} &\mbox{ a.e. in }\Omega,\\
-\div(\underline{A}\nabla\psi)+w\cdot\nabla\psi+\underline{V}\psi= \mu\psi &\mbox{ a.e. in }\Omega,
\end{array}\right.
\]
and Lemma \ref{lemcomparison} yields $\psi=\underline{\varphi}$ and therefore $\mu=\underline{\lambda}$ and $\underline{v}\cdot \nabla\underline{\varphi}=w\cdot \nabla\underline{\varphi}=-\left\vert \underline{v}\right\vert\times\left\vert \nabla\underline{\varphi}\right\vert$ a.e. in~$\Omega$.\par
As far as assertion $(c)$ is concerned, define, for a.e. $x\in \Omega$,
\[
W(x)=\left\{
\begin{array}{ll}
\underline{V}(x) &\mbox{ if }\underline{V}(x)\leq 0,\\
0 & \mbox{ if }\underline{V}(x)>0.
\end{array}
\right.
\]
Observe that $|W|\le|\underline{V}|$ a.e. in $\Omega$, whence $(\underline{A},\underline{v},W)\in{\mathcal{A}}_{\Omega,M,\Lambda,p,q,w_{1,p},w_{1,\infty},w_{2,q},w_{2,\infty},\tau_{1,p},\tau_{1,\infty},\tau_{2,q},\tau_{2,\infty}}$. If $\mu=\lambda_1(\Omega,\underline{A},\underline{v},W)$ and $\psi=\varphi_{\Omega,\underline{A},\underline{v},W}$, one therefore has $\underline{\lambda}\leq \mu$ and 
\[
-\div(\underline{A}\nabla\psi)+\underline{v}\cdot \nabla\psi+W\psi=\mu\psi\mbox{ a.e. in }\Omega,
\]
while, since $W\underline{\varphi}\leq \underline{V}\underline{\varphi}$ a.e. in $\Omega$,
\[
-\div(\underline{A}\nabla\underline{\varphi})+\underline{v}\cdot \nabla\underline{\varphi}+W\underline{\varphi}\leq \mu\underline{\varphi}\mbox{ a.e. in }\Omega.
\]
Lemma \ref{lemcomparison} therefore shows that $\psi=\underline{\varphi}$, and it follows that $W\underline{\varphi}=\underline{V}\underline{\varphi}$ a.e. in $\Omega$, which implies $W=\underline{V}$ a.e. in $\Omega$ since $\underline{\varphi}>0$ in $\Omega$, and this is assertion $(c)$.\par
Assume now that $\left\Vert w_{1,p}\underline{v}\right\Vert_p<\tau_{1,p}$ and  $\left\Vert w_{1,\infty}\underline{v}\right\Vert_{\infty}<\tau_{1,\infty}$, and define, for a.e. $x\in \Omega$,
\[
w(x)=\left\{
\begin{array}{ll}
\displaystyle -\left(\left\vert \underline{v}(x)\right\vert+\varepsilon\right) \frac{\nabla\underline{\varphi}(x)}{\left\vert \nabla\underline{\varphi}(x)\right\vert} &\mbox{ if }\nabla\underline{\varphi}(x)\neq 0,\\
0 & \mbox{ if }\nabla\underline{\varphi}(x)=0,
\end{array}
\right.
\]
where $\varepsilon>0$ is choosen in such a way that $\left\Vert  w_{1,p}w\right\Vert_p\leq \tau_{1,p}$ and $\left\Vert w_{1,\infty}w\right\Vert_{\infty}\leq \tau_{1,\infty}$, so that, if $\mu=\lambda_1(\Omega,\underline{A},w,\underline{V})$, one has $\underline{\lambda}\leq \mu$. Let $\psi=\varphi_{\Omega,\underline{A},w,\underline{V}}$. Observe that
\[
-\div(\underline{A}\nabla\underline{\varphi})+w\cdot\nabla\underline{\varphi}+\underline{V}\underline{\varphi}\leq \underline{\lambda}\underline{\varphi}\leq \mu\underline{\varphi}\mbox{ a.e. in }\Omega,
\]
since $w\cdot\nabla\underline{\varphi}=-\left(\left\vert \underline{v}\right\vert+\varepsilon\right) \left\vert \nabla\underline{\varphi}\right\vert\le-|\underline{v}|\times|\nabla\underline{\varphi}|\le\underline{v}\cdot\nabla\underline{\varphi}$ a.e. in $\Omega$, while
\[
-\div(\underline{A}\nabla\psi)+w\cdot\nabla\psi+\underline{V}\psi= \mu\psi \mbox{ a.e. in }\Omega.
\]
Another application of Lemma \ref{lemcomparison} yields that $\psi=\underline{\varphi}$ and therefore $w\cdot\nabla\underline{\varphi}=\underline{v}\cdot\nabla\underline{\varphi}$, so that $-\varepsilon\left\vert \nabla\underline{\varphi}\right\vert=0$ a.e. in $\Omega$, which is impossible. One argues similarly for $\underline{V}$, using the fact that $\underline{V}\leq 0$ a.e. in $\Omega$. \hfill\fin\break

\noindent{\bf Proof of Lemma \ref{exist}.} Write
$$\underline{\lambda}=\underline{\lambda}(\Omega,M,\Lambda,p,q,w_{1,p},w_{1,\infty},w_{2,q},w_{2,\infty},\tau_{1,p},\tau_{1,\infty},\tau_{2,q},\tau_{2,\infty}).$$
There exist a sequence $(A_k)_{k\in\N}\in W^{1,\infty}(\Omega,{\mathcal S}_n(\R))$ with $\left\Vert A_k\right\Vert_{W^{1,\infty}(\Omega,{\mathcal S}_n(\R))}\leq M$ and $A_k\geq \Lambda\Id$ a.e. in $\Omega$, a sequence $(v_k)_{k\in\N}\in L^{\infty}(\Omega,\R^n)$ with $\left\Vert w_{1,p}v_k\right\Vert_{p}\leq \tau_{1,p}$ and $\left\Vert w_{1,\infty}v_k\right\Vert_{\infty}\leq \tau_{1,\infty}$, and a sequence $(V_k)_{k\in\N}\in L^{\infty}(\Omega)$ with $\left\Vert w_{2,q}V_k\right\Vert_{q}\leq \tau_{2,q}$ and $\left\Vert w_{2,\infty}V_k\right\Vert_{\infty}\leq \tau_{2,\infty}$, such that
$$\lambda_k:=\lambda_1(\Omega,A_k,v_k,V_k)\rightarrow\underline{\lambda}\hbox{ as }k\to+\infty.$$
For each $k\in\N$, call $\varphi_k=\varphi_{\Omega,A_k,v_k,V_k}$, so that
\[
-\div(A_k\nabla\varphi_k)+v_k\cdot\nabla\varphi_k+V_k\varphi_k=\lambda_k\varphi_k\mbox{ in }\Omega\mbox{ and }\varphi_k=0\mbox{ on }\partial\Omega.
\]
By the usual elliptic estimates, the $\varphi_k$'s are uniformly bounded in $W^{2,r}(\Omega)$ for all $1\leq r<+\infty$, and therefore in $C^{1,\alpha}(\overline{\Omega})$ for all $0\leq \alpha<1$. Therefore, up to a subsequence, one may assume that, for some nonnegative function $\underline{\varphi}\in W^{2,r}(\Omega)$ for all $1\leq r<+\infty$, $\varphi_k\rightharpoonup \underline{\varphi}$ weakly in $W^{2,r}(\Omega)$ for all $1\le r<+\infty$ and $\varphi_k\to\underline{\varphi}$ strongly in $C^{1,\alpha}(\overline{\Omega})$ for all $0\leq \alpha<1$, as $k\to+\infty$. Similarly, there exists $\underline{A}\in W^{1,\infty}(\Omega,{\mathcal S}_n(\R))$ such that (up to extraction), $A_k\rightarrow \underline{A}$ uniformly in $\overline{\Omega}$ and, for each $1\leq j\leq n$, $\partial_jA_k\rightharpoonup \partial_j\underline{A}$ in $\sigma(L^{\infty}(\Omega),L^1(\Omega))$ componentwise. In particular, $\|\underline{A}\|_{W^{1,\infty}(\Omega,{\mathcal{S}}_n(\R))}\le M$ and $\underline{A}\ge\Lambda	\Id$ a.e. in $\Omega$. Finally, up to extraction again, there exists $\omega\in L^{\infty}(\Omega)$ such that $\left\vert v_k\right\vert \rightarrow\omega\ge 0$ in $\sigma(L^r(\Omega),L^{r'}(\Omega))$ for all $1<r\le+\infty$ (where $1/r+1/r'=1$) and there exists $\underline{V}\in L^{\infty}(\Omega)$ such that $V_k\rightarrow \underline{V}$ in $\sigma(L^r(\Omega),L^{r'}(\Omega))$ for all $1<r\le+\infty$. Since, for all $k\geq 1$, by Cauchy-Schwarz,
\[
-\div(A_k\nabla\varphi_k)-\left\vert v_k\right\vert \left\vert \nabla\varphi_k\right\vert+V_k\varphi_k\leq \lambda_k\varphi_k\mbox{ a.e. in }\Omega,
\]
one has
\[
-\div(\underline{A}\nabla\underline{\varphi})-\omega\left\vert \nabla\underline{\varphi}\right\vert+\underline{V}\underline{\varphi}\leq \underline{\lambda}\underline{\varphi}\mbox{ a.e. in }\Omega.
\]
Define now, for a.e. $x\in\Omega$,
\[
\underline{v}(x)=\left\{
\begin{array}{ll}
\displaystyle -\omega(x)\frac{\nabla\underline{\varphi}(x)}{\left\vert \nabla\underline{\varphi}(x)\right\vert} & \mbox{if }\nabla\underline{\varphi}(x)\neq 0,\\
0 & \mbox{if } \nabla\underline{\varphi}(x)=0,
\end{array}
\right.
\]
so that $\underline{v}\cdot \nabla\underline{\varphi}=-\omega\left\vert \nabla\underline{\varphi}\right\vert$ a.e. in $\Omega$. One therefore has
\be\label{underlines}
-\div(\underline{A}\nabla\underline{\varphi})+\underline{v}\cdot\nabla\underline{\varphi}+\underline{V}\underline{\varphi}\leq \underline{\lambda}\underline{\varphi}\mbox{ a.e. in }\Omega\mbox{ and }\underline{\varphi}=0\mbox{ on }\partial\Omega.
\ee
Observe that $(\underline{A},\underline{v},\underline{V})\in {\mathcal A}_{\Omega,M,\Lambda,p,q,w_{1,p},w_{1,\infty},w_{2,q},w_{2,\infty},\tau_{1,p},\tau_{1,\infty},\tau_{2,q},\tau_{2,\infty}}$. Define now $\mu=\lambda_1(\Omega,\underline{A},\underline{v},\underline{V})$ and $\psi=\varphi_{\Omega,\underline{A},\underline{v},\underline{V}}$, so that $\underline{\lambda}\leq \mu$. It follows from (\ref{underlines}) that
\[
\left\{\begin{array}{ll}
-\div(\underline{A}\nabla\psi)+\underline{v}\cdot \nabla\psi+\underline{V}\psi=\mu\psi & \mbox{ a.e. in }\Omega,\\
-\div(\underline{A}\nabla\underline{\varphi})+\underline{v}\cdot\nabla\underline{\varphi}+\underline{V}\underline{\varphi}\leq \mu\underline{\varphi} & \mbox{ in }\Omega.
\end{array}\right.
\]
Moreover, $\underline{\varphi}\ge 0$ in $\Omega$, $\psi>0$ in $\Omega$, $\underline{\varphi}=\psi=0$ on $\partial\Omega$ and $\left\Vert \underline{\varphi}\right\Vert_{\infty}=\left\Vert \psi\right\Vert_{\infty}=1$. Lemma \ref{lemcomparison} therefore yields $\underline{\varphi}=\psi$, hence $\underline{\lambda}=\mu$. This ends the proof of Lemma \ref{exist}.\hfill\fin

\begin{rem}\label{remlp}{\rm What happens in Theorem \ref{fixedlp} if one drops the $L^{\infty}$ bounds for $v$ or $V$ ? Even if one still assumes that $v$ and $V$ are qualitatively in $L^{\infty}$ (so that the principal eigenvalue of $-\div(A\nabla)+v\cdot\nabla+V$ is well-defined by \cite{bnv}), it turns out that the infimum or the supremum considered there may not be achieved. For instance, fix $1<p<n$, $\Omega\in {\mathcal C}$, $\tau>0$, $A=\hbox{Id}$ and $V=0$ in $\overline{\Omega}$, and define
\[
\underline{\lambda}(\Omega,\tau)=\inf_{v\in L^{\infty}(\Omega,\R^n),\ \left\Vert v\right\Vert_p\leq \tau} \lambda_1(\Omega,\hbox{Id},v,0).
\]
Since the operator $-\Delta+v\cdot\nabla$ satisfies the maximum principle in $\Omega$, its principal eigenvalue is positive, for each $v\in L^{\infty}(\Omega,\R^n)$, and therefore $\underline{\lambda}(\Omega,\tau)\ge 0$. We claim that $\underline{\lambda}(\Omega,\tau)=0$. Indeed, fix $\rho_0>0$ such that there exists a ball $B_0$ with radius $\rho_0$ included in $\Omega$. Call $x_0$ its center. For all $A>0$ large enough, define $\rho_A\in (0,\rho_0)$ such that $A(\alpha_n\rho_A^n)^{1/p}= \tau$ (recall that $\alpha_n$ is the Lebesgue measure of the Euclidean unit ball in $\R^n)$, let $B_A=B(x_0,\rho_A)$ be the ball with the same center $x_0$ as $B_0$ and with radius $\rho_A$, and set $v=A\times{\bf 1}_{B_A}e_r(\cdot-x_0)$ in $\Omega$, so that $\left\Vert v\right\Vert_{L^p(\Omega,\R^n)}=\tau$. One has $\lambda_1(\Omega,\hbox{Id},v,0)\leq \lambda_1(B_A,\hbox{Id},Ae_r(\cdot-x_0),0)$ since $B_A\subset\Omega$. But
$$\lambda_1(B_A,\hbox{Id},Ae_r(\cdot-x_0),0)=\frac{\mu_A}{\rho_A^2},$$
where $\mu_A=\lambda_1(\widetilde{B},\hbox{Id},A\rho_Ae_r(\cdot-x_0),0)$ is the principal eigenvalue of $-\Delta+A\rho_Ae_r(\cdot-x_0)\cdot \nabla$ on the ball $\widetilde{B}$ with center $x_0$ and with radius $1$, under Dirichlet boundary condition. Notice that $A\rho_A\to+\infty$ as $A\to+\infty$ since $1<p<n$. Furthermore, $\mu_A=\lambda_1(B^n_1,A\rho_Ae_r)$ under the notation (\ref{lambda1Omegav}), where $B^n_1$ is the Euclidean ball of $\R^n$ with center $0$ and radius $1$. It then follows immediately from Appendix~\ref{sec73} that $\log\mu_A\sim -A\rho_A$ when $A\rightarrow +\infty$ (see also \cite{fr} for related results under stronger regularity assumptions). As a consequence, there exists $A_0>0$ such that, for all $A>A_0$, 
\[
\lambda_1(B_A,\hbox{Id},Ae_r(\cdot-x_0),0)\leq \frac{e^{-A\rho_A/2}}{\rho_A^2}=\alpha_n^{2/n}\tau^{-2p/n}A^{2p/n}e^{-(\alpha_n^{-1/n}\tau^{p/n}A^{1-p/n})/2},
\]
and this expression goes to $0$ when $A\rightarrow +\infty$, which proves the claim.\par
Similarly, one can show that, if $\Omega\in {\mathcal C}$, $\tau>0$ and $1<q<n/2$ are fixed,
\be\label{claimV}
\inf_{V\in L^{\infty}(\Omega);\ \left\Vert V\right\Vert_q\leq \tau} \lambda_1(\Omega,\hbox{Id},0,V)=-\infty.
\ee
Indeed, fix $\rho_0$ as before, and, for all $A$ large enough, let $\rho_A\in(0,\rho_0)$ such that $A^q\alpha_n\rho_A^n=\tau^q$, and set $V=-A\times{\bf 1}_{B_A}$ where $B_A=B(x_0,\rho_A)$ is defined as previously, so that $\left\Vert V\right\Vert_{L^q(\Omega)}=\tau$. One has 
\be\label{VA}\baa{rcl}
\lambda_1(\Omega,\hbox{Id},0,V)\ \leq\ \lambda_1(B_A,\hbox{Id},0,-A) & = & \lambda_1(B_A,\hbox{Id},0,0)-A\\
& = & \displaystyle{\frac{C}{\rho_A^2}}-A=C\alpha_n^{2/n}\tau^{-2q/n}A^{2q/n}-A,\eaa
\ee
where $C=\lambda_1(\widetilde{B},\hbox{Id},0,0)>0$. The right-hand side of (\ref{VA}) goes to $-\infty$ when $A\rightarrow +\infty$, due to the choice of $q$. This ends the proof of the claim (\ref{claimV}).}
\end{rem}


\subsubsection{The case of $L^{\infty}$ constraints}

When solving optimization problems for $\lambda_1(\Omega,A,v,V)$ if $A$ is fixed and $v,V$ vary and satisfy $L^{\infty}$ bounds, we can precise the conclusions of Theorem \ref{fixedlp}. Fix $\Omega\in {\mathcal C}$ and $A\in W^{1,\infty}(\Omega,{\mathcal S}_n(\R))$ such that $A\ge\gamma\hbox{Id}$ in $\Omega$ for some positive real number $\gamma>0$. Given $\tau_1,\tau_2\geq 0$ and $w_1\in L^{\infty}_+(\Omega)$, define
\[
\displaystyle {\mathcal A}_{\Omega,A,w_1,\tau_{1},\tau_{2}}=\displaystyle \left\{(v,V)\in L^{\infty}(\Omega,\R^n) \times L^{\infty}(\Omega);\ \left\Vert w_1v\right\Vert_{\infty}\leq \tau_{1},\ \left\Vert V\right\Vert_{\infty}\leq \tau_{2},\right\},
\]
and
\[
\left\{\begin{array}{lll}
\displaystyle \underline{\lambda}(\Omega,A,w_1,\tau_{1},\tau_{2}) & = & \displaystyle \inf_{(v,V)\in {\mathcal A}_{\Omega,A,w_1,\tau_{1},\tau_{2}}} \lambda_1(\Omega,A,v,V),\\
\displaystyle \overline{\lambda}(\Omega,A,w_1,\tau_{1},\tau_{2}) & = & \displaystyle \sup_{(v,V)\in {\mathcal A}_{\Omega,A,w_1,\tau_{1},\tau_{2}}} \lambda_1(\Omega,A,v,V).
\end{array}\right.
\]
The optimization results under these constraints are the following ones:

\begin{theo} \label{fixedlinfty}
Let $\Omega\in {\mathcal C}$, $A\in W^{1,\infty}(\Omega,{\mathcal S}_n(\R))$ such that $A\ge\gamma{\rm{Id}}$ in $\Omega$ for some $\gamma>0$, $\tau_1,\tau_2\geq 0$ and $w_1\in L^{\infty}_+(\Omega)$ be given.Then, 
\begin{itemize}
\item[$(1)$]
there exist a unique vector field $\underline{v}\in L^{\infty}(\Omega,\R^n)$ with $\left\Vert w_1\underline{v}\right\Vert_{\infty}\leq \tau_1$ and a unique function $\underline{V}\in L^{\infty}(\Omega)$ with $\left\Vert \underline{V}\right\Vert_{\infty}\leq \tau_2$, such that
$$\underline{\lambda}(\Omega,A,w_1,\tau_1,\tau_2)=\lambda_1(\Omega,A,\underline{v},\underline{V}).$$
Moreover, if $\underline{\varphi}=\varphi_{\Omega,A,\underline{v},\underline{V}}$, one has
\begin{itemize}
\item[$(a)$]
$\underline{v}\cdot\nabla\underline{\varphi}=-\tau_1w_1^{-1} \left\vert \nabla\underline{\varphi}\right\vert\mbox{a.e. in }\Omega$,
\item[$(b)$]
$\left\vert \underline{v}(x)\right\vert w_1(x)=\tau_{1}\mbox{ a.e. in }\Omega$,
\item[$(c)$]
$\underline{V}(x)=-\tau_2\mbox{ a.e. in }\Omega$,
\end{itemize} 
Furthermore, $\nabla\underline{\varphi}(x)\neq 0$ a.e. in $\Omega$ and $\underline{v}(x)=-\tau_1w_1(x)^{-1}\nabla\underline{\varphi}(x)/|\nabla\underline{\varphi}(x)|$ a.e. in $\Omega$.
\item[$(2)$]
there exist a unique vector field $\overline{v}\in L^{\infty}(\Omega,\R^n)$ with $\left\Vert w_1\overline{v}\right\Vert_{\infty}\leq \tau_1$ and a unique function $\overline{V}\in L^{\infty}(\Omega)$ with $\left\Vert \overline{V}\right\Vert_{\infty}\leq \tau_2$, such that
$$\overline{\lambda}(\Omega,A,w_1,\tau_1,\tau_2)=\lambda_1(\Omega,A,\overline{v},\overline{V}).$$
Moreover, if $\overline{\varphi}=\varphi_{\Omega,A,\overline{v},\overline{V}}$, one has
\begin{itemize}
\item[$(a)$]
$\overline{v}\cdot\nabla\overline{\varphi}=\tau_1w_1^{-1}\left\vert \nabla\overline{\varphi}\right\vert\mbox{a.e. in }\Omega$,
\item[$(b)$]
$\left\vert \overline{v}(x)\right\vert w_1(x)=\tau_{1}\mbox{ a.e. in }\Omega$,
\item[$(c)$]
$\overline{V}(x)=\tau_2\mbox{ a.e. in }\Omega$,
\end{itemize} 
Furthermore, $\nabla\overline{\varphi}(x)\neq 0$ a.e. in $\Omega$ and $\overline{v}(x)=+\tau_1w_1(x)^{-1}\nabla\overline{\varphi}(x)/|\nabla\overline{\varphi}(x)|$ a.e. in $\Omega$.
\end{itemize}
\end{theo}

\noindent{\bf Proof. }As in the proof of Theorem \ref{fixedlp}, let us focus on the minimization problem. The existence of $\underline{v}$ and $\underline{V}$ such that $\lambda_1(\Omega,A,\underline{v},\underline{V})=\underline{\lambda}(\Omega,A,w_1,\tau_1,\tau_2)$ is obtained in the same way as in Lemma \ref{exist}, except that one has to define, for almost every $x\in \Omega$,
\[
\underline{v}(x)=\left\{
\begin{array}{ll}
\displaystyle -\tau_1w_1(x)^{-1}\frac{\nabla\underline{\varphi}(x)}{\left\vert \nabla\underline{\varphi}(x)\right\vert} &\mbox{ if }\nabla\underline{\varphi}(x)\neq 0,\\
0 &\mbox{ if }\nabla\underline{\varphi}(x)=0,
\end{array}
\right.
\]
and $\underline{V}(x)=-\tau_2$ for all $x\in \Omega$, and that we do not need to introduce the vector field $w$. To prove the uniqueness of $\underline{V}$, proceed as in the proof of assertion $(c)$ in Lemma \ref{step1}. We are now left with the task of proving the uniqueness of $\underline{v}$  and the fact that $w_1(x)\left\vert \underline{v}(x)\right\vert=\tau_1$ for almost every $x\in \Omega$.\par
First, arguing as in the proof of Lemma \ref{step1}, one shows that, if $v$  and $V$ are such that $\lambda_1(\Omega,A,v,V)=\underline{\lambda}(\Omega,A,w_1,\tau_1,\tau_2)$ and if $\varphi=\varphi_{\Omega,A,v,V}$, then $v\cdot \nabla\varphi=-\tau_1w_1^{-1}\left\vert \nabla\varphi\right\vert$ a.e. in $\Omega$.\par
To conclude, we need the following lemma:

\begin{lem} \label{uniquenessbis}
Let $\lambda\in\R$ and $\psi\in W^{2,r}(\Omega)$ for all $1\le r<+\infty$, be such that $\psi=0$ on $\partial\Omega$, $\psi>0$ in 
$\Omega$, $\left\Vert \psi\right\Vert_{\infty}=1$ and
\[
-{\rm{div}}(A\nabla\psi)-\tau_1w_1^{-1}\left\vert 
\nabla\psi\right\vert-\tau_2\psi=\lambda\psi\mbox{ in }\Omega.
\]
Let $v\in L^{\infty}(\Omega,\R^n)$ be such that $\left\Vert 
w_1v\right\Vert_{\infty}\leq\tau_1$ and 
$\lambda_1(\Omega,A,v,-\tau_2)=\underline{\lambda}(\Omega,A,w_1,\tau_1,\tau_2)$. Then $\lambda=\underline{\lambda}(\Omega,A,w_1,\tau_1,\tau_2)$ and $\psi=\varphi_{\Omega,A,v,-\tau_2}$.
\end{lem}

\noindent{\bf Proof of Lemma \ref{uniquenessbis}.} Let $v$ be as above and set $\varphi=\varphi_{\Omega,A,v,-\tau_2}$, so that
\[
-\div(A\nabla\varphi)-\tau_1w_1^{-1}\left\vert \nabla\varphi\right\vert-\tau_2\varphi=-\div(A\nabla\varphi)+v\cdot \nabla\varphi-\tau_2\varphi=\underline{\lambda}(\Omega,A,w_1,\tau_1,\tau_2)\varphi\mbox{ in }\Omega
\]
by what we have just seen. Define also
\[
w(x)=
\left\{
\begin{array}{ll}
\displaystyle -\tau_1w_1(x)^{-1}\frac{\nabla\psi(x)}{\left\vert \nabla 
\psi(x)\right\vert} & \mbox{if }\nabla\psi(x)\neq 0,\\
\\
0 & \mbox{if }\nabla\psi(x)=0.
\end{array}
\right.
\]
One has $\left\Vert w_1w\right\Vert_{\infty}=\tau_1$ and
\[
-\div(A\nabla\psi)+w\cdot \nabla\psi-\tau_2\psi=-\div(A\nabla\psi)-\tau_1w_1^{-1}\left\vert 
\nabla \psi\right\vert-\tau_2\psi=\lambda\psi\mbox{ in }\Omega,
\]
so that, since $\psi>0$ in $\Omega$ and $\psi=0$ on $\partial\Omega$, by the characterization of the principal eigenfunction and the normalization $\|\psi\|_{\infty}=1$, one has $\psi=\varphi_{\Omega,A,w,-\tau_2}$ and
$$\lambda=\lambda_1(\Omega,A,w,-\tau_2)\geq\underline{\lambda}(\Omega,A,w_1,\tau_1,\tau_2).$$
As a consequence,
\[
-\div(A\nabla\psi)+v\cdot \nabla\psi-\tau_2\psi\geq -\div(A\nabla\psi)-\tau_1w_1^{-1}\left\vert 
\nabla\psi\right\vert-\tau_2\psi=\lambda\psi\geq 
\underline{\lambda}(\Omega,A,w_1,\tau_1,\tau_2)\psi\mbox{ in }\Omega,
\]
while
$$-\div(A\nabla\varphi)+v\cdot\nabla\varphi-\tau_2\varphi
=\lambda_1(\Omega,A,v,-\tau_2)\varphi
=\underline{\lambda}(\Omega,A,w_1,\tau_1,\tau_2)\varphi\hbox{ in }\Omega$$
by assumption. Since $\psi>0$ in $\Omega$, another application of Lemma~\ref{lemcomparison} shows that 
$\psi=\varphi=\varphi_{\Omega,A,v,-\tau_2}$, and that 
$\lambda=\underline{\lambda}(\Omega,A,w_1,\tau_1,\tau_2)$.\hfill\fin\break\par

With the help of Lemma \ref{uniquenessbis}, we conclude the proof of Theorem \ref{fixedlinfty}. Let $v_1\in L^{\infty}(\Omega,\R^n)$ and $v_2\in L^{\infty}(\Omega,\R^n)$ with $\left\Vert w_1v_1\right\Vert_{\infty}\leq \tau_1$ and $\left\Vert w_1v_2\right\Vert_{\infty}\leq \tau_1$ be such that
$$\underline{\lambda}(\Omega,A,w_1,\tau_1,\tau_2)=\lambda_1(\Omega,A,v_1,-\tau_2)=\lambda_1(\Omega,A,v_2,-\tau_2),$$
and set
$$\varphi_1=\varphi_{\Omega,A,v_1,-\tau_2}\hbox{ and }\varphi_2=\varphi_{\Omega,A,v_2-\tau_2}.$$
Since $v_1\cdot\nabla\varphi_1=-\tau_1w_1^{-1}\left\vert \nabla\varphi_1\right\vert$ and $v_2\cdot\nabla\varphi_2=-\tau_1w_1^{-1}\left\vert \nabla\varphi_2\right\vert$ a.e. in $\Omega$, one has
\[
\left\{\begin{array}{ll}
-\div(A\nabla\varphi_1)-\tau_1w_1^{-1}\left\vert \nabla\varphi_1\right\vert-\tau_2\varphi_1=\underline{\lambda}(\Omega,A,w_1,\tau_1,\tau_2)\varphi_1 & \mbox{ in }\Omega,\\
-\div(A\nabla\varphi_2)-\tau_1w_1^{-1}\left\vert \nabla\varphi_2\right\vert-\tau_2\varphi_2=\underline{\lambda}(\Omega,A,w_1,\tau_1,\tau_2)\varphi_2 & \mbox{ in }\Omega,
\end{array}\right.
\]
with $\varphi_1,\varphi_2\in\bigcap_{1\le r<+\infty}W^{2,r}(\Omega)$, $\varphi_1,\varphi_2>0$ in $\Omega$, $\|\varphi_1\|_{\infty}=\|\varphi_2\|_{\infty}=1$, and $\varphi_1=\varphi_2=0$ on $\partial\Omega$. Lemma \ref{uniquenessbis} shows that $\varphi_1=\varphi_{\Omega,A,v_2,-\tau_2}=\varphi_2:=\varphi$, so that $v_1\cdot\nabla\varphi=v_2\cdot\nabla\varphi=-\tau_1w_1^{-1}\left\vert \nabla\varphi\right\vert$. It follows that $v_1=v_2$ and $\left\vert v_1\right\vert=\left\vert v_2\right\vert=\tau_1w_1^{-1}$ a.e. on the set $\left\{x\in \Omega;\ \nabla\varphi(x)\neq 0\right\}$. It remains to be observed that $\nabla\varphi(x)\neq 0$ a.e. in $\Omega$. Indeed, if $E=\left\{x\in\Omega;\ \nabla\varphi(x)=0\right\}$, one has $\div(A\nabla\varphi)=0$ a.e. in $E$, so that $-\tau_2\varphi=\lambda_1(\Omega,A,v_1,-\tau_2)\varphi$ in $E$, and since $\lambda_1(\Omega,A,v_1,-\tau_2)>-\tau_2$, one has $\left\vert E\right\vert=0$. \hfill\fin\break

If, in Theorem \ref{fixedlinfty}, we specialize to the case when $\Omega$ is a ball and the diffusion matrix $A$ is equal to $\Lambda\Id$ with $\Lambda$ radially symmetric, we obtain a more complete description of the unique minimizer and maximizer. More precisely, we have:

\begin{theo} \label{fixedlinftyball}
Assume that $\Omega$ is a Euclidean ball centered at $0$ with radius $R>0$, let $\Lambda\in L^{\infty}_+(\Omega)\cap W^{1,\infty}(\Omega)$ be radially symmetric, set $A=\Lambda\!\!\Id$ and use the same notations as in Theorem \ref{fixedlinfty}, under the extra assumption that $w_1$ is radially symmetric Then, $\underline{v}=\tau_1w_1^{-1}e_r$, $\overline{v}=-\tau_1w_1^{-1}e_r$ a.e. in $\Omega$, and $\underline{\varphi}$ and $\overline{\varphi}$ are radially symmetric and decreasing.
\end{theo}

\noindent{\bf Proof. }Let $\varphi=\underline{\varphi}=\varphi_{\Omega,\Lambda\!\!\Id,\underline{v},-\tau_2}$ where $\underline{v}$ is given in Theorem~\ref{fixedlinfty}. One has
\begin{equation} \label{eqphi}
-\div(\Lambda\nabla\varphi)-\tau_1w_1^{-1}\left\vert \nabla\varphi\right\vert-\tau_2\varphi=\underline{\lambda}(\Omega,A,w_1,\tau_1,\tau_2)\varphi
\end{equation}
in $\Omega$. If $S$ is any orthogonal transformation in $\R^n$ and $\psi=\varphi\circ S$, then $\psi\in W^{2,r}(\Omega)$ for all $1\le r<+\infty$, satisfies (\ref{eqphi}), vanishes on $\partial\Omega$ and is positive in $\Omega$. Lemma \ref{uniquenessbis} therefore yields $\varphi=\psi$, which means that $\varphi$ is radially symmetric, so that there exists a function $u:\left[0,R\right]\rightarrow \R$ such that $\varphi(x)=u\left(\left\vert x\right\vert\right)$ for all $x\in \overline{\Omega}$, and $u$ is $C^{1,\alpha}(\left[0,R\right])$ for all $0\leq \alpha<1$. Let $0\leq r_1<r_2<R$. Remind that $\underline{\lambda}(\Omega,A,w_1,\tau_1,\tau_2)=\lambda_1(\Omega,A,\underline{v},-\tau_2)>-\tau_2$. Since
\begin{equation} \label{maxprincphi}
-\div(\Lambda\nabla\varphi)+\underline{v}\cdot\nabla\varphi=(\underline{\lambda}(\Omega,A,w_1,\tau_1,\tau_2)+\tau_2)\varphi>0
\end{equation}
in $\Omega$,  the maximum principle applied to $\varphi$ in $B_{r_2}$ yields that $\varphi\geq u(r_2)$ in $B_{r_2}$, which means that $u(r_1)\geq u(r_2)$. Moreover, if $u(r_1)=u(r_2)$, the strong maximum principle implies that $\varphi$ is constant in $B_{r_2}$, which is impossible because of (\ref{maxprincphi}). Therefore, $u(r_1)>u(r_2)$. Finally, if $0\leq r_1<r_2=R$, one has immediately $u(r_1)>u(r_2)=0$. Thus, $u$ is decreasing in $\left[0,R\right]$, and this yields at once $\underline{v}=\tau_1w_1^{-1}e_r$ from Theorem~\ref{fixedlinfty}.\par
The arguments for $\overline{\varphi}$ and $\overline{v}$ are entirely analogous and this completes the proof of Theorem~\ref{fixedlinftyball}.\hfill\fin


\subsection{Faber-Krahn inequalities}\label{sec62}

{\bf Proof of Theorem~\ref{faberkrahn}.} First, since $\|V\|_{\infty}\le\tau_2$, it follows from \cite{bnv} that
$$\lambda_1(\Omega,A,v,V)\ge\lambda_1(\Omega,A,v,-\tau_2)=-\tau_2+\lambda_1(\Omega,A,v,0).$$
But $\lambda_1(\Omega,A,v,0)>0$. Theorem~\ref{th1bis} then yields the existence of a positive constant $\theta=\theta(\Omega,n,\overline{M}_A,\underline{m}_{\Lambda},\tau_1)>0$ depending only on $(\Omega,n,\overline{M}_A,\underline{m}_{\Lambda},\tau_1)$, and the existence of two radially symmetric $C^{\infty}(\overline{\Omega^*})$ fields $\Lambda^*>0$, $\omega^*\ge 0$ such that, for $v^*=\omega^*e_r$ in $\Omega^*$,
$$\left\{\baa{l}
\displaystyle{\mathop{\rm{ess}\ \rm{inf}}_{\Omega}}\ \Lambda \le \displaystyle{\mathop{\min}_{\overline{\Omega^*}}}\ \Lambda^*\le \displaystyle{\mathop{\max}_{\overline{\Omega^*}}}\ \Lambda^*\le \displaystyle{\mathop{\rm{ess}\ \rm{sup}}_{\Omega}}\ \Lambda,\ \|(\Lambda^*)^{-1}\|_{L^1(\Omega^*)}=\|\Lambda^{-1}\|_{L^1(\Omega)},\\
\|v^*\|_{L^{\infty}(\Omega^*,\R^n)}\le\|v\|_{L^{\infty}(\Omega,\R^n)},\ \|\ |v^*|^2\times(\Lambda^*)^{-1}\|_{L^1(\Omega^*)}=\|\ |v|^2\times\Lambda^{-1}\|_{L^1(\Omega)},\eaa\right.$$
and
$$\lambda_1(\Omega,A,v,0)\ge\lambda_1(\Omega^*,\Lambda^*{\rm{Id}},v^*,0)\times(1+\theta).$$
Observe that $\left\Vert v^{\ast}\right\Vert_{\infty}\leq \tau_1$. It follows from Theorem~\ref{fixedlinftyball} (with $w_1=1$) that
\[
-\tau_2+\lambda_1(\Omega^{\ast},\Lambda^{\ast}\hbox{Id},v^{\ast},0)=\lambda_1(\Omega^{\ast},\Lambda^{\ast}\hbox{Id},v^{\ast},-\tau_2)\ge\lambda_1(\Omega^{\ast},\Lambda^{\ast}\hbox{Id},\tau_1e_r,-\tau_2)
\]
and $\lambda_1(\Omega^*,\Lambda^*{\rm{Id}},v^*,0)\ge\lambda_1(\Omega^*,\Lambda^*{\rm{Id}},\tau_1e_r,0)$. Therefore,
$$\lambda_1(\Omega,A,v,V)\ge\lambda_1(\Omega^{\ast},\Lambda^{\ast}\hbox{Id},\tau_1e_r,-\tau_2)+\theta\times\lambda_1(\Omega^*,\Lambda^*{\rm{Id}},\tau_1e_r,0)$$
since $\theta$ and $\lambda_1(\Omega^*,\Lambda^*{\rm{Id}},\tau_1e_r,0)$ are positive.\par
Let us now estimate $\lambda_1(\Omega^*,\Lambda^*{\rm{Id}},\tau_1e_r,0)$ from below. Call $\varphi=\varphi_{\Omega^*,\Lambda^*{\rm{Id}},\tau_1e_r,0}$, $\lambda=\lambda_1(\Omega^*,\Lambda^*{\rm{Id}},\tau_1e_r,0)>0$ and
$$U(x)=\tau_1\int_0^{|x|}\Lambda^*(re)^{-1}dr\hbox{ for all }x\in\overline{\Omega^*},$$
where $e$ is an arbitrary unit vector. Multiply the equation
$$-\div(\Lambda^*\nabla\varphi)+\tau_1e_r\cdot\nabla\varphi=\lambda\varphi\hbox{ in }\Omega^*$$
by $\varphi e^{-U}\in H^1_0(\Omega^*)\cap W^{1,\infty}(\Omega^*)$ and integrate by parts over $\Omega^*$. It follows from the definition of $U$ that
$$\int_{\Omega^*}\Lambda^*|\nabla\varphi|^2e^{-U}=\lambda\int_{\Omega^*}\varphi^2e^{-U}\le\lambda\int_{\Omega^*}\varphi^2.$$
The last inequality holds since $\lambda>0$, and $U\ge 0$ in $\Omega^*$. But $\Lambda^*\ge\hbox{ess inf}_{\Omega}\ \Lambda\ge\underline{m}_{\Lambda}>0$, whence $U\le\tau_1\underline{m}_{\Lambda}^{-1}R$ in $\overline{\Omega^*}$, where $R=\alpha_n^{-1/n}|\Omega|^{1/n}>0$ is the radius of $\Omega^*$. Finally,
$$\lambda\int_{\Omega^*}\varphi^2\ge\underline{m}_{\Lambda}e^{-\tau_1\underline{m}_{\Lambda}^{-1}\alpha_n^{-1/n}|\Omega|^{1/n}}\int_{\Omega}|\nabla\varphi|^2,$$
whence
$$\lambda\ge\underline{m}_{\Lambda}e^{-\tau_1\underline{m}_{\Lambda}^{-1}\alpha_n^{-1/n}|\Omega|^{1/n}}\times|\Omega|^{-2/n}\alpha_n^{2/n}(j_{n/2-1,1})^2=:\kappa>0$$
from (\ref{RFK}) and (\ref{variational}). The conclusion of Theorem~\ref{faberkrahn} follows with the choice
$$\eta=\eta(\Omega,n,\overline{M}_A,\underline{m}_{\Lambda},\tau_1)=\theta\times\kappa>0.$$\hfill\fin\break

\noindent{\bf{Proof of Corollary~\ref{corFK}.}} Assume first that $\Omega$ is not a ball. Write
$$\lambda_1(\Omega,A,v,V)\ge\lambda_1(\Omega,A,v,0)+\mathop{\hbox{ess inf}}_{\Omega}V.$$
Under the notations of Corollary~\ref{corFK}, then Theorem~\ref{faberkrahn} applied to $\lambda_1(\Omega,A,v,0)$ with $\Lambda=\gamma_A$ clearly gives $\Lambda^{\ast}=\gamma_A$ in $\overline{\Omega^{\ast}}$, so that
$$\lambda_1(\Omega,A,v,0)>\lambda_1(\Omega^{\ast},\gamma_A{\hbox{Id}},\|v\|_{\infty}e_r,0),$$
whence
$$\lambda_1(\Omega,A,v,V)>\lambda_1(\Omega^{\ast},\gamma_A{\hbox{Id}},\|v\|_{\infty}e_r,0)+\mathop{\hbox{ess inf}}_{\Omega}\ V=\lambda_1(\Omega^*,\gamma_A{\hbox{Id}},\|v\|_{\infty}e_r,\mathop{\hbox{ess inf}}_{\Omega}\ V).$$
Assume now that $\Omega$ is a ball. From Theorem~\ref{th1} applied to $\lambda_1(\Omega,A,v,0)$ with $\Lambda=\gamma_A$, there exists $v^*\in L^{\infty}(\Omega^*,\R^n)$ such that $\|v^*\|_{L^{\infty}(\Omega^*,\R^n)}\le\|v\|_{L^{\infty}(\Omega,\R^n)}$ and
$$\lambda_1(\Omega^*,\gamma_A\hbox{Id},v^*,0)\le\lambda_1(\Omega,A,v,0).$$
But $\lambda_1(\Omega^*,\gamma_A\hbox{Id},\|v\|_{L^{\infty}(\Omega,\R^n)}e_r,0)\le\lambda_1(\Omega^*,\gamma_A\hbox{Id},v^*,0)$ from Theorem~\ref{fixedlinftyball} with $w_1=1$. Therefore,
$$\lambda_1(\Omega,A,v,V)\ge\lambda_1(\Omega^*,\gamma_A\hbox{Id},\|v\|_{L^{\infty}(\Omega,\R^n)}e_r,\mathop{\hbox{ess inf}}_{\Omega}V).$$
The conclusion of Corollary~\ref{corFK} follows immediately.\hfill\fin

\begin{rem}\label{interpretation}{\rm If, in Theorem~\ref{faberkrahn}, we specialize to the case when $A=\gamma\hbox{Id}$ and $\Lambda=\gamma>0$ is a given constant, then we have immediately
$$\lambda_1(\Omega,\gamma\hbox{Id},v,V)>\lambda_1(\Omega^{\ast},\gamma{\hbox{Id}},\|v\|_{\infty}e_r,-\|V\|_{\infty})$$
provided that $\Omega\in{\mathcal{C}}$ is not a ball. Furthermore, if $\Omega$ is a ball, say with center $x_0$, the uniqueness statement in Theorem~\ref{fixedlinftyball} shows that $\lambda_1(\Omega,\gamma{\hbox{Id}},v,V)\ge\lambda_1(\Omega,\gamma{\hbox{Id}},\|v\|_{\infty}e_r(\cdot-x_0),-\|V\|_{\infty})$, where the inequality is strict if $v\neq\|v\|_{\infty}e_r(\cdot-x_0)$ or $V\neq-\|V\|_{\infty}$. Finally, we obtain that, if $\Omega\in {\mathcal C}$, then
\be\label{fkgamma}
\lambda_1(\Omega,\gamma{\hbox{Id}},v,V)\geq \lambda_1(\Omega^{\ast},\gamma{\hbox{Id}},\|v\|_{\infty}e_r,-\|V\|_{\infty})
\ee
and the equality holds if and only if, up to translation, $\Omega=\Omega^{\ast}$, $v=\|v\|_{\infty}e_r$ and $V=-\|V\|_{\infty}$.\par
A rough parabolic interpretation of inequality (\ref{fkgamma}) can be as follows: consider the evolution equation $u_t=\gamma\Delta u-v\cdot\nabla u-Vu$ in $\Omega$, for $t>0$, with Dirichlet boundary condition on $\partial\Omega$, and with an initial datum at $t=0$. Roughly speaking, minimizing $\lambda_1(\Omega,\gamma{\hbox{Id}},v,V)$ with given measure $|\Omega|$ and with given $L^{\infty}$ constraints $\|v\|_{\infty}\le\tau_1$ and $\|V\|_{\infty}\le\tau_2$ can be interpreted as looking for the slowest exponential time-decay of the solution $u$. The best way to do that is: 1)~to try to minimize the boundary effects, namely to have the domain as round as possible, 2)~to have $-V$ as large as possible, that is $V$ as small as possible, and 3) it is not unreasonable to say that the vector field $-v$ should as much as possible point inwards the domain to avoid the drift towards the boundary. Of course diffusion, boundary losses, transport and reaction phenomena take place simultaneously, but these heuristic arguments tend to lead to the optimal triple $(\Omega,-v,-V)=(\Omega^*,-\tau_1e_r,\tau_2)$ (up to translation).\par
In Theorem~\ref{faberkrahn}, it follows from the above proofs that the inequality
$$\lambda_1(\Omega^*,\Lambda^*\hbox{Id},\tau_1e_r,\tau_3)\le\lambda_1(\Omega,A,v,V)-\eta$$
holds if the assumption $\|V\|_{\infty}=\|V\|_{L^{\infty}(\Omega)}\le\tau_2$ is replaced by: $\hbox{ess inf}_{\Omega}V\ge\tau_3$. Furthermore, since $\lambda_1(\Omega,A,v,V)\ge\lambda_1(\Omega,A,v,\hbox{ess inf}_{\Omega}V)=\lambda_1(\Omega,A,v,0)+\hbox{ess inf}_{\Omega}V$ with a strict inequality if $V$ is not constant (see \cite{bnv}), it is then immediate to check that formula (\ref{fkgamma}) still holds when $-\|V\|_{\infty}$ is replaced by $\hbox{ess inf}_{\Omega}V$, and that the case of equality can be extended similarly. The parabolic interpretation is the same as above if the condition $\|V\|_{\infty}\le\tau_2$ is replaced by: $\hbox{ess inf}_{\Omega}V\ge\tau_3$.}
\end{rem}

\begin{rem}\label{nonsmooth3}{\rm If $\Omega$ is a general open subset of $\R^n$ with finite measure, and if $A\in W^{1,\infty}(\Omega,{\mathcal{S}}_n(\R))$, $v\in L^{\infty}(\Omega,\R^n)$, $V\in L^{\infty}(\Omega)$ are such that $A\ge\gamma\hbox{Id}$ in $\Omega$ for some constant $\gamma>0$, then we claim that
\be\label{comparns3}\baa{rcl}
\lambda_1(\Omega,A,v,V) & \ge & \lambda_1(\Omega^*,\gamma\hbox{Id},\|v\|_{L^{\infty}(\Omega,\R^n)}e_r,\displaystyle{\mathop{\hbox{ess inf}}_{\Omega}}\ V)\\
& \ge & \lambda_1(\Omega^*,\gamma\hbox{Id},\|v\|_{L^{\infty}(\Omega,\R^n)}e_r,-\|V\|_{L^{\infty}(\Omega)}).\eaa
\ee
Indeed, given $\epsilon>0$, as in Remark~\ref{nonsmooth1}, there exists a non-empty set $\Omega'=\Omega'_{\epsilon}\in{\mathcal{C}}$ such that $\Omega'\subset\subset\Omega$ and $\lambda_1(\Omega',A,v,V)\le\lambda_1(\Omega,A,v,V)+\epsilon$. Then the arguments used in the proof of Corollary~\ref{corFK} (with $\gamma$ instead of $\gamma_A$) imply that
$$\lambda_1(\Omega',A,v,V)\ge\lambda_1(\Omega^*_{\epsilon},\gamma\hbox{Id},\|v\|_{L^{\infty}(\Omega,\R^n)}e_r,\mathop{\hbox{ess inf}}_{\Omega}\ V),$$
where $\Omega^*_{\epsilon}$ is the ball centered at the origin with the same measure as $\Omega'$. 
Therefore, 
$$\lambda_1(\Omega,A,v,V)+\epsilon\ge\lambda_1(\Omega^*_{\epsilon},\gamma\hbox{Id},\|v\|_{L^{\infty}(\Omega,\R^n)}e_r,\mathop{\hbox{ess inf}}_{\Omega}\ V)\ge\lambda_1(\Omega^*,\gamma\hbox{Id},\|v\|_{L^{\infty}(\Omega,\R^n)}e_r,\mathop{\hbox{ess inf}}_{\Omega}\ V)$$
since $\Omega^*_{\epsilon}\subset\Omega^*$, and (\ref{comparns3}) follows immediately. Notice that (\ref{comparns3}) holds in particular with $\gamma=\hbox{ess inf}_{\Omega}\ \Lambda[A]$.}
\end{rem}


\SE{Appendix}\label{secapp}


\subsection{Proof of the approximation lemma \ref{lemapprox}}

Fix $k\in\N$. Call
$$r_{i,k}=\frac{iR}{k+1}\ \hbox{ for }i=0,\ldots,k+1$$
and
$$r_{i+1/2,k}=\left(\frac{r_{i,k}^n+r_{i+1,k}^n}{2}\right)^{1/n}\in(r_{i,k},r_{i+1,k})\ \hbox{ for } i=0,\ldots,k.$$
Remember that
$$|\Omega_{\rho^{-1}(r_{i+1,k}),\rho^{-1}(r_{i,k})}|=\alpha_n(r_{i+1,k}^n-r_{i,k}^n)=|S_{r_{i,k},r_{i+1,k}}|\hbox{ for all }i\in\{0,\ldots,k\}.$$
Let us first define the function $g_k$ almost everywhere in $\Omega^*$: for $i\in\{0,\ldots,k\}$ and $x\in S_{r_{i,k},r_{i+1,k}}$ such that $|x|\neq r_{i+1/2,k}$, set
$$g_k(x)=G_k(|x|),$$
where
$$G_k(r)=\sup\left\{a\in\R;\ \displaystyle{\frac12}\left|\left\{x\in\Omega_{\rho^{-1}(r_{i+1,k}),\rho^{-1}(r_{i,k})},\ g(x)>a\right\}\right|\ge\alpha_n\left|r_{i+1/2,k}^n-r^n\right|\right\}$$
for all $r\in(r_{i,k},r_{i+1/2,k})\cup(r_{i+1/2,k},r_{i+1,k})$.\footnote{In each shell $S_{r_{i,k},r_{i+1,k}}$, the function $g_k$ is then a kind of Schwarz decreasing rearrangement of the function $g$ in $\Omega_{\rho^{-1}(r_{i+1,k}),\rho^{-1}(r_{i,k})}$, with respect to the inner radius $r_{i+1/2,k}$.} It then follows by definition that $g_k$ is radially symmetric, nondecreasing with respect to $|x|$ in $S_{r_{i,k},r_{i+1/2,k}}$ and nonincreasing with respect to $|x|$ in $S_{r_{i+1/2,k},r_{i+1,k}}$, and that
$$\left\{\baa{l}
\left|\left\{x\in\Omega_{\rho^{-1}(r_{i+1,k}),\rho^{-1}(r_{i,k})},\ g(x)> t\right\}\right|=\left|\left\{x\in S_{r_{i,k},r_{i+1,k}},\ g_k(x)>t\right\}\right|\\
\\
\left|\left\{x\in\Omega_{\rho^{-1}(r_{i+1,k}),\rho^{-1}(r_{i,k})},\ g(x)\ge t\right\}\right|=\left|\left\{x\in S_{r_{i,k},r_{i+1,k}},\ g_k(x)\ge t\right\}\right|\eaa\right.$$
for all $i\in\{0,\ldots,k\}$ and $t\in\R$. As a consequence, the restriction of $g_k$ to $S_{r_{i,k},r_{i+1,k}}$ is in $L^{\infty}(S_{r_{i,k},r_{i+1,k}})$,
$$\displaystyle{\mathop{\hbox{ess inf}}_{\Omega}}\ g\le\displaystyle{\mathop{\hbox{ess inf}}_{\Omega_{\rho^{-1}(r_{i+1,k}),\rho^{-1}(r_{i,k})}}}g=\displaystyle{\mathop{\hbox{ess inf}}_{S_{r_{i,k},r_{i+1,k}}}}g_k\le\displaystyle{\mathop{\hbox{ess sup}}_{S_{r_{i,k},r_{i+1,k}}}}g_k=\displaystyle{\mathop{\hbox{ess sup}}_{\Omega_{\rho^{-1}(r_{i+1,k}),\rho^{-1}(r_{i,k})}}}g\le\displaystyle{\mathop{\hbox{ess sup}}_{\Omega}}\ g$$
and
\be\label{gkh}
\int_{S_{r_{i,k},r_{i+1,k}}}g_k=\int_{\Omega_{\rho^{-1}(r_{i+1,k}),\rho^{-1}(r_{i,k})}}g=\int_{S_{r_{i,k},r_{i+1,k}}}h
\ee
for all $i=0,\ldots,k$, by assumption (\ref{gh}). Therefore, $g_k\in L^{\infty}(\Omega^*)$ with
$$\displaystyle{\mathop{\hbox{ess inf}}_{\Omega}}\ g=\displaystyle{\mathop{\hbox{ess inf}}_{\Omega^*}}\ g_k\le\displaystyle{\mathop{\hbox{ess sup}}_{\Omega^*}}\ g_k=\displaystyle{\mathop{\hbox{ess sup}}_{\Omega}}\ g$$
and
$$\left\{\baa{l}
|\{x\in\Omega,\ g(x)>t\}|=|\{x\in\Omega^*,\ g_k(x)>t\}|\\
\\
|\{x\in\Omega,\ g(x)\ge t\}|=|\{x\in\Omega^*,\ g_k(x)\ge t\}|\eaa\right.$$
for all $t\in\R$.\par
Let us now define the sequence of functions $(\underline{g}_k)_{k\in\N}$. Fix $k\in\N$. For each $i\in\{0,\ldots,k\}$, the function $G_k$ is by construction nondecreasing in the interval $(r_{i,k},r_{i+1/2,k})$ and nonincreasing in the interval $(r_{i+1/2,k},r_{i+1,k})$. Furthermore,
$$\mathop{\hbox{ess inf}}_{(r_{i,k},r_{i+1/2,k})}G_k=\mathop{\hbox{ess inf}}_{(r_{i+1/2,k},r_{i+1/2,k})}G_k=\mathop{\hbox{ess inf}}_{\Omega_{\rho^{-1}(r_{i+1,k}),\rho^{-1}(r_{i,k})}}g\ge\mathop{\hbox{ess inf}}_{\Omega}\ g.$$
Therefore, in each of the intervals $(r_{i,k},r_{i+1/2,k})$ and $(r_{i+1/2,k},r_{i+1,k})$, the function $G_k$ can be approximated uniformly and from below by piecewise constant functions which are larger than or equal to $\hbox{ ess inf}_{\Omega}\ g$. As a consequence, there exists a piecewise constant function $\widetilde{G}_k$ defined in $[0,R]$ such that
\be\label{tildeGk}
\mathop{\hbox{ess inf}}_{\Omega}\ g\le\widetilde{G}_k(r)\le G_k(r)\ (\le\mathop{\hbox{ess sup}}_{\Omega}\ g)\ \hbox{ for all }r\in(0,r_{1/2,k})\cup\cdots\cup(r_{k+1/2,k},R)
\ee
and
\be\label{tildeGkbis}
\|G_k-\widetilde{G}_k\|_{L^{\infty}(0,R)}\le\frac{1}{k+1}.
\ee
Let $0=\rho_{0,k}<\rho_{1,k}<\cdots<\rho_{N_k+1,k}=R$ be a subdivision adapted to $\widetilde{G}_k$ (with $N_k\in\N$), namely $\widetilde{G}_k$ is equal to a constant $m_{j,k}\in[\mathop{\hbox{ess inf}}_{\Omega}g,\mathop{\hbox{ess sup}}_{\Omega}g]$ in each interval $(\rho_{j,k},\rho_{j+1,k})$ for $j=0,\ldots,N_k$. Choose a real number $\underline{\rho}_k$ such that
$$0<\underline{\rho}_k<\min\left(\min_{0\le j\le N_k}\frac{\rho_{j+1,k}-\rho_{j,k}}{2},\frac{1}{(N_k+1)(k+1)}\right).$$
Let $\zeta$ be a fixed $C^{\infty}(\R,\R)$ function such that $0\le \zeta\le 1$ in $\R$, $\zeta=0$ in $(-\infty,1/3]$ and $\zeta=1$ in $[2/3,+\infty)$. Denote $\underline{G}_k$ the function in $[0,R]$ by:
$$\left\{\baa{ll}
\underline{G}_k(r)=m_{0,k} & \hbox{for }r\in[0,\rho_{1,k}-\underline{\rho}_k],\\
\underline{G}_k(r)=m_{j,k} & \hbox{for }r\in[\rho_{j,k}+\underline{\rho}_k,\rho_{j+1,k}-\underline{\rho}_k]\hbox{ and }1\le j\le N_k-1\hbox{ (if }N_k\ge 1\hbox{)},\\
\underline{G}_k(r)=m_{N_k,k} & \hbox{for }r\in[\rho_{N_k,k}+\underline{\rho}_k,R],\eaa\right.$$
$$\baa{l}\forall\ 0\le j\le N_k-1,\\ \forall\ r\in[\rho_{j+1,k}-\underline{\rho}_k,\rho_{j+1,k}],\eaa\ 
\underline{G}_k(r)=\left\{\baa{l}m_{j,k}\qquad\qquad\qquad\qquad\qquad\qquad\ \hbox{if }m_{j,k}\le m_{j+1,k},\\
m_{j,k}+(m_{j+1,k}-m_{j,k})\times\zeta\left(\displaystyle{\frac{r-\rho_{j+1,k}+\underline{\rho}_k}{\underline{\rho}_k}}\right)\\
\qquad\qquad\qquad\qquad\qquad\qquad\qquad\hbox{if }m_{j+1,k}<m_{j,k},\eaa\right.$$
and
$$\baa{l}\forall\ 1\le j\le N_k,\\ \forall\ r\in[\rho_{j,k},\rho_{j,k}+\underline{\rho}_k],\eaa\ 
\underline{G}_k(r)=\left\{\baa{l}m_{j,k}\qquad\qquad\qquad\qquad\qquad\qquad\hbox{if }m_{j-1,k}\ge m_{j,k}\\
m_{j,k}+(m_{j-1,k}-m_{j,k})\times\zeta\left(\displaystyle{\frac{\rho_{j,k}+\underline{\rho}_k-r}{\underline{\rho}_k}}\right)\\
\qquad\qquad\qquad\qquad\qquad\qquad\qquad\hbox{if }m_{j-1,k}<m_{j,k}.\eaa\right.$$
The function $\underline{G}_k$ is well-defined and $C^{\infty}$ in $[0,R]$ and
$$\mathop{\hbox{ess inf}}_{(0,R)}\ \widetilde{G}_k=\min_{0\le j\le N_k}m_{j,k}\le\underline{G}_k(r)\le\widetilde{G}_k(r)\hbox{ for all }r\in(0,\rho_{1,k})\cup\cdots\cup(\rho_{N_k,k},R),$$
whence
$$\mathop{\hbox{ess inf}}_{\Omega}\ g\le\underline{G}_k\le\widetilde{G}_k\le G_k\le\mathop{\hbox{ess sup}}_{\Omega}\ g\ \hbox{ almost everywhere in }[0,R]$$
by (\ref{tildeGk}). The function defined by
$$\underline{g}_k(x)=\underline{G}_k(|x|)\hbox{ for all }x\in\overline{\Omega^*}$$
is radially symmetric and of class $C^{\infty}(\overline{\Omega^*})$ and it satisfies
$$\mathop{\hbox{ess inf}}_{\Omega}\ g\le\underline{g}_k\le g_k\le\mathop{\hbox{ess sup}}_{\Omega}\ g\ \hbox{ almost everywhere in }\Omega^*.$$\par
Fix now $q\in[1,+\infty)$ and let us check that $g_k-\underline{g}_k\to 0$ in $L^q(\Omega^*)$ as $k\to+\infty$. One has
$$\baa{rcl}
\|g_k-\underline{g}_k\|_{L^q(\Omega^*)} & \le & \|G_k(|\cdot|)-\widetilde{G}_k(|\cdot|)\|_{L^q(\Omega^*)}+\|\widetilde{G}_k(|\cdot|)-\underline{G}_k(|\cdot|)\|_{L^q(\Omega^*)}\\
& \le & \displaystyle{\frac{(\alpha_nR^n)^{1/q}}{k+1}}+\|\widetilde{G}_k(|\cdot|)-\underline{G}_k(|\cdot|)\|_{L^q(\Omega^*)}\eaa$$
by (\ref{tildeGkbis}). On the other hand, the definition of $\underline{G}_k$ and formula (\ref{tildeGk}) imply that
$$\|\widetilde{G}_k-\underline{G}_k\|_{L^{\infty}(0,R)}=\max_{0\le j\le N_k}|m_{j,k}-m_{j+1,k}|\le 2\|\widetilde{G}_k\|_{L^{\infty}(0,R)}\le 2\|g\|_{L^{\infty}(\Omega)}.$$
Using once again the definition of $\underline{G}_k$, it follows that
$$\baa{rcl}
\|\widetilde{G}_k(|\cdot|)-\underline{G}_k(|\cdot|)\|_{L^q(\Omega^*)} & = & \left[n\alpha_n\displaystyle{\sum_{j=0}^{N_k}}\left(\displaystyle{\int_{\rho_{j,k}}^{\rho_{j,k}+\underline{\rho}_k}}(\widetilde{G}_k(r)-\underline{G}_k(r))^qr^{n-1}dr\right.\right.\\
& & \qquad\qquad\left.\left.+\displaystyle{\int_{\rho_{j+1,k}-\underline{\rho}_k}^{\rho_{j+1,k}}}(\widetilde{G}_k(r)-\underline{G}_k(r))^qr^{n-1}dr\right)\right]^{1/q}\\
& \le & \left[2n\alpha_n(N_k+1)\underline{\rho}_kR^{n-1}(2\|g\|_{L^{\infty}(\Omega)})^{q}\right]^{1/q}\\
& \le & 2\|g\|_{L^{\infty}(\Omega)}\times\left(\displaystyle{\frac{2n\alpha_nR^{n-1}}{k+1}}\right)^{1/q}\eaa$$
from the choice of $\underline{\rho}_k$. Thus,
$$\|g_k-\underline{g}_k\|_{L^q(\Omega^*)}\le\frac{(\alpha_nR^n)^{1/q}}{k+1}+2\|g\|_{L^{\infty}(\Omega)}\times\left(\frac{2n\alpha_nR^{n-1}}{k+1}\right)^{1/q}\to 0\hbox{ as }k\to+\infty.$$\par
Finally, let us check that the sequences $(g_k)_{k\in\N}$ and $(\underline{g}_k)_{k\in\N}$ converge to $h$ as $k\to+\infty$ in $L^p(\Omega^*)$ weak for all $1<p<+\infty$ and in $L^{\infty}(\Omega^*)$ weak-*. Let $\phi$ be in $C(\overline{\Omega^*},\R)$ and fix $\epsilon>0$. Since the unit sphere $\mathbb{S}^{n-1}$ is compact and $\phi$ is uniformly continuous in $\overline{\Omega^*}$, there exists $k_0\in\N$ and a finite family of measurable pairwise disjoint subsets $U_1,\ldots,U_q$ of $\mathbb{S}^{n-1}$ with positive area, such that $\mathbb{S}^{n-1}=U_1\cup\cdots\cup U_q$, and
\be\label{kq}
|\phi(x)-\phi(y)|\le\epsilon\hbox{ for all }x,y\in\overline{\Omega^*}\backslash\{0\}\hbox{ such that }\left\{\baa{l}\left|\ |x|-|y|\ \right|\le\displaystyle{\frac{1}{k_0+1}},\\
\exists\ j\hbox{ such that }\displaystyle{\frac{x}{|x|}},\ \displaystyle{\frac{y}{|y|}}\in U_j.\eaa\right.
\ee
Fix any $k$ such that $k\ge k_0$. Use the notation $x=r\theta$ with $r=|x|$ and $\theta=x/|x|$ for the points of $\Omega^*\backslash\{0\}$. Denote by $d\sigma$ the surface measure on $\mathbb{S}^{n-1}$. For all $i\in\{0,\ldots,k\}$ and $j\in\{1,\ldots,q\}$, call
$$\phi_{i,j}=\frac{\displaystyle{\int_{r_{i,k}}^{r_{i+1,k}}}\displaystyle{\int_{U_j}}\phi(r\theta)\ d\sigma(\theta)\ dr}{\displaystyle{\int_{r_{i,k}}^{r_{i+1,k}}}\displaystyle{\int_{U_j}}\ d\sigma(\theta)\ dr}.$$
Since $g_k$ and $h$ are radially symmetric and satisfy (\ref{gkh}) for all $i=0,\ldots,k$, it follows that
$$\int_{r_{i,k}}^{r_{i+1,k}}\int_{U_j}(g_k(r\theta)-h(r\theta))\phi_{i,j}\ d\sigma(\theta)\ dr=0$$
for all $i\in\{0,\ldots,k\}$ and $j\in\{1,\ldots,q\}$. Thus,
$$\int_{\Omega^*}g_k\phi-\int_{\Omega^*}h\phi=\sum_{i=0}^k\ \sum_{j=1}^q\int_{r_{i,k}}^{r_{i+1,k}}\int_{U_j}(g_k(r\theta)-h(r\theta))(\phi(r\theta)-\phi_{i,j})\ d\sigma(\theta)\ dr$$
and then
$$\baa{rcl}
\left|\displaystyle{\int_{\Omega^*}}g_k\phi-\displaystyle{\int_{\Omega^*}}h\phi\right| & \le & \displaystyle{\sum_{i=0}^k}\ \displaystyle{\sum_{j=1}^q}\displaystyle{\int_{r_{i,k}}^{r_{i+1,k}}}\displaystyle{\int_{U_j}}(\|g\|_{L^{\infty}(\Omega)}+\|h\|_{L^{\infty}(\Omega^*)})\ \epsilon\ d\sigma(\theta)\ dr\\
& = & \alpha_nR^n(\|g\|_{L^{\infty}(\Omega)}+\|h\|_{L^{\infty}(\Omega^*)})\ \epsilon\eaa$$
for all $k\ge k_0$ (remember that $\|g_k\|_{L^{\infty}(\Omega^*)}=\|g\|_{L^{\infty}(\Omega)}$). Since $\epsilon>0$ was arbitrary, one concludes that
$$\int_{\Omega^*}g_k\phi\to\int_{\Omega^*}h\phi\hbox{ as }k\to+\infty.$$
Since this is true for every $\phi\in C(\overline{\Omega^*},\R)$, standard density arguments imply then that
\be\label{weak}
\int_{\Omega^*}g_k\phi\to\int_{\Omega^*}h\phi\hbox{ as }k\to+\infty,\hbox{ for all }\phi\in L^{p'}(\Omega^*)\hbox{ and for all }p'\in[1,+\infty),
\ee
namely $g_k\rightharpoonup h$ as $k\to+\infty$ in $L^p(\Omega^*)$ weak for all $p\in(1,+\infty)$ and in $L^{\infty}(\Omega^*)$ weak-*. Lastly, since $g_k-\underline{g}_k\to 0$ as $k\to+\infty$ in $L^p(\Omega^*)$ for all $p\in[1,+\infty)$ and since the functions $\underline{g}_k$ are uniformly bounded in $L^{\infty}(\Omega^*)$, it follows from (\ref{weak}) and standard density arguments that
$$\int_{\Omega^*}\underline{g}_k\phi\to\int_{\Omega^*}h\phi\hbox{ as }k\to+\infty,\hbox{ for all }\phi\in L^{p'}(\Omega^*)\hbox{ and for all }p'\in[1,+\infty).$$
Thus, $\underline{g}_k\rightharpoonup h$ as $k\to+\infty$ in $\sigma(L^p(\Omega^*),L^{p'}(\Omega^*))$ for all $1<p\le+\infty$.\par
The construction of the functions $\overline{g}_k$ is similar to that of the functions $\underline{g}_k$, but they approximate the functions $g_k$ from above.\hfill\fin


\subsection{A remark on distribution functions} \label{distribution}

Let $\alpha\le\beta\in \R$ and $m>0$ be fixed. We extend a definition which we used just before Corollary~\ref{cor1}: ${\mathcal F}_{\alpha,\beta}(m)$ stands for the set of right-continuous non-increasing functions $\mu:\R\rightarrow \left[0,m\right]$ such that
\[
\mu(t)=m\mbox{ for all } t<\alpha\mbox{ and }\mu(t)=0\mbox{ for all }t\geq \beta.
\]
In this appendix, we prove the following fact:

\begin{pro} \label{givendistri}
Let $\alpha\le\beta\in \R$, $m>0$, $\mu\in {\mathcal F}_{\alpha,\beta}(m)$ and $\Omega\in {\mathcal C}$ such that $|\Omega|=m$. Then, there exists $V\in L^{\infty}(\Omega)$ such that $\mu=\mu_V$.
\end{pro}

\noindent{\bf Proof. }This fact is rather classical, but we give here a quick proof for the sake of completeness. Let $\varphi$ be the solution of
\[
\left\{
\begin{array}{rcll}
-\Delta\varphi & = & 1 & \mbox{in }\Omega,\\
\varphi & = & 0 & \mbox{on }\partial\Omega.
\end{array}
\right.
\]
Observe that the function $\varphi$\ belongs to $W^{2,p}(\Omega)$ for all $1\leq p<+\infty$, to $C^{1,\gamma}(\overline{\Omega})$ for all $0\leq \gamma<1$ and is analytic and positive in $\Omega$. Let $M=\max_{\overline{\Omega}} \varphi$ and, for all $0\leq a<M$, define (as in Section \ref{rearrangement})
\[
\Omega_a=\left\{x\in \Omega;\ \varphi(x)>a\right\}.
\]
Set also $\Omega_M=\emptyset$. Remember that, for all $0\leq a\leq M$, $\left\vert\partial\Omega_a\right\vert=0$.\par

Define now, for all $x\in \Omega$,
\[
V(x)=\sup\left\{s\in \R;\ \mu(s)>\left\vert \Omega_{\varphi(x)}\right\vert\right\}.
\]
Notice first that this supremum is well-defined for all $x\in \Omega$. Indeed, if $x\in \Omega$, one has $\varphi(x)>0$, therefore $0\leq \left\vert \Omega_{\varphi(x)}\right\vert<\left\vert \Omega\right\vert$.\par

We now claim that $V$ is measurable and bounded in $\Omega$ and that $\mu_V=\mu$. Indeed, let $t\in \R$. By definition of $V$, for all $x\in \Omega$,
\[
V(x)>t\Leftrightarrow \left(\exists s>t\mbox{ such that }\mu(s)>\left\vert \Omega_{\varphi(x)}\right\vert\right)\Leftrightarrow \left\vert \Omega_{\varphi(x)}\right\vert<\mu(t),
\]
where the last equivalence follows from the right-continuity of $\mu$ and the fact that this function is non-increasing. Define now, for all $0\leq a\leq M$, $F(a)=\left\vert \Omega_a\right\vert$. The previous equivalence yields
\[
\mu_V(t)=\left\vert \left\{x\in \Omega;\ F(\varphi(x))<\mu(t)\right\}\right\vert.
\]
Since the function $F:\left[0,M\right]\rightarrow \left[0,\left\vert \Omega\right\vert\right]$ is decreasing, one-to-one and onto, one obtains that $\left\{x\in \Omega;\ V(x)>t\right\}$ is measurable for all $t\in\R$, and that
\[
\mu_V(t)=\left\vert \left\{x\in \Omega;\ \varphi(x)>F^{-1}(\mu(t))\right\}\right\vert=\left\vert \Omega_{F^{-1}(\mu(t))}\right\vert=\mu(t),
\]
where the last equality uses the definition of $F$. Finally, $\left\vert \left\{x\in \Omega;\ V(x)>\beta\right\}\right\vert=\mu(\beta)=0$ and, for all $s<\alpha$,
\[
\left\vert \left\{x\in \Omega;\ V(x)\leq s\right\}\right\vert=\left\vert \Omega\right\vert-\mu(s)=0,
\]
which shows that $V\in L^{\infty}(\Omega)$. \hfill\fin


\subsection{Estimates of $\lambda_1(B^n_R,\tau e_r)$ as $\tau\to+\infty$}\label{sec73}

We recall that $\lambda_1(\Omega,v)$ is defined as $\lambda_1(\Omega,\hbox{Id},v,0)$ for $v\in L^{\infty}(\Omega,\R^n)$. We call $B^n_R$ the open Euclidean ball of $\R^n$ with center $0$ and radius $R>0$, and we set
$$G_n(m,\tau)=\lambda_1(B^n_{(m/\alpha_n)^{1/n}},\tau e_r)$$
for all $n\in\N\backslash\{0\}$, $m>0$ and $\tau\ge 0$. Notice that $G_n(m,\tau)$ is always positive.\par
Our goal here is to discuss the behavior of $G_n(m,\tau)$ for large $\tau$. Indeed, if, in Theorem~\ref{faberkrahn}, $\Lambda$ is a constant $\gamma>0$, then, with the same notations as in Theorem~\ref{faberkrahn},
$$\baa{rcl}
\lambda_1(\Omega,A,v,V)\ \ge\  \lambda_1(\Omega^*,\gamma\hbox{Id},\tau_1e_r,-\tau_2)
& =& \gamma\lambda_1(\Omega^*,\hbox{Id},\tau_1\gamma^{-1},0)-\tau_2\\
& = & \gamma G_n(|\Omega|,\tau_1\gamma^{-1})-\tau_2.\eaa$$
The constants $\gamma$ and $\tau_2$ appear as multiplicative and additive constants in the previous formula. The function $[0,+\infty)\ni\tau\mapsto G_n(m,\tau)>0$ is obviously continuous, and decreasing (as a consequence of Theorem~\ref{fixedlinftyball}). However, the behaviour when $\tau\to+\infty$ is not immediate. It is the purpose of the following lemma, which was used in Remark~\ref{remlp}. When $\Lambda$ is not constant in Theorem~\ref{faberkrahn} but still satisfies some given lower and upper bounds, the following lemma provides some bounds of $\lambda_1(\Omega^*,\Lambda^*\hbox{Id},\tau_1e_r,-\tau_2)$ when $\tau_1\to+\infty$.

\begin{lem}\label{lemGn} For all $m>0$, $\tau^{-2}e^{\tau m/2}G_1(m,\tau)\to 1$ as $\tau+\infty$, and one even has
\be\label{F1}
\exists\ C(m)\ge 0,\ \exists\ \tau_0\ge 0,\ \forall\ \tau\ge\tau_0,\quad|\tau^{-2}e^{\tau m/2}G_1(m,\tau)-1|\le C(m)\tau e^{-\tau m/2}.
\ee
Moreover, for all $n\geq 2$ and $m>0$, $G_n(m,\tau)>G_1(2(m/\alpha_n)^{1/n},\tau)$ for all $\tau\ge 0$, and
\be\label{Fntau}
-\tau^{-1}\log G_n(m,\tau)\to m^{1/n}\alpha_n^{-1/n}\hbox{ as }\tau\rightarrow +\infty.
\ee
\end{lem}

In \cite{fr}, with probabilistic arguments, Friedman proved some lower and upper logarithmic estimates, as $\varepsilon\to 0^+$, for the first eigenvalue of general elliptic operators $-a_{ij}\varepsilon^2\partial_{ij}+b_i\partial_i$ with $C^1$ drifts $-b=-(b_1,\ldots,b_n)$ pointing inwards on the boundary (see also \cite{v}). Apart from the fact that the vector field $e_r$ is not $C^1$ at the origin, the general result of Friedman would imply the asymptotics (\ref{Fntau}) for $\log G_n(m,\tau)=\log \lambda_1(B^n_{(m/\alpha_n)^{1/n}},\tau e_r)$. For the sake of completeness, we give here a proof of (\ref{Fntau}) with elementary analytic arguments. Lemma \ref{lemGn} also provides the precise equivalent of $G_1(m,\tau)$ for large $\tau$. However, giving an equivalent for $G_n(m,\tau)$ when $\tau$ is large and $n\geq 2$ is an open question.\hfill\break

\noindent{\bf{Proof of Lemma \ref{lemGn}.}} First, to prove (\ref{F1}), fix $m>0$ and $\tau\geq 0$, set $\Omega=(-R,R)$ with $2R=m$, and
denote
$$\lambda=\lambda_1(\Omega, \tau e_r)$$
and $\varphi=\varphi_{\Omega,\hbox{Id},\tau e_r,0}$. Theorem \ref{fixedlinftyball} ensures that $\varphi$
is an even function, decreasing in $[0,R]$ and that
\[
-\varphi^{\prime\prime}(r)+\tau\varphi^{\prime}(r)=\lambda\varphi(r)\mbox{ for all }0\leq r\le R,
\]
with $\varphi(R)=0$, $\varphi>0$ in $(-R,R)$ and $\varphi^{\prime}(0)=0$ (in particular, the above equality holds in the classical sense in $[0,R]$). For all $s\in [0,\tau R]$, define $\psi(s)=\varphi(s/\tau)$, so that $\psi$
satisfies the equation
\[
-\psi^{\prime\prime}(s)+\psi^{\prime}(s)=\frac{\lambda}{\tau^2}\psi(s)\mbox{ for all }0\leq s\le\tau R,
\]
with $\psi(\tau R)=0$ and $\psi^{\prime}(0)=0$. Notice that $\lambda$ depends on $\tau$, but since, for all $\tau\geq 0$,
$0<\lambda\leq \lambda_1((-R,R),0)$, there exists $\tau_0>0$ such that $\tau^2\geq 4\lambda$ for all $\tau\geq \tau_0$, and we
will assume that $\tau\geq \tau_0$ in the sequel. The function $\psi$ can be computed explicitly: there exist $A,B\in \R$ such
that, for all $0\leq s\leq \tau R$,
\[
\psi(s)=Ae^{\mu_+r}+Be^{\mu_-r},
\]
where $\displaystyle \mu_{\pm}=(1\pm \sqrt{1-4\lambda/\tau^2})/2$. Using the boundary values of $\psi$ and
$\psi^{\prime}$, one obtains after straightforward computations:
\[
\lambda=\frac{\tau^2}4 \left(1+\sqrt{1-\frac{4\lambda}{\tau^2}}\right)^2 e^{-\sqrt{1-\frac{4\lambda}{\tau^2}}\tau R}.
\]
Since $\lambda$ remains bounded when $\tau\rightarrow +\infty$, it is then straightforward to check that $\lambda\sim \tau^2e^{-\tau R}$ when $\tau\rightarrow +\infty$, and that (\ref{F1}) follows.\par
We now turn to the proof of assertion (\ref{Fntau}). Let $n\geq 2$, $m>0$, $\tau\geq 0$ and $\Omega=B^n_R$ be such that $\left\vert \Omega\right\vert=m$, so that one has $R=(m/\alpha_n)^{1/n}$ and $G_n(m,\tau)=\lambda_1(\Omega,\tau e_r)$. We first claim that
\[
G_n(m,\tau)> G_1(2R,\tau).
\]
Indeed, write
$$\lambda=\lambda_1(\Omega,\tau e_r)\ \hbox{ and }\varphi_n=\varphi_{\Omega,\hbox{Id},\tau e_r,0}.$$
Similarly,
$G_1(2R,\tau)=\lambda_1((-R,R),\tau e_r)$, and we denote $\mu=\lambda_1((-R,R),\tau e_r)$ and $\varphi_1=\varphi_{(-R,R),\hbox{Id},\tau
e_r,0}$ (where $\hbox{Id}$ is then understood as the $1\times 1$ identity matrix). As before, define $\psi_n(y)=\varphi_n(y/\tau)$ for all $y\in \tau\overline{\Omega}=\overline{B^n_{\tau R}}$ and $\psi_1(r)=\varphi_1(r/\tau)$ for all
$r\in[-\tau R,\tau R]$. Finally, since $\psi_n$ is radially symmetric, let $u_n:\left[0,\tau R\right]\rightarrow \R$ such that
$\psi_n(y)=u_n(\left\vert y\right\vert)$ for all $y\in \tau\overline{\Omega}=\overline{B^n_{\tau R}}$. One has
\begin{equation} \label{onedim}
\left\{
\begin{array}{ll}
\displaystyle -u_n^{\prime\prime}(r)-\frac{n-1}ru_n^{\prime}(r)+u_n^{\prime}(r)=\frac{\lambda}{\tau^2}u_n(r)& \mbox{ in }(0,\tau
R],\\
\\
\displaystyle -\psi_1^{\prime\prime}(r)+\psi_1^{\prime}(r)=\frac{\mu}{\tau^2}\psi_1(r) &\mbox{ in }[0,\tau R],
\end{array}
\right.
\end{equation}
with $u_n^{\prime}(0)=u_n(\tau R)=0$, $\psi^{\prime}_1(0)=\psi_1(\tau R)=0$. 

Assume that $\lambda\leq \mu$. Since $u_n^{\prime}<0$ in $(0,\tau R]$ and $u_n\ge 0$, one obtains
\be\label{onedimbis}
\left\{
\begin{array}{ll}
\displaystyle -u_n^{\prime\prime}(r)+u_n^{\prime}(r)\leq \frac{\mu}{\tau^2}u_n(r) &\mbox{ in }[0,\tau R],\\
\\
\displaystyle -\psi_1^{\prime\prime}(r)+\psi_1^{\prime}(r)=\frac{\mu}{\tau^2}\psi_1(r) &\mbox{ in }[0,\tau R].
\end{array}
\right.
\ee
Since $\psi_1^{\prime}(\tau R)<0$ by Hopf lemma, while $\psi_1(r)>0$ in $[0,\tau R)$, $u_n(r)>0$ in $[0,\tau R)$ and the functions
$u_n$ and $\psi_1$ belong (at least) to $C^1([0,\tau R])$, there exists then $\gamma>0$ such that
$\gamma\psi_1(r)>u_n(r)$ for all $0\le r<\tau R$. Define $\gamma^{\ast}$ as the infimum of all the $\gamma>0$ such that
$\gamma\psi_1>u_n$ in $[0,\tau R)$, observe that $\gamma^{\ast}>0$ and define $z=\gamma^{\ast}\psi_1-u_n$ which is non-negative
in $[0,\tau R]$ and satisfies
\begin{equation} \label{ineqdimone}
-z^{\prime\prime}(r)+z^{\prime}(r)-\frac{\mu}{\tau^2}z(r)\ge 0
\end{equation}
for all $0\le r\le\tau R$ and $z(\tau R)=0$. 

Assume that there exists $0<r<\tau R$ such that $z(r)=0$. The strong maximum principle shows that $z$ is identically zero in
$[0,\tau R]$, which means that $\gamma^{\ast}\psi_1=u_n$ in $[0,\tau R]$, and even that $\psi_1=u_n$ because
$\psi_1(0)=u_n(0)=1$. But this is impossible according to (\ref{onedim}) and (\ref{onedimbis}).

Thus, $z>0$ in $(0,\tau R)$. Furthermore, $z'(0)=0$, hence $z(0)>0$ from Hopf lemma. Another application of Hopf lemma shows that $z^{\prime}(\tau R)<0$. Therefore, there exists $\kappa>0$ such that $z>\kappa u_n$ in $[0,\tau R)$, whence
$$\frac{\gamma^*}{1+\kappa}\psi_1>u_n\ \hbox{ in }[0,\tau R),$$
which is a
contradiction with the definition of $\gamma^{\ast}$. 

Finally, we have obtained that $\mu<\lambda$, which means that $G_n(m,\tau)>G_1(2R,\tau)$.

We now look for a reverse inequality. To that purpose, let $\varepsilon\in(0,1)$ and $R_0>0$ such that $\displaystyle
\frac{n-1}{R_0}<\varepsilon$. In the following computations, we always assume that $\tau R>R_0$. Define $u_n$ and $\lambda$ as
before. Let
$$\displaystyle \mu^{\prime}=\lambda_1\left(\left(-\left(R-\frac{R_0}{\tau}\right),\left(R-\frac{R_0}{\tau}\right)\right), \tau(1-\varepsilon)e_r\right)$$
and $w$ the normalized corresponding eigenfunction, so that
\[
\left\{
\begin{array}{l}
\displaystyle -w^{\prime\prime}(r)+\tau (1-\varepsilon) w^{\prime}(r)=\mu^{\prime}w(r) \displaystyle \mbox{ in }\left[0,R-\frac{R_0}{\tau}\right],\\
\\
\displaystyle w^{\prime}(0)=0,\ w>0\hbox{ in }\displaystyle\left[0,R-\frac{R_0}{\tau}\right),\ w\left(R-\frac{R_0}{\tau}\right)=0.
\end{array}
\right.
\] 
For all $R_0\leq x\leq \tau R$, define $\displaystyle v(x)=w\left(\frac{x-R_0}{\tau}\right)$, which satisfies
\[
\left\{
\begin{array}{l}
\displaystyle -v^{\prime\prime}(r)+(1-\varepsilon)v^{\prime}(r)=\frac{\mu^{\prime}}{\tau^2}v(r)\ \mbox{ in }\left[R_0,\tau
R\right],\\
\\
\displaystyle v^{\prime}(R_0)=0,\ v>0\hbox{ in }[R_0,\tau R),\ v(\tau R)=0.
\end{array}
\right.
\] 
Assume that $\lambda\geq \mu^{\prime}$. Since $(n-1)/R_0<\varepsilon$ and $u_n'(r)<0$ in $(0,\tau R]$, one therefore has
\[
\left\{
\begin{array}{ll}
\displaystyle -u_n^{\prime\prime}(r)+(1-\varepsilon)u_n^{\prime}(r)\geq \frac{\mu^{\prime}}{\tau^2}u_n(r) &\mbox{ in }[R_0,\tau
R],\\
\\
\displaystyle -v^{\prime\prime}(r)+(1-\varepsilon)v^{\prime}(r)=\frac{\mu^{\prime}}{\tau^2}v(r) &\mbox{ in }[R_0,\tau R].
\end{array}
\right.
\]
Arguing as before, we see that there exists $\gamma>0$ such that $\gamma u_n>v$ in $[R_0,\tau R)$. Define $\gamma^{\ast}$ ($>0$) as the
infimum of all such $\gamma$'s and define $z=\gamma^{\ast}u_n-v$, which is nonnegative in $[R_0,\tau R]$ and satisfies
$\displaystyle -z^{\prime\prime}+(1-\varepsilon)z^{\prime}-(\mu^{\prime}/\tau^2)z\ge 0$ in $[R_0,\tau R]$. 

Assume that $z(r)=0$ for some $r\in (R_0,\tau R)$. The strong maximum principle ensures that $z$ is $0$ in $[R_0,\tau R]$, which
means that $u_n=v$ in $[R_0,\tau R]$, which is impossible because $u_n^{\prime}(R_0)<0=v^{\prime}(R_0)$. 

Therefore, $z>0$ everywhere in $(R_0,\tau R)$. Furthermore, $z'(R_0)<0$, thus $z(R_0)>0$. On the other hand, by Hopf lemma, $z^{\prime}(\tau R)<0$. Thus, there exists $\kappa>0$ such that $z>\kappa v$ in $[R_0,\tau R)$, whence $(\gamma^*/(1+\kappa))u_n>v$ in $[R_0,\tau R)$. This contradicts the definition of $\gamma^{\ast}$.

Thus, we have established that $\lambda<\mu^{\prime}$. Straightforward computations (similar to those of the proof of (\ref{F1})) show
that
\[
\lambda<\mu^{\prime}=\frac{\tau^2}4\left(1-\varepsilon+\sqrt{(1-\varepsilon)^2-\frac{4\mu'}{\tau^2}}\right)^2
e^{-\sqrt{(1-\varepsilon)^2-\frac{4\mu'}{\tau^2}}(\tau R-R_0)},
\]
and, since $\lambda>G_1(2R,\tau)$, formula (\ref{F1}) and the fact that $m=\alpha_n R^n$ end the
proof of (\ref{Fntau}).\hfill\fin


\end{document}